\documentclass[11pt]{article}

\usepackage{amsmath,amssymb,amsfonts,enumitem,amsthm,accents,pstricks,pst-node,pstricks-add,mathrsfs,xy}
\usepackage{multicol}
\xyoption{all}

\usepackage{color}
\usepackage{makeidx}
\usepackage{graphicx}
\usepackage{setspace}
\usepackage{extarrows}
\numberwithin{equation}{section}

\theoremstyle{plain}
\newtheorem{theorem}{Theorem}
\numberwithin{theorem}{subsection}
\newtheorem{lemma}[theorem]{Lemma}
\newtheorem{conjlemma}[theorem]{Conjectural Lemma}
\newtheorem{corollary}[theorem]{Corollary}
\newtheorem{proposition}[theorem]{Proposition}
\newtheorem{conjecture}[theorem]{Conjecture}
\newtheorem{notation}[theorem]{Notation}
\newtheorem{definition}[theorem]{Definition}
\newtheorem{condition}[theorem]{Condition}

\theoremstyle{definition}
\newtheorem{example}[theorem]{Example}
\newtheorem{exercise}[theorem]{Exercise}

\theoremstyle{remark}
\newtheorem{remark}[theorem]{Remark}

\def\ov#1{\overline{#1}}
\def\tn#1{\textnormal{#1}}
\def\mf#1{\mathfrak{#1}}
\def\wt#1{\widetilde{#1}}
\def\wh#1{\widehat{#1}}
\def\vr{\varrho}
\def\ll{\left\langle}
\def\rr{\right\rangle}
\def\mc{\mathcal}
\def\lra{\longrightarrow}
\def\dbar{\bar\partial}
\def\dpst{\displaystyle}
\def\ve{\varepsilon}
\def\hb{\hbar}

\newcommand{\abs}[1]{\left\vert #1 \right\vert}
\newcommand{\lrp}[1]{\left( #1 \right)}
\newcommand{\lrc}[1]{\left\{ #1 \right\}}
\newcommand{\one}{\tn{\tiny{(1)}}}
\newcommand{\two}{\tn{\tiny{(2)}}}

\newcommand{\zero}{\tn{\tiny{(0)}}}

\newcommand{\tb}{\tn{\tiny{(b)}}}
\newcommand{\tc}{\tn{\tiny{(c)}}}

\newcommand{\sk}{\tn{\tiny{(k)}}}
\newcommand{\skk}{\tn{\tiny{(k+1)}}}

\def\bEqu#1{\begin{equation}\label{#1}}
\def\eEqu{\end{equation}}

\def\bsEq{\begin{equation*}}
\def\esEq{\end{equation*}}

\def\bDef#1{\begin{definition}\label{#1}}
\def\eDef{\end{definition}}

\def\bThm#1{\begin{theorem}\label{#1}}
\def\eThm{\end{theorem}}

\def\bCon#1{\begin{conjecture}\label{#1}}
\def\eCon{\end{conjecture}}

\def\bLem#1{\begin{lemma}\label{#1}}
\def\eLem{\end{lemma}}

\def\bCLm#1{\begin{conjlemma}\label{#1}}
\def\eCLm{\end{conjlemma}}

\def\bRem#1{\begin{remark}\label{#1}}
\def\eRem{\end{remark}}

\def\bExa#1{\begin{example}\label{#1}}
\def\eExa{\end{example}}

\def\bExe#1{\begin{exercise}\label{#1}}
\def\eExe{\end{exercise}}

\def\bPro#1{\begin{proposition}\label{#1}}
\def\ePro{\end{proposition}}

\def\bCor#1{\begin{corollary}\label{#1}}
\def\eCor{\end{corollary}}

\def\bFig#1{\begin{figure}\label{#1}}
\def\eFig{\end{figure}}

\def\bProof{\begin{proof}}
\def\eProof{\end{proof}}

\def\bItem{\begin{itemize}[leftmargin=*]}
\def\eItem{\end{itemize}}

\def\bEnum{\begin{enumerate}[label=$(\arabic*)$,leftmargin=*]}
\def\eEnum{\end{enumerate}}

\def\pr{\tn{pr}}
\def\nd{\tn{d}}

\def\id{\tn{id}}
\def\vir{\tn{vir}}
\def\vfc{\tn{VFC}}
\def\st{\tn{st}}
\def\ev{\tn{ev}}
\def\ker{\tn{ker}}
\def\coker{\tn{coker}}
\def\neck{\tn{neck}}

\def\multi{\tn{multi}}
\def\GW{\tn{GW}}
\def\top{\tn{top}}
\def\Def{\tn{Def}}
\def\Obs{\tn{Obs}}

\def\cA{\mc{A}}
\def\cB{\mc{B}}
\def\cJ{\mc{J}}
\def\cM{\mc{M}}
\def\cO{\mc{O}}
\def\cE{\mc{E}}

\def\cP{\mc{P}}
\def\cK{\mc{K}}
\def\cN{\mc{N}}

\def\cY{\mc{Y}}
\def\cZ{\mc{Z}}
\def\cW{\mc{W}}
\def\cS{\mc{S}}
\def\cL{\mc{L}}
\def\cV{\mc{V}}
\def\cU{\mc{U}}
\def\cC{\mc{C}}
\def\cH{\mc{H}}
\def\cW{\mc{W}}
\def\cT{\mc{T}}

\def\R{\mathbb R}
\def\C{\mathbb C}
\def\Z{\mathbb Z}
\def\Q{\mathbb Q}
\def\P{\mathbb P}
\def\H{\mathbb H}

\def\N{\mathbb N}

\def\mfi{\mf{i}}
\def\mfj{\mf{j}}

\def\mfq{\mf{q}}
\def\aut{\mf{Aut}}
\def\mfD{\mf{D}}
\def\mfB{\mf{B}}
\def\mfd{\mf{d}}
\def\mfs{\mf{s}}
\def\mft{\mf{t}}
\def\mfz{\mf{z}}
\def\mfL{\mf{L}}

\def\la{\lambda}

\def\ep{\epsilon}
\def\De{\Delta}
\def\de{\delta}
\def\om{\omega}
\def\Om{\Omega}
\def\si{\sigma}
\def\Si{\Sigma}
\def\al{\alpha}

\def\Ga{\Gamma}
\def\ze{\zeta}

\makeindex
\begin{document}
\title{Gromov-Witten Theory\\ via Kuranishi Structures}
\author{Mohammad F.~Tehrani and Kenji Fukaya\thanks{Partially supported by NSF Grant No. 1406423 and Simons Collaboration on Homological Mirror Symmetry.} }
\date{\today}
\maketitle

\section*{Preface}\label{sec:preface}

In the past few decades, moduli spaces have been a major tool to study families of geometric objects and to define geometric invariants of various kind for these spaces. 
Among various moduli spaces that have been extensively studied, moduli space of pseudoholomorphic maps have been a corner stone of various subjects in symplectic geometry and string theory. For example, through these moduli spaces we can define Gromov-Witten invariants; these are rational numbers that, in certain situations, count pseudoholomorphic curves meeting prescribed conditions in a given symplectic manifold.  

\noindent
In physics, more specifically in topological sigma models, GW~invariants were first studied by Witten \cite{Wi}.
Conjectural predictions of mirror symmetry \cite{COGP} relating the GW~invariants of a symplectic manifold on one side, and the solutions of Picard-Fuchs equation over some ``mirror" family of complex manifolds on the other side, raised the efforts and interests to rigorously define these numbers mathematically. 
Naively, GW~invariants are obtained by integrating cohomology classes against the fundamental class of the moduli space; equally, they are given by intersecting cycles on the moduli space. 

\noindent
However, a priori, every moduli space of pseudoholomorphic maps of some fixed topological type does not carry the structure of a (compact) differentiable orbifold in general; this is mainly due to the transversality issues of the corresponding Fredholm operators. 
Nevertheless, one would like to at least associate a ``virtual'' fundamental class to the moduli space in question.

\noindent
In past twenty years, several different\footnote{But to some extent similar.}  approaches to the construction of virtual fundamental class have emerged, each of which has its own advantages, drawbacks, and motivations.
Among all proposed analytical/topological approaches, we can mention the works of  Li-Tian~\cite{LT}, Fukaya-Ono~\cite{FO}, Siebert~\cite{Si}, Hofer~\cite{Ho},  and many other recent works.

\noindent
Construction of virtual fundamental class for moduli spaces via Kuranishi models is a traditional technique for describing the local structure of moduli spaces cut out by non-linear equations whose linearization is Fredholm.  
First in the ``semi-positive" case, the transversality issue was handled by Ruan-Tian \cite{RT} via  global inhomogeneous perturbations of the Cauchy-Riemann equation. 
The Kuranishi method can be regarded as a local and multi-valued version of the perturbation method.

\noindent
A more elaborate version of Kuranishi structure is used by Fukaya, Oh, Ohta, and Ono \cite{FOOO}  over moduli spaces of pseudoholomorphic discs and strips to define and study the Lagrangian intersection Floer homology of two Lagrangian submanifolds. Such moduli spaces can also be used to define GW-type invariants counting pseudoholomorphic discs (and higher genus analogues of that) with boundary on Lagrangian submanifolds.

\noindent
In this manuscript, we review the construction of GW~virtual fundamental class via Kuranishi structures for moduli spaces of pseudoholomorphic maps defined on closed Riemann surfaces. We consider constraints coming from the ambient space and Deligne-Mumford moduli, called primary insertions, as well as intrinsic classes such as $\psi$-classes and Hodge classes. Readers interested in the Floer theoretic aspects of Kuranishi structures should consult  \cite{FOOO}.

\noindent
The first part of this article, i.e. Sections 2-4, is about abstract and topological aspects of Kuranishi structures and is fairly self-contained. The second part, i.e. Sections 5 and 6, is about construction of a natural class of Kuranishi structures for moduli spaces of interest. In part 2, some of the main analytical steps, such as the gluing theorem, are stated without proof. For the proof we refer to \cite{FOOO-detail} and \cite{FOOO-detail3}.

\noindent
\textbf{Disclaimer:}
The construction provided in this article is a repackaging of the relevant materials in \cite{FO, FOOO-detail,FOOO-detail2,FOOO-detail3}, and the seminar series delivered by Fukaya at the Simons center during Spring semester of 2014, which are available at SCGP's video portal \cite{F-SCGP}. It also relies on our discussions in the past two years to simplify or rephrase some of the statements and definitions. For example, we use ``dimensionally graded systems" in place of ``good coordinate systems" to simplify\footnote{A dimensionally graded system is a good coordinate system where the index set is just $\Z$.}  the notation and avoid complications of working with partially ordered sets.
Our main priority has been to collect all the concepts and ideas which are relevant to Gromov-Witten theory in the work of FOOO in one relatively short article. 
There is no doubt that some of the techniques, details, and statements in this note do not directly come  from the original article of Fukaya and Ono \cite{FO} or the book of Fukaya-Oh-Ohta-Ono \cite{FOOO}.  
In \cite{FOOO-detail}, motivated by feedback to and questions about \cite{FO} and \cite{FOOO}, FOOO provided more details and further arguments to explain their construction in \cite{FO} and \cite{FOOO}.
We have adapted some of this material here.
It is beyond the scope and goal of this text to cover the long history and all the contributions of various people that have led to our current understanding of virtual fundamental classes. For that, we refer the interested readers to \cite{FOOO-detail} and \cite{MW2}. In Sections~1-3 of \cite{MW2}, McDuff and Wehrheim give a detailed description of their approach and the motivations behind it. They compare their approach to that of FOOO and various others, and list some issues that they believe have not been addressed properly before. In \cite[Part 6]{FOOO-detail}, FOOO discuss the questions that were raised\footnote{Originally, as a Google group discussion.} about the details of their approach, provide some history, and address these questions with details. An interested reader may read both to get  a sense of how things have evolved to their current form.
Last but not least, some of the arguments in this note have benefited from many helpful and exciting discussions  Mohammad Tehrani had with Dusa McDuff, John Morgan, and Gang Tian. We thank them especially for their interest and support.

\begin{spacing}{.7}
\tableofcontents
\end{spacing}

\section{Introduction}\label{sec:intro}
Enumerative geometry, originally a branch of algebraic geometry with a history of more than a century, concerns certain counts of algebraic objects, mainly algebraic curves, in generic situations. 
This has been a challenging field of study for various reasons; the category of algebraic (more generally holomorphic) functions does not admit partitions of unity, algebraic objects are fairly rigid, and a precise formulation of the necessary intersection theory needs elaborate techniques (etc.). 
On the other hand, these tools are easily accessible and play a fundamental role in topology and smooth category.

\noindent
Symplectic manifolds are a softened version of K\"ahler manifolds, such as complex projective varieties, where we drop the integrality condition of the complex structure and mainly focus on the underlying non-degenerate two form which we call the symplectic form. 
At a quick glance, symplectic manifolds and complex projective varieties
seem very far apart. 
Surprisingly, in \cite{G}, Gromov combined the rigidity of algebraic geometry with the flexibility of the smooth category and initiated the use of ``pseudoholomorphic curves" as a generalization of holomorphic curves, or better said ``parametrized holomorphic curves''.
These allow the formulation of symplectic analogues of enumerative questions from algebraic geometry as well-defined invariants of symplectic manifolds.
Moreover, from a purely symplectic perspective, these invariants give us a tool to distinguish isotopy classes of symplectic manifolds which are indistinguishable as smooth manifolds, see \cite{Ruold}.

\noindent
Other important applications of pseudoholomorphic curves include (now classical) Floer homology for giving lower bounds for the number of fixed points of Hamiltonian diffeomorphisms, initiated by Floer \cite{F1986} and further enhanced to prove Arnold conjecture in some cases; e.g. see \cite{HS} and \cite{On}. 
This idea was also extended to study the intersection of a pair of Lagrangians, known as Lagrangian intersection Floer theory \cite{FOOO}, in an arbitrary symplectic manifold.  

\noindent
The invariants of symplectic manifolds obtained by counting the number of pseudoholomorphic maps has also relations 
to topological sigma models in string theory physics. 
Among various works in this direction, the 1988 paper of Witten \cite{Wi} laid a background for further studies. 
Since then, the invariants obtained by the counts of pseudoholomorphic maps into symplectic manifolds  are called \textbf{Gromov-Witten} invariants.

\noindent
Let us begin with a review of ``classical" Gromov-Witten theory, as a theory of counting pseudoholomorphic maps; we refer interested reader unfamiliar with the basics of this subject to \cite{MS2004}.
Let $(X^{2n},\om)$ be a symplectic manifold, which we assume to be closed  (i.e. compact without boundary) throughout this article. An almost complex structure $J$ on $TX$ is an endmorphism $J\colon TX\!\lra\! TX$ such that $J^2\!=\!-\id$. We say $J$ is $\om$-compatible if $g_J\!\equiv\! \om(\cdot,J\cdot)$ is a metric; in this case the action of $J$ is an isometry (with respect to $g_J$). More generally, we say $J$ is $\om$-tame if 
$$
\om(v,Jv)\!>\!0\qquad \forall~v\!\in\! TX.
$$ 
Let $\cJ(X,\om)$ and $\cJ_\tau(X,\om)$ denote the space\footnote{With its natural topology.} of smooth almost complex structures that are $\om$-compatible and $\om$-tame, respectively. Fortunately, each of these spaces is contractible; therefore, any number associated to a tuple $(X,\om,J)$ which is invariant under the deformation of $J$ within these spaces is an invariant of $(X,\om)$.   

\noindent
In the K\"ahler category, the complex structure gives an almost complex structure on the underlying real manifold which is compatible with the K\"ahler symplectic form. However, over a  general symplectic manifold,  the Nijenhuis tensor 
\bEqu{equ:Nij}
\aligned
&N_J \!\in\! \Gamma(X,\Om^2_X\otimes TX),\\
&N_J(u,v)\!\equiv\![u,v]\!+\!J[u,Jv]\!+\!J[Ju,v]\!-\![Ju,Jv]\quad\forall u,v\!\in\! TX,
\endaligned
\eEqu
calculates the deviation of a compatible triple $(X,\om,J)$ from being K\"ahler; i.e. $(X,\om,J)$ is K\"ahler if and only if $N_J\!\equiv \!0$; see \cite[Section 2.1]{MS2004}.

\subsection{Moduli space of pseudoholomorphic maps}\label{sec:basicJhol}
Given a fixed $J\!\in\!\cJ_\tau(X,\om)$,  $A\!\in\!\! H_2(X,\Z)$, and $g,k\!\in\! \!\Z^{\geq 0}$, a $k$-marked genus $g$ degree $A$ (i.e. the image in $X$ represents the homology class $A$) 
\textbf{pseudoholomorphic} map, or $J$-holomorphic map, is a differentiable map\footnote{More precisely, a $W^{\ell,p}$-smooth map for some $p\!>\!1$ with $\ell p\!>\!2$, where $W^{\ell,p}(\Si,X)$ is the space of $\ell$-times differentiable maps whose $\ell$-th derivative has finite $L_p$-norm.} $u\colon \Si\!\lra\! X$ from a closed genus $g$ Riemann surface (complex curve) $(\Si,\mfj)$ with $k$ distinct ordered points 
$$
\vec{z}\!=\!(z^1,\ldots,z^k)\!\subset\! \Si
$$ 
to $X$ that satisfies the \textbf{Cauchy-Riemann equation}
\bEqu{equ:J-holo}
\Gamma(\Sigma,\Omega_{(\Sigma,\mfj)}^{0,1}\otimes_\C u^* TX) \ni \dbar u= \frac{1}{2}\lrp{\nd u+J\circ \nd u \circ \mfj}\equiv 0.
\eEqu
It follows from an elliptic bootstrapping argument, cf. \cite[Theorem B.4.1]{MS2004}, that every such map is automatically smooth (more precisely, it is as smooth as the almost complex structure). 
For every positive integer $N$ let 
$$
[N]\!\equiv\!\{1,2,\ldots,N\}.
$$
Two $k$-marked\footnote{We write an ordered set of marked points as $\vec{z}\!\equiv\! (z^1,\ldots,z^k)\! \subset\! \Si$. Thus, in our notion, upper indices indicate the ordering of the marked points and lower indices, if present, take care of other labelings of the marked points.} $J$-holomorphic maps $u_i: (\Si_i,\mfj_i,\vec{z}_i)\!\lra\! X$, with $i\!=\!1,2$, are said to be \textbf{equivalent} if there is a holomorphic diffeomorphism 
\bEqu{equ:autoh}
\varphi\colon(\Si_1,\mfj_1)\lra(\Si_2,\mfj_2) ~~\tn{s.t.}~~ u_2\!=\!u_1\circ \varphi\quad\tn{and}\quad\varphi(z_1^a)\!=\!z_2^a\quad \forall a\!\in\! [k].
\eEqu
Let $\aut \big(\Si,\mfj\big)$ be the group of biholomorphisms of $(\Si,\mfj)$, $\aut \big(\Si,\mfj,\vec{z}\big)$ be the subgroup fixing the marked points, and $\aut \big(u,\Si,\mfj,\vec{z}\big)$ be the subgroup that does not change the map and the marked points; we have
$$
\aut \big(u,\Si,\mfj,\vec{z}\big)\subset \aut \big(\Si,\mfj,\vec{z}\big)\subset \aut \big(\Si,\mfj \big).
$$
We denote by $\cM_{g,k}(X,A,J)$, or $\cM_{g,k}(X,A)$ for short, to be the space of equivalence classes  of $k$-marked genus $g$ degree $A$ $J$-holomorphic maps. This space carries a natural Haussdorf topology
and is in fact metrizable\footnote{The Gromov topology is paracompact, Hausdorff, and locally metrizable; hence, by Smirnov's theorem, it is metrizable. See Remark~\ref{rem:top-remark} for a comparison of various metrizability theorems.}. 
However, it is normally not compact and usually very singular. 
A priori, the Cauchy-Riemann equation (\ref{equ:J-holo}) is independent of the symplectic structure. Nevertheless, the symplectic form plays an important role; it provides an energy function which is essential for compactifying $\cM_{g,k}(X,A)$.

\noindent
In order to compactify the moduli space, we need to add ``proper" limits  of  sequences of $J$-holomorphic maps. 
By the celebrated compactness theorem of Gromov \cite[Theorem 15.B]{G}, a sequence of $J$-holomorphic maps in an almost complex manifold with a uniform area bound (with respect to some $J$-invariant metric $g$) has a subsequence which limits to a ``stable" connected set of $J$-holomorphic maps attached to each other along some nodes; see Theorem~\ref{thm:gromov}. 

\noindent
A marked (smooth) complex curve $C\!=\!(\Si,\mfj,\vec{z})$ is said to be \textbf{stable} if $\aut(C)$ is finite. This is the case whenever $2g+k\!\geq\! 3$. Similarly, a marked $J$-holomorphic map $f\!\equiv\! (u,C)$ is said to be stable if $\aut(f)$ is finite. A \textbf{nodal} marked complex curve (or Riemann surface) is a connected space obtained by considering a disjoint union of smooth curves and identifying them along some pairs of (unordered) distinguished points disjoint from the marked points; i.e.
$$
C\!=\!(\Si,\mfj,\vec{z})\!\equiv\!\bigg(\coprod_{i=1}^N \big(\Si_i,\mfj_i,\vec{z}_i,\{w_i^a\}_{a\in[k_i]}\big)\bigg)/\sim,\quad w_i^a\sim w_{\si(i,a)}^{\de(i,a)},
$$
where $(i,a)\!\lra\!(\si(i,a),\de(i,a))$ is a fixed point free involution on the set of pairs 
$$
\{(i,a)\colon i\!\in\! [N],~a\!\in\! [k_i]\};
$$
see Section~\ref{sec:stable} or \cite[Section 22.2]{mirror}. The ordered set of marked points on the union corresponds to an ordering on  
$$
 \vec{z}=\bigcup_{i=1}^N \vec{z}_i
 $$ 
 which respects the ordering of each $\vec{z}_i$. The genus of a nodal curve is the sum of the genera of the connected components and the genus of its dual graph; see Equation~\ref{equ:genus}. 
\begin{figure}
\begin{pspicture}(-4,-5)(10,-4)
\psset{unit=.3cm}
\psellipse[linewidth=.08](10,-14)(5,2)\pscircle[linewidth=.08](17,-14){2}  
\pscircle*(18.41,-12.59){.25}\pscircle*(18.41,-15.41){.25}\pscircle*(15,-14){.25}\pscircle*(10.5,-12.7){.25}
\rput(19.2,-12.2){$z^1$}\rput(19.4,-15.8){$z^3$}\rput(9.4,-12.7){$z^2$}
\psarc[linewidth=.08](10,-17){3.16}{70}{110}
\psarc[linewidth=.08](10,-11){3.16}{245}{295}
\end{pspicture}
\caption{A $3$-marked $1$-nodal curve of genus $1$.}
\label{fig:311}
\end{figure}
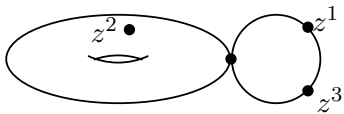
Figure~\ref{fig:311} illustrates a $3$-marked $1$-nodal curve of genus $1$ made of two smooth curves. A {nodal} $J$-holomorphic map is union of $J$-holomorphic maps on $\{\Si_i\}_{i\in[N]}$ that have identical values on each pair of identified points. 
Similarly, a nodal $J$-holomorphic map is said to be stable if its automorphism group is finite; see (\ref{equ:auto-nodal}) and (\ref{equ:auto-nodal-map}).
The degree of a nodal map is the sum of degrees of its components. 
Then, we denote by $\ov\cM_{g,k}(X,A)$ to be the space of equivalence classes of all $k$-marked genus $g$ degree $A$ $J$-holomorphic maps. The Deligne-Mumford space of stable curves, denoted by $\ov\cM_{g,k}$, corresponds to the special case where $X$ is just a  point. In order to make our notation consistent, in the unstable range $2g\!+\!k\!<\!3$, we define $\ov\cM_{g,k}$ to be a point!\footnote{In the original approach of \cite{FO}, $\ov\cM_{1,0}$ is taken to be the $j$-line; however, for the sake of consistency in Theorem~\ref{thm:VFC}, we take it to be a point in this article.}

\begin{figure}
\begin{pspicture}(-2,-5)(10,-1)
\psset{unit=.3cm}
\psellipse[linewidth=.08](0,-4)(5,2)\pscircle[linewidth=.08](7,-4){2}  
\pscircle*(5,-4){.25}\pscircle*(9,-4){.25}\rput(10,-4){$z^1$}
\psarc[linewidth=.08](0,-7){3.16}{70}{110}
\psarc[linewidth=.08](0,-1){3.16}{245}{295}
\rput(6,0){$u\neq  \tn{constant}$}
\psline[linewidth=.08]{->}(5,-.6)(6,-2)
\rput(15,-3){$\st$}
\psline[linewidth=.12]{->}(12,-4)(18,-4)
\psellipse[linewidth=.08](25,-4)(5,2)\pscircle*(30,-4){.25}\rput(31,-4){$z^1$}
\psarc[linewidth=.08](25,-7){3.16}{70}{110}
\psarc[linewidth=.08](25,-1){3.16}{245}{295}

\psellipse[linewidth=.08](0,-14)(5,2)\pscircle[linewidth=.08](7,-14){2}  
\pscircle[linewidth=.08](7,-10){2} 
\pscircle*(5,-14){.25}
\rput(1,-11){$u\neq  \tn{constant}$}
\psline[linewidth=.08]{->}(4,-11.6)(5.5,-12.5)
\pscircle*(7,-12){.25}\pscircle*(5.59,-8.59){.25}\pscircle*(8.41,-8.59){.25}  
\rput(4.7,-8.3){$z^1$}\rput(9.3,-8.3){$z^2$}
\psarc[linewidth=.08](0,-17){3.16}{70}{110}
\psarc[linewidth=.08](0,-11){3.16}{245}{295}
\rput(15,-13){$\st$}
\psline[linewidth=.12]{->}(12,-14)(18,-14)
\psellipse[linewidth=.08](25,-14)(5,2)\pscircle[linewidth=.08](32,-14){2}  
\pscircle*(30,-14){.25}\pscircle*(33.41,-12.59){.25}\pscircle*(33.41,-15.41){.25}
\rput(34.3,-12.2){$z^1$}\rput(34.4,-15.8){$z^2$}
\psarc[linewidth=.08](25,-17){3.16}{70}{110}
\psarc[linewidth=.08](25,-11){3.16}{245}{295}
\end{pspicture}
\caption{Examples of the stabilization process.}
\label{fig:st311}
\end{figure}
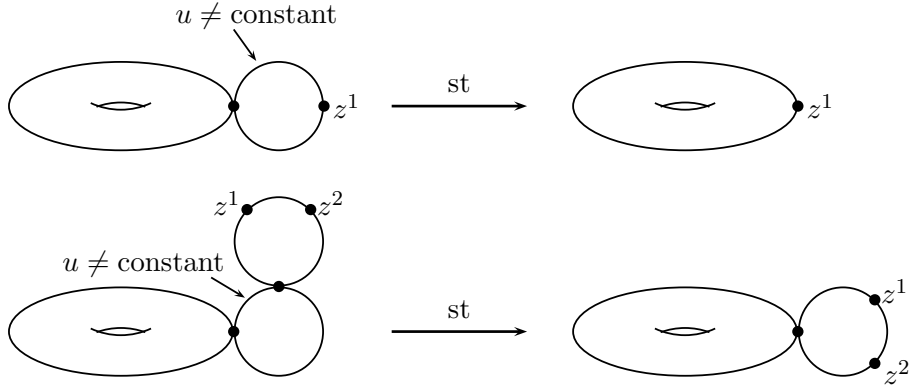

\noindent
Let 
\bEqu{equ:eval}
\ev_i\colon \ov\cM_{g,k}(X,A) \lra X, \quad \ev_i \big([u,\Si,\mfj,\vec{z}]\big)\!=\!u(z^i)\!\in\! X,
\eEqu
be the natural \textbf{evaluation} map at the marked points.
For every (equivalence class of) stable map $[u,\Si,\mfj,\vec{z}]$, forgetting about $u$, we obtain a nodal curve which a priori may not be stable; it may have spherical components with at most two distinguished and marked points, or $\Si$ could be a smooth genus $1$ curve with no marked point. In this situation, we can collapse every unstable component into a node or a marked point and obtain a  stable curve 
$[\ov{\Si},\ov{\mfj},\vec{z}]$ with possibly fewer components but the same number of marked points; see \cite[Section 26.3]{mirror} for a more through description.
This process gives us a well-defined \textbf{forgetful} map from the space of stable maps into the Deligne-Mumford space,
\bEqu{equ:st} 
\st\colon \ov\cM_{g,k}(X,A) \lra \ov\cM_{g,k}, \quad \st \big([u,\Si,\mfj,\vec{z}]\big) = [\ov{\Si},\ov\mfj,\vec{z}]. 
\eEqu
Figure~\ref{fig:st311} top shows an example of stabilization process where an spherical component with one marked point and one distinguished point is collapsed into a marked point on the next component. Figure~\ref{fig:st311} bottom shows an example where an spherical component with two distinguished points is collapsed into a node.
By a similar procedure, for every $I\!\subset\! [k]$ we obtain similar forgetful maps (this time we remove some of the marked points)
\bEqu{equ:pi-I} 
\aligned
&\pi_I\colon \ov\cM_{g,k}(X,A) \lra \ov\cM_{g,k-|I|}(X,A), \\ 
& \pi_I ([u,\Si,\mfj,(z^j)_{j\in[k]}]) = [u,\ov{\Si},\ov\mfj,(z^j)_{j\in [k]-I} ], 
\endaligned
\eEqu
that removes $(z^j)_{j\in I}$, collapses the possible unstable components, and relabels the remaining points by $1,\ldots,k\!-\!|I|$, preserving the ordering; cf.  \cite[Section 23.4]{mirror}.

\noindent
With this setup, the Gromov's Compactness Theorem~\ref{thm:gromov}, states that every sequence in $\ov\cM_{g,k}(X,A)$ has a subsequence with a well-defined limit in $\ov\cM_{g,k}(X,A)$ and this sequential convergence provides a compact Hausdorff  topology on $\ov\cM_{g,k}(X,A)$. Moreover, with respect to this topology (\ref{equ:eval}), (\ref{equ:st}), and (\ref{equ:pi-I})  are continuous.

\subsection{GW invariants}\label{sec:GW}
If  $ \ov\cM_{g,k}(X,A)$ has the structure of a smooth (in some suitable sense) oriented orbifold of the expected dimension 
\bEqu{equ:exp-dim}
d\!=\!2((n\!-\!3)(1\!-\!g)\!+\!c_1(A)\!+\!k)
\eEqu 
and the evaluation and forgetful maps are smooth, Gromov-Witten invariants are defined via integrals of the form
\bEqu{equ:GW}
\GW^{X}_{g,A}(\theta_1,\ldots,\theta_k; \kappa)= \int_{[\ov\cM_{g,k}(X,A)]^{FC}} \ev_1^*\theta_1\wedge \cdots\wedge \ev_k^*\theta_k \wedge \st^*\kappa,
\eEqu
where  $\theta_i$ are de Rham cohomology classes on $X$, $\kappa$ is a cohomology class on $\ov\cM_{g,k}$, and $[\ov\cM_{g,k}(X,A)]^{FC}$ is the fundamental class of the oriented orbifold structure (see Remark~\ref{rem:final-remarks}). 
Similarly, we may consider the push-forward of the singular homology fundamental class $[\ov\cM_{g,k}(X,A)]^{FC}\in H_*(\ov\cM_{g,k}(X,A),\Q)$ under the map 
\bEqu{equ:evst}
\ev\times \st\equiv \ev_1\times\cdots \times \ev_k\times \st\colon \ov\cM_{g,k}(X,A)\lra X^k\times \ov\cM_{g,k}
\eEqu
and define the image 
$$
\GW^{X}_{g,A}= (\ev\times \st)_*\big([\ov\cM_{g,k}(X,A)]^{FC}\big)\in H_d(X^k\times \ov\cM_{g,k},\Q)
$$ to be the \textbf{Gromov-Witten fundamental class} of $\ov\cM_{g,k}(X,A)$. 
In this sense, GW~invariants are defined via intersection of homology/cohomology classes in $X^k\!\times\! \ov\cM_{g,k}$ with the GW fundamental class.

\noindent
Except in some very special cases (e.g. $X\!\cong\! \C\P^n$ and $g\!=\!0$) $\ov\cM_{g,k}(X,A)$ does not carry a nice orbifold structure of the right dimension; it may have many components which are singular, intersect non-trivially, and have different dimensions! Treating these difficulties is the main point of this manuscript.
\subsection{Semi-positive case}\label{sec:positive}
There is still one case, the case of semi-positive manifolds, where the GW invariants can be defined by a direct count of geometric objects. Let us review this case; see the original article of Ruan-Tian \cite{RT1997} or the book of McDuff-Salamon \cite{MS2004} for more details.
A marked $J$-holomorphic map with smooth domain 
$$
f\!\!=\![u,\Si,\mfj,\vec{z}]\in \cM_{g,k}(X,A)
$$ 
is said to be \textbf{somewhere injective} if there is a point $z\!\in\! \Si$ such that 
$$
\nd u(z)\!\neq\! 0,\qquad u^{-1}(u(z))\!=\!\{z\}.
$$
Let $\cM_{g,k}^*(X,A)\subset \cM_{g,k}(X,A)$ denote the subspace of somewhere injective maps. 
By \cite[Theorem A]{MS2004}, over a dense subset of second category 
$$
\cJ^{\tn{reg}}(X,\om,g,A) \subset  \cJ_\tau(X,\om)
$$ 
(such $J$ is called \textbf{regular} with respect to $g$ and $A$), $\cM_{g,k}^*(X,A)$ is a naturally oriented smooth manifold of the expected real dimension~(\ref{equ:exp-dim}).
\bDef{def:pos}
A closed $2n$-dimensional symplectic manifold $(X,\om)$  is called \textbf{semi-positive} if $c_1(A)\!\geq\! 0$ for all $A\!\in\! H_2(M)$ such that 
\bEqu{equ:pos}
\om(A)\!>\!0\quad\tn{and}\quad c_1(A)\!\geq\! 3\!-\!n.
\eEqu
It is called \textbf{positive}\footnote{Other terms such as {monotone}, {strongly semi-positive}, {convex}, etc. have also been used in the literature for this notion.} if $c_1(A)\!> \!0$ for all $A\!\in\! H_2(M)$ such that (\ref{equ:pos}) holds.
\eDef

\noindent
If $(X,\om)$ is positive and $J$ is regular with respect to all  $(g',A')$ such that $g'\!\leq\! g$ and $\om(A')\!\leq\! \om(A)$, then a basic dimension counting argument shows that the image of boundary strata
$$
\partial\ov\cM_{g,k}(X,A,J)= \ov\cM_{g,k}(X,A,J)\setminus \cM_{g,k}^*(X,A,J)
$$ 
under $\ev\!\times\! \st$ in $X^k\!\times\! \ov\cM_{g,k}$ is a set of at least real codimension $2$ and can be ignored; i.e. the inclusion $\cM_{g,k}^*(X,A)\!\subset\! \ov\cM_{g,k}(X,A)$ gives rise to a GW \textbf{pseudocycle}
$$
\GW^{X}_{g,A}\subset H_d(X^k\times \ov\cM_{g,k},\Q)
$$
independent of the choice of the regular $J$; see \cite[Theorem 6.6.1]{MS2004} for more details. This way, we can extend the definition of GW fundamental class to the case of positive manifolds. 

\noindent
Given a (sufficiently) smooth map $u\colon\!(\Si,\mfj)\!\lra\! (X,J)$, let 
\bEqu{equ:E01}
\cE_f = \Gamma(\Sigma,\Om_{\Si,\mfj}^{0,1}\otimes_\C u^*TX),\quad \tn{with}~~f\!=\!(u,(\Si,\mfj)),
\eEqu
be the space of (sufficiently) smooth $u^*TX$-valued $(0,1)$-forms with respect to $\mfj$ on $\Si$ and $J$ on $TX$. Given a \textbf{perturbation} term $\nu\!\in\! \cE_f$, we say $u$ is $(J,\nu)$-holomorphic if it satisfies the perturbed Cauchy-Riemann equation  
\bEqu{equ:Jnu-holo}
\dbar u= \nu.
\eEqu
Given a ``smooth'' family of perturbations over the Deligne-Mumford space (see \cite{RT1997} for the precise description) we define $\ov\cM_{g,k}(X,A,J,\nu)$ to be the space of equivalence classes of perturbed pseudoholomorphic maps in a similar way.  
If $(X,\om)$ is semi-positive, by \cite[Corollary 3.9]{RT1997}, the space $\ov\cM_{g,k}(X,A,J,\nu)$ is Hausdorff and compact with respect to a similarly defined Gromov convergence topology. By \cite[Theorem 3.16]{RT1997}, for generic $J$ and sufficiently small $\nu$, the main stratum $\cM_{g,k}^*(X,A,J,\nu)$ consisting of somewhere-injective maps is a smooth manifold of the expected dimension (\ref{equ:exp-dim}). 
By \cite[Theorem 3.11]{RT1997}, it has a canonical orientation. By \cite[Proposition 3.21]{RT1997}, the image of  $\partial\ov\cM_{g,k}(X,A,J,\nu)$ under $\ev\times \st$
is contained in images of maps from smooth even-dimensional manifolds of at least $2$ real dimension less than the main stratum. Thus, similar to the positive case, 
the inclusion 
$$
\cM_{g,k}^*(X,A,J,\nu)\!\subset\! \ov\cM_{g,k}(X,A,J,\nu)
$$ 
provides a GW pseudocycle 
independent of the choice of the admissible almost complex structure $J$ or the perturbation $\nu$. 

\noindent
In the positive/semi-positive situations, the resulting GW invariants are \textbf{enumerative}, in the sense that they can be interpreted as a finite $\Q$-weighted count of $(J,\nu)$-holomorphic maps of fixed degree and genus meeting some prescribed cycles at the marked points.

\subsection{Virtual Fundamental Class}\label{sec:VFC-intro}
Extending the definition of GW invariants to more general symplectic manifolds occupies the rest this manuscript. Let us briefly indicate how we shall proceed. Locally, the argument is analogous to the semi-positive case.
Roughly speaking, we consider local perturbations in a consistent way and patch the solution sets of the perturbed equations (\ref{equ:Jnu-holo}) to obtain a so called \textbf{Virtual Fundamental Class}, or \textbf{VFC} for short.

\noindent
To this end, we construct a \textbf{Kuranishi structure} on the moduli space. Abstractly speaking, a Kuranishi structure on a compact metrizable   space $M$ is made of \textbf{Kuranishi charts} on open sets of $M$ (similar to manifolds) which are ``compatible" along intersections. An abstract Kuranishi chart of dimension $d$ over an open subset $F\!\subset\! M$ (called the support or footprint of the chart) is made of  an orbifold vector bundle $\pr\colon\cU\!\lra\!\cV$, an orbifold section $s\colon\cV\!\lra \!\cU$, and a homeomorphism from the underlying topological space of the zero set of $s$ onto $F$, $\ov\psi\colon\ov{s}^{-1}(0)\!\lra\! F\subset M$, such that 
$$
\dim \cV\!-\! \tn{rank}~\cU\! =\!d.
$$ 
See Section~\ref{sec:prelim-orbi} for the basics of orbifolds and Section~\ref{sec:prelim-kur} for the precise definition of Kuranishi charts.

\noindent
One can put a Kuranishi structure on a fixed topological space $M$ in many different ways. 
However, there is a natural way of constructing Kuranishi charts over moduli spaces of stable maps. Given a (family of) $J$-holomorphic map 
$$
u\colon\!(\Si,\mfj)\!\lra\!(X,J),
$$
instead of considering a fixed perturbation $\nu \!\in \!\cE_f$ as in (\ref{equ:Jnu-holo}), we consider a finite dimensional $\aut(f)$-invariant subspace $E_f \!\subset \!\cE_f$ of perturbations and look at the space of all maps 
$$
V_f= \{u'\colon(\Si,\mfj)\lra(X,J)\colon~\exists~\nu \in E_f ~\tn{s.t}~\dbar u'= \nu~~\tn{and}~~\mfd(u,u')<\ep\};
$$
here $\ep\!>\!0$ is a small positive number and $\mfd(u,u')\!<\!\ep$ means $u'$ is $\ep$-close to $u$ with respect to some specified distance function $\mfd$.
For an ``appropriate" and sufficiently large choice of $E_f$ (see Definition~\ref{def:obs-space}), $V_f$ is $\aut(f)$-invariant, $[V_f/\aut(f)]$ is a smooth orbifold chart, 
$$
\cU\!=\![U_f\!=\!(V_f \times E_f)/\aut(f)] \lra \cV\!=\![V_u/\aut(f)]
$$ 
is an orbibundle chart, and the map 
$$
s\colon \cV_f\lra \cU_f,\quad s([u'])=[u',\dbar u']
$$
is an orbifold section. The zero set of this section is the set of actual $J$-holomorphic maps (on the fixed domain $(\Si,\mfj)$) sufficiently close to $u$. More generally, in Section~\ref{sec:main} we will consider the case where the complex structure $\mfj$ on $\Si$ can change and the domain has marked points. We extend this construction to nodal maps and cover the moduli space with Kuranishi charts of this specific type. In Section~\ref{sec:VFC}, we construct a VFC for every abstract oriented Kuranishi structure. In the case of natural Kuranishi structures on a moduli space of pseudoholomorphic maps, this gives us a VFC which we use to define GW invariants. The following theorem is the main result of this article in its simplest form; see Theorem~\ref{thm:Nat-Kur} for the precise statement.

 \begin{theorem}[{\cite[Theorem 1.3]{FO}}]\label{thm:VFC}
Let $(X^{2n},\om)$ be a closed symplectic manifold, $J\!\in\! \cJ(X,\om)$ be an arbitrary compatible\footnote{One may consider tame almost complex structures, instead.} almost complex structure, $A\!\in\! H_2(X,\Z)$, and $g,k\!\in\! \Z^{\geq 0}$. With the dimension $d$ as in (\ref{equ:exp-dim}), 
there exist a class of natural oriented cobordant $d$-dimensional Kuranishi structures on $\ov{\cM}_{g,k}(X,A,J)$ and an associated rational homology class $$
[\ov{\cM}_{g,k}(X,A,J)]^\vfc\in H_d(X^k\!\times\! \ov{\cM}_{g,k},\Q)
$$ 
which is independent of $J$ and which only depends on the isotopy class of $(X,\om)$.
\eThm

\noindent
The VFC provided by Theorem~\ref{thm:VFC} only allows us to define GW invariants with \textbf{primary insertions}; i.e. constraints coming from the ambient space $X$ or the Deligne-Mumford space $\ov{\cM}_{g,k}$. A more complete list of GW invariants include invariants involving $\psi$-classes, Hodge classes, and other similarly defined intrinsic classes. First, let us briefly recall the definition of these cohomology classes. 

\noindent
For every genus $g$ $k$-marked degree $A$ $J$-holomorphic map $f\equiv(u,\Si,\mfj,\vec{z})$ (nodal or not), let $\cL_i|_f$ be the holomorphic cotangent space of $(\Si,\mfj)$ at the $i$-th marked point, 
$$
\cL_i|_f=\mc{T}^*_{z^i}\Si_\mfj.
$$
Since every $z^i$ is fixed with respect to $\aut(f)$, $\aut(f)$ acts on $\cL_i|_f$. If $\ov\cM_{g,k}(X,A)$ is an orbifold, then the quotient spaces $[\cL_i|_f/\aut(f)]$ for different $[f]\!\in\! \ov\cM_{g,k}(X,A)$ form an orbifold complex line bundle $\cL_i$ over $\ov\cM_{g,k}(X,A)$. This line bundle is called the $i$-th \textbf{tautological} line bundle. In this situation, the $i$-th $\psi$-class, denoted by $\psi_i$, is the orbifold first Chern class of the $i$-th tautological line bundle; see \cite[Definition 25.2.1]{mirror}. In the general case, for every natural Kuranishi structure on $\ov\cM_{g,k}(X,A)$, these one-dimensional complex vector spaces form a ``Kuranishi complex line bundle" in the sense of Definition~\ref{def:Kur-bundle}; also see \cite[Section 25.2]{mirror}. In this case, $\psi_i$ is the ``virtual" first Chern class of $\cL_i$ against VFC of $\ov\cM_{g,k}(X,A)$; see the next paragraph. Similarly, for every $g\geq 1$, we can define a rank $g$ Kuranishi vector bundle $\cE_g$, called \textbf{Hodge bundle}, whose fiber at every $[f]\equiv[u,\Si,\mfj,\vec{z}]$ is equal to 
\bEqu{Hodge-fiber_e}
[H^0(\Sigma,\om_{\Si,\mfj})/\aut(f)],
\eEqu
where $\om_{\Si,\mfj}$ is the dualizing sheaf of $(\Si,\mfj)$; see \cite[Section 25.3.1]{mirror}. For smooth $\Si$, $\om_{\Si,\mfj}$ is simply the sheaf of antiholomorphic 1-forms $\Om^{0,1}_{(\Si,\mfj)}$. Then, the Hodge class $\la_g$ is the ``virtual" top Chern class of $\cE_g$ against VFC of $\ov\cM_{g,k}(X,A)$.

\noindent
In order to define a VFC involving the top Chern class of some Kuranishi vector bundle $\cL$, we change the original Kuranishi structure of Theorem~\ref{thm:VFC} into some augmented Kuranishi structure by increasing the obstruction bundle. For example, in the case of GW invariants,  at every $J$-holomorphic map $f$, we replace the obstruction bundle $E_f$ with $E_f\oplus L_f$, where $L_f$ is the fiber of $\cL$ (before taking quotient w.r.t $\aut(f)$) at $f$, and extend the Kuranishi map by zero on the second factor.  Then, we define the ``virtual" top Chern class of $\cL$ against VFC of $\ov\cM_{g,k}(X,A)$ to be the VFC of the augmented system. We can then pair the VFC of the augmented system with constraints from $X$ and $\ov{\cM}_{g,k}$ to define invariants involving a mixture of primary and non-primary insertions. See Section~\ref{sec:bundle-kuranishi} and Section~\ref{sec:GW-VFC} for more details. 

\bRem{rmk:family-version}
Given a compact smooth oriented submanifold $M_\cJ\!\subset\! \cJ(X,\om)$ of real dimension $d_\cJ$, let 
\bEqu{equ:famlyM}
\ov{\cM}_{g,k}(X,A,M_\cJ)= \bigcup_{J\in M_\cJ}\ov{\cM}_{g,k}(X,A,J)
\eEqu
be the family of moduli spaces over $M_\cJ$ (with the Gromov topology) and 
$$
\pi_\cJ\colon \ov{\cM}_{g,k}(X,A,M_\cJ)\lra M_\cJ
$$ 
be the projection map. 
Then a family version of Theorem~\ref{thm:VFC} holds: there exists a natural family of oriented $(d\!+\!d_\cJ)$-dimensional Kuranishi structures on (\ref{equ:famlyM}) and an associated relative homology class (obtained via $\ev\!\times\! \st\!\times\!\pi_\cJ$)
$$
[\ov{\cM}_{g,k}(X,A,M_\cJ)]^\vfc\in H_{d+d_\cJ}(X^k\!\times \!\ov{\cM}_{g,k}\!\times\! M_\cJ,X^k\!\times \!\ov{\cM}_{g,k}\!\times\! \partial M_\cJ,\Q)
$$ 
which only depends on the isotopy class of $(X,\om)$ and $M_\cJ$. For example, this family version of Theorem~\ref{thm:VFC}, with $M_\cJ\cong S^2$, has been used to define GW-type invariants of certain $S^1$-family of $\tn{K}3$-surfaces; see \cite{YZ}.
\eRem

\bRem{rmk:funt-significance}
It is important to mention that unlike in algebraic geometry, all we need in the symplectic approach to Gromov-Witten theory is a space on which we can do intersection theory. So the main functorial aspects of studying moduli spaces in algebraic geometry, e.g. the existence of universal families, are absent in this approach. The result of this flexibility is a spectrum of different analytical/topological approaches  \cite{LT, FO, MW3, Ho, Jo, CM} to formalizing such intersection theory in symplectic geometry. 
\eRem

\subsection{Outline}\label{sec:outline}

In Sections~\ref{sec:orbifold} and \ref{sec:orbibundle}, we will review some basic definitions and facts about orbifolds, orbibundles, and sections of orbibundles. We will mostly focus on effective orbifolds, however, at certain points throughout the paper, we will consider some basic non-effective orbifolds as well. In Section~\ref{sec:multisection}, we introduce the notion of a multisection.  A multisection $\mfs$ of an orbibundle $\pr\colon \cE\lra \cM$ is locally given by a family of sections of an orbifold chart invariant under the group action.  As we will show in Section~\ref{sec:perturbation}, multisections give us a way of deforming possibly non-transversal orbifold sections into a compatible collection of transversal multisections. Otherwise, it is not always feasible to perturb an arbitrary orbifold section into a transversal one. In Section~\ref{sec:resolution}, we introduce the notion of a resolution of a multisection. Given a multisection $\mfs$, we decompose the base orbifold into a finite union of  open sets labeled by positive integers $\ell\!\in \N$, such that in the region with label $\ell$, $\mfs$ is given by a set of $\ell$ distinct branches (possibly with different multiplicities). Then, a resolution is a compatible collection of \textit{branched orbifold $\ell$-covering spaces} over the corresponding regions such that on each piece $\mfs$ is a weighted push-forward of an orbifold section $s_\ell$ of the pull back bundle.  In Section~\ref{sec:euler}, we use these resolutions to define the orbifold Euler class of the corresponding relatively oriented orbibundle. This construction lays out a foundation for the more complicated construction of VFC in Section~\ref{sec:VFC}. 

\noindent
In Sections~\ref{Introductory remarks} and \ref{sec:basics}, we introduce an abstract notion of a Kuranishi structure $\cK$ for a compact metrizable space $M$. It consists of a set of compatible Kuranishi charts $\cU_p$ indexed by $p\!\in\!M$.  Each $\cU_p$ is a local orbibundle with a section and a homeomorphism, called footprint map, from the zero set of the section onto an open subset $F_p$ of $M$. The open sets $\{F_p\}_{p\in M}$ cover $M$, and the Kuranishi charts are compatible on the overlaps in a suitable sense. In Section~\ref{sec:level}, we assemble this uncountable data into a finite set $\mfL$ of compatible orbibundles with sections and footprint maps  
$$
\aligned
&\pr_i\colon \cE(i)\lra \cY(i), \quad s_i\colon \cY(i)\lra \cE(i), \quad \tn{rank}~\cE(i)\!=\!i,\\
 &\psi\colon s_i^{-1}(0)\lra F(i)\!\subset \!M,\quad  \forall i\! \in\! \mfB\! \subset\! \N,
 \endaligned
$$
called  a \textbf{dimensionally graded system} (or \textbf{DGS}), such that the set of footprints $\{F(i)\}$ gives us a ``nice" covering of $M$. 
The compatibility conditions are given by orbibundle embeddings  $\cE(i)|_{\cY(i,j)}\!\lra \cE(j)$
from the restriction to open sub-orbifolds of the orbibundles  of the lower rank into the orbibundles of the higher rank. In Section~\ref{shrinking}, we introduce the notation of a shrinking of a dimensionally graded system. A shrinking allows us to replace orbifold pieces of an DGS with smaller orbifolds to accommodate certain perturbations and gluing-pasting arguments. Our construction of a VFC for $(M,\cK)$ in this article depends on many auxiliary choices, including a particular choice of a Kuranishi structure itself. For this reason, in Section~\ref{sec:cobordism}, we briefly review the notion of a cobordism between two Kuranishi structures, and more generally, between two sets of auxiliary data. In Section~\ref{sec:def}, we extend the perturbation results of Section~\ref{sec:perturbation} to dimensionally graded systems. Finally, in Section~\ref{sec:VFC}, we give two different constructions of a virtual fundamental class associated to a cobordism class of Kuranishi structures. In the simplest case of an DGS made of a single orbibundle, the resulting VFC is simply the Euler class of the orbibundle. We give three examples of this situation in Section~\ref{sec:GW-VFC}.

\noindent
In Section~\ref{sec:moduli}, we will review some well-known facts about the structure of the Deligne-Mumford space \cite{DM,RS} of stable curves and moduli spaces of stable pseudoholomorphic maps. 

\noindent
In Section~\ref{sec:main}, we show that the moduli spaces of stable pseudoholomorphic maps have a natural cobordism class of  Kuranishi structures.  In order to construct Kuranishi charts around different types of maps in such a moduli space, in Section~\ref{sec:canonical-smooth}, we start with the simplest case of pseudoholomorphic maps with smooth stable domain. In Section~\ref{sec:canonical-nodal}, we extend this construction to maps with stable nodal domain, and in Section~\ref{sec:natural}, we address the case of maps with unstable domain. The collection of Kuranishi charts constructed in these three sections, which we call \textbf{primary charts}, are not necessarily compatible along the overlaps. For this reason,  in Section~\ref{sec:natural-KUR}, we modify the construction and construct  a new collection of Kuranishi charts, which we call \textbf{induced charts}, satisfying all the requirements of a Kuranishi structure. A collection of induced charts, and thus the resulting Kuranishi structure, depends on a finite collection of primary charts. Different such collections result in different but cobordant Kuranishi structures.
Together with the results of Section~\ref{sec:VFC}, this gives the GW virtual fundamental class of Theorem~\ref{thm:VFC}.

\noindent 
The most difficult step in the construction of Kuranishi charts around maps with nodal domain is establishing the gluing argument of smoothing the nodes and its exponential decay properties. 
Due to the lack of space, we will state and use the gluing theorem as a black box 
and refer the interested reader to the dedicated article of FOOO \cite{FOOO-detail3}.
The main purpose of this article is to give a fairly detailed, yet efficient, description of Kuranishi structures and of the resulting VFC, and to show that how they can be used to define GW invariants.


\section{Preliminaries}\label{sec:prelim-orbi}
In \cite{Sa}, Satake introduced the notion of orbifold (which he called $V$-manifold) and orbibundle to generalize the notion of manifold and vector bundle, in order to include local quotients of smooth manifolds and vector bundles by finite group actions; see also the book of Thurston \cite{Th}. 
In this chapter, we briefly recall the definitions and prove necessary statements for the construction of VFC. 
We mainly consider effective orbifolds.  More general case of non-effective orbifolds can be better studied in the context of stacks or groupoids; see \cite{Fa} and \cite{Me} for an introduction to stacks, and \cite[Section 1.4]{ALR}, \cite[Section 2]{Mc2006}, \cite{Mo}, and the works of Haefliger \cite{Ha1971,Ha1984,Ha2001} to learn about groupoids.

\subsection{Orbifolds}\label{sec:orbifold}
\bDef{def:local-orbi}
An \textbf{effective orbifold chart} of real dimension $n$ on a  Haussdorff topological space $M$ is a tuple $\cV\!=\!(V,G,\psi)$ such that $V$ is a connected open manifold of real dimension $n$, $G$ is a finite group acting effectively and smoothly on $V$, and $\psi\colon V\! \lra\! F\! \subset\! M$  is a continuous map of the form $\psi\!=\ov\!\psi\!\circ\! \pi$, where $\pi\colon V\!\lra\! V/G$ is the quotient map and $\ov\psi\colon V/G\! \lra\! F$ is a homeomorphism (with respect to the quotient topology on $V/G$) onto an open subset of $M$. 
\eDef

\noindent
Given an effective orbifold chart $(V,G,\psi)$, for every $x\in V$ let 
$$
G_x= \{g\in G\colon g\cdot x = x\}
$$
be the \textbf{isotropy} group at $x$.

\bDef{def:refinement}
Given an $n$-dimensional effective orbifold chart $\cV\!=\!(V,G,\psi)$ on $M$, a \textbf{refinement} of $\cV$ is another  $n$-dimensional effective orbifold chart $\cV'\!=\!(V',G',\psi')$ together with a group homomorphism $h\!\colon\!G'\!\lra\! G$ and an $h$-equivariant smooth embedding $\phi \colon\!V'\!\lra\! V$ such that $\psi'\!=\!\psi\!\circ\! \phi$.
Two charts $\cV_1$ and $\cV_2$ on $M$ are said to be \textbf{compatible} if they  admit a common refinement around every point of $\psi_1(V_1)\cap \psi_2(V_2)$.
\eDef
\noindent
For the group homomorphism of refinements, it follows from the effectiveness of the group action that $h$ is an injective homomorphism.

\noindent
More generally, if $\ov{\phi}\colon\! M' \!\lra\! M$ is a continuous map between two topological spaces, given two effective orbifold charts $(V',G',\psi')$ and $(V,G,\psi)$ on $M'$ and $M$, respectively, a \textbf{smooth lift} 
\bEqu{equ:orbi-lift}
(\phi,h)\colon (V',G',\psi')\lra (V,G,\psi)
\eEqu
of $\ov{\phi}|_{\psi'(V')}\!\colon\!  \psi'(V')\!\lra\! \psi(V)$  consists of a group homomorphism $h\colon\! G' \lra G$  and an $h$-equivariant smooth map $\phi \colon\! V' \lra V$ such that 
\bEqu{equ:pi2pi1}
\psi\circ \phi =\ov{\phi} \circ \psi'.
\eEqu
In this situation, for every $x\!\in\! V'$ we get an induced homomorphism 
\bEqu{equ:stalk-hom}
h_x\colon G'_x\lra G_{\phi(x)}.
\eEqu
With this definition, a refinement is a smooth lift of the inclusion map over an open embedding $F'\!\subset\! F\! \subset\! M$.
The group homomorphism of a general smooth lift may not be injective or surjective.
In the case of refinements, every two such smooth embedding lifts $(\phi_i,h_i)\colon\!\cV'\!\lra\! \cV$, with $i\!=\!1,2$, differ by an element of $G$, i.e. there exists $g\!\in\! G$ such that $\phi_2\!=\!g \circ \phi_1$. If $g \big(\tn{Im}(\phi_1))\cap \tn{Im}(\phi_2)\big)\!\neq \!\emptyset$, then $g\! \in\! \tn{Im}(h)$;
therefore, in this case, (\ref{equ:stalk-hom}) is an isomorphism and the conjugacy class of the isotropy group for different lifts of $\ov{x}\!=\!\psi(x)\in M$ is independent of the particular choice. We denote by $\tn{I}_{\ov{x}}$ to be the conjugacy class of the isotropy groups over different lifts of $\ov{x}\in M$ (w.r.t. the orbifold structure $\cM$) and define
\bEqu{equ:isotropy-order}
\mf{I}_{\cM}\colon M\lra \Z, \quad \mf{I}_{\cM}(\ov{x})=|\tn{I}_{\ov{x}}|.
\eEqu
We use the function $\mf{I}_{\cM}$ in Section~\ref{sec:resolution} and Section~\ref{sec:euler} in the Definition of Euler class of an orbibundle.

\bDef{def:orbi-atlas}
Given a metrizable topological space $M$, an $n$-dimensional smooth  \textbf{ effective orbifold atlas} on $M$ is a countable collection 
\bEqu{equ:atlas}
\cA\equiv\{\cV_\al=(V_\al,G_\al,\psi_\al)\}_{\al\in S}
\eEqu
of compatible $n$-dimensional orbifold charts covering $M$. Two effective orbifold atlases $\cA$ and $\cA'$ are said to be \textbf{equivalent} if they admit a common refinement\footnote{An atlas $\cA'$ is a refinement of another atlas $\cA$, if every chart of $\cA'$ is the refinement of some chart of $\cA$.}. An \textbf{effective orbifold structure} $\cM$ on $M$ is an equivalence class of such atlases. 
\eDef

\noindent
It actually requires a little bit of work to show that ``refinement'' defines an equivalence relation on the set of orbifold atlases. We refer to \cite[Section 1]{ALR} for a more detailed discussion of orbifold atlases, their refinements, and the refinement equivalence relation between atlases. By replacing diffeomorphisms with homeomorphisms in the defining equation of the orbifold charts and the lift maps, we can define topological orbifolds, similarly. By Definition~\ref{def:local-orbi}, the orbifolds considered in this article, unless specifically mentioned (e.g. in Section~\ref{sec:cobordism}), do not have boundary. 

\begin{notation}\label{not:calconv}
In the rest of this article, we let calligraphic letters $\cM$, $\cY$, $\cZ$, and $\cW$ to denote for orbifold structures on topological spaces (or manifolds) $M$, $Y$, $Z$, and $W$, respectively. We also use $\cM$ for moduli spaces (as they are orbifolds in some nice cases). Similarly, after the introduction of orbibundles in Section~\ref{sec:orbibundle}, we let  $\cE$ denote an orbibundle with underlying topological space $E$. 
\end{notation}

\noindent
\bRem{rmk:effectivity}
\emph{Unless explicitly mentioned, e.g. in comparison of our results and framework with similar works and in the construction of Euler class, the orbifolds considered in this manuscript are effective. Therefore, we omit the word ``effective'' from our terminology for simplicity. We refer to \cite[Section 4]{FOOO2013} for some discussion on the complexity of non-effective orbifold quotients.}
\eRem

\noindent
The simplest examples of orbifolds are global quotients of manifolds by discrete groups (which indeed can be infinite).
If $M$ is a smooth manifold and $G$ is a discrete group acting effectively and properly discontinuously on
$M$, then $M/G$ has the structure of an orbifold which we denote by $[M/G]$. 
More generally, if $G$ is a compact Lie group acting smoothly, effectively, and almost freely (i.e.~with finite stabilizers) on a smooth manifold $M$, the quotient space $M/G$ has the structure of an orbifold. For example, all the weighted projective spaces are orbifolds of this sort with $M\!=\!\C^N\!-\!\{0\}$ and $G\!=\!\C^*$. In the terminology of \cite[Definition 1.8]{ALR}, the former is called an \textbf{effective global quotient} and the latter is called an \textbf{effective quotient}.

\noindent
As another example, if $\Si$ is a Riemann surface of genus $g$, hence a manifold, by removing finitely many disjoint open discs $\{D_i\}_{i=1}^k$ from $\Si$ and gluing back copies of $[D_i/\Z_{m_i}]$, where $\Z_{m_i}\!=\! \Z/(m_i \Z)$ acts by multiplication with $m_i$-th roots of unity, we obtain new orbifold structures on the underlying space $\Si$. Thurston \cite[Theorem 13.3.6 and p. 304]{Th} shows that the result is an effective global quotient unless $g\!=\!0$ and $k\!=\!1$ or $k\!=\!2$ with $m_1\!\neq\! m_2$. 
Examples of this form naturally arise in the study of moduli spaces; see Section~\ref{subsec:elliptic}.
On the other hand, by \cite[Theorem 1.23]{ALR}, every (effective) orbifold is an effective quotient of its frame bundle. 

\noindent
As the previous example shows, a fixed topological space can carry many different orbifold structures. 
It is possible to adopt most of the invariants of manifolds, e.g. fundamental group and Euler characteristic, to orbifolds and these invariants are usually different from the corresponding invariants  of the underlying topological space. 

\bExa{exa:M11}
Among various examples of orbifolds, we are particularly interested in those arising from moduli problems. For example, the moduli space $\cM_{1,1}$ of isomorphism classes of smooth elliptic curves with $1$ marked point is isomorphic to $[\H/ \tn{SL}_2(\Z)]$, where
$$ 
\H=\{ \tau\in \C\colon \tn{im}(\tau) >0\},\quad A\tau =\frac{a\tau+b}{c\tau + d}\quad \forall A=\left[ \begin{array}{cc}
a& b \\
c& d  \end{array} \right]\!\in\!\tn{SL}_2(\Z).
$$
Since the action of $-\tn{Id}$ is trivial, $\cM_{1,1}$ is not an effective orbifold. In order to get an effective orbifold, we need to reduce the group action to 
$$
\tn{PSL}_2(\Z)\!=\!\tn{SL}_2(\Z)/\ll -\tn{Id} \rr.
$$ 
Then, the reduced orbifold structure
$$
\cM_{1,1}^{\tn{red}}\equiv [\H/ \tn{PSL}_2(\Z)]
$$ 
is the parameter  space of isomorphism classes of smooth elliptic curves with $1$ marked point and $\cM_{1,1}$ is some (stacky) $\Z_2$-quotient of that.
An alternative way to approach the issue is to increase the number of marked points. For example, the moduli space $\cM_{1,2}$ of isomorphism classes of smooth elliptic curves with $2$ marked points is an open subset of the effective global quotient 
$$
[(\H \times \C)/(\tn{SL}_2(\Z) \ltimes  \Z^2)]
$$ 
where
$$
(\tau,z)\sim \bigg( A\tau, (c\tau + d)^{-1}(z+m+\tau n)\bigg)\quad\forall 
A\!=\!\left[ \begin{array}{cc}
a& b \\
c& d  \end{array} \right]\!\in\! \tn{SL}_2(\Z),~(m,n)\!\in\! \Z^2.
$$
In our construction of VFC, we often use the trick of adding more marked points to avoid encountering stacks and more complicated groupoids.
\eExa

\noindent
\bDef{def:orbi-oriented}
Given an orbifold atlas $\cA$
on $M$ as in (\ref{equ:atlas}), we say $\cA$ is \textbf{oriented} if each $V_\al$ is oriented, the action of $G_\al$ is orientation-preserving, and the orientations agree on the overlaps. An \textbf{oriented} orbifold structure on $M$ is a maximal oriented atlas.
\eDef 

\bDef{def:orbi-maps}
An \textbf{orbifold smooth map}, $\phi\colon\cM_1\!\lra\!\cM_2$, between two orbifolds $\cM_1$ and $\cM_2$ is an underlying continuous map $\ov{\phi}\colon\!M_1\lra M_2$ that admits local smooth lifts as in (\ref{equ:orbi-lift}). 
An \textbf{orbifold morphism} is an orbifold smooth map for which the induced homomorphisms $h_x\colon G_x\!\lra\! G_{\phi(x)}$ as in (\ref{equ:stalk-hom}) are isomorphism.  
An orbifold \textbf{diffeomorphism} is an orbifold morphism with an inverse.
An \textbf{orbifold embedding} is an orbifold morphism such that the lift maps (\ref{equ:orbi-lift}) are smooth embeddings. Finally, a \textbf{sub-orbifold} $\cM'$ of an orbifold structure $\cM$ on $M$ is a subspace $M'\!\subset\! M$ together with an orbifold structure $\cM'$ on $M'$ such that the inclusion map $\ov\iota\colon M'\lra M$ gives an orbifold embedding $\iota\colon \cM'\lra \cM$.
\eDef

\noindent
The isomorphism condition on (\ref{equ:stalk-hom}) for orbifold embeddings is not standard across the literature but is the one we use in our construction of Kuranishi structures in this article. 

\begin{notation}\label{not:phi}
In Definition~\ref{def:orbi-maps}, the local group homomorphisms $h$ of (\ref{equ:orbi-lift}) are embedded in the notation $\phi$, i.e. in Definition~\ref{def:orbi-maps} $\phi$ denotes a global map between orbifolds and at the level of charts the map $\phi$ also comes with a group homomorphism $h$ as in (\ref{equ:orbi-lift}).
\end{notation}

\bExa{exa:morphism-vs-map}
Let $\cM\!=\![\R^2/\Z_2]$ where $\Z_2$ acts by $(x,y)\!\lra\! (-x,y)$. 
Then the map $\phi\colon \R\!\lra\! \cM$, given by $\phi(y)\!=\!(0,y)$, is a smooth orbifold map which is \textit{not} an orbifold embedding. 
\eExa

\bRem{rem:final-remarks}
We conclude this section with some final remarks on the structure and properties of orbifolds. 
In what follows, $\cM$ is a closed (compact without boundary) orbifold .
\bEnum
\item\label{l:singM} If $\cM$ is oriented, then the singular locus 
$$
M^{\tn{sing}}\!\equiv \!\{\ov{x}\in M\colon \tn{I}_{\ov{x}} \neq 1\}
$$ 
of $\cM$ has real codimension greater than or equal to $2$. Let $M^{\tn{sm}}$ be the complement of $M^{\tn{sing}}$ in $M$; $M^{\tn{sm}}$ has the structure of an oriented smooth  manifold. Then, the inclusion map $M^{\tn{sm}}\!\lra\! M$ provides a pseudocycle, and thus a fundamental class for $M$; see \cite{Zi} for the definition and properties of pseduocycles. In the case of semi-positive symplectic manifolds in Section~\ref{sec:positive}, we use this idea to define GW fundamental class of $\ov\cM_{g,k}(X,A)$.

\item\label{l:OrbiStrata} More generally, given a group $H$, let $M^H\!\subset\! M$ be the set of points where the isotropy group is isomorphic to $H$. For example, for the trivial group $\ll1\rr$, $M^{\ll 1\rr}\!=\! M^{\tn{sm}}$. By the argument in \cite[Section 5]{FOOO2013}, every $M^{H}$ has the structure of a smooth manifold (induced from the orbifold structure on $\cM$). Let $\cH$ be the set of all $H$ such that $M^{H}\!\neq\! \emptyset$; if $M$ is compact, then $\cH$ is finite and $M\!=\! \bigcup_{H\in \cH} M^{H}$ gives a stratification of $M$ into smooth manifolds. 

\item For every point $\ov{x}\!\in\! M$, there exists an orbifold chart $\cV\!=\!(V,G,\psi)$ \textbf{centered at} $\ov{x}$ in the sense that $\psi^{-1}(\ov{x})\!=\!\{x\}$ and $G\!=\!G_x$. 

\item An orbifold chart $(V,G,\psi)$ is said to be \textbf{linear}, whenever $V\!\subset\! \R^n$ is a connected open set and the (restriction to $V$ of the) action of $G$ is linear. As 
\cite[Page 2]{ALR} indicates, since smooth actions are locally linear, any orbifold has an atlas consisting of linear charts. Moreover, we can choose the linear charts so that $G\!\subset\! O(n)$, and if $\cM$ is orientable, we can choose them so that $G\!\subset\! SO(n)$. This is the point of view in the original approach of \cite{FO}; all the local structures taken there are linear. In the case of orbifold charts centered at $\ov{x}$, we can choose the linear chart so that the unique preimage $x$ of $\ov{x}$ is the origin.

\item\label{l:groupoid} As we mentioned earlier, orbifolds can be approached  via the frame work of \textbf{Lie groupoids}. A Lie groupoid is a category in which all morphisms are invertible (i.e. they are isomorphisms), the space of objects $\wt{M}$ and the space of morphisms $R$ (relations) are smooth manifolds (possibly with infinitely many, but countable, connected components), and the following set of structure maps are all smooth. For $x,y\!\in\! \wt{M}$, let $R(x,y)$ be the set of morphisms from $x$ to $y$. Since all morphisms are invertible, there is an \textbf{inversion} involution $\iota\colon R\!\lra \!R$ which takes $R(x,y)$ to $R(y,x)$. From this set-up, we get the \textbf{source} and \textbf{target} maps 
$$
s\colon R\lra \wt{M}\times \wt{M}~~~\tn{and}~~~t\colon R\lra \wt{M}\times \wt{M}
$$ 
which are defined on every element $\gamma\!\in\! R(x,y)$ by $s(\gamma)\!=\!x$ and $t(\gamma)\!=\!y$, respectively. Note that $s\circ \iota\! =\!t$. 
By composing morphisms, whenever possible, we  get a multiplication map 
$$
m \colon R{}_t\!\times_{s}\!R \lra R, \qquad m(\gamma_1,\gamma_2)\!=\!\gamma_2\circ \gamma_1\quad \forall
\gamma_1\times \gamma_2\!\in\! R(x,y)\!\times\! R(y,z).
$$ 
The \textbf{orbit space} of $(\wt{M},R)$, denoted by $M\!\equiv\!\wt{M}/R$, is the quotient space of equivalence classes of objects in $\wt{M}$ with respect to relations in $R$.

\noindent
A \textbf{Lie groupoid orbifold} is a lie groupoid where 
$$
s\!\times\! t\colon \!R\!\to\! \wt{M}\!\times\! \wt{M}
$$ 
is proper and $s,t\colon R\!\lra\! \wt{M}$ are local diffeomorphisms; see \cite[Example 1.33]{ALR} for more details.

\noindent
Let $\cA$ be an orbifold atlas on $M$ as in (\ref{equ:atlas}). Associated to $\cA$, there is a Lie groupoid orbifold with orbit space $M$ given by
$$
\wt{M}=\coprod_{\al \in S} V_\al,\quad R= \coprod_{\al,\beta \in S} \Gamma_{\al\beta},
$$
where $\Gamma_{\al\beta}$ is the disjoint union of the irreducible components of the fiber product $V_\al\times_{M}V_\beta$.
Therefore, all (effective) orbifolds are special cases of Lie groupoid orbifolds but the latter includes non-effective orbifolds as well. The most relevant example to us is the Deligne-Mumford space of stable marked curves which is worked out in details in \cite{RS}. By \cite[Theorem 6.5]{RS}, for every $g,k\!\in\! \Z^{\geq 0}$ such that $2g\!+\!k\!\geq\! 3$, $\ov\cM_{g,k}$ has the structure of a Lie groupoid orbifold. If $(g,k)\!\neq\! (1,1)$ or $(2,0)$, then the Lie groupoid orbifold structure comes from an effective orbifold as above, but in those two special cases, one either has to mildly generalize the notion of orbifold or to use the language of Lie groupoid orbifold.

\eEnum
\eRem

\subsection{Orbibundles}\label{sec:orbibundle}
Orbibundles are the analogue of vector bundles for orbifolds. They show up as obstruction bundles, tangent bundle of orbifolds, and normal  bundle of orbifold embeddings in the construction of VFC. 

\bDef{def:orbibundlechart}
An \textbf{orbibundle chart} of rank $r$ over an $n$-dimensional orbifold chart $\cV\!=\!(V,G,\psi)$ is a triple 
\bEqu{equ:orbi-bundle-chart}
\cU\!=\!(\pr\colon \!U\!\lra\! V,G,\psi)
\eEqu
where $\pr\colon \!U\!\lra\! V$ is a smooth $G$-equivariant\footnote{i.e. the $G$-action on $U$ and $V$ commutes with the projection map, and for every 
$g\!\in\! G$ and $x\!\in\! V$, $g\colon U_x \!\lra\! U_{g\cdot x}$ is linear.} vector bundle of rank $r$ over $V$.
\eDef

\noindent
Similarly, given an orbibundle chart (\ref{equ:orbi-bundle-chart}), a \textbf{refinement} of $\cU$ is another orbibundle chart $\cU'\!=\!(\pr'\colon \!U'\!\lra\! V',G',\psi')$ of the same rank over a refinement $\cV'\!=\!(V',G',\psi')$ of $\cV$ which is defined via an embedding $(\phi,h)$,
and an $h$-equivariant fiber-wise linear smooth embedding $\mfD\phi\colon U'\lra U$ such that the following diagram commutes,
\bEqu{bundle-lifts}
\xymatrixcolsep{2pc}\xymatrix{
U' \ar[r]^{\mfD\phi} \ar[d]^{\pr'} & U \ar[d]^{\pr} & \\
V' \ar[r]^{\phi}  & V~. &
  }
\eEqu

\begin{notation}\label{fD_n}
We use the notation $\mfD \phi$ for the map between vector bundles of orbibundle charts to emphasize on its linearity and compatibility with $\phi$. In the example of tangent bundle below and few other instances, $\mfD \phi$ is actually the derivative $\phi$, but otherwise, the letter $\mfD$ behind $\phi$ is just a  notation.
\end{notation}

\noindent
Two orbibundle chart enhancements $\cU_1$ and $\cU_2$ of the orbifold charts $\cV_1$ and $\cV_2$ are said to be \textbf{compatible}, whenever they  admit a common refinement around every point of $\psi_1(V_1)\cap \psi_2(V_2)$. 
 
 \bDef{def:orbibundle}
 Given an orbifold atlas $\cA_\cM$ on $M$ as in (\ref{equ:atlas}), an \textbf{orbibundle atlas} of rank $r$ over $\cA_\cM$ is a compatible enhancement of the orbifold charts of $\cA_\cM$ to orbibundle charts of rank $r$,
 \bEqu{equ:atlasE}
 \cA_\cE\equiv\{\cU_\al\!=\!(\pr_\al\colon \!U_\al\!\lra\! V_\al,G_\al,\psi_\al)\}_{\al\in S}.
 \eEqu
 We say two orbibundle atlases $\cA_\cE$ and $\cA_\cE'$ over $\cM$ are \textbf{equivalent} if they admit a common refinement. An \textbf{orbibundle} $\cE$ on $\cM$ is an equivalence class of such atlases. 
\eDef

\noindent
Equivalently, we can think of the total space of an orbibundle as an orbifold, itself, which admits certain type of orbifold projection map.

\bDef{def:orbibundle2}(alternative definition)
Let $\cE$ and $\cM$ be two orbifolds, and 
\bEqu{equ:orbi-projection-map}
\pr\colon \cE\lra \cM
\eEqu
 be a surjective orbifold smooth map\footnote{Surjective means the underlying continuous map $\ov{\pr}\colon E\lra M$ is surjective.}. We say (\ref{equ:orbi-projection-map}) is an \tn{orbibundle} of rank $r$, if for every sufficiently small orbifold chart $\cV\!=\!(V,G,\psi_V)$ of $\cM$  there is an orbifold chart $\cU\!=\!(U,G,\psi_U)$ over $\ov{\pr}^{-1}(\psi_V(V))\!\subset\! E$ such that 
 $\pr\colon U\!\lra\! V$ is a smooth ($G$-equivariant by definition) vector bundle of rank $r$ over $V$ and for every $g\!\in\! G$ and $x\!\in\! V$, $g\colon U_x \lra U_{g\cdot x}$ is linear.
\eDef
\bExe{two-dfn}
Prove that Definition~\ref{def:orbibundle} and Definition~\ref{def:orbibundle2} are equivalent.
\eExe

\noindent
\bRem{rmk:OBL}
Similar to Remark~\ref{rem:final-remarks}.(4), for every $\ov{x}\!\in\! M$, we can take a sufficiently small orbibundle chart  $\cU$ \textbf{centered at} $\ov{x}$  which is linear and trivial in the sense that
$$
V\!\subset\! \R^n,~~x\!=\!\psi^{-1}(\ov{x})\!=\!0\!\in\! V,~~ U\!=\!V\!\times\! \R^{r}, 
$$
the group $G$ acts linearly on both $V$ and $\R^r$, and $\pr$ is the projection map onto the first component.
This is the point of view in the original approach of \cite{FO}.
\eRem

\bExa{exa:OrbiTangent}
Let $\cM$ be a smooth orbifold of dimension $n$, then its tangent bundle is a well defined orbibundle of rank $n$. 
Given a smooth orbifold chart $\cV\!=\!(V,G,\psi)$ for $\cM$, the action of $G$ on $V$ lifts to $U\!=\!TV$ by 
$$
\rho \colon G\times U \lra U,\quad \rho(g,(x,v_x))=(g\cdot x, (\nd g)_x (v_x)).
$$
\eExa

\bRem{rem:pull-back}
In general, if $\phi\colon \cM'\!\lra\! \cM$ is an orbifold smooth map and $\pr\colon\cE\!\lra\! \cM$ is an orbibundle, then there is no well-defined notion of  pull-back orbibundle similar to the manifold case. 
For example, with $\cM$ and $\phi$ as in Example~\ref{exa:morphism-vs-map} and $\cE=T\cM$, the fiber of underlying space $E$ over $[(0,y)]$ is a $\Z_2$-quotient of $\R^2$, but every orbibundle over $\R$ is a vector bundle. Nevertheless, the pull-back bundle is well-defined whenever the homomorphisms $h_x$ that come with $\phi$ as in Definition~\ref{def:orbi-maps} are  monomorphism; see \cite[Theorem 2.43]{ALR}.
We will use this fact in Section~\ref{sec:resolution}.
\eRem

\bExa{eax:OrbiNormal}
Let $\phi\colon\! \cM'\lra \!\cM$ be a smooth embedding of orbifolds as in Definition~\ref{def:orbibundle}. 
Then there exists a well-defined normal orbibundle 
$$
\cN_{\cM'}\cM\cong \frac{T\cM|_{M'}}{T\cM'} \lra \cM'
$$
generalizing the notion of normal bundle for smooth embedding of manifolds.  
\eExa

\noindent
\bDef{def:orbibundle-maps}
An \textbf{orbibundle smooth map}, $(\mfD\phi,\phi)\colon(\cE_1,\cM_1)\!\lra\!(\cE_2,\cM_2)$, between two orbibundles $\cE_1$ and $\cE_2$ over $\cM_1$ and $\cM_2$, respectively, is given by an underlying pair of continuous maps $(\ov{\mfD\phi},\ov{\phi})\colon\!(E_1,M_1)\lra(E_2,M_2)$ that commute with the projection maps $\ov{\pr}_1$ and $\ov{\pr}_2$, admit local smooth lifts as in (\ref{bundle-lifts}), and $\mfD\phi$ is linear on each fiber.
An \textbf{orbibundle morphism} is an orbibundle smooth map for which the underlying orbifold smooth map $\phi\colon\cM_1\lra \cM_2$ is a morphism. An \textbf{orbibundle embedding} is an orbibundle morphism such that the lift maps (\ref{bundle-lifts}) are smooth embeddings. Finally, a \textbf{sub-orbibundle} $\cE'$ of an orbibundle $\cE$ over $\cM$ is a subspace $(E',M')\!\subset\! (E,M)$ together with an orbibundle structure $\cE'\!\lra\!\cM'$ on $(E',M')$  such that the inclusion map $\ov\iota\colon\! (E',M')\!\lra\!(E,M)$ gives an orbibundle embedding $\iota\colon\! (\cE',\cM')\!\lra\!(\cE,\cM)$.
\eDef
\noindent 
We can similarly extend the notions of quotient bundle, tensor product, pull-back bundle, etc. to the context of orbibundles. 

\noindent
From the perspective of Definition~\ref{def:orbibundle}, an \textbf{orbifold section} $s$ of an orbibundle $\pr\colon\!\cE\!\lra\! \cM$ is an orbifold map $s\colon\! \cM\!\lra\! \cE$ such that $\pr\circ s\!=\!\tn{id}$. The simplest orbifold section is the zero section. It is clear from the definition of orbibundle that the zero section is an orbifold embedding whose normal bundle is canonically isomorphic to $\cE$ itself. We say that a section is \textbf{transverse} if all local lifts are transverse to the zero section. 

\bExa{exa:Pmn}
For $m,n\!>\!0$, with $\gcd(m,n)\!=\!1$, consider $M\!=\!S^2$  and let  $\cM\!=\!\P^1_{m,n}$ denote the orbifold structure of type $(m,n)$ obtained by replacing neighborhoods of $0$ and $\infty$ in  $\P^1\!\cong\! S^2$ with $\Z_m$ and $\Z_n$ quotients of a disk, respectively. 
We can cover $\P^1_{m,n}$ with two orbifold charts $\cV_0$ and $\cV_\infty$ where $V_0\!\cong\! V_\infty\! \cong \!\C$, $\Z_m$ and $\Z_n$ act by the corresponding roots of unity, $\psi_0(z)\!=\!z^m$, and $\psi_\infty(w)\!=\!w^n$ ($w^{-n}\!=\!z^{m}$ on the overlap). 
Let $\cE=T\cM$. Then it is easy to see that 
$$
s_0\colon V_0\lra TV_0,~ s_0(z)=z\frac{\partial}{\partial{z}},~~~\tn{and}~~~s_\infty\colon V_\infty\lra TV_\infty,~ s_\infty(w)=\frac{-m}{n}w\frac{\partial}{\partial{w}},
$$
define a holomorphic orbifold section of the orbibundle $\cE$ which is transverse to the zero section. 
\eExa

\noindent
For an orbibundle $\cE\!\lra\! \cM$, let 
$$
\tn{det}~\cE \lra \cM
$$
be the determinant (top exterior power) orbi-line bundle of $\cE$.

\bDef{def:rel-oriented}
We say an orbibundle $\cE\lra \cM$ is \textbf{relatively orientable} if 
\bEqu{equ:relative-trivialization}
\tn{det}~\cE\otimes \tn{det}~T\cM\! \cong \!\cM\!\times \!\R.
\eEqu
A relative orientation is a choice of trivialization in (\ref{equ:relative-trivialization}). 
\eDef
\noindent
If $s$ is a transverse section of a relatively oriented orbibundle $\cE\lra\cM$, then 
$$
Z(s)=\ov{s}^{-1}(0)\subset M
$$
inherits a possibly non-effective but oriented orbifold structure $\cZ(s)$ from the pair $(\cE,\cM)$. In the case of Example~\ref{exa:Pmn}, $\cZ(s)$ is a disjoint union of the two points $0$ and $\infty$ with trivial (thus non-effective) actions of $\Z_m$ and $\Z_n$, respectively. For simplicity, we will allow these basic non-effective orbifold structures in Section~\ref{sec:euler}.

\subsection{Multisections}\label{sec:multisection}
In this section, we introduce the notion of multisection. We are forced to consider multisections in order to achieve transversality. 
Let $s\colon \cM\!\lra\! \cE$ be an orbifold section, $\cU$ as in (\ref{equ:orbi-bundle-chart}) be an orbibundle chart of $\cE$, $x\!\in\! V$ such that $G_x\!\neq\! 1$, and assume $G_x$ acts non-trivially on $U_x$. Then $s(x)\!\in\! U_x^{G_x}$ and if this invariant subspace is zero, every section has to vanish at that point. Therefore, it is quite likely (e.g. if the set of such points has larger dimension than $\dim V\!-\! \tn{rank}~U$) that $s$ does not admit any transversal deformation.

\noindent
This leads us to consider multisections, i.e. a $G$-invariant set of possibly non $G$-equivariant sections in place of a single equivariant section. Given a vector bundle $\pr\colon U \!\lra\! V$, we denote by $S^\ell(U)$ to be the quotient of $U^{\times \ell}$  
($\ell$-times fiber-product of $U$ over $V$) with respect to the action of the symmetric group of order $\ell$,
$$
v^1\times \cdots \times v^\ell \sim v^{\si(1)} \times \cdots \times v^{\si(\ell)}
\qquad \forall x\in V,~ v^1,\ldots,v^\ell \in U_x,~\tn{and}~\si \in S_\ell.
$$
The projection map $\pr$ extends to a similarly denoted map $\pr\colon S^\ell(U)\! \lra\! V$; however, the fibers are -usually- singular. 
An $\ell$-\textbf{section} is then a (continuous) section of $S^\ell(U)$. 

\noindent
Similarly, given an orbibundle chart $\cU$ as in (\ref{equ:orbi-bundle-chart}), we can lift the $G$-action to a component-wise $G$-action on $U^{\times \ell}$ and it descends to a $G$-action on $S^\ell(U)$ commuting with $\pr$. 
Then, by definition, an \textbf{orbifold} $\ell$-\textbf{section} for $\cU$ is a $G$-equivariant (continuous) section of $S^\ell(U)$. 
If $G\!=\!\{ g_1, \ldots, g_\ell\}$, every arbitrary section $s\colon V\lra U$ can be naturally enhanced to an orbifold $\ell$-section $\mf{s}\!=\![G\cdot s]\in S^\ell(U)$ defined by
\bEqu{equ:nat-enhance}
\mf{s}=[g_1\cdot s, \ldots, g_\ell\cdot s],\quad (g_i\cdot s) (x)= g_i (s(g_i^{-1}(x)))~~\forall i\!\in\![\ell],~x\!\in\! V.
\eEqu

\noindent
Every $\ell$-section can be canonically extended  to orbifold $\ell m$-sections, with $m\!\in\!\N$: we replace $[v^1, \ldots , v^\ell]$ with
$$ 
m\cdot [v^1, \ldots , v^\ell]=[\underbrace{v^1, \ldots , v^\ell~,\ldots,~v^1, \ldots , v^\ell}_\text{\tn{repeated}~m~\tn{times}} ].
$$
We say two orbifold multisections $\mfs_1 \!\in\! S^{\ell_1}(U)$ and $\mfs_2\!\in\! S^{\ell_2}(U)$ are \textbf{equivalent}, and write $\mfs_1\!\cong\! \mfs_2$, if there exists some common multiple $\ell$ of $\ell_1$ and $\ell_2$ such that
 \bEqu{equ:common-mult}
 \frac{\ell}{\ell_2}\cdot  \mfs_1= \frac{\ell}{\ell_1}\cdot \mfs_2;
 \eEqu 
 the equality means that the un-ordered set of of values of $\frac{\ell}{\ell_2}\cdot  \mfs_1$ and $\frac{\ell}{\ell_1}\cdot \mfs_2$ at every point $x\!\in\! V$ are the same in $U_x^\ell$.
The map
\bEqu{equ:multi-sum}
S^{\ell_1}(U)\times S^{\ell_2}(U)\lra S^{\ell_1\ell_2}(U),\qquad \big([v^i]_{i\in[\ell_1]}, [u^j]_{j\in[\ell_2]}\big) \lra [v^i+u^j]_{\substack{i\in[\ell_1]\\~j\in[\ell_2]}}
\eEqu
extends the definition of addition to the equivalence classes of orbifold multisections. The extended sum operation is commutative and associative but does not have an inverse; hence, the space of equivalence classes of multisections together with this addition is just a commutative monoid.

\noindent
Given an orbibundle $\pr\colon \cE\!\lra\! \cM$ with an orbibundle atlas $\cA_\cE$ as in (\ref{equ:atlasE}), an \textbf{orbifold multisection} with respect to $\cA_\cE$ is a family $\mfs\!\equiv\! \{\mfs_\al\}_{\al\in S}$ of (continuous) $G_\al$-equivariant $\ell_\al$-sections  $\mfs_\al \!\in\! S^{\ell_\al}(U_\al)$ which are equivalent along intersections; c.f. \cite[Definition 4.13]{Mc2006}  for a global (coordinate independent) functorial definition. Note that if $S$ is finite, e.g. whenever $M$ is compact, we can take a large enough $\ell$ such that all the local representatives $\mfs_\alpha$ are orbifold $\ell$-sections and $\mfs_\alpha\!=\!\mfs_\beta$ on the overlap. A multisection with respect to $\cA_\cE$ induces a multisection with respect to any refinement of $\cA_\cE$. Given two multisections $\mfs$ and $\mfs'$ over  $\cA_\cE$ and $\cA'_\cE$, respectively, $\mfs+\mfs'$ is a multisection defined over any common refinement of $\cA_\cE$ and $\cA'_\cE$.

\bDef{def:liftable}
Let $\cU\!=\!(\pr\colon \!U\!\!\lra\!\! V,G,\psi)$ be an orbibundle chart and $\mfs$ be an orbifold $\ell$-section of that. We say $\mfs$ is \textbf{liftable} (over this chart) if there is a $G$-invariant family of $\ell$ (not necessarily $G$-equivariant) sections 
$$
s^i\colon V\lra U,\quad \forall i\!\in\! [\ell],
$$
such that $\mfs\!=\![s^1,\ldots,s^\ell]$. Similarly,  with notation as in the paragraph before Example~\ref{exa:Pmn}, let $\pr\colon\! \cE\!\lra\! \cM$ be an orbibundle and $\mfs\!=\!\{\mfs_\al\}_{\al\in S}$ be an orbifold multisection with respect to some orbibundle atlas $\cA_\cE$.
We say $\mfs$ is liftable over $\cA_\cE$ if every $\mfs_\al$ is liftable. A multisection is called \textbf{locally liftable} if it is liftable over some orbibundle atlas.
\eDef

\noindent
Given two multisections $\mfs$ and $\mfs'$  which are liftable over  $\cA_\cE$ and $\cA'_\cE$, respectively, the multisection sum $\mfs+\mfs'$ is liftable over any common refinement of $\cA_\cE$ and $\cA'_\cE$.

\bExa{crazy-lift}
With $\cM=\P^1_{m,n}$ as in Example~\ref{exa:Pmn}, consider the trivial line bundle 
$$
\cE=\P^1_{m,n}\times \C.
$$
Over $\cV_0$, consider the section $s_0(z)= \rho_0(|z|) z$ where $\rho_0$ is a bump function supported near $z\!=\!0$. 
The extension of this section to $\cV_\infty$ is a multisection of the form 
$$
\mfs_0(w)= [\rho_0(|w|^{-n/m})w^{-n/m}]\in S^{m}(\cV_\infty\times \C).
$$ 
Similarly, 
consider the section $s_\infty(w)= \rho_\infty(|w|) w$ and extend it to a multisection of the form 
$$
\mfs_\infty(z)= [\rho_\infty(|z|^{-m/n})z^{-m/n}]\in S^{n}(\cV_0\times \C)
$$ 
over $\cV_0$.
Let $\mfs=\mfs_0+\mfs_\infty$. Since both $\mfs_0$ and $\mfs_\infty$ are liftable on their defining charts, $\mfs$ is a locally liftable multisection; however,  it is not liftable over the orbibundle atlas made  of $\cV_0\times \C$ and  $\cV_\infty\times \C$, i.e. $\cV_0\times\C$ and $\cV_\infty\times \C$ do not form a sufficiently refined atlas for $\mfs$.
\eExa

\bRem{rem:branches}
If $\mfs$ is a liftable orbifold $\ell$-section of an orbibundle chart $\cU$, the decomposition $
\mfs\!=\![s^1,\ldots,s^\ell]$
into a $G$-invariant set of sections of $U$ is not unique. For example, assume $U=\R\times \R$, $V=\R$, $\pr$ is the projection onto the first factor, and $G=\tn{id}$. For each pair $(\ep_1,\ep_2)=(\pm ,\pm )$, define 
$$
s_{\ep_1,\ep_2}(t)=\begin{cases} \ep_1 e^{-1/|t|} &\mbox{if } t\leq 0 \\ 
\ep_2 e^{-1/|t|} & \mbox{if } t\geq 0 \end{cases}
$$
to be the corresponding section of $U$. Then, the equations
$$
\mfs_1= [s_{+,+},s_{-,-}], \quad \mfs_2=[s_{+,-},s_{-,+}] ,
$$
describe different decompositions of the same $2$-section $\mfs\!=\! (\mfs_1\!\cong\! \mfs_2)$ into  a set of two smooth branches.
\eRem
\noindent
\bRem{rem:lifted}
In light of the preceding remark, it will be important in the construction of Euler class of an orbibundle in Section~\ref{sec:euler}, and its generalization to the construction of VFC of a Kuranishi structure in Section~\ref{sec:VFC}, that we consider a  consistent choice of the liftings of a multisection over different local charts. For this reason, we consider lifted multisections instead of just liftable multisections; see definition below. This change of perspective, compared to the original definitions in \cite{FO,FOOO,FOOO-detail}, simplifies some of the steps in the construction of perturbations and VFC in Sections~\ref{sec:prelim-kur} and~\ref{sec:VFC}.
\eRem

\bDef{def:lifted}
Let $\cE\!\lra\!\cM$ be  an orbibundle and $\mfs$ be a locally liftable multisection of $\cE$. A \textbf{lifting} of $\mfs$ consists of an orbibundle atlas $\cA_\cE$ as in (\ref{equ:atlasE}) with respect to which $\mfs\!=\!\{\mfs_\al\}_{\al\in S}$ is liftable and a choice of lifting of every $\mfs_\al$ into an unordered $G_\al$-invariant set of sections (called \textbf{branches}),
$$
\mfs_\al=[s_\al^1,\ldots,s_\al^{\ell_\al}],
$$
which are compatible on the overlaps in the following sense. For every common refinement $\Phi_\al,\Phi_\beta\colon \cU_{\al,\beta}\lra \cU_\al,\cU_\beta$, as in the paragraph following Definition~\ref{def:orbibundlechart}, the pull-back branches 
$$
\mfD\phi_\al^* (\ell_\beta\cdot \{s_\al^1,\ldots,s_\al^{\ell_\al}\})~~~\tn{and}~~~\mfD\phi_\beta^* (\ell_\al\cdot \{s_\beta^1,\ldots,s_\beta^{\ell_\beta}\})
$$
are equal as an unordered set with multiplicities. By a \textbf{lifted} multisection, we mean a locally liftable multisection together with a choice of lifting.
\eDef
\noindent
Given a lifting over some orbibundle atlas $\cA_\cE$, by restricting to sub-charts, we obtain a lifting over any refinement $\cA'_\cE$ of $\cA_\cE$. Therefore, the multisection sum of two lifted multisections is lifted as well. As we showed in Example~\ref{rem:branches}, a multisection can have many different liftings; thus, the notion of ``equivalence" between lifted multisections should be enhanced as well. 

\bDef{def:equivalence-of-lifted}
We say two lifted multisections $\mfs$ and $\mfs'$ are \textbf{equivalent}, and write $\mfs\cong \mfs'$, if they are equivalent as in (\ref{equ:common-mult}) and up to some multiple have the same germ of branches around every point.
\eDef

\bDef{def:liftable-transverse}
An orbifold multisection is said to be \textbf{smooth} (resp. $C^m$, with $m\!\geq\! 0$),  if it is lifted\footnote{This is a prerequisite.} and each local branch of the chosen lift is smooth (resp. $C^m$). An orbifold smooth (or $C^m$ for some $m\!>\!0$) multisection is said to be \textbf{transversal} if each local branch of the chosen smooth lift is transversal to the zero section.
\eDef

\noindent
We denote the set of $C^m$-smooth (hence lifted) orbifold multisections of $\cE$ by $C^m_{\multi}(\cE)$, the set of $C^m$-smooth orbifold sections of $\cE$ by $C^m(\cE)$, and the set of $C^m$-smooth functions on $\cM$ by $C^m(\cM)$. 
\bDef{def:limit}
Let $\cE\!\lra\!\cM$ be an orbibundle, $(\mfs_i)_{i=1}^\infty$ be a sequence of lifted orbifold multisections, and $||\cdot ||$ be some specified norm. We say this sequence $||\cdot ||$-\textbf{converges} to the lifted orbifold multisection $\mfs_\infty$, and write 
$$
\lim_{i\lra \infty} \mfs_i = \mfs_\infty,
$$ 
if there exist a finite atlas $\cA_\cE$ as in (\ref{equ:atlasE}) and a family of positive integers $\{\ell_\al\}_{\al\in S}$ such that the following conditions hold.
\bItem
\item For every $i\in \N\cup\{\infty\}$,  the given lifting of $\mfs_i$ is equivalent to a lifting $\{\mfs_{i,\al}\}_{\al\in S}$ of the form
\bEqu{equ:components2}
\mfs_{i,\al}=[s_{i,\al}^1,\ldots, s_{i,\al}^{\ell_\al}],\quad s_{i,\al}^j\colon V_\al\lra U_\al~~\forall\al\! \in\! S,~ j\!\in\![\ell_\alpha].
\eEqu
\item Every $s_{i,\al}^j$ has finite $||\cdot ||$-norm, and (for some choice of ordering of (\ref{equ:components2}))
\bEqu{equ:limitmulti}
\lim_{i\lra \infty} ||s_{i,\alpha}^j \!- \!s_{\infty,\al}^j||=0, \quad \forall\al\! \in\! S,~ j\!\in\![\ell_\alpha].
\eEqu
\eItem
\eDef

\subsection{Perturbations}\label{sec:perturbation}
\noindent
The purpose of this section is to introduce  a method which allows us to perturb any smooth orbifold section to a smooth (and hence lifted) multisection transverse to the zero section.

\noindent
Assume $\cE\!\lra\! \cM$ is an orbibundle and $s$ is an orbifold $C^m$-section of $\cE$.
Let 
\bEqu{equ:sequence-inc}
K_1\!\subset\! W_1 \!\subset\! K_2 \!\subset\! M
\eEqu
such that $K_1$ and $K_2$ are compact and $W_1$ is open.
Let 
$$
\mf{t}_1\!\in\! C^m_{\tn{multi}}(\cE|_{W_1})
$$ 
be a given lifted multisection of $\cE$ over $W_1$ such that $s|_{W_1}+\mf{t}_1$ is a transversal orbifold multisection of $\cE|_{W_1}$.

\noindent
Let $K_1^c$ and $W^c_1$  be the complements of $K_1$ and $W_1$ in $M$, respectively. 
By (\ref{equ:sequence-inc}), $W^c_1$ is closed and $K_1^c$ is an open neighborhood of $W_1^c$. By compactness of $K_2$ and Remark~\ref{rmk:OBL}, there exists a finite orbibundle atlas $\cA_\cE$ as in (\ref{equ:atlasE}) for $\cE$ over a neighborhood $W_2$ of $K_2$ with the following properties.
\emph{
\bItem
\item For every $\al\!\in\! S$, $V_\al$ is isomorphic to some small ball in $\R^{n+n_\al}$ and $U_\al$ is trivial; 
we fix a choice of smooth trivialization $U_\al\cong V_\al\times \R^{n_\al}$ and an identification
\bEqu{nalpha}
V_\al\cong B_{\de_\al}(0)\subset \R^{n+n_\al}.
\eEqu
\item There exists a decomposition $S\!=\!S_1 \sqcup S_{2}$ such that 
\bEqu{equ:cov-property}
K_1\! \subset\! \big(O_1\!=\! \bigcup_{\al\in S_1} \psi_\al(V_\al)\big) \!\subset\! W_1,\qquad  W_1^c\cap K_2\!\subset\! \big(O_2\!=\! \bigcup_{\al\in S_2} \psi_\al(V_\al)\big) \!\subset\! K_1^c\cap W_2.
\eEqu
\item Restricted to $O_1$, the given lifting of $\mft_1$ is equivalent to a lifting $\{\mft_{1,\al}\}_{\al\in S_1}$ of the form
\bEqu{equ:components}
\mft_{1,\al}=[t_{\al}^1,\ldots, t_{\al}^{\ell_\al}],\quad t_{\al}^j\colon V_\al\lra U_\al\quad \forall j\!\in\! [\ell_\al],~\al\! \in\! S_1.
\eEqu
\item The underlying orbifold atlas 
$$
\cA_{\cW_2}\!\equiv\!\{\cV_\al=(V_\al,G_\al,\psi_\al)\}_{\al\in S}
$$ 
covering $W_2$ admits an orbifold smooth \textbf{partition of unity}  $\{ \chi_\al\}_{\al\in S}$ with respect to $K_2$ (c.f. \cite[p.~361]{Sa}), i.e. each $\chi_\al\colon V_\al\lra [0,1]$ is a compactly supported $G_\al$-invariant smooth function on $V_\al$ and
$$
\ov\chi(x)= \hspace{-.1in}\sum_{\al\in S\colon\ov{x}\in \psi_\al(V_\al)} \ov{\chi}_\al(\ov{x})=1\quad\forall\ov{x}\!\in\! K_2.
$$
\eItem
}
\noindent
Let 
$$
W_3=\{\ov{x}\in W_2: \ov\chi(\ov{x})\neq 0\},\quad K_2\!\subset\!W_3\!\subset\!W_2.
$$
\noindent
We fix a choice of such atlas. Define
$$
\ov\eta\colon W_2\lra \R, \qquad \ov\eta(\ov{x})\!=\!\sum_{\al\in S_1}\ov\chi_\al(\ov{x});
$$
this is the underlying continuous function of an orbifold smooth function supported in $O_1$ with  $\ov\eta|_{K_1}\!\equiv \!1$. Via the cut-off function $\ov\eta$, we extend $\mft_1$ to the well-defined lifted $C^m$-multisection $\eta\mft_1$ over entire $W_2$.

\noindent
Let $\nu$ be an arbitrary set of not necessarily equivariant local $C^m$-sections 
$$
\nu\equiv \{\nu_\al\colon V_\al\lra U_\al\}_{\al\in S_2}.
$$ 
With respect to the chosen set of trivializations of each $U_\al$ and the identifications (\ref{nalpha}), we define the $C^m$-norm of $\nu$ by
$$
||\nu||^2_{C^m}= \sum_{\al\in S_2} || \nu_\al||^2_{C^m}.
$$
Let $\tn{Pert}(\cA_\cE)$ be the Banach space of such $\nu$ with finite norm.  Corresponding to every $\nu\!\in\! \tn{Pert}(\cA_\cE)$, we build a $C^m$-smooth lifted multisection $\mft(\nu)$ of $\cE|_{W_2}$ such that $\mft(\nu)|_{K_1}\!=\!\mft_1|_{K_1}$ in the following way.

\noindent
For every $\beta\!\in\! S_2$, if $G_\beta \!=\!\{g_1,\ldots, g_\ell\}$, let 
$$
\mft_\beta= [g_1\cdot\nu_\beta, \ldots, g_\ell\cdot\nu_\beta]
$$ 
be the corresponding natural lifted multisection in (\ref{equ:nat-enhance}). Let
\bEqu{equ:sum-of-t}
\mft_2\!=\!\sum_{\beta\in S_2} \chi_\beta \mft_\beta.
\eEqu
Since $\chi_\beta$ is supported in $V_\beta$, $\chi_\beta \mft_\beta$ naturally extends to a lifted multisection over $W_2$. The sum considered in (\ref{equ:sum-of-t}) is the sum of multisections as in (\ref{equ:multi-sum}). Since the multisection sum of a finite set of lifted multisections is lifted, $\mft_2$ is a lifted $C^m$-multisection supported in $O_2$ (hence it is supported away from $K_1$). Finally, we put 
$$
\mft(\nu) =\eta \mft_1 \!+\!  \mft_2;
$$ 
this is again a sum of multisections as in (\ref{equ:multi-sum}). It is lifted over an appropriate finite refinement of $\cA_\cE$.
More precisely, for every $\ov{x}\in W_3$, let
$$
S(\ov{x})=\{\al \in S\colon \ov{x}\in \psi_\al(V_\al)\},
$$
and $\cU_{x}\!=\!(U_x,V_x,G_x,\psi_x,\pr_x)$ be an orbibundle chart centered at $\ov{x}$ which is a refinement of $\cU_\al$, for every $\al\!\in\! S(\ov{x})$, and 
$$
 \psi_x(V_x)\!\subset\!W_3,\quad \ov\chi_\al|_{\psi_x(V_x)}\!\equiv\! 0 \qquad \forall\al\!\notin \!S(\ov{x}).
$$
The latter condition can be satisfied because support of $\ov\chi_\al$, for every $\al\!\notin\! S(\ov{x})$, is a compact subset of $V_\al$ and $\ov{x}\!\notin\! \psi_\al(V_\al)$. By the choice of $\cU_x$, for every $\ov{x}\!\in\! W_3$, the multisection $\mft(\nu)$ and hence $s+\mft(\nu)$ is lifted over $\cU_x$; moreover, the set of branches is unique up to ordering. Let 
\bEqu{equ:decomp-nu}
(s+\mft(\nu))|_{V_x}=[s+t(\nu)^1,\ldots,s+t(\nu)^{\ell_x}]
\eEqu
be the set of branches. 
By the compactness of $K_2$, there exists a finite set $T$ of such charts covering a neighborhood $W_4\!\subset\!W_3$ of $K_2$; moreover, fixing an ordering of $S$, an ordering of $G_\beta$, for every $\beta\!\in\! S_{2}$, and an ordering of the set branches of $\mft_{1,\al}$, for all $\al\!\in\! S_{1}$, determines\footnote{Although it is clear from the definition of $\mft(\nu)$, the precise argument for this claim is a bit tedious to write.} an ordering of branches in (\ref{equ:decomp-nu}), for all $\ov{x}\in T$.
Then, for every $\nu \!\in\!\tn{Pert}(\cA_\cE)$ let 
$$
\cZ(\nu) \!\equiv\! \coprod_{\ov{x}\in T} \coprod_{i=1}^{\ell_x} \{ y\in V_x\colon (s+t(\nu)^i)(y)=0\} \subset \coprod_{\ov{x}\in T} \coprod_{i=1}^{\ell_x} V_x .
$$
This completes the construction of extension and its zero set.
Every $\nu\!\in\! \tn{Pert}(\cA_\cE)$ is a \textbf{perturbation term} and the construction above gives us a systematic way of constructing a multisection $\mft(\nu)$ out of a perturbation term. Next, we are going to show that for generic $\nu$, the perturbed multisection $s\!+\!\mft(\nu)$ is transverse to the zero section on $W_4$. Note that by definition, $\mft(\nu)$ is already transverse to the zero section along $K_1$.

\vskip.1in
\noindent
A set 
$$
\Omega\!=\!\{\nu_1,\ldots,\nu_N\}\! \subset\! \tn{Pert}(\cA_\cE)
$$ 
is called a \textbf{globally generating set relative to $K_1$}, if for every $\ov{x}\!\in\! T$, every $\!i\in\! [\ell_x]$, and every $y\!\in\! V_x$, either $\ov\eta(\ov{y})\!=\!1$ or for every $u\!\in\! U_y$ there exists a linear combination $\nu\!=\!\sum_{j\in [N]} a_j \nu_j$ such that $t(\nu)^i(y)\!=\!u$. 

\bExa{exa:GGS}
With $n_\al$ as in (\ref{nalpha}), let 
$$
I=\{j\equiv(j_\al)_{\al\in S_2}\colon j_\al \in [n_\al]\}.
$$
For every $j\!\in\! I$, let
$$
\nu_j \equiv ((\nu_j)_\al)_{\al\in S_{2}},~~ (\nu_j)_\al\colon V_\al \lra U_\al\cong V_\al \times \R^{n_\al},~~ (\nu_j)_\al= \!\!\underbrace{(0,\ldots,0,1,0,\ldots,0)}_\text{\tn{only} 1 \tn{at the} $j_\al$-th \tn{position}}\!.
$$ 
Then the set $\Omega=\{\nu_j\}_{j\in I}$ is a globally generating set relative to $K_1$.
\eExa

\noindent
Given a globally generating set $\Omega=\{\nu_1,\ldots,\nu_N\} $, let 
$$
\R^\Om=\{  \nu\!=\!\sum_{j=1}^N a_j \nu_j,~a_j\in \R\}
$$
and 
$$
\cZ=\bigcup_{\nu \in \R^\Om} \cZ(\nu)\!\times\! \{\nu\},\quad \cZ\subset \bigg(\coprod_{\ov{x}\in T} \coprod_{i=1}^{\ell_x} V_x\bigg) \times \R^\Om.
$$
Let $\pi\colon \cZ\lra \R^\Om$ be the projection map. By Sard's theorem and transversality of $s+\mft(\nu)$ on the complement of support of $(1-\ov\eta)$, for large enough $m$ (i.e. if the sections are differentiable enough), generic $\nu\! \in\! \R^\Om$ is a regular value of $\pi$. This implies that for generic $\nu\!\in\! \R^\Om$,  $s\!+\!\mft(\nu)$ is a transverse multisection on the complement of $\eta\!\equiv\! 1$. Since $s+\mft(\nu)$ is already transverse on $O_1$, this implies that $s+\mft(\nu)$ is transverse everywhere on $W_4$. We summarize the outcome of this argument as the following proposition. This proposition is a relative (and more general) version of \cite[Lemma 3.14]{FO} which is used (but not explicitly stated) in \cite[Section 6]{FOOO-detail}.

\bPro{pro:estimate} Assume $\cE\!\lra\!\cM$ is an orbibundle and $s$ is an orbifold $C^m$-section of $\cE$, for some sufficiently large $m\!\in \N$. Let 
$$
K_1\!\subset\! W_1 \!\subset\! K_2 \!\subset\! M
$$
such that $K_1$ and $K_2$ are compact and $W_1$ is open.
Let 
$$
(\mft_{1,i})_{i\in \N}\subset C^m_{\tn{multi}}(\cE|_{W_1})
$$ 
be a given sequence of multisections of $\cE$ over $W_1$ such that $ (s|_{W_1}\!+\!\mft_{1,i})_{i\in \N}$ is a sequence of transversal orbifold multisections of $\cE|_{W_1}$ and $\lim_{i\to \infty} \mft_{1,i} \!=\!0$ in $C^m$-convergence norm over $W_1$. Then, there exists a neighborhood\footnote{This is $W_4$ of the construction above.} $W$ of $K_2$ in $M$ and a sequence of a $C^m$-multisections $(\mft_i)_{i\in \N}\!\subset\! C^m_{\tn{multi}}(\cE|_{W})$ such that $\mft_i|_{K_1}\cong \mft_{1,i}|_{K_1}$, 
$(s+\mft_i)_{i\in \N}$ is a sequence of transversal orbifold multisections of $\cE|_{W}$, and $\lim_{i\to \infty} \mft_i \!=\!0$ in $C^m$-convergence norm over $W$.
\ePro

\bRem{rem:smoothness}
The orbifolds, Kuranishi structures, and gluing maps considered in this article are assumed to be (sufficiently) smooth, only to satisfy of the smoothness condition in Sard's theorem we used above. Nevertheless, $C^1$-smoothness is enough for the rest of arguments throughout the article.
In genus $0$ (and assumably higher genus), a middle ground between the $C^1$-smoothness and smoothness condition is the 
$\tn{C}^1\tn{SS}$ ($C^1$ Stratified Smooth) condition considered in \cite{Rob}. According to \cite[Lemma 5.2.2]{Rob}, Sard's theorem holds in the category of $\tn{C}^1\tn{SS}$-smooth manifolds as well.
\eRem

\subsection{Resolution of multisections}\label{sec:resolution}
Let $E\!\lra\! M$ be a relatively oriented vector bundle and $s\colon M\lra E$ be a transversal section with compact support. Then the zero set of $s$ is a smooth oriented manifold which defines the Euler class of $E$ as a singular homology class
in $M$; see the beginning of Section~\ref{sec:euler}. Similarly, we can define the Euler class of orbibundles, if we can find a transversal orbifold section.

\noindent
As we argued at the beginning of Section~\ref{sec:multisection}, it is quite likely that an orbibundle does not admit any transversal orbifold section. In the light of Proposition~\ref{pro:estimate}, we consider a transversal multisection, in place of an orbifold section, and  construct the Euler class from the ``zero set" of the multisection. To this end, we introduce the notion of 	``resolution" of multisections which allows us, after passing to some covering spaces, to reduce the problem back to the of case of orbifold sections.

\noindent
Throughout this section, assume $\pr\colon\! \cE\!\lra\! \cM$ is a smooth (or $C^m$ for some $m\!>\!0$) relatively oriented orbibundle and $\mfs$ is a lifted orbifold multisection of $\cE$ with compact support.  By assumption, there is a finite orbibundle atlas $\cA_\cE$ as in (\ref{equ:atlasE}) covering a neighborhood of the zero set of $\mfs$  such that the given lifting of $\mfs$ is equivalent to a lifting $\{\mfs_{\al}\}_{\al\in S}$ of the form
\bEqu{equ:local-presentation}
\mfs_{\al}=[s_{\al}^1,\ldots, s_{\al}^{\ell_\al}],\quad s_{\al}^j\colon V_\al\lra U_\al\quad\forall\al\! \in\! S,~ j\!\in\![\ell_\alpha].
\eEqu
Let  
\bEqu{equ:sero-set}
Z(\mfs)= \bigcup_{\substack{\al\in S\\ j\in [\ell_\al]}} \psi_\al((s_\al^j)^{-1}(0))\subset M
\eEqu
be the zero set of $\mfs$.

\noindent
In general, $Z(\mfs)$ can be a complicated subspace of $M$, i.e. it may not be a sub-orbifold or  it may not even admit a  triangulation. In order to resolve this issue, we introduce the notion of \textbf{resolution} of $\mfs$. In short, a resolution of $\mfs$ consist of a covering of a neighborhood $W$ of $Z(\mfs)$ by a finite number of open sets, and a compatible system of orbifold covering spaces over these open sets, such that the multisection $\mfs$ on each piece is the push forward of some orbifold section on the covering. In the next section, we start with a resolution of $\mfs$, compatibly triangulate the zero sets of the associated orbifold sections, and push forward these triangulations (with some $\Q$-coefficients) into $Z(\mfs)$ to define the Euler class of $\cE$.
We begin the discussion by recalling the notion of covering space for orbifolds; see \cite{Th} for more details.

\bDef{def:k-cover}
An orbifold smooth map $q\colon\cM_1\lra \cM_2$ is called a $c$-\textbf{covering map}, if every $\ov{x}\!\in\! M_2$ has an orbifold chart $(V,G,\psi)$ centered around that, over which the orbifold structure on $\ov{q}^{-1}(\psi(V))$ is of the form
\bEqu{equ:cover-charts}
\coprod_{a=1}^N V \lra \coprod_{a=1}^N V/G_a,
\eEqu
where $G_a$ (for every $a\!\in\! [N]$) is a subgroup of $G$, 
\bEqu{equ:sum-of-terms}
 \sum_{a=1}^N \abs{G/G_a}=c,
\eEqu
and the map $q$ on each component of (\ref{equ:cover-charts}) is given by the identity map on $V$ as in the following diagram,
$$
\xymatrix{
V \ar[d]^{\pi} \ar[r]^{\overset{q}{\cong}} & V \ar[d]^{\pi} \\
V/G_a \ar[r]^{\ov{q}}        & V/G\;.       }
$$
\eDef
\noindent
For example, if $\cM_2$ is an effective global quotient $[M/ G]$, then the manifold $\cM_1\!=\! M$ is a $|G|$-covering of $\cM_2$. In fact, for every $x\!\in\! M$ and some sufficiently small $G_x$-invariant open neighborhood $V$ of $x$, $[V/ G_x]$ is an orbifold chart of $[M/ G]$ centered at $\ov{x}$, $\{ [g]\cdot V \}_{[g]\in G/G_x}$ are the manifold charts around different preimages of $[V/ G_x]$, and (\ref{equ:sum-of-terms}) is simply equal to 
$$
\sum_{i=1}^{|G/G_x|} |G_x| = |G|.
$$ 
If $\cM_2$ is a manifold, then $\cM_1$ is a covering space in the usual sense. 
It is clear from Definition~\ref{def:k-cover} that every covering map is an orbifold map such that all the homomorphisms $h_x$ as in Definition~\ref{def:orbi-maps} are monomorphism.
Therefore, by Remark~\ref{rem:pull-back}, if $q\colon\! \cM_1\!\lra\! \cM_2$ is a $c$-covering map, every orbibundle $\cE_2$ over $\cM_2$ naturally pulls back to an orbibundle $\cE_1\!=\! q^*\cE_2$ over $\cM_1$. If $\cE_2\!\lra\! \cM_2$ is relatively oriented, then so is $\cE_1\!\lra\!\cM_1$.

\vskip.1in
\noindent 
Let $\cM_{1,o}$, with $1\!\leq\! o\! \leq\! t$,
be the connected components of $\cM_1$ and suppose $q|_{\cM_{1,o}}$ is a $c_o$-covering; we have
$$
\sum_{o=1}^t c_o=c.
$$
Let $\{w(o)\}_{o=1}^t$ be a set of positive integers, called \textbf{weights} below, and set 
$$
\ell=\sum_{o=1}^t c_o~w(o).
$$
Then, every orbifold section $s$ of $\cE_1$ descends to a lifted $\ell$-section $\mfs$ of $\cE_2$, called the \textbf{weighted push forward multisection}, in the following way. With notation as in Definition~\ref{def:k-cover}, let $\pr\colon U\!\lra\! V$ be an orbibundle chart of $\cE_2$ and 
$$
s^{a}\colon V\lra U, \quad 1\!\leq\! a\! \leq\! N,
$$
be the family of local $G_a$-equivariant representatives of $s$ over different components of $q^{-1}[V/ G]$ in (\ref{equ:cover-charts}). Let $c'_a\!=\!|G/G_a|$, for all $1\!\leq\! a\! \leq\! N$. Each orbifold chart $[V/ G_a]$ for $\cM_1$ belongs to some connected component $\cM_{1,o(a)}$; we define $w(a)\!\in\! \Z^{>0}$ to be $w(o(a))$. We define the weighted push forward multisection 
$$
\mfs\!= \!q^{w}_*s
$$ 
to be the multisection given by the following set of branches on the orbifold chart $[V/ G]$ for $\cM_2$,
\bEqu{equ:push-forward}
\mfs|_V= [\underbrace{g_{11}s^1,\ldots, g_{1c'_1}s^1,\dots}_\text{\tn{repeated}~$w(1)$~\tn{times}},\ldots\ldots,
\underbrace{g_{N1}s^N,\ldots, g_{Nc'_N}s^N,\ldots}_\text{\tn{repeated}~$w(N)$~\tn{times}}],
\eEqu
where $\{g_{a1},\ldots,g_{ac'_a}\}$ forms a set of representatives of $G/G_a$. In simple words, $\mfs$ is a lifted multisection whose branches over $V$ are give by the section $s$ over different lifts of $V$ in (\ref{equ:cover-charts}), and the number of times each branch appears is equal to the weight of the component it comes from.

\bRem{rem:ration}
Replacing the weight function $w$ with a uniform  multiple of that does not change the equivalence class of the resulting weighted push forward multisection. We define the \textbf{weight ratio} by 
\bEqu{equ:weights}
\ov{w}(o)= \frac{w(o)}{\ell} = \frac{w(o)}{\sum_{o'=1}^t c_{o'}~ w(o')}.
\eEqu
We can think of (\ref{equ:weights}) as a locally constant rational function 
$$
\ov{w}\colon M_1 \lra \Q.
$$
We use the weight ratio function (\ref{equ:weights}) in  the formulation of Euler class singular cycle in (\ref{equ:Euler-class}). Note that $\ov{w}\equiv \frac{1}{c}$, if all the weights are equal to $1$.
\eRem

\vskip.1in
\noindent
Let us return to the construction of resolution. Given an orbibundle chart $\cU\!=\!(\pr\colon \!U\!\lra\! V,G,\psi)$ and a lifted $\ell$-section 
$
\mfs\!=\![s^1,\ldots,s^\ell]
$
on $\cU$, we define the \textbf{multiplicity} of every $s^i$ \textbf{over} $V$ to be the number of branches $s^j$ which are equal to $s^i$ over $V$. For a push forward multisection of the form (\ref{equ:push-forward}), for generic choices of $s^1,\ldots,s^N$, $\mfs$ has $c\!=\!\sum_{a=1}^N c'_a$ distinct branches where each $g_{ij}s^i$ has multiplicity $w(i)$.

\noindent
For every $\ov{x}\!\in\! M$, let 
\bEqu{equ:val-function}
\tn{val}_\mfs(\ov{x})\!\in\! \N
\eEqu
be the number of  distinct germs of branches of the given transversal multisection $\mfs$ around $\ov{x}$. In other words, if $\tn{val}_\mfs(\ov{x})=c$, then there exists some open set $W\!\subset\! M$ around $\ov{x}$ such that for every orbibundle chart $\cU$ with $\psi(V)\!\subset\! W$, $\mfs|V$ has exactly $c$ distinct branches.  The function $\tn{val}_\mfs$ is well-defined and upper semi-continuous; see Figure~\ref{fig:val}. Therefore, for every $c\!\in\! \N$, the set
$$
M(\geq\! c)= \{\ov{x}\!\in\! M\colon\tn{val}_\mfs(\ov{x})\!\geq \!c\}
$$
is closed. 
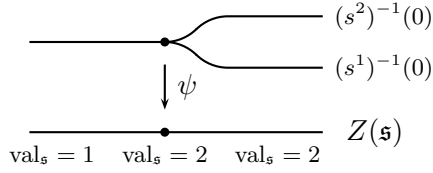
\begin{figure}
\begin{pspicture}(-6,-6.1)(10,-5)
\psset{unit=.3cm}
\psline(-6,-16)(0,-16)
\psarc(0,-14){2}{-90}{-45}
\psarc(2.8,-16.86){2}{90}{135}
\psarc(0,-18){2}{45}{90}
\psarc(2.8,-15.14){2}{-135}{-90}
\psline(2.8,-17.14)(7,-17.14)
\psline(2.8,-14.86)(7,-14.86)
\psline(-6,-20)(7,-20)
\psline{->}(0,-17)(0,-19)
\rput(1,-18){$\psi$}
\rput(9.7,-17.14){\footnotesize{$(s^1)^{-1}(0)$}}
\rput(9.7,-14.86){\footnotesize{$(s^2)^{-1}(0)$}}
\rput(9.3,-20){$Z(\mfs)$}
\pscircle*(0,-20){0.2}\rput(0,-21){\footnotesize{$\tn{val}_\mfs=2$}}
\pscircle*(0,-16){0.2}
\rput(5,-21){\footnotesize{$\tn{val}_\mfs=2$}}
\rput(-5,-21){\footnotesize{$\tn{val}_\mfs=1$}}
\end{pspicture}
\caption{The $\tn{val}_\mfs$ function for a multisection with two branches.}
\label{fig:val}
\end{figure}
Let
$$
M(c)= \{\ov{x}\in M\colon~\tn{val}_\mfs(\ov{x})\!=\! c\}\quad \forall c\!\in\!\N.
$$

\bDef{def:resolution}
For any $c\!\in\! \N$ and any compact set $K\!\subset\! M(c)$, by a \textbf{resolution of $\mfs$ around $K$} we mean an open neighborhood $W$ of $K$ and an orbifold $\wt\cW$ with a $c$-covering map $q\colon \wt\cW \lra \cW= \cM|_{W}$, such that the lifted multisection $\mfs$ restricted to $W$ is equivalent to the weighted push forward of some orbifold section $s$ of $\wt\cE=q^*(\cE|_{W})$.
\eDef

\bLem{lem:cover}
For every $c\!\in\! \N$ and every compact set $K\!\subset\! M(c)$ such that
$$
K\cap M(>\! c)\!=\!\emptyset,
$$ 
there exists a resolution of $\mfs$ around $K$.
\eLem

\begin{proof}
For every $\ov{x}\!\in\! K$ choose an orbibundle chart $\cU_x\!=\!(\pr_x\colon \!U_x\!\lra\! V_x,G_x,\psi_x)$ centered at $\ov{x}$ such that 
$$
\mfs|_{V_{x}}\!=\![s^1,\ldots,s^{\ell_x}],
$$
for some $\ell_x\!\geq\!c$, has exactly $c$ distinct branches and $\psi_x(V_x)\cap M(> \!c)\!=\!\emptyset$. Since $K$ is compact, there exists a finite subset of points $I\!\subset\! K$ such that  
$$
W= \bigcup_{\ov{x}\in I} \psi_x(V_{x})\subset M
$$
is an open neighborhood of $K$ in $M$ disjoint from $M(>\! c)$. 
Over each $V_{x}$, if $s^{j_1},\ldots, s^{j_{c}}$ forms a set of distinct branches of $\mfs|_{V_{x}}$,  let 
$$
J_x\!=\!\{j_1,\ldots, j_{c}\}\subset [\ell_x]
$$ 
and define the weight $w(j_a)$ of $j_a$ to be the multiplicity of $s^{j_a}$.
Since the set of branches of $\mfs|_{V_{x}}$ is $G_x$-invariant, $G_x$ acts on $J_x$ and hence on the product
$$
\wt{V}_{x}= V_{x}\times J_x.
$$
The projection map into the first component, $q_x\colon\wt{V}_{x}\!\lra\! V_{x}$, is a $G_x$-equivariant $c$-covering map. 
If $|I|\!=\!1$, say $I\!=\!\{\ov{x}\}$, the effective global quotient orbifold $\wt\cW\!=\! [\wt{V}_{x}/ G_x]$ has the desired property. In fact, every connected component of $\wt{V}_{x}/G_x$ is of the form 
\bEqu{equ:comp-o}
(V_x\times o)/G_{x}\cong V_x/G_{x,o},
\eEqu
where $o\!=\!G\cdot j_a$ is some $G_x$-orbit in $J_x$ and $G_{x,o}$ is the isotropy group of $j_a$ (Changing the choice of $j_a$ in $o$ changes everything by conjugation). Since the weight function defined above is constant in each orbit, we define the weight $w(o)$ over the connected component  (\ref{equ:comp-o}) to be the weight of $j_a$.  
Restricted to such component, the $G_{x,o}$-equivariant orbifold section 
$
s_o\colon V_x \lra U_x
$ 
can be chosen to be $s^{j_a}$. If $|I|\!>\!1$, the conclusion follows from Lemma~\ref{lem:cover2} below and induction on $|I|$.
\eProof
 
\bLem{lem:cover2}
For some $c\!\in\! \N$, let $K_1,K_2\!\subset\! M(c)$ be two compact subsets such that 
$$
K_1\cap M(>\!c)\!=\!K_2\!\cap\! M(> \!c)\!=\!\emptyset
$$ 
and both $K_1$ and $K_2$ admit a resolution of $\mfs$ around them.
Then $K\!=\!K_1\!\cup\! K_2$ admits a resolution of $\mfs$ around it.
\eLem

\bProof
By assumption, for $j\!=\!1,2$, there exist a neighborhood $W_j$ of $K_j$, disjoint from $M(>\!c)$, and an orbifold $\wt\cW_j$ with a $c$-covering map 
$$
q_j\colon \wt\cW_j\lra \cW_j= \cM|_{W_j},
$$ 
such that the multisection $\mfs$ restricted to $W_j$ is equivalent to the weighted push forward, with weight $w_j$, of an orbifold section $s_j$ of $\wt\cE_j=q_j^*(\cE|_{W_j})$.

\noindent
For every $\ov{x}\!\in\! K_{12}\!=\! K_1\cap K_2$, we can choose an orbibundle chart $\cU_{x}$
centered at $\ov{x}$, such that
$\mfs|_{V_{x}}$ is lifted and has  exactly $c$ distinct branches, 
$\psi_x(V_x)$ is a subset of $W_1\cap W_2$, 
the covering map $q_j$ (for $j\!=\!1,2$) lifts over $V_x$ as in Definition~\ref{def:k-cover}, and  $\mfs={q_{j*}^{w_j}}s_j$ as in (\ref{equ:push-forward}). 
Moreover, we can choose these charts sufficiently small so that the following two conditions hold.
\bEnum
\item If $\psi_x(V_x)\cap \psi_y(V_y)\neq \emptyset$, there exists some orbibundle chart $\cU\!=\!(\pr\colon \!U\!\lra\! V,G,\psi)$ with the same properties such that 
\bEqu{equ:whenever0}
\psi_x(V_x)\cup \psi_y(V_y)\subset \psi(V).
\eEqu
\item For every $\ov{x},\ov{y}\in K_{12}$, 
\bEqu{equ:whenever1}
\psi_x(V_x)\cap \psi_y(V_y)\cap M(c)\neq \emptyset,~\tn{whenever}~\psi_x(V_x)\cap \psi_y(V_y)\neq \emptyset.
\eEqu
\eEnum

\noindent
Since $K_{12}$ is compact, a finite set $I_{12}$ of such charts covers $K_{12}$. Let 
$$
W_{12}= \bigcup_{\ov{x}\in I_{12}} \psi_x(V_x),\quad K_{12}\subset W_{12}\subset W_1\cap W_2.
$$

\noindent
Since $K_1\!\setminus\! W_{12}$ and $K_2\!\setminus\! W_{12}$ are compact and 
$$
(K_1\setminus W_{12})\cap (K_2\setminus W_{12})=\emptyset,
$$
by possibly reducing $W_1$ and $W_2$ to smaller open sets, we can assume that 
$$
W_1\cap W_2=W_{12}.
$$
For every $\ov{x}\!\in\! I_{12}$, let 
$$
q_x\colon\!\wt{V}_x=V_x\times J_x\!\lra\! V_x
$$ 
be as in the proof of Lemma~\ref{lem:cover}. 
On the other hand, for any of $j\!=\!1,2$ and by Definition~\ref{def:k-cover}, $q_j^{-1}[V_x/ G_x]$ is also of the form $(V_x\times J^j_x)/G_x$, 
where 
$$
J^j_x=\coprod_{a=1}^{N^j} G_x/G_{x,a}^j, \quad |J^j_x|=c~.
$$
Here $\{G_{x,a}^j\}_{a=1}^N$ is a set of subgroups of $G_x$ and $G_x$ acts on $J^j_x$ by acting on each coset $G_x/G_{x,a}^j$ from the left.
For either of $j\!=\!1,2$, by definition, both $J_x$ and $J^j_x$ correspond to the set of distinct branches of the multisection $\mfs$ on $V_x$; therefore, there is a unique identification between the set of orbits of the action of $G_x$ on $J^1_x$ and $J^2_x$. This gives an identification of $\wt\cW_1$ and $\wt\cW_2$ along $q_1^{-1}[V_x/ G_x]$ and $q_2^{-1}[V_x/ G_x]$. The natural identifications corresponding to distinct $\ov{x},\ov{y}\in I_{12}$ are compatible along the intersections for the following two reasons. 
By (\ref{equ:whenever1}), we can identify the set of distinct branches over $V_x$ and $V_y$ by their agreement near a point 
\bEqu{equ:zglue}
\ov{z}\!\in\! \psi_x(V_x)\cap \psi_y(V_y)\cap M(c).
\eEqu
By (\ref{equ:whenever0}), every two intersecting charts lie in a bigger chart over which there are exactly $c$ distinct branches; this condition assures us that the previous set of identifications is independent of the choice of the point $\ov{z}$ in (\ref{equ:zglue}). We conclude that, $\ov{q}_1^{-1}W_{12}$ and $\ov{q}_2^{-1}W_{12}$ are isomorphic as topological spaces and the orbifold structure induced on $\ov{q}_1^{-1}W_{12}$ by $\wt{\cW}_1$ is the same as the orbifold structure induced on $\ov{q}_2^{-1}W_{12}$ by $\wt\cW_2$. Identifying $\wt\cW_1$ and $\wt\cW_2$ along this common intersection, we obtain a $c$-covering $\wt\cW$ of $\cW|_{W_1 \cup W_2}$. Finally, notice that the weight functions $w_1$ and $w_2$, up to multiplication by a constant, are determined by the multiplicity of distinct branches of the multisection $\mfs$. Therefore,  after possibly multiplying each of $w_1$ and $w_2$ by a constant, ${q_{1*}^{w_1}}s_1$ and ${q_{2*}^{w_2}}s_2$ are equal along $W_{12}$ and we obtain a weight function $w$ on the union such that $\mfs|_{W_1\cup W_2}=q^{w}_*s$.
\end{proof}

\bRem{rem:John}
Let $K$ be a compact subspace of $M$. More generally, we can define the notion of the \textbf{germ of an orbifold $c$-covering} of $\cM$ along $K$: simply the direct limit of orbifold $c$-coverings of neighborhoods of $K$ under restriction. In order to define such a germ, it suffices for each point of $\ov{x}$ in $K$ to admit an orbifold chart centered at $\ov{x}$ and a $c$-covering as in Definition~\ref{def:k-cover} with the property that if 
$\ov{y}\!\in\! \psi_x(V_x)$, then restricted to some smaller open set $V'_y\! \subset\! V_y$ such that $V'_y$ is a refinement of $V_x$, the orbifold covering of $V_y$ restricted to $V'_y$ is ``naturally" identified with the restriction of the orbifold covering of $V_x$ restricted to $V'_y$. For example, in the case of Lemma~\ref{lem:cover}, the natural identifications are given by the set of distinct branches of $\mfs$.
\eRem

\noindent
By assumption, $Z(\mfs)$ is compact; therefore, $\tn{val}_\mfs|_{Z(\mfs)}$ takes only finitely many values
$$
c_1\!<\!\cdots\!<\!c_N.
$$
We use Lemma~\ref{lem:cover} to build a sequence of resolutions of $\mfs$ around some compact subsets of 
$$
Z(\mfs)_i= M(c_i)\cap Z(\mfs)\quad \forall 1\!\leq\!  i \!\leq \! N.
$$
We then use these coverings to build the overall resolution of $\mfs$.

\noindent
Starting with $i\!=\!\!N$, since $M(\geq \!c_N)\!=\!M(c_N)$ is closed and $Z(\mfs)$ is compact, we conclude that $Z(\mfs)_N\!= \!M(c_N)$ is compact. Therefore, by Lemma~\ref{lem:cover},  there exist~a neighborhood $W_N$ of $Z(\mfs)_N$, and an orbifold $\wt\cW_N$ with a $c_N$-covering map 
$$
q_N\colon \wt\cW_N\lra \cW_N= \cM|_{W_N},
$$ 
such that the multisection $\mfs$, restricted to $W_N$, is equivalent to the weighted push forward, with weight $w_N$, of some orbifold section $s_N$ of 
$$
\wt\cE_N= q_N^*(\cE|_{W_N}).
$$

\noindent
We proceed inductively. Having constructed $q_i\colon\!\wt\cW_{i}\lra \cW_{i}$, for all $N\!\geq\! i\! \geq\! t$, such that 
$$
W_{\geq t} = \bigcup_{t\leq i \leq N} W_i 
$$
is a neighborhood of $Z(\mfs)\cap M(\geq\! c_{t})$, let
$$
K= (Z(\mfs) \cap M(c_{t-1})) \setminus W_{\geq t}\;.
$$
By the upper semicontinuity of $\tn{val}_\mfs$ and compactness of $Z(\mfs)$, $K$ is a compact subset of $M(c_{t-1})$ disjoint from $M(c_i)$, for all $i\!\geq\! t$.
Therefore, again by Lemma~\ref{lem:cover},  there exist a neighborhood $W_{t-1}$ of $K$ disjoint from $M(\geq\!c_t)$, and an orbifold $\wt\cW_{t-1}$ with a $c_{t-1}$-covering map 
$$
q_{t-1}\colon \wt\cW_{t-1}\lra \cW_{t-1}= \cM|_{W_{t-1}}
$$
such that the multisection $\mfs$, restricted to $W_{t-1}$, is equivalent to the weighted push forward, with weight $w_{t-1}$, of some orbifold section $s_{t-1}$ of 
$$
\wt\cE_{t-1}\!=\!q_{t-1}^*(\cE|_{W_{t-1}}).
$$
We continue this induction process until we exhaust $Z({\mfs})$ at $c_1$.

\noindent
It also follows from the proof of Lemma~\ref{lem:cover2} that on the overlaps 
$$
W_{i,j}= W_i\cap W_{j}\quad \forall 1\!\leq \!j\!<\!i \!\leq\! N,
$$ 
we have $q_i\!=\!q_{j}\circ q_{i,j}$, where $q_{i,j}$ is a component-wise\footnote{For example if $c_i\!=\!3$ and $c_{j}\!=\!2$, $q_{i,j}$ is a $2$-covering map on some connected components and an isomorphism everywhere else.} covering map
$$
q_{i,j}\colon \wt\cW_i|_{\ov{q}_i^{-1}(W_{i,j})}\lra \wt\cW_{j}|_{\ov{q}_{j}^{-1}(W_{i,j})}.
$$
Finally, the weight ration functions (\ref{equ:weights}) are related to each other on the overlaps in the following way.
\bLem{lem:weight-relation}
For every $i\!\in\![N]$, let $\mf{I}_i=\mf{I}_{\wt{\cW}_i}$ be the integer valued function\footnote{Which is equal to the size of isotropy group at each point.} in (\ref{equ:isotropy-order}) and $\ov{w}_i$ be the weight ratio rational function in (\ref{equ:weights}). For every  $1\!\leq\! j\!<\!i\leq\! N$, the push forward via $\ov{q}_{i,j}$ of the twisted weight ratio function 
$(\ov{w}_i/\mf{I}_i)_{\ov{q}_i^{-1}(W_{i,j})}$  is equal to $(\ov{w}_j/\mf{I}_j)_{\ov{q}_j^{-1}(W_{i,j})}$, i.e. 
$$
\ov{q}_{i,j*} (\ov{w}_i/\mf{I}_i)|_{\ov{q}_i^{-1}(W_{i,j})}=(\ov{w}_j/\mf{I}_j)|_{\ov{q}_j^{-1}(W_{i,j})}.
$$
In other words, for every $\ov{x}\in \ov{q}_j^{-1}(W_{i,j})\subset \wt{W}_{j}$, we have
\bEqu{equ:weight-relation}
\frac{\ov{w}_j(\ov{x})}{\mf{I}_j(\ov{x})}=\sum_{\ov{y}\in \ov{q}_{i,j}^{-1}(\ov{x})\subset \ov{q}_i^{-1}(W_{i,j})} \frac{\ov{w}_i(\ov{y})}{\mf{I}_i(\ov{y})}.
\eEqu
\eLem

\bProof
The equality (\ref{equ:weight-relation}) is a consequence of the following explicit calculation.
Let $\cU\!=\!(\pr\colon \!U\!\lra\! V,G,\psi)$ be an orbibundle chart of $\cE$ over which $\mfs$ is equal to a lifted $\ell$-section
$$
\mfs\cong[s^1,\ldots,s^\ell].
$$ 
Suppose $\mfs$ restricted to $V$ has $c$ distinct branches given by the index set
$$
J=\{{j_1},\ldots,{j_{c}}\}\subset [\ell].
$$ 
For every $\ov{x}\!\in \!\psi(V)$ with $\tn{val}_\mfs(\ov{x})\!=\!c_x\!\leq\! c$, let $\cU_x\!=\!(\pr|_{U_x}\colon \!U_x\!\lra\! V_x,G_x,\psi|_{V_x})$ be a sub-chart of $\cU$ centered around $\ov{x}$ over which $\mfs$ has $c_x$ distinct branches given by the index set
$$
J_x=\{{i_1},\ldots,i_{c_x}\}\subset [\ell].
$$ 
There exists an obvious projection map $\pi\colon\! J\lra J_x$  which sends $j_a$ to $i_b$ whenever 
$$
s^{j_a}|_{V_x}=s^{i_b}|_{V_x}.
$$ 
This projection map commutes with the subgroup action of $G_x$ on $J$ and the action of $G_x$ on $J_x$. In this projection, a $G_x$-orbit $O$ in $J$ can be larger than its image $G_x$-orbit $o=\pi(O)$ in $J_x$.
As in the proof of Lemma~\ref{lem:cover}, corresponding to $\cU$ we obtain a $c$-covering  
$$
q_2\colon \cW_2= [(V\times J)/ G]\lra  [V/ G]
$$  
whose connected components are of the form 
$$
(V\times \mf{O})/G\cong V/G_{\mf{O}},
$$
where $\mf{O}$ is an orbit of the action of $G$ on $J$ and $G_{\mf{O}}<G$ is the isotropy group of some element of $\mf{O}$. Let $\mf{I}_2\colon W_2\lra \Z$ be the corresponding function as in (\ref{equ:isotropy-order}).
Similarly, corresponding to $\cU_x$ we obtain a $c_x$-covering 
 $$
 q_1\colon \cW_1= [(V_x\times J_x)/ G_x] \lra [V_x/ G_x]
 $$  
 whose connected components are of the form 
$$
(V_x\times o)/G_x\cong V_x/(G_{x})_{o},
$$
where $o$ is an orbit of the action of $G_x$ on $J_x$ and $(G_{x})_{o}\!<\!G_x$ is the isotropy group of some element of $o$. The function $\mf{I}_1\colon W_1\lra \Z$ is defined similarly.
Let $w_1$ and $w_2$ be the weight functions of these two coverings; $w_2(\mf{O})$ is the multiplicity of $s^{j_a}$ in $\{s^1,\ldots,s^\ell\}$ for any $j_a\!\in\! \mf{O}$, and $w_1(o)$ is the multiplicity of $s^{i_b}|_{V_x}$ in $\{s^1|_{V_x},\ldots,s^\ell|_{V_x}\}$ for any $i_b \in o$. Let $q_{2,1}$ be the projection map from $q_2^{-1}([V_x/ G_x])$ to $\cW_1$. We have the following commutative diagram
$$
\xymatrix{
 & [(V_x\times J)/ G_x] \ar[ld]^{q_{2,1}}\ar[ldd]^{q_2} \ar@{^{(}->}[r]^{\tn{~~refine}} &[(V\times J)/ G]\ar[dd]^{q_2}\\
[(V_x\times J_x)/ G_x] \ar[d]^{q_1} &         &               \\
[V_x/ G_x] \ar@{^{(}->}[rr]^{\tn{refine}}  &  & [V/ G], }
$$
where $q_{2,1}$ sends $[y\times O]_{G_x}$ to $[y\times \pi(O)]_{G_x}$.
Fix  a $G_x$-orbit $o$ in $J_x$, say $o\!=\!G_x\cdot \{i_1\}$, and let $\ov{x}_o\!=\![x\times o]_{G_x}$ be the unique element in $\ov{q}_1^{-1}(\ov{x})$ which belongs to the connected component $(V_x\times o)/G_x$. 
Let $\{O_{a}\}_{a=1}^t$ be the set of all $G_x$-orbits in $J$ where $\pi(O_a)\!=\!o$, and denote by $\ov{x}_a\!=\![x\times O_a]_{G_x}$ to be the unique element in $\ov{q}_2^{-1}(\ov{x})$ which belongs to $(V_x\times O_a)/G_x$. 
We have
$$
\ov{q}_{2,1}^{-1}(\ov{x}_o)=\{\ov{x}_a\}_{a=1}^t,
$$
and 
\bEqu{equ:w-decomp}
\lrp{(\ov{q}_{2,1})_*\lrp{\frac{\ov{w_2}}{\mf{I}_2}}}(\ov{x}_o)=\sum_{a=1}^t \frac{\ov{w}_2(\ov{x}_a)}{|G_{x,a}|}, 
\eEqu
where $G_{x,a}<G_x$ is the isotropy subgroup of some element of $O_a$. By definition, if $O_a= G_x\cdot \{j_a\}$, where $s^{j_a}|_{V_x}=s^{i_1}|_{V_x}$, then 
\bEqu{equ:compare-w}
\aligned
\ov{w}_1(\ov{x}_o)&=\frac{\tn{multiplicity}|_{V_x} (s^{i_1})}{\ell}=\sum_{a=1}^t \frac{\tn{multiplicity}|_V (s^{j_a})}{\ell}\frac{|O_a|}{|o|}=\\
&\sum_{a=1}^t \frac{\tn{multiplicity}|_V (s^{j_a})}{\ell}\frac{|G_{x,o}|}{|G_{x,a}|}=\sum_{a=1}^t \ov{w}_2(\ov{x}_a)\frac{|G_{x,o}|}{|G_{x,a}|},
\endaligned
\eEqu
where $G_{x,o}\!<\!G_x$ is the isotropy subgroup of some element of $o$. Moving $|G_{x,o}|$ in (\ref{equ:compare-w}) to the left-hand side, together with (\ref{equ:w-decomp}), we obtain
$$
\frac{\ov{w}_1(\ov{x}_o)}{\mf{I}_1(\ov{x}_o)}= \frac{\ov{w}_1(\ov{x}_o)}{|G_{x,o}|}= \sum_{a=1}^t \frac{\ov{w}_2(\ov{x}_a)}{|G_{x,a}|}=\lrp{(\ov{q}_{2,1})_*\lrp{\frac{\ov{w_2}}{\mf{I}_2}}}(\ov{x}_o);
$$
this finishes the proof of (\ref{equ:weight-relation}).

\eProof

\noindent
We have thus proved that every multisection $\mfs$ with compact support admits a resolution in the following sense.

\bDef{def:resolution-conditions}
Let $\cE\!\lra\!\cM$ be an orbibundle and $\mfs$ be a lifted multisection with compact support.
A resolution of $\mfs$ consist of a set of $c_i$-covering maps
\bEqu{equ:seq-covering}
q_i\colon\wt\cW_i\lra \cW_i,\quad i=1,\ldots,N,
\eEqu
over open subsets of $\cM$ such that the following properties hold.
\bEnum
\item We have 
$$
Z(\mfs)\subset \bigcup_{i=1}^N W_i, \quad W_t\cap M(>\!c_t)\!=\!\emptyset,\quad \tn{and}\quad  M(\geq\! c_t)\subset \bigcup_{i=t}^N W_i\quad \forall t\!\in\![N].
$$
\item For every $i\!\in\! [N]$, the multisection $\mfs$ restricted to $W_i$ is equivalent to the weighted push forward with weight $w_i$ of some orbifold section $s_i$ of $\wt\cE_i$. 
\item\label{l:Wij-qij} On the overlaps 
$$
W_{i,j}=W_{j,i}= W_i\cap W_j \quad  \forall 1\!\leq\! j\!<\!i\!\leq\! N,
$$ 
we have $q_i=q_j\circ q_{i,j}$, where $q_{i,j}$ is a component-wise covering map
$$
q_{i,j}\colon \wt\cW_i|_{\ov{q}_i^{-1}(W_{i,j})}\lra \wt\cW_j|_{\ov{q}_j^{-1}(W_{i,j})};
$$
moreover,
\bEqu{equ:third-condition}
q_{i,j}\circ q_{j,k}\!=\! q_{i,k}\quad\mbox{on}\quad W_i\cap W_j\cap W_k \quad \forall 1\!\leq \!k\!<\!j\!<\!i\!\leq \!N. 
\eEqu
\item For every  $1\!\leq\! j\!<\!i\leq\! N$, 
\bEqu{equ:wIiwIj}
(\ov{w}_j/\mf{I}_{\cW_j})|_{\ov{q}_j^{-1}(W_{i,j})}= (\ov{q}_{i,j})_* (\ov{w}_i/\mf{I}_{\cW_i})|_{\ov{q}_i^{-1}(W_{i,j})}.
\eEqu
\eEnum
\eDef

\bRem{rem:manifold-covering}
For a generic choice of $\mfs$, the orbifolds $\wt\cW_i$ are all manifolds. In fact, this happens if for every sufficiently small orbibundle chart $(\pr\colon \!U\!\lra\! V,G,\psi)$, in the local presentation 
$$
\mfs|_V=[s^1,\ldots,s^\ell], \quad s^i\colon V\lra U~~\forall~i\in [\ell],
$$
none of the $s^i$ is equivariant with respect to some subgroup of $G$; see Equation~\ref{equ:comp-o}. This is the case if we choose $\ell$ large enough and the branches $s^i$ generically\footnote{We need to consider a sufficiently large $\ell$, with respect to $|G|$, to achieve this.}.
In this situation, the triangulation argument in the next section simply involves manifolds, vector bundles, and triangulation of smooth submanifolds, which is standard. For the specific type of transverse perturbations constructed in Section~\ref{sec:perturbation}, if we choose the perturbation terms $\nu_\al$ generically, the resulting multisection has this property.
We continue without this assumption so that our formula for the resulting rational cycle applies to every arbitrary choice of transversal multisection in the upcoming examples.
\eRem

\subsection{Euler class}\label{sec:euler}
\noindent
Let $E\!\lra\! M$ be a relatively oriented vector bundle over a manifold without boundary and $s\colon \!M\!\lra\! E$ be a transverse section with compact support. Then $s^{-1}(0)\!\subset\! M$ is a compact oriented submanifold of $M$; thus it defines a singular homology class in 
$$
H_{\dim M-\tn{rank} E}(M,\Z).
$$
Two such sections $s$ and $s'$ are called deformation equivalent if there exists a one-parameter family of sections $(s_t)_{t\in [0,1]}$ with compact support such that $s\!=\!s_0$ and $s'\!=\!s_1$. Let $\hb$ be a deformation equivalence class of sections of $E$ with compact support and $s\!\in\!\hb$ be transverse. Then the homology class of $s^{-1}(0)$, which we call the Euler class of the pair $(E,\hb)$ and denote it by $e(E,\hb)$, is a homological invariant of the pair $(E,\hb)$. If $M$ is oriented, this homology class is  Poincar\'e dual to a cohomology class with compact support
$$
\tn{PD}(e(E,\hb))\in H^{\dim M-\tn{rank} E}_c(M).
$$  
\noindent
\bExa{exa:Mathoverflow}
Let $E\!\lra M$ be the trivial complex line bundle over the open unit disk in the plane.  For every $k\geq 0$, let $\hb_k$ be the homotopy class of the sections with compact support containing $s(z)\!=\!z^k$. Then
$$
 \langle e_k, M\rangle =k,\quad e_k=\tn{PD}(e(E,\hb_k))\in H^2_c(M).
$$
\eExa

\noindent
Similarly, after some simple generalization of this argument, if $s\colon \cM \lra \cE$ is a transverse orbifold section of a relatively oriented orbibundle, $\ov{s}^{-1}(0)\!\subset\! M$ inherits a possibly non-effective oriented orbifold structure from the pair $(\cE,\cM)$ and its homology class, with some $\Q$-coefficients that come from the orbifold structure,  defines a homology class in $H_{\dim M-\tn{rank} E}(M,\Q)$ invariant under deformation of $s$.

\noindent
As we argued before, it is quite likely that an orbibundle does not admit any transverse orbifold section; thus this construction does now work in all cases. In the light of Proposition~\ref{pro:estimate}, we consider a transversal multisection deformation $\mfs$ of $s$ and  construct an Euler class via the zero set $Z(\mfs)$. To this end, we use a resolution of the multisection to reduce the construction back to the case of orbibundles and orbifold sections. 

\noindent
In the context of Gromov-Witten theory, an enhanced version of the Euler class for the obstruction bundle of an associated Kuranishi structure defines the Gromov-Witten VFC of the underlying moduli space.\\

\noindent
Recall that a triangulation of some topological subspace $Z\!\subset\! M$ is a  simplicial complex $\De$, 
together with a homeomorphism $\si \colon\De\! \lra\! Z$.
A pure simplicial $n$-complex $\De$ is a simplicial complex where every simplex of dimension
less than $n$ is a face of some $n$-simplex in $\De$.
If $\De$ is a cycle and if the top simplexes  of $\De$ are compatibly oriented, then $\si(\De)$ defines a class in $H_n(M,\Z)$.
In general, the zero set $Z(\mfs)$ (\ref{equ:sero-set}) of a transverse multisection $\mfs$ may not admit any triangulation because different branches may intersect wildly. 
However, in what follows, we consider  a ``triangulation" of the zero set of the sections of some resolution of $\mfs$ and push that into $Z(\mfs)$, with some rational weights, to obtain the desired singular $n$-cycle with $\Q$-coefficients.
Our construction is inductive, explicit, but somewhat sketchy at certain points; for a more functorial and detailed construction, we refer to the definition and construction of branched coverings in \cite{Mc2006}. 

\vskip.1in
\noindent
For the rest of this section, we fix a relatively oriented orbibundle $\cE\!\lra\!\cM$, a transverse orbifold multisection with compact support $\mfs$, and a resolution 
$$
\{q_i\colon\wt\cW_i\lra \cW_i,\quad s_i\colon\wt\cW_i\lra \wt\cE_i\}_{i\in [N]}
$$
of $\mfs$ as in Definition~\ref{def:resolution-conditions}.
Since $\mfs$ is assumed to be transverse, for each $i\!=\!1,\ldots,N$, the orbifold section $s_i$
is transverse to the zero section. Let 
\bEqu{equ:zero-orbifold}
\wt{Z}_i\!=\!Z(s_i)\!\subset\! \wt{W}_i
\eEqu
be the zero set of $s_i$. If $\wt\cW_i\!=\!\wt{W}_i$ is a manifold, see Remark~\ref{rem:manifold-covering}, $\wt{Z}_i$ is an oriented submanifold of $\wt{W}_i$ and can be triangulated.
More generally, in the orbifold case, over each local chart $(V,G,\psi)$ of $\wt\cW_i$, since 
$$
{s}_{i,V}\!=\! {s}_i|_V\colon V\lra U
$$ 
is $G$-equivariant and transverse to the zero section, ${s}_{i,V}^{-1}(0)$ is a $G$-invariant $n$-dimensional oriented submanifold of $V$. The $G$-action on (any component of) ${s}_{i,V}^{-1}(0)$ can be non-effective, i.e. the orbifold structure $\wt\cW_i$ induces a possibly non-effective orbifold structure $\wt\cZ_i$ on $\wt{Z}_i$. More precisely, for each connected component ${s}_{i,V}^{-1}(0)_o$ of ${s}_{i,V}^{-1}(0)$, let $G_o$ be the isotropy group of a generic point\footnote{If ${s}_{i,V}^{-1}(0)_o$ is a single point $\{y\}$, then $G_o\!=\!G_{y}$. } $y\!\in\! {s}_{i,V}^{-1}(0)_o$. This is well-defined because by the orientability of the action on ${s}_{i,V}^{-1}(0)_o$, the isotropy group is fixed outside a set of codimension $2$ in ${s}_{i,V}^{-1}(0)_o$; see Remark~\ref{rem:final-remarks}\ref{l:singM}.
In other words,
$G_o\!<\!G$ is the largest subgroup acting as identity on ${s}_{i,V}^{-1}(0)_o$.  Thus, $G_o\!<\! G$ is normal and the action of $G/G_o$ on ${s}_{i,V}^{-1}(0)_o$ is effective. 
Therefore, 
\bEqu{equ:effectivization}
[{s}_{i,V}^{-1}(0)_o/(G/G_o)]
\eEqu 
is an effective orbifold chart on $\wt{Z}_i$. We denote this reduced effective orbifold structure by $\wt{\cZ}_i^{\tn{red}}$.
Globally, for each connected component $\wt{Z}_{i,o}$ of $\wt{Z}_i$, let $\tn{I}_o$ be the conjugacy class of the isotropy group of a generic point of $\wt{Z}_{i,o}$.  We conclude that $\wt\cZ_{i,o}$ is possibly a non-effective orbifold with global isotropy group $\tn{I}_o$ such that 
\bEqu{equ:reduced-structure}
\wt{\cZ}_{i,o}^{\tn{red}}=\wt{\cZ}_{i,o}/\tn{I}_o.
\eEqu
For this reason, in the construction of Euler class rational cycle below, the isotropy group $\tn{I}_o$ contributes by $\frac{1}{|\tn{I}_o|}$ to the weight  of each simplex mapped into $ \wt{Z}_{i,o} $. 

\noindent
Similar to case of manifolds (see \cite{Wh} for the manifold case) every orbifold admits a ``good" triangulation in the following sense\footnote{We do not prove the existence of such triangulation in this article. We do not know of a reference that discusses the construction either. However, existence of such triangulation is a straightforward generalization of \cite{Wh}.  As Remark~\ref{rem:manifold-covering} indicates, we can avoid orbifolds by considering a generic choice of transverse multisection.}. 

\bDef{def:good-tri}
A \textbf{good} triangulation of an $n$-dimensional orbifold\footnote{For this definition, the orbifold structure does not really need to be effective} $\cM$ is a triangulation  $\ov\si\colon \De \lra M$ with a pure simplicial $n$-complex $\De=\{\De^b\}_{b\in \Om}$ such that the following conditions hold.
\bEnum
\item For every $n$-simplex $\De^b\in \De$, there exists an orbifold chart $(V,G,\psi)$ of 
 $\cM$ such that $\ov\si(\De^b)\subset \psi(V)$ and $\ov\si|_{\De^b}$ lifts to a continous embedding $\si\colon\De^b\lra V$.
 \item The interior of every subface $\de$ of every simplex $\De^b$ is mapped into a strata $M^{H}$ of $M$; see Remark~\ref{rem:final-remarks}\ref{l:OrbiStrata}.
\eEnum
\eDef

\noindent
If $\cM$ is oriented, by the second condition in Definition~\ref{def:good-tri}, the interior of every $n$-dimensional simplex and the interior of every $(n-1)$-dimensional face of that are mapped into the smooth locus $M^{\tn{sm}}$ (which is a manifold); see Remark~\ref{rem:final-remarks}\ref{l:singM}. Therefore, up to codimension 2 faces, a triangulation of $\cM$ gives a triangulation of the smooth manifold $M^{\tn{sm}}$ and this is enough to conclude that $\ov\si\colon \De \lra M$ defines a cycle. This cycle is the fundamental cycle of $\cM$.
If $q\colon \cM_1\!\lra\!\cM_2$ is a $c$-covering, it is clear from Definition~\ref{def:k-cover} that a good triangulation of $\cM_2$, after possibly subdividing into smaller triangles\footnote{So that the image of each simplex sits inside an orbifold chart as in Definition~\ref{def:k-cover}}, canonically lifts to a good triangulation of $\cM_1$.

\bDef{def:triangulation}
With notation as above, given a resolution 
$$
\mf{R}\equiv \lrp{
\{q_i\colon \wt\cW_i\lra \cW_i\}_{i\in [N]}, \{q_{i,j}\colon \wt\cW_i|_{\ov{q}_i^{-1}(W_{i,j})}\lra \wt\cW_j|_{\ov{q}_j^{-1}(W_{i,j})}\}_{1\leq j < i \leq N}
}
$$
of $\mfs$ as in Definition~\ref{def:resolution-conditions}, an \textbf{admissible} triangulation of $\mf{R}$ consists of a good triangulation 
$$
\ov\si_i\colon \De_i \lra T_i$$ 
of some compact subset $T_i\!\subset \!\wt{Z}_i$, for every $i\!=\!1,\ldots, N$, where $\De_i$ is a pure simplicial $n$-complex $\De_i=\{\De_{i}^b\}_{b\in \Om_i}$, such that the following conditions hold.
\bEnum
\item The union of images of $T_j$ in each $\wt{Z}_i$ cover $\wt{Z}_i$, i.e.  for every $1\!\leq \!i\! \leq \!N$, 
$$
T_i \cup \bigcup_{j>i} \ov{q}_{j,i}(T_j\cap \ov{q}_j^{-1}(W_{i,j})) \cup \bigcup_{j<i} \ov{q}_{i,j}^{-1}(T_j)= \wt{Z}_i.
$$
\item For every $i,j\!\in\![N]$, with $i\!\neq\! j$, there exist a sub-triangulation
$$ 
\ov\si_{i}\colon \De_{i,j}\lra T_{i,j} \subset T_i
$$
for some pure simplicial $n$-complex $\De_{i,j}=\{\De_{i,j}^b\}_{b\in \Om_{i,j}}$ such that
\bEnum
\item $T_{i,j}\subset \ov{q}_i^{-1}(W_{i,j})$,
\item  $\ov{q}_{i,j}^{-1}(T_{j,i})=T_{i,j}$, i.e. $T_{i,j}$ is the lift of $T_{j,i}$ via $\ov{q}_{i,j}$, whenever $i\!>\!j$,
\item  and 
$$
\ov{q}_{i,j}(T_i \cap \ov{q}_j^{-1}(W_{i,j}))\cap T_j= T_{j,i}\quad \forall i\!>\!j,
$$ 
i.e. the overlap of $T_i$ and $T_j$ via $\ov{q}_{i,j}$ is equal to $T_{i,j}$ and $T_{j,i}$ in the corresponding spaces. 
\eEnum

\eEnum
\eDef
\noindent
Note that we are considering a triangulation of $\wt{Z}_i$ and not the whole space $\wt{W}_i$. However, we can assume that the triangulation of each $\wt{Z}_i$ is the restriction of some good triangulation of $\wt{W}_i$ and these extended triangulations are compatible on the overlaps as well (as in Definition~\ref{def:triangulation}). 
We do not go into the details of existence of such an admissible triangulation in this article. The proof of existence is by induction on $i\!=\!1,\ldots,N$, and uses the fact that every $\wt{Z}_i$ is triangulatable and $M$ is metrizable. For the induction, we would start from $i\!=\!1$ (the smallest value) and consider a triangulation of $\wt{Z}_1$. After passing to some barycentric subdivision of the triangulation, we can compatibly extend it to the next level $\wt{Z}_2$ and continue until  $i\!=\!N$.

\vskip.1in
\noindent
Given an admissible triangulation of a resolution of $\mfs$, as above, we construct the Euler class singular $n$-cycle of $(\cE,\mfs)$ in the following way. 

\noindent
For every $i\!=\!1,\ldots,N$, and every $n$-simplex $\De_i^b \in \De_i$, we say $\De_i^b$ is \textbf{primary}, if there is no $j\!<\!i$ such that $\De_i^b \in T_{i,j}$. Let $\De^{\tn{prim}}_i\subset \De_i$ be the subset of primary $n$-simplexes and 
$$
\De^{\tn{prim}}=\coprod_{i=1}^N\De^{\tn{prim}}_i.
$$
For every $\De_i^b \!\in\! \De^{\tn{prim}}_i$, by Definition~\ref{def:good-tri}, there exists an orbibundle chart 
$$
(U,V,G,\psi,\pr)
$$ 
of 
 $\wt\cE_i\lra \wt\cW_i$ such that $\ov\si_i(\De_i^b)\subset \psi(V)$ and $\ov\si_i|_{\De_i^b}$ lifts to 
 \bEqu{equ:lifted-simplex}
 \si_i\colon\De_i^b\lra V
 \eEqu
By assumption, the restriction ${s}_{i,V}\colon V\!\lra\! U$ of the orbifold section ${s}_i$ on $\wt\cW_i$ to $V$ is transverse to the zero section. Since the bundle $U\!\lra\! V$ is relatively oriented, 
${s}_{i,V}^{-1}(0)$ is oriented. We orient $\De_i^b$ such that (\ref{equ:lifted-simplex}) is orientation preserving. 
Let $\mf{I}_i(\De_i^b)$ be equal to $\mf{I}_i(\ov{x})$ for some generic $\ov{x}\!\in \!\ov\si_i(\De_i^b)$. In other words, $\mf{I}_i(\De_i^b)$ is equal to the order of isotropy group of the connected component $\wt{Z}_{i,o}$ containing  $\ov\si_i(\De_i^b)$; see (\ref{equ:effectivization}), (\ref{equ:reduced-structure}), and the argument there.
Let $\ov{w}_i(\De_i^b)$ be the value of the weight ratio function $\ov{w}_i\colon \wt{W}_i \lra \Q$ along $\ov\si_i(\De_i^b)$.

\bPro{pro:chain-cycle}
With notation as in Lemma~\ref{lem:weight-relation} and Definition~\ref{def:triangulation},
the singular $n$-chain 
\bEqu{equ:Euler-class}
e(\cE,\mfs)=\sum_{\De_i^b\in \De^{\tn{prim}}} \frac{\ov{w}_i(\De_i^b)}{\mf{I}_i(\De_i^b)} ~ \ov{q}_i( \ov\si_i(\De_i^b))\subset Z(\mfs)\subset M
\eEqu
is a rational $n$-cycle. Moreover, its homology class is independent of the choice of resolution, admissible triangulation, and deformation equivalence class of the transversal multisection $\mfs$; thus we denote it by $e(\cE,\mfs)$.
\ePro

\bExa{exa:Euler-of-global}
Let $M$ be a connected closed manifold, $G=\{g_1,\ldots,g_c\}$ be a finite subgroup of diffeomorphisms of $M$,  and $\cM=[M/ G]$ be the resulting effective  global quotient orbifold. In this case, the manifold $M$ is a $c$-covering of the orbifold $\cM$ and $T\cM\lra \cM$ is relatively oriented. A generic\footnote{If $s$ is not generic, $\tn{val}_\mfs$ may not be a constant function in which case $M$ is not a resolution of $\mfs$ as we built above.} transverse section $s$ of $TM$ gives a push forward transverse multisection 
$$
\mfs= [g_1\cdot {s}, \ldots, g_c\cdot {s}]
$$ 
of $T\cM$ for which $\tn{val}_\mfs\!\equiv\! c$.
Let $s^{-1}(0)\!=\!\{p_1,\cdots, p_N\}$. For every $i\!\in\! [N]$, let $\ep_i\in \{\pm 1\}$ be the intersection number of $p_i$.  The weight ratio of every $p_i$ is simply $\frac{1}{|G|}$ and $\mf{I}_{M}\colon M\!\lra\! \Q$ is the constant function $1$. Then, by (\ref{equ:Euler-class}), the orbifold Euler characteristic (class) of $\cM$ is equal to
$$
\chi(\cM)=e(T\cM)=\sum_{i=1}^\ell  \frac{\ep_i}{|G|}= \frac{\chi(M)}{|G|} \in \Q\cong H_0(M/G,\Q).
$$
\eExa

\bExa{exa:Euler-of-PMN}
With notation as in Example~\ref{exa:Pmn}, since $s$ is already a transverse orbifold section, we do not need to pass to any non-trivial cover. We have $Z(s)=\{0,\infty\}$, isotropy group of $0$ has order $m$, and isotropy group of $\infty$ has order $n$; therefore, by (\ref{equ:Euler-class}), $\chi(T\P^1_{m,n})=\frac{1}{m}+\frac{1}{n}$.
\eExa

\newtheorem*{proofofpro:chain-cycle}{Proof of Proposition~\ref{pro:chain-cycle}}
\begin{proofofpro:chain-cycle}
We just show that (\ref{equ:Euler-class}) is a cycle. Proof of independence of its homology class from various choices is by constructing a cobordism between different choices and we skip that.

\noindent
Assume $n\!=\!\dim \cM\!-\!\tn{rank}~\cE\!>\!0$, otherwise, $\partial e(\cE,\mfs)\!=\!0$ for dimensional reason. Let  $\partial \De^{\tn{prim}}$ and $\partial \De^{\tn{prim}}_i$  be the set of all $(n-1)$-dimensional faces -or simply faces- of primary simplexes in $\De^{\tn{prim}}$ and $\De^{\tn{prim}}_i$, respectively.

\bDef{def:int-ext}
We call $\de\!\in\! \partial \De^{\tn{prim}}_i$ \textbf{internal} if there exists another $\de'\!\in\! \partial \De^{\tn{prim}}_i$ such that $\ov\si_i(\de)\!=\!\ov\si_i(\de')$. We call $\de\!\in\! \partial \De^{\tn{prim}}_i$ \textbf{transitive} if there exists $j\!\neq\! i$ and $\de'\!\in\! \partial \De^{\tn{prim}}_j$ such that 
\bEqu{equ:face-relation}
\ov\si_i(\de)=\ov{q}_{j,i} \ov\si_j(\de')~~ \mbox{if } j\!>\!i,~\tn{or}\quad\ov\si_j(\de)=\ov{q}_{i,j} \ov\si_i(\de')~~ \tn{if } i\!>\!j;
\eEqu
see Figure~\ref{fig:PFT}. In this situation, $j\!\in\! [N]$ is called an \textbf{associated} index to $\de$.
\eDef

\begin{figure}
\begin{pspicture}(-6,-6.5)(10,-4)
\psset{unit=.3cm}
\psline(-8,-16)(0,-16)  
\psline(-7.5,-16.2)(-7.5,-15.8)
\psline(-5.5,-16.2)(-5.5,-15.8)
\psline(-3.5,-16.2)(-3.5,-15.8)
\psline(-1.5,-16.2)(-1.5,-15.8)
\psline(-4,-15)(0,-15)
\psline(-4,-17)(0,-17)
\psline(-1.5,-17.2)(-1.5,-16.8)
\psline(-1.5,-15.2)(-1.5,-14.8)
\psline(-3.5,-17.2)(-3.5,-16.8)
\psline(-3.5,-15.2)(-3.5,-14.8)
\psline(0.4,-17.2)(0.4,-16.8)
\psline(0.4,-15.2)(0.4,-14.8)

\psarc(0,-13){2}{-90}{-45}
\psarc(2.8,-15.86){2}{90}{135}
\psarc(0,-19){2}{45}{90}
\psarc(2.8,-16.14){2}{-135}{-90}
\psline(2.8,-18.14)(5.5,-18.14)
\psline(2.8,-18.34)(2.8,-17.96)\psline(5.2,-18.34)(5.2,-17.96)
\psline(2.8,-14.04)(2.8,-13.66)\psline(5.2,-14.04)(5.2,-13.66)
\psline(2.8,-13.86)(5.5,-13.86)
\psline(-8,-20)(5.5,-20)
\psline{->}(-5,-17.3)(-5,-19.3)
\rput(-4,-18.3){$\ov{q}_1$}
\psline{->}(4.5,-18.5)(4.5,-19.7)
\rput(3.5,-19){$\ov{q}_2$}
\psarc[linestyle=dashed,dash=1pt]{->}(-5,-15.8){2}{0}{70}
\psline[linestyle=dashed,dash=1pt]{->}(-5,-15.8)(-5,-14)
\psarc[linestyle=dashed,dash=1pt]{<-}(-5,-15.8){2}{110}{180}
\rput(-6,-13){\tiny{Elements of $\De_1^{\tn{prim}}$}}

\psarc[linestyle=dashed,dash=1pt]{->}(2,-13.7){2}{0}{70}
\psline[linestyle=dashed,dash=1pt]{->}(2,-13.7)(2,-11.7)
\psarc[linestyle=dashed,dash=1pt]{<-}(3.5,-14.8){4}{130}{180}
\rput(2,-11){\tiny{Elements of $\De_2^{\tn{prim}}$}}

\psline[linestyle=dashed,dash=1pt]{->}(-7.5,-18)(-7.5,-16.3)
\psline[linestyle=dashed,dash=1pt]{->}(-7.5,-18)(-5.5,-16.3)
\psline[linestyle=dashed,dash=1pt]{->}(-7.5,-18)(-3.5,-16.3)
\rput(-10.5,-18.7){\tiny{Internal faces in $\partial \De_1^{\tn{prim}}$}}

\psline[linestyle=dashed,dash=1pt]{->}(5.3,-16)(5.2,-14.2)
\psline[linestyle=dashed,dash=1pt]{->}(5.3,-16)(2.8,-14.2)
\psline[linestyle=dashed,dash=1pt]{->}(5.3,-16)(0.6,-15.1)
\rput(7.3,-16.5){\tiny{Internal faces in $\partial \De_2^{\tn{prim}}$}}

\psellipse[linestyle=dashed,dash=1pt](-1.5,-16.2)(.5,2)
\psline[linestyle=dashed,dash=1pt]{->}(-1.5,-18.2)(-1.5,-21)
\rput(-1.5,-21.5){\tiny{Transitive faces}}

\end{pspicture}
\caption{A schematic description of primary simplexes, internal faces, and transitive faces.}
\label{fig:PFT}
\end{figure}
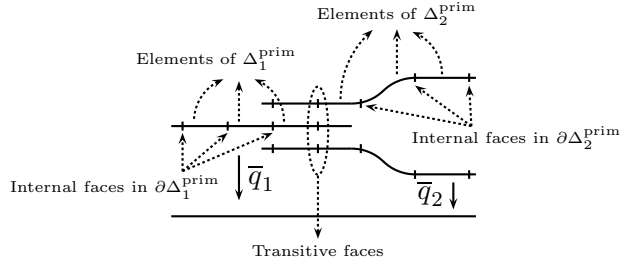

\bLem{lem:unique-pairs}
Every $\de\in \partial \De^{\tn{prim}}$ is either internal or transitive, but not both.
For every internal $\de$, the corresponding face $\de'$ is unique. For every transitive face $\de\in \partial \De^{\tn{prim}}_i$ there is only one index $j\in [N]$ associated to $\de$.
\eLem

\bProof
For every $\De_i^b\! \in\!  \De^{\tn{prim}}_i$ and  every $\de\!\in\! \partial \De_i^b$, let $\de^o$ be the interior of $\de$. By the orientability of $\wt\cZ_i$ and the second condition in Definition~\ref{def:good-tri}, $\ov\si_i(\de^o)$ lies in the smooth locus $\wt{Z}_i^{\tn{sm}}$ of the reduced orbifold structure $\wt\cZ_i^{\tn{red}}$. Therefore, as in the manifold case, either there exists another $\de'\!\in\! \partial \De^{\tn{prim}}_i$ (which could be a different face of the same simplex $\De_i^b$) such that $\ov\si_i(\de')\!=\!\ov\si_i(\de)$, or $\ov\si_i(\De_i^{\tn{prim}})$ does not cover any sufficiently small enough neighborhood of any point $\ov{x}\!\in\! \ov\si_i(\de^o)\in \wt{Z}_i^{\tn{sm}}$.  
In the first case, $\de$ is an internal face. Moreover, similar to the manifold case, since $\wt\cZ_i$ is oriented, the contribution of $\de$ and $\de'$ to $\partial e(\cE,\mfs)$ have equal absolute values and different orientations  (the weight ratio and the order of isotropy group over the simplexes contacting $\de$ and $\de'$ are the same, because they belong to the same connected component of $\wt{Z}_i$); therefore, they cancel out.

\noindent
In the second case, for every $\ov{x}\!\in\! \ov\si_i(\de^o)$, it follows from the first condition in Definition~\ref{def:triangulation} that there exists another $\De_j^a\!\in\! \De^{\tn{prim}}_j$, for some $j\!\in\! [N]$ different from $i$, and $\de'\!\in\! \partial \De_j^a$ such that 
$$
\begin{cases} 
\ov{q}_{i,j}(\ov{x})\in \ov\si_j(\de') & \mbox{if } i>j,\\
\ov{x}\in \ov{q}_{j,i}(\ov\si_j(\de')) & \mbox{if } j>i .
\end{cases}
$$ 
In either case, it follows from Definition~\ref{def:triangulation}.2(c) that $\ov{x}\!\in\! T_{i,j}$. Finally, it follows from Condition 2(b) that
$$
\begin{cases} 
\ov{q}_{i,j}(\ov\si_i(\de))= \ov\si_j(\de') & \mbox{if } i>j,\\
\ov{q}_{j,i}(\ov\si_j(\de'))= \ov\si_i(\de) & \mbox{if } j>i .
\end{cases}
$$ 

\noindent
An internal face can not be transitive because by Definition~\ref{def:triangulation}.2(c), if 
$$
\de,\de'\in \partial\De_i^{\tn{prim}},~~\de\in \partial\De_i^b,~~\de'\in \partial\De_i^{b'},~~\ov\si_i(\de)=\ov\si_i(\de'),~~\tn{and}~\ov\si_i(\de)=\ov{q}_{j,i}\ov\si_j(\de''),
$$ 
for some $\de''\!\in\! \partial\De_j^{b''}$, where $\De_j^{b''}\!\in\!\De_j^{\tn{prim}}$ and $j\!<\!i$, since $\ov\si_i(\de)$ locally divides $\wt{Z}^{\tn{sm}}_i$ into two components, either 
$\ov\si_i(\De_i^b)\!=\!\ov{q}_{j,i}(\ov\si_j(\De_j^{b''}))$ or $\ov\si_i(\De_i^{b'})\!=\!\ov{q}_{j,i}(\ov\si_j(\De_j^{b''}))$. This is a contradiction because it implies that $\De_j^{b''}$ is not primary. The case of $j\!<\!i$ is similar.

\noindent
Suppose there exist $i\!<\!j\!<\!k$ and 
$$
\de\in \partial\De_i^b,\quad \de'\in \partial\De_j^{b'},\quad\tn{and}~~\de\in \partial\De_k^{b''}
$$
where $\De_i^b$, $\De_j^{b'}$, and $\De_k^{b''}$ are primary and 
$$
\ov\si_i(\de)=\ov{q}_{j,i}(\ov\si_j(\de')),\quad \ov\si_i(\de)=\ov{q}_{k,i}(\ov\si_k(\de'')).
$$
First, it follows from (\ref{equ:third-condition}) that 
\bEqu{equ:3rd-one}
\ov\si_j(\de)=\ov{q}_{k,j}(\ov\si_k(\de'')).
\eEqu
Since $\ov\si_i(\de)$ locally divides $\wt{Z}^{\tn{sm}}_i$ into two components, both 
$$
\ov{q}_{j,i}(\ov\si_j(\De_j^{b'}))\quad\tn{and}\quad \ov{q}_{k,i}(\ov\si_k(\De_k^{b''}))
$$
are mapped to the same side of $\ov\si_i(\de)$. It follows from (\ref{equ:3rd-one}) and Definition~\ref{def:triangulation}.2(c) that $\ov\si_j(\De_j^{b'})\!=\!\ov{q}_{k,j}(\ov\si_k(\De_k^{b''}))$. This is a contradiction to the assumption that both $\De_j^{b'}$ and $\De_k^{b''}$ are primary. This establishes the last claim of Lemma~\ref{lem:unique-pairs}.
\eProof
\noindent
From Lemma~\ref{lem:unique-pairs}, we conclude that internal faces come in pairs with identical images and that the contribution of each pair to $\partial e(\cE,\mfs)$ is zero. As for transitive faces, proof of Lemma~\ref{lem:unique-pairs} shows that the elements of $\partial \De^{\tn{prim}}$ can be grouped into tuples of the form
\bEqu{equ:comparable-pairs}
(\de; \de'_1,\ldots,\de'_t), \quad \de\in \partial \De^{\tn{prim}}_i, \quad \de'_1,\ldots,\de'_t\in \partial \De^{\tn{prim}}_j \quad\tn{for some} ~j\!>\!i,
\eEqu
where $j$ is the unique associated index to $\de$ and  $\de'_1,\ldots,\de'_t$ are the set of all faces where (\ref{equ:face-relation}) holds. Finally, by (\ref{equ:weight-relation}), the contribution of
$\de$ and the sum of contributions of $\de'_1,\ldots,\de'_t$ to $\partial e(\cE,\mfs)$ cancel out each other, i.e.
\bEqu{equ:contibutions}
\sum_{a=1}^t \frac{ \ov{w}_j(\De_j^a)}{\mf{I}_j(\De_j^a)}=\frac{\ov{w}_i(\De_i^b)}{\mf{I}_i(\De_i^b)} ,
\eEqu
where $\De_i^b$ is the simplex containing $\de$ and $\De_j^a$, for each $a\!=\!1,\ldots,t$, is the simplex containing $\de'_a$. More precisely, by choosing a generic point $\ov{x}\!\in\! \ov\si_i(\de^o)$ such that $\mf{I}_i(\De_i^b)=\mf{I}_i(\ov{x})$ (this is possible because $\wt{Z}_i\setminus \wt{Z}_i^{\tn{sm}}$ is of codimension at least $2$),
$$
\frac{\ov{w}_i(\De_i^b)}{\mf{I}_i(\De_i^b)} = \frac{\ov{w}_i(\ov{x})}{\mf{I}_i(\ov{x})}=
\sum_{\ov{y}\in \ov{q}_{j,i}^{-1}(\ov{x})} \frac{\ov{w}_j(\ov{y})}{\mf{I}_j(\ov{y})}= \sum_{a=1}^t \frac{ \ov{w}_j(\De_j^a)}{\mf{I}_j(\De_j^a)}.
$$
\qed
\end{proofofpro:chain-cycle}

\section{Abstract Kuranishi structures}\label{sec:prelim-kur}
\subsection{Introductory remarks}\label{Introductory remarks}
\noindent
In differential geometry, various moduli spaces appear as the set of solutions of nonlinear partial differential equations. 
The zero set of the defining equation of the moduli space can be very singular. The idea of using obstruction bundle to resolve this issue goes back to Kuranishi \cite{Ku}:
Let $(X,J)$ be a closed holomorphic manifold with the group of holomorphic diffeomorphisms $\aut(X,J)$. In order to describe the local structure of the moduli space of complex structures on $X$ near $J$, Kuranishi proved that there exist
\bItem
\item a finite dimensional complex manifold $V$ on which $\aut(X,J)$ acts homomorphically, 
\item a complex representation $E$ of $\aut(X,J)$,
\item and an $\aut(X,J)$-equivariant holomorphic map $s\colon V\lra E$,
\eItem   
such that the set of complex manifold $(X,J')$ with $J'$ sufficiently 
close to $J$ has a canonical bijection to the quotient space  $s^{-1}(0)/\aut(X,J)$.
Roughly speaking, we can exchange the possibly singular moduli space $s^{-1}(0)$ with a larger smooth space $V$ and an obstruction space $E$ such that the moduli space is the zero set of a section of the vector bundle $V\!\times\! E$. 
After his work, such a local structure is called a ``Kuranishi model" or a ``Kuranishi chart".
For example, motivated by this construction, in \cite{Tau}, Taubes used this idea to study the space of solutions of the Yang-Mills equations near bubbling in way that led to the modern use.

\noindent
In 1980's, first by Donaldson in gauge theory \cite{D} and then by Gromov \cite{G} in the theory
of pseudoholomorphic curves, the idea of  using the fundamental homology class of moduli spaces to obtain invariants in various topological contexts  was introduced. In the theory of pseudoholomorphic curves, whenever the moduli space in question has a nice orbifold structure, Gromov-Witten invariants were first obtained in this way. For example,
such a theory was first rigorously built in the case of semi-positive symplectic manifolds
by Ruan~\cite{Ru96} and Ruan-Tian \cite{RT} and is now a straight-forward tool in these cases, see McDuff-Salamon \cite{MS2004}.

\noindent
However, moduli spaces of pseudoholomorphic maps of arbitrary genus over arbitrary symplectic manifolds can be quite singular and may have components of larger than expected dimension. In this situation, the moduli space does not have a fundamental homology class and a precise definition of Gromov-Witten invariants, as some sort of intersection theory on the moduli space, needs more complicated tools. In order to solve the issue, we use the idea of Kuranishi to define a well-defined intersection theory on the compactified  moduli space of pseudoholomorphic maps. We cover the moduli space with Kuranishi charts as above; these charts are obtained via finite dimensional reductions of the Cauchy-Riemann equation defining the moduli space and generalize the orbifold structure of the moduli space in the semi-positive cases. A suitably defined ``compatible" covering of the moduli space with Kuranishi charts is called a Kuranishi structure. After adding some extra assumptions to the definition of Kuranishi structure\footnote{e.g. existence of ``tangent bundle", Hausdorffness, $\ldots$}, we define an intersection theory which results in Gromov-Witten invariants and generalizes the definition of these invariants in the semi-positive case.

\noindent
In Section~\ref{sec:basics}, we introduce the notion of abstract Kuranishi chart and Kuranishi structure over a compact metrizable topological space so that a version of basic concepts like tangent space, orientability, and smoothness make sense. In Section~\ref{sec:level}, in order to simplify the calculations, we assemble the uncountable data of Kuranishi structure into a finite set of triples of orbifolds, orbibundles, and orbifold sections, with transition maps from lower rank orbibundles to higher ones. The result is called a dimensionally graded system.
In Section~\ref{sec:def}, similar to Section~\ref{sec:perturbation}, we consistently perturb the orbifold sections of the dimensionally graded system into a consistent set of transversal but multivalued sections of the obstruction bundles. 

\bRem{rem:category}
It is yet unknown how to build a category from the set of spaces with Kuranishi structure; thus, a so called Kuranishi category is the dream of this theory. The problem is that we do not know how to define morphisms. 
By a Kuranishi category, we mean a category that includes morphisms between Kuranishi structures and such that basic operations like fiber product are well-defined.
For some attempts in this direction see \cite{Jo}.
 In algebraic geometry, theory of 2-categories with 2-morphisms and 2-fiber products is one such big category that includes all moduli spaces of holomorphic curves in complex projective varieties. 
\eRem

\subsection{Kuranishi structures}\label{sec:basics}

In this section, we continue to use the notation conventions of Notation~\ref{not:calconv},~\ref{not:phi}, and~\ref{fD_n}.

\bDef{def:kur-chart} 
A \textbf{Kuranishi chart}  of real dimension $n$ on a  compact metrizable topological space $M$ is a tuple $\cU\!\equiv\! (\pr\colon\!U\!\lra\!V,G, s, \psi)$ made up of
\bEnum
\item an orbibundle chart $(\pr\colon\!U\!\lra\!V,G)$ of rank $\dim V \!-\! n$ as\footnote{With the exception of the map $\psi$ which in this case is only defined on the zero set of Kuranishi map.} in Definition~\ref{def:orbibundlechart},
\item\label{l:K-map} a smooth $G$-equivariant section $s\colon \!V\!\lra\! U$, called \textbf{Kuranishi map} or section,
\item  and a continuous map $\psi\colon\! s^{-1}(0)\!\lra\! M$ defined by $\psi\!=\!\ov\psi\circ \pi$, where 
$$
\pi\colon s^{-1}(0)\lra \ov{s}^{-1}(0)= s^{-1}(0)/G
$$ 
is the quotient map and $\ov\psi\colon\! \ov{s}^{-1}(0)\!\lra\! F\!\subset\! M$ is a homeomorphism into an open subset $F$, called \textbf{footprint}, of $M$.
\eEnum
\eDef
\noindent
Similar to Section~\ref{sec:orbifold}, for $x\!\in\! s^{-1}(0)$ and $\ov{x}\!=\!\psi(x)\!\in\! M$, let $G_x$ and $\tn{I}_{\ov{x}}$ denote the isotropy group and the conjugacy class of the isotropy group, respectively. We say $\cU$ is a Kuranishi neighborhood \textbf{centered at} $\ov{x}\!\in\! M$ if $G_x\!=\!G$; in this case there is a unique $x\!\in\! s^{-1}(0)$ such that $\psi(x)\!=\!\ov{x}$. 

\noindent
Given a Kuranishi chart $\cU$ as above, a subgroup $G'\!< \!G$, and a $G'$-invariant connected open subset $V'\!\subset\! V$ such that 
$$
g(V')\cap V'=\emptyset \qquad \forall g\!\in\! G\!\setminus\! G',
$$
by restricting to $V'$, we obtain a Kuranishi chart 
\bEqu{Sub-Chart_dfn}
\cU'\!=\! \big(\pr\colon U'\!=\!\pr^{-1}(V')\lra V',G',s'\!=\!s|_{V'},\psi'\!=\!\psi|_{s^{-1}(0)\cap V'},\big)
\eEqu
which we call a \textbf{sub-chart}. For every $\ov{x}\in \psi(s^{-1}(0))$, by reducing to a sub-chart, we can always obtain a Kuranishi neighborhood centered at $\ov{x}$. Similar to the argument before Example~\ref{exa:OrbiTangent}, every sufficiently small Kuranishi (sub-) chart centered at any $\ov{x}\!\in\! M$ is isomorphic to a linear chart where $U$ is a trivial bundle; this is the point of view in the original approach of \cite{FO}.

\bDef{def:intersection}
Let $\cU_i\!=\! (\pr_i\colon \!U_i\!\lra\!V_i,G_i,s_i,\psi_i)$,  with $i\!=\!1,2$, be a pair of $n$-dimensional Kuranishi charts 
 with intersecting footprints $F_i$, such that 
 $$
 \ov{p}\in F_{12}= F_1\cap F_2~~\tn{and}~~\dim V_1\leq \dim V_2.
 $$ 
A \textbf{coordinate change} from $\cU_1$ to $\cU_2$ \textbf{centered at} $\ov{p}$ consists of a sub-chart\footnote{Any object labeled by $1,2$ is the restriction to $\cU_{1,2}$ of the corresponding object in $\cU_1$.} $\cU_{1,2}\!\subset\! \cU_1$  
centered at $\ov{p}$ and an orbibundle embedding 
$$
\Phi_{12}\equiv(\mfD\phi_{12}, \phi_{12})\colon \cU_{1,2}\lra\cU_2
$$
as in Definition~\ref{def:orbibundle-maps} such that  
\bEqu{equ:cond-on-spsih}
\mfD\phi_{12}\circ s_{1,2} = s_{2} \circ \phi_{12}, \qquad  \psi_{1,2}=\psi_{2} \circ \phi_{12}.
\eEqu 
\eDef
\noindent
Recall that by definition of orbibundle embedding, the group homomorphism 
\bEqu{equ:isom-stabilizers}
h_{12}\colon G_{1,2}\!=\!(G_{1,2})_{p}\lra (G_2)_{\phi(p)}
\eEqu
that comes with $\phi_{12}$ (see Notation~\ref{not:phi}) is an isomorphism.

\noindent
\bExa{exa:simple-kur}
Let $M\!=\!(0,3)$ and cover $M$  with two $1$-dimensional Kuranishi charts $\cU_1$ and $\cU_2$ such that  $U_2\!=\!(0,2)\!\times\! (-1,1) \!\times\! \R$, $V_2\!=\!(0,2)\!\times\! (-1,1)$, $\pr_2$ is the projection onto the first two components, $s_2(x,y)\!=\!y$, $V_1\!=\!(1,3)$, $U_1$ is trivial, both isotropy groups are trivial, and  $\psi_1$ and $\psi_2$ are the obvious maps to $(0,2)$ and $(1,3)$, respectively. 
In this situation, for any $x\!\in\! (1,2)$, by restricting $\cU_1$  to $(1,2)$, we obtain a coordinate change centered at $x$ from $U_1$ to $U_2$. 
\noindent
If we replace the trivial group $G_2$ with $G_2\!=\!\Z_2$, acting on $U_2$ by $(x,y,z)\lra (x,-y,-z)$, the resulting charts do not admit a coordinate change anymore; the isomorphism condition (\ref{equ:isom-stabilizers}) can never be satisfied. 
\eExa
\noindent
Given  a coordinate change as above, since (\ref{equ:isom-stabilizers}) is an isomorphism and $\phi_{12}$ is embedding, we conclude that for every $x_1\!\in\! s_{1,2}^{-1}(0)$ and $x_2\!\in\! s_2^{-1}(0)$ with the same image $\ov{x}$ in $M$, $h_{12}$ induces an isomorphism between the isotropy groups $(G_1)_{x_1}$ and $(G_2)_{x_2}$. Therefore, similarly to orbifolds, the conjugacy class of the isotropy group $\tn{I}_{\ov{x}}$ is well-defined.


\bDef{def:cocycle}
Let $\cU_i$, with $i\!=\!1,2,3$, be a triple of $n$-dimensional Kuranishi charts centered at $\ov{p}_i\!\in\! M$ with intersecting footprints $F_i$, such that 
$$
\ov{p}_1\!\in\! F_{123}\!=\! F_1\!\cap\! F_2\!\cap F_3,\quad\ov{p}_2\!\in\! F_{23}\!= \!F_2\cap F_3,\quad\tn{and}~~\dim V_1\!\leq\! \dim V_2 \!\leq\! \dim V_3.
$$
Given a triple set of coordinate change maps $\Phi_{ij}\colon \cU_{i,j}\! \lra\! \cU_{j}$ centered at $\ov{p}_i$, 
with $1\!\leq\! i\!<\!j\leq\!3$, we say they satisfy \textbf{cocycle condition} if 
\bEqu{equ:cocycle-condition}
\big(\Phi_{13}\!=\!\Phi_{23}\!\circ\!\Phi_{12}\big)|_{\ov{V}_{1,2,3}},\quad\tn{where}~~ \ov{V}_{1,2,3}\!= \!\ov\phi_{12}^{-1}(\ov{V}_{2,3})\!\cap\! \ov{V}_{1,3}\!\subset\! \ov{V}_{1,2}.
\eEqu
\eDef

\noindent
More explicitly, the equality (\ref{equ:cocycle-condition}) means that for every connected component $V_{1,2,3;\al}$ of 
$
V_{1,2,3}\!= \!\phi_{12}^{-1}(V_{2,3})\!\cap\! V_{1,3}\!\subset\! V_{1,2},
$
there exists $g_\al  \!\in\! G_3$ such that
\bEnum
\item $\phi_{13}|_{V_{1,2,3;\al}}\!= \!g_\al\circ \phi_{23}\circ\phi_{12} |_{V_{1,2,3;\al}}$,
\item $\mfD\phi_{13}|_{U_{1,2,3;\al}}\!=\! g_\al\circ \mfD\phi_{23}\circ \mfD\phi_{12} |_{U_{1,2,3;\al}}$, where $U_{1,2,3;\al}\!=\!\pr_1^{-1}(V_{1,2,3;\al})$,
\item and $\tn{ad}(g_\al)\circ h_{13}\!=\!h_{23}\circ h_{12}$, where for every $g\!\in \!G_3$, $\tn{ad}(g)\colon G_3\!\lra\! G_3$ is the adjoint action $ g' \lra g g' g^{-1}$.
\eEnum

\bRem{rem:necesity}
Since the coordinate change  maps commute with the section and footprint maps,
the composite maps on the right-hand sides of the first and second condition
and the coordinate change map $\Phi_{13}$ automatically satisfy the cocycle condition along the graph of $s_1|_{V_{1,2,3}}$.
Otherwise, the cocycle condition away from the graph of $s_1$ has to be considered as a separate condition\footnote{Because we will  perturb the Kuranishi maps and we need to preserve their compatibility throughout the deformation process.}. In the definition of orbifolds, where there is no Kuranishi map, the cocycle condition is automatic. Also, in the natural construction of Kuranishi structure over moduli space of pseudoholomorphic maps in Section~\ref{sec:natural-KUR}, the cocycle condition naturally follows from the construction.
\eRem

\bDef{def:Kur-structure}
A Kuranishi structure $\cK$ of real dimension $n$ on a compact metrizable topological space $M$ consists of  a Kuranishi chart $\cU_{p}$ of real dimension $n$ centered at $\ov{p}$ with footprint $F_{p}$, for every $\ov{p}\!\in\! M$, and a coordinate change map $\Phi_{qp}\colon\cU_{q,p}\lra \cU_{p}$, for every $\ov{q}\in F_{p}$, 
such that whenever $\ov{r}\!\in\! F_{p}\cap F_{q}$ and $\ov{q}\!\in\! F_{p}$, then $\Phi_{qp}$, $\Phi_{rq}$, and $\Phi_{rp}$ satisfy the cocycle condition.
\eDef
\noindent
\begin{notation}
Ideally, we should label the coordinate charts $\cU_p$ as $\cU_{\ov{p}}$. However, since by definition the preimage $p\!=\!\psi^{-1}(\ov{p})$ is unique, we label the Kuranishi charts and the coordinate change maps between them by $p,q, r,$ etc, instead of $\ov{p}, \ov{q},\ov{r}$, etc, respectively.
\end{notation}

\noindent
Given a Kuranishi structure
$$
 \cK\!=\!\big\{(\cU_{p})_{\ov{p}\in M}, (\Phi_{qp}\colon\cU_{q,p}\lra \cU_{p})_{\ov{p}\in M,\ov{q}\in F_p}\big\}
 $$ 
 as in Definition~\ref{def:Kur-structure}, define the \textbf{level map}
\bEqu{equ:d-function}
L\colon M \lra \Z^{\geq 0},\qquad L(\ov{p})=\tn{rank}~U_{p}\;.
\eEqu
This is an upper semi-continuous function and gives us a stratification
\bEqu{L-starata}
M\!=\! \coprod_{k\in \Z^{\geq 0}} M(k),\qquad M(k)\!=\!L^{-1}(k),
\eEqu
where each \textbf{level set} $M(k)$ is an open subset of $M(\geq\! k)\!=\!\bigcup_{\ell\geq k} M(\ell)$. Let 
\bEqu{equ:def-of-B}
\mfB= \mfB(M,\cK)\!=\!\{k\in \Z^{\geq 0}\colon~M(k)\!\neq\! \emptyset\}.
\eEqu
It follows from the compactness of $M$ and the upper semi-continuity of $L$ that $\mfB$ is a finite set and $M(\geq\!k)$ is a closed, hence compact, subset of $M$.

\bExa{exa:global-kur}
The simplest examples of (abstract) Kuranishi structures come from orbibundles over orbifolds. 
Given a non-negative integer $n$, let $\pr\colon\!\cE\!\lra \!\cY$ be an orbibundle of rank $\dim \cY\!-\!n$ and $s\colon\!\cY\!\lra \!\cE$ be an orbifold section with compact zero locus $M\!=\! \ov{s}^{-1}(0)\!\subset\! Y$. For every $\ov{x}\!\in\! M$, we can take the $n$-dimensional Kuranishi neighborhood $\cU_x$ centered at $\ov{x}$ to be an orbibundle chart of $\cE$ centered at $\ov{x}$ together with the Kuranishi map $s$. The coordinate change maps are simply the refinements and the cocycle condition is a consequence of the compatibility condition of orbibundle charts. In this situation, we say $M$ has a \textbf{pure orbibundle} Kuranishi structure. 
We may also consider a disjoint union 
$$
\coprod_{i=1}^m(\pr_i\colon\cE_i\lra\cY_i)
$$ 
of orbibundles of rank 
$$
\tn{rank}~\cE_i=\dim \cM_i-n.
$$
Given a set of orbifold sections $\{s_i\}_{i=1}^m$ as above, the zero set 
$$
M\!=\!\coprod_{i=1}^m \big(M_i\equiv\ov{s}_i^{-1}(0)\big)
$$ 
inherits a Kuranishi structure as above that we still denote by a pure orbibundle Kuranishi structure.
The examples in Section~\ref{sec:examples} are moduli spaces of pseudoholomorphic maps which admit Kuranishi structures of this type. 
\eExa
\noindent
In the next section, using the level function $L$ above,  we assemble the local charts of an arbitrary Kuranishi structure into a finite union of pure orbibundle Kuranishi structures compatible along intersections. We call the result a \textbf{dimensionally graded system}, or \textbf{DGS} for short, associated to the Kuranishi structure. In \cite{FO} and the related followups, this is denoted by ``good coordinate system" or ``dimensionally graded good coordinate system". In an DGS, the pure orbibundles are indexed by their rank, while in good coordinate systems, they are indexed by some partially ordered set; thus, as D. Joyce \cite[Definition 3.5]{Jo2} names it, an DGS is an \textit{excellent} good coordinate system. This notion of ``reduction" in \cite{MW2} plays the same role as DGS.

\noindent
For every arbitrary smooth vector bundle $\pr\colon U\!\lra\! V$ and a smooth section $s$, we can canonically decompose the tangent bundle of $U$ along the zero set of $s$ to 
\bEqu{decomp}
TU|_{s^{-1}(0)}= U|_{s^{-1}(0)}\oplus TV|_{s^{-1}(0)}.
\eEqu
With respect to this decomposition, and by projection onto the the first factor in the right-hand side of (\ref{decomp}), we obtain the linearization map
$$
\nd s\colon TV|_{s^{-1}(0)} \lra U_{s^{-1}(0)}.
$$
Therefore, for a coordinate change embedding map as in Definition~\ref{def:intersection}, we obtain a commutative diagram
\bEqu{comm-tangent}
\xymatrix{
0\ar[r] & TV_{1,2}|_{s_{1,2}^{-1}(0)}\ar@{^{(}->}[r]^{d\phi_{12}} \ar[d]^{\nd s_{1,2}} & TV_2 |_{\phi_{12}(s_{1,2}^{-1}(0))}   \ar[r]\ar[d]^{\nd s_2}     & \cN_{\phi_{12}(V_{1,2})}{V_2}|_{\phi_{12}(s_{1,2}^{-1}(0))}  \ar[d]^{\nd s_{2/1}= \nd s_2/\nd s_{1,2}}  \ar[r] & 0\\
0\ar[r] & U_{1,2}|_{s_{1,2}^{-1}(0)}    \ar@{^{(}->}[r]^{\mfD\phi_{12}}    &  U_{2}|_{\phi_{12}(s_{1,2}^{-1}(0))} \ar[r]   &(U_2/\mfD\phi_{12}(U_{1,2}))|_{\phi_{12}(s_{1,2}^{-1}(0))} \ar[r]& 0.
}
\eEqu
where the last vertical map, denoted by $\nd s_{2/1}$, is defined in the following way. 
For every $v\in TV_2|_{\phi_{12}(s_{1,2}^{-1}(0))}$, 
$$
\nd s_2 (v)\in U_2|_{\phi_{12}(s_{1,2}^{-1}(0))}.
$$ 
By the left-hand side of (\ref{equ:cond-on-spsih}), if $v\!\in\! \nd \phi_{12} (TV_{1,2})\subset TV_2$, then 
$$
\nd s_2 (v)\in \mfD\phi_{12}(U_{1,2})\subset U_2;
$$
therefore, we obtain a well-defined quotient map $\nd s_{2/1}$ from the normal bundle of the embedding into the quotient bundle (restricted to $\phi_{12}(s_{1,2}^{-1}(0))$).
Note that $G_{1,2}\cong h_{12}(G_{1,2})\subset G_2$ acts on both 
$$
\cN_{\phi_{12}(V_{1,2})}{V_2}|_{\phi_{12}(s_{1,2}^{-1}(0))}~~~\tn{and}~~~(U_2/\mfD\phi_{12}(U_{1,2}))|_{\phi_{12}(s_{1,2}^{-1}(0))},
$$ 
and $\nd s_{2/1}$ is equivariant
with respect to this action. 

\bDef{def:Kuranishi-tangent}
We say that a  Kuranishi structure $\cK$ on a topological space $M$ has \textbf{tangent bundle} if for every coordinate change map $\Phi_{qp}\colon\cU_{q,p}\lra \cU_p$ in $\cK$ the map 
\bEqu{dspq}
\nd s_{p/q}\colon\cN_{\phi_{qp}(V_{q,p})}{V_p}|_{\phi_{qp}(s_{q,p}^{-1}(0))}\lra \big(U_p/\mfD\phi_{qp}(U_{q,p})\big)|_{\phi_{qp}(s_{q,p}^{-1}(0))},
\eEqu
given by (\ref{comm-tangent}) is an orbibundle isomorphism. 
In other words, we require the restriction of $s_p$ to $\phi_{qp}(s_{q,p}^{-1}(0))$ to be transverse in the normal direction modulo $U_{q,p}$. By a \textbf{Kuranishi space}, we mean a pair $(M,\cK)$ where $\cK$ is a Kuranishi structure with tangent bundle on $M$.
\eDef

\bRem{rem:trivial-case}
If the two charts $\cU_1$ and $\cU_2$ used in (\ref{comm-tangent}) have the same dimension, $\nd s_{2/1}$ is trivial and the tangent bundle condition is automatically satisfied. Therefore, an equidimensional set of Kuranishi neighborhoods automatically satisfies the tangent bundle condition. This is for example the case when $\cK$ is pure orbibundle Kuranishi structure.
\eRem

\noindent
We say a Kuranishi chart $\cU\!=\!(\pr\colon\!U\!\lra\!V,G,s,\psi)$ is \textbf{orientable} if the orbibundle $[U/ G]\lra [V/ G]$ is relatively orientable; i.e. there is a $G$-equivariant isomorphism
\bEqu{index-bundle}
\det U^* \otimes \det TV\cong V\!\times\! \R,
\eEqu
where the action of $G$ on $\R$ is the trivial action.
An orientation is a choice of (homotopy class of) such trivialization.
Given a $G_{q,p}$-equivariant isomorphism (\ref{dspq}), it induces a canonical $G_{q,p}$-equivariant isomorphism
\bEqu{iso-ori}
\det \nd s_{p/q}\colon\big( \det U_{q,p}^* \otimes \det TV_{q,p}\big)|_{s_{q,p}^{-1}(0)} \lra \big(\det U_p^*\otimes \det TV_p\big)|_{\phi_{qp}(s_{q,p}^{-1}(0))}.
\eEqu

\bDef{def:Kuranishi-orientable}
We say that a  Kuranishi space $(M,\cK)$ is \textbf{orientable} if every Kuranishi chart $\cU_p$ can be oriented in a way that following conditions hold. For every coordinate change map $\Phi_{qp}\colon\cU_{q,p}\lra \cU_p$ in $\cK$, the canonical isomorphism (\ref{iso-ori}) is orientation-preserving. Moreover, if $L(\ov{p})\!=\!L(\ov{q})$, then the isomorphism 
$$
\big( \det U_{q,p}^* \otimes \det TV_{q,p}\big)\stackrel{\cong}{\lra} \big(\det U_p^*\otimes \det TV_p\big)|_{\phi_{qp}(V_{q,p})}
$$
is orientation preserving. 
\eDef

\bRem{rem:trivial-case-2}
A pure orbibundle Kuranishi space as in Example~\ref{exa:global-kur} is orientable if and only if each orbibundle $\cE_i\!\lra\!\cY_i$ is relatively orientable.
\eRem
\subsection{Dimensionally graded systems}\label{sec:level}
Similar to the orbifold case of Section~\ref{sec:euler}, our idea for constructing a VFC for an arbitrary Kuranishi space is to perturb the Kuranishi maps $s_p$, in a compatible way with coordinate change maps, and make them transverse to the zero section. We construct these perturbations inductively; we start from the points with the lowest value of $L$ in (\ref{equ:d-function}) and move up. 
For this induction to work properly, we build a pure orbibundle Kuranishi structure -in the sense of Example~\ref{exa:global-kur}- around a compact subset of each level set $M(k)$ with orbibundle embeddings from (sub-orbibundles of) lower levels to higher levels. In this section, we introduce the notion of dimensionally graded system for Kuranishi structures which makes sense of this outline. This notion is a special case of the notion of ``good coordinate system" in \cite{FO} where the orbibundles are graded by their rank instead of an arbitrary partially ordered set.

\noindent
\begin{notation}
For a topological subspace $A\!\subset\! B$, let $\tn{cl}(A)\!\subset\! B$ denote the closure of $A$ inside $B$. We are adopting this notation, instead of the common notion $\ov{A}$, in order to avoid confusion with the quotient space of group actions. However, there is one special case, the case of compactified moduli spaces of maps and curves, $\ov\cM_{g,k}(X,A)$ and $\ov\cM_{g,k}$, respectively, where the over-line has its usual meaning.
\end{notation}

\bDef{level-system}
Let $\cK$ be an $n$-dimensional Kuranishi structure on $M$. A \textbf{dimensionally graded system  (or DGS)} for $(M,\cK)$ consist of an orbibundle (possibly empty!)
$$
\pr_i\colon \cE(i) \lra \cY(i)
$$ 
and an orbifold section 
$s_i\colon\! \cY(i)\!\lra\! \cE(i)$, for every $i\!\in\! \mfB\!=\!\mfB(M,\cK)\!\subset\!\Z$ with notation as in (\ref{equ:def-of-B}), such that the following conditions hold (See the notation convention in Notation~\ref{not:calconv}).
\bEnum
\item\label{l:dim} For every $i\!\in\! \mfB$, (if $Y(i)\neq \emptyset$) $\tn{rank}~\cE(i)\!=\!i$ and $\dim \cY(i)\!=\!n\!+\!i$.

\item\label{l:footprint} There are open embeddings $\ov{\psi}_i\colon\ov{s}_i^{-1}(0)\lra M$ such that
the set of footprints 
$$
\{F(i)= \ov{\psi}_i(\ov{s}_i^{-1}(0))\}_{i\in \mfB}
$$ 
is an open covering of $M$ and 
$$
\tn{cl}\!\lrp{F(k)}\cap M(>\! k)=\emptyset\quad \forall k\in \mfB;
$$ 
see Figure~\ref{F(i)}.

\begin{figure}
\begin{pspicture}(10,1)(11,2)
\psset{unit=.3cm}

\pscircle*(50,5){.25}\rput(49,4.5){\tiny{$M(2)$}}
\pscircle[linestyle=dashed,dash=2pt](50,5){2.7}

\psline(50,5)(50,10)\psline(50,5)(60,5)\rput(46,3.5){\tiny{$M(1)$}}
\psline[linearc=1]{->}(46.3,3.9)(48,6)(50,7)
\psline[linearc=1]{->}(46.9,3.5)(50,3.8)(51.5,5)

\psline[linestyle=dashed,dash=2pt](50.5,12)(50.5,6)
\psarc[linestyle=dashed,dash=2pt](51,6){.5}{180}{270}
\psline[linestyle=dashed,dash=2pt](51,5.5)(62,5.5)

\psline[linestyle=dashed,dash=2pt](53,6)(60,6)
\psline[linestyle=dashed,dash=2pt](53,4)(60,4)
\psarc[linestyle=dashed,dash=2pt](60,5){1}{-90}{90}
\psarc[linestyle=dashed,dash=2pt](53,5){1}{90}{-90}

\psline[linestyle=dashed,dash=2pt](51,7)(51,10)
\psline[linestyle=dashed,dash=2pt](49,7)(49,10)
\psarc[linestyle=dashed,dash=2pt](50,7){1}{-180}{0}
\psarc[linestyle=dashed,dash=2pt](50,10){1}{0}{180}

\rput(58,10){$F(0)$}
\rput(53,7.2){\tiny{$M(0)$}}

\rput(53,12){$F(1)$}
\psline{->}(52.8,11.5)(50,10.3)
\rput(61,7){$F(1)$}
\psline{->}(60.5,6.5)(60.5,5)
\rput(54,3){$F(2)$}
\psline{->}(52.8,3)(50.8,3)

\end{pspicture}
\caption{A possible configuration of footprints of different level.}
\label{F(i)}
\end{figure}
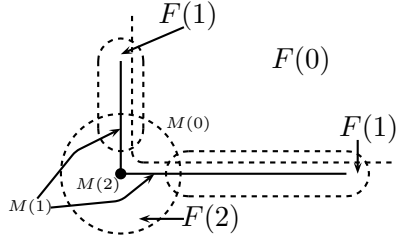

\item\label{l:intersection-2} For every $i,j\!\in\!\mfB$, with $i\!<\!j$, there exists an open subset $Y(i,j)\!\subset\!Y(i)$ and an embedding of  orbibundles commuting with the section and footprint maps 
$$
\aligned
&\Phi_{ij}=(\mfD \phi_{ij},\phi_{ij})\colon(\cE(i,j),\cY(i,j))\lra(\cE(j),\cY(j)),\\
& \cY(i,j)\!=\cY(i)|_{Y(i,j)},~~\cE(i,j)\!=\cE(i)|_{\cY(i,j)},
\endaligned
$$ 
as in Definition~\ref{def:orbibundle-maps}, such that\footnote{Anything indexed by $i,j$ denotes the restriction of the corresponding object indexed by $i$ to $\cE(i,j)\lra\cY(i,j)$.} 
\bEqu{equ:Fij=FicapFj}
F(i,j)= \ov{\psi}_{i,j}(\ov{s}_{i,j}^{-1}(0))=F(i)\cap F(j),\quad\tn{where}~~s_{i,j}\!=\!s_i|_{\cY(i,j)}.
\eEqu

\item\label{l:intersection-source} For $i,j,k\!\in\! \mfB$, with $i\!<\!j\!<\!k$,  
$$
\ov\phi_{ij}^{-1}(Y(j,k))=Y(i,k)\cap Y(i,j).
$$
We denote the intersection on right by $Y(i,j,k)$ and the restriction of $\cE(i,j)$ to $Y(i,j,k)$ by  $\cE(i,j,k)$.
 \item\label{l:intersection-target} For $i,j,k\!\in \mfB$, with $i\!<\!j\!<\!k$,  
 $$\ov\phi_{ik}(Y(i,k))\cap \ov\phi_{jk}(Y(j,k))=\ov\phi_{ik}(Y(i,j,k)).$$ 
 
\item\label{l:intersection-cycle}  For every $i,j,k\!\in\!\mfB$, with $i\!<\!j\!<\!k$, the cocycle condition
\bEqu{equ-of-maps}
(\Phi_{jk}\circ \Phi_{ij})|_{\cY(i,j,k)}=\Phi_{ik}|_{\cY(i,j,k)}
\eEqu
of Definition~\ref{def:cocycle} holds.

\item\label{l:compatible} For every $i\!\in\! \mfB$ and $\ov{p}\!\in\! F(i)$, there exists a sub-chart $\cU_{p,i}\!\subset\! \cU_p$ centered at $\ov{p}$ and an orbibundle embedding
\bEqu{equ:Kur-to-level}
\Phi_{pi}=(\mfD\phi_{pi},\phi_{pi})\colon ([U_{p,i}/ G_{p,i}],[V_{p,i}/ G_{p,i}])\lra (\cE(i),\cY(i))
\eEqu
that commutes with the projection, section, and footprint maps, and the following compatibility conditions (with $\cK$) hold.
\begin{enumerate}[label=(\alph*)]
\item Whenever $\ov{q}\!\in\! F_p\cap F(i)$ and $\ov{p}\!\in\! F(i)$, restricted to $\ov\phi_{qp}^{-1}(\ov{V}_{p,i})\cap \ov{V}_{q,i}$ we have $\ov\Phi_{pi}\circ \ov\Phi_{qp}\!=\!\ov\Phi_{qi}$, where $\ov{V}_{p,i}$ is the underlying topological space of the orbifold $[V_{p,i}/ G_{p,i}]$.
\item Whenever $\ov{p}\!\in\! F(i)\!\cap\! F(j)$ and $i\!<\!j$, restricted to $\ov\phi_{pi}^{-1}(Y_{i,j})\cap \ov{V}_{p,j}$ we have $\ov\Phi_{ij}\circ\ov\Phi_{pi}\!=\!\ov\Phi_{pj}$.
\end{enumerate}
\item\label{l:tangent} If in addition $(\cK,M)$ is a Kuranishi space, i.e. if $\cK$ has tangent bundle, then for every $i\!\in\! \mfB$ and every $\ov{p}\!\in \!F(i)$, with $\Phi_{pi}$ as in (\ref{equ:Kur-to-level}), the quotient map $\nd s_{i/p}$ as in (\ref{comm-tangent}) is an equivariant isomorphism of vector bundles. 
\eEnum
\eDef

\bRem{inclusion-property}
It follows from Definition~\ref{level-system}.\ref{l:footprint}\ref{l:compatible} that 
$$
\quad M(\geq\! k)\subset \hspace{-.1in}\bigcup_{i\in\mfB(\geq k)} F(i),\qquad F(k)\subset \hspace{-.1in}\bigcup_{\ov{p}\in M(k)} F_p.
$$
\eRem

\bRem{rem:explanation-of-conditions}
Note that by Definition~\ref{level-system}.\ref{l:footprint}, the necessary condition 
$$\tn{rank } \cU_{p,i} \leq \tn{rank }\cU_i$$ for~\ref{l:compatible} to be plausible is automatically satisfied.
The first and second conditions in Definition~\ref{level-system} together with \ref{l:compatible} simply assert that $(\cE(i),\cY(i))$ is made of Kuranishi charts centered at points of level $i$. Conditions~\ref{l:intersection-2} and \ref{l:intersection-cycle} allow us to glue different $(\cE(i),\cY(i))$ along their common intersections and form a thickening of the moduli space; see Lemma~\ref{sim-thick} and the argument preceding that. The equality (\ref{equ:Fij=FicapFj}) in~\ref{l:intersection-2} insures that every coordinate change sub-orbibundle covers the whole intersection of the corresponding footprints. 
Conditions~\ref{l:intersection-2} and \ref{l:intersection-source} are necessary for~\ref{l:intersection-cycle} to make sense.
We need~\ref{l:intersection-target} as well as~\ref{l:intersection-cycle} for the relation defined in Lemma~\ref{sim-thick} to be an equivalence relation. Condition~\ref{l:compatible} indicates how the level structure is related to the underlying Kuranishi structure. Finally, Condition~\ref{l:tangent} (as stated) needs some further clarification. The orbibundle $\cE(i)$ itself, by definition, is made of local orbibundle charts. Then~\ref{l:tangent} means that for each $\ov{p}\!\in\! F(i)$, with $\Phi_{pi}$ as in (\ref{equ:Kur-to-level}), there exist a sub-chart $\cU_1\!=\!(\pr_1\colon\!U_1\!\lra\!V_1,G_1,s_1,\psi_1)$ of $\cU_p$ centered at $\ov{p}$ and a defining orbibundle chart (together with the restriction of the footprint map and section maps) $\cU_2\!=\!(\pr_2\colon\!U_2\!\lra\!V_2,G_2,s_2,\psi_2)$ of $\cE(i)$
centered at 
$$
\ov{p}\!\cong\! \ov\psi_i^{-1}(\ov{p})\!\in \!\ov{s}_i^{-1}(0)\!\subset\! Y(i),
$$ 
over which $\Phi_{pi}$ is given by an equivariant embedding of vector bundles 
$$
(\mfD\phi,\phi)\colon\cU_1\lra \cU_2,\quad  \mfD\phi \circ s_1=s_2\circ \phi,\quad \ov\phi_1\circ \psi=\psi_2\circ \phi,
$$
and the quotient map $\nd s_{2/1}$ as in (\ref{comm-tangent}) is an equivariant isomorphism of vector bundles. 
\eRem

\bRem{global-case}
In the case of pure orbibundle Kuranishi structures introduced in Example~\ref{exa:global-kur}, we simply set
$$
\aligned
&\big(M(i),\cE(i),\cY(i)\big)=\hspace{-.2in}\coprod_{j\in[m]\colon \tn{rank}~\cE_j=i} \hspace{-.2in}(M_j,\cE_j,\cY_j)\\
&\tn{and}~~
\big(\cE(i,j),\cY(i,j)\big)=\emptyset\qquad \forall i,j\!\in\! \mfB,~i\!<\!j.
\endaligned
$$
Therefore, there are no overlap conditions for pure orbibundle Kuranishi structures and we are essentially in the realm of orbifold and orbibundles. This also illuminates the idea behind dimensionally graded systems; we assemble the uncountable data of Kuranishi structure into a finite set of orbibundles, graded by their ranks, with transition embeddings from lower grades to higher grades. We then use this finite set to inductively build up the VFC.
\eRem

\bLem{lem:level-tangent}
Assume
\bEqu{level-tuple}
\mfL=\mfL(\cK)\equiv\big(\{\cE(i),\cY(i),s_i,\psi_i\}_{i\in \mfB},\{\cE(i,j),\cY(i,j),\Phi_{ij}\}_{i<j,~i,j\in\mfB}\big)
\eEqu
 is a dimensionally graded system for  the Kuranishi space $(M,\cK)$.
Then $\mfL$ also has tangent bundle in the sense that for every $i,j\!\in\!\mfB$, with $i\!<\!j$, the map
\bEqu{global-dsij}
\nd s_{j/i} \colon \cN_{\phi_{ij}(\cY(i,j))}\cY(j)|_{\ov\phi_{ij}(\ov{s}_{i,j}^{-1}(0))} \lra \big(\cE(j)/\mfD\phi_{ij} \cE(i,j)\big)|_{\ov{\phi}_{ij}(\ov{s}_{i,j}^{-1}(0))}.
\eEqu
defined locally by the right column of (\ref{comm-tangent}) is an isomorphism of orbibundles.
\eLem

\bProof The point is that the condition (\ref{global-dsij}) is local, and by the compatibility conditions \ref{l:tangent} and \ref{l:compatible} of Definition~\ref{level-system}, local charts of $(\cE(i),\cY(i))$ and that of Kuranishi structure are essentially the same. More precisely, for every $\ov{p}\!\in \!F(i,j)$, let $\cU_0\!=\!(\pr_0\colon\!U_0\!\lra\!V_0,G_0,s_0,\psi_0)$ be a sub-chart of $\cU_p$ centered at $\ov{p}$, 
$\cU_1\!=\!(\pr_1\colon\!U_1\!\lra\!V_1,G_1,s_1,\psi_1)$ be a sub-chart of $\cE(i)$ centered at 
$$
\ov{p}\!\cong\! \ov\psi_i^{-1}(\ov{p})\!\in\! \ov{s}_i^{-1}(0)\!\subset Y(i) ,$$
and $\cU_2\!=\!(\pr_2\colon\!U_2\!\lra\!V_2,G_2,s_2,\psi_2)$ be a sub-chart of $\cE(j)$ centered at $\ov{p}$, such that $\Phi_{pi}$, $\Phi_{pj}$ and $\Phi_{ij}$ are given by equivariant embeddings of vector bundles 
$$
(\mfD\phi_{ab},\phi_{ab})\colon\cU_a\lra \cU_b\qquad \forall a\!=\!0,1,~b\!=\!1,2,~a\!<\!b,
$$
and the quotient maps $\nd s_{1/0}$ and $\nd s_{2/0}$ are as in (\ref{comm-tangent}). To distinguish between preimages of $\ov{p}$ in $V_0$, $V_1$, and $V_2$, we denote them by $p_0$, $p_1$, and $p_2$, respectively. We then obtain a commutative diagram
$$
    \xymatrix{  
        0 \ar[d]
            && 
   0 \ar[d]\\
    \cN_{\phi_{01}(V_0)}V_1|_{p_1=\phi_{01}(p_0)} \ar[rr]^{\nd s_{1/0}}\ar[d]^{\nd \phi_{12}}
    && 
    \big(U_{1}/\mfD\phi_{01}(U_{0})\big)|_{p_1} \ar[d]^{\nd \mfD\phi_{12}}\\
    \cN_{\phi_{02}(V_0)}V_2|_{p_2} \ar[rr]^{\nd s_{2/0}} \ar[d]^{\pi}
    &&
   \big(U_{2}/\mfD\phi_{02}(U_{0})\big)|_{p_2} \ar[d]^{\pi}\\  
       \cN_{\phi_{12}(V_1)}V_2|_{p_2} \ar[rr]^{\nd s_{2/1}}\ar[d]
    &&
   \big(U_{2}/\mfD\phi_{12}(U_{1})\big)|_{p_2}\ar[d]\\
          0 
            && 
  		 0 
   }
$$
where $\pi$ is the obvious projection map.
Since the columns are exact and the first two rows are isomorphism, we conclude that the bottom row is also an isomorphism.
\eProof

\noindent
\bDef{def:Kuranishi-orientable-2}
Let $\mfL$ be an DGS with tangent bundle as in Lemma~\ref{lem:level-tangent}. We say $\mfL$ is \textbf{orientable} (resp. oriented) if every orbibundle $\cE(i)\!\lra\!\cY(i)$ is relatively orientable (resp. oriented) and there exists a (resp. the given) set of trivializations 
\bEqu{equ:index-bundle-i-0}
\det \cE(i)^* \otimes \det T\cY(i)\cong \cY(i)\!\times\! \R \quad \forall i\!\in\!\mfB
\eEqu 
(resp. are) compatible with isomorphisms
\bEqu{equ:iso-det-bundles}
\det \nd s_{j/i}\colon\!\big(\! \det \cE(i,j)^* \otimes \det T\cY(i,j)\big)|_{\ov{s}_{i,j}^{-1}(0)} \!\lra\! \big(\!\det \cE(j)^*\otimes \det T\cY(j)\big)|_{\ov\phi_{ij}(\ov{s}_{i,j}^{-1}(0))}
\eEqu
induced by (\ref{global-dsij}).
\eDef

\noindent
Similarly to Lemma~\ref{lem:level-tangent}, if $(M,\cK)$ is orientable (resp. oriented) and $\mfL$ is an DGS for $(M,\cK)$, then a shrinking of $\mfL$, which we will introduce in Section~\ref{shrinking}, is orientable (resp. oriented).

\noindent
In the rest of this section, we introduce an equivalence relation on the disjoint union 
\bEqu{dis-uion}
\coprod_{i\in \mfB} Y(i)
\eEqu
such that the quotient space $Y(\mfL)$ can be seen as  a global \textbf{thickening} of the space $M$ in which a perturbation of the zero set $M\!\subset\! Y(\mfL)$ will reside.  

\bLem{sim-thick}
The relation $\sim$ on (\ref{dis-uion}) defined by
\bEqu{equ:sim-relation}
\ov{x}\in Y(i)\sim \ov{y}\in Y(j) \Leftrightarrow \begin{cases} i=j~~\tn{and}~~\ov{x}=\ov{y},~~\tn{or}\\ 
 i< j,~\ov{x}\in Y(i,j),~\tn{and}~~\ov\phi_{ij}(\ov{x})=\ov{y},~~\tn{or}\\
 j<i,~\ov{y}\in Y(j,i),~\tn{and}~~\ov\phi_{ji}(\ov{y})=\ov{x},
 \end{cases}  
\eEqu
is an equivalence relation.
\eLem

\bProof
Among reflexivity, symmetry, and transitivity, only the last one is somewhat non-trivial. Suppose $\ov{x}_1\!\sim\! \ov{x}_2$ and $\ov{x}_2\!\sim \!\ov{x}_3$, $\ov{x}_a \!\in \!Y(j_a)$.  If two of $j_a$ are equal, the desired conclusion is trivial. If $j_1\!<\!j_2\!<\!j_3$, $\ov{x}_1\!\sim\! \ov{x}_3$ follows from  Definition~\ref{level-system}.\ref{l:intersection-cycle}. If $j_2 \!<\!j_1,j_3$, we may assume $j_1\!<j_3$, then again 
$$
\ov{x}_2\!\in\! \ov\phi_{j_2j_1}^{-1}(Y(j_1,j_3))=Y(j_2,j_1)\cap Y(j_2,j_3)
$$ 
and the result follows from the cocycle condition. If $j_1,j_3\!<\! j_2$, we may assume $j_1\!<\!j_3$, then by  Definition~\ref{level-system}.\ref{l:intersection-target}, 
$$
\ov{x}_2\!\in\! \ov\phi_{j_1,j_2}(Y(j_1,j_2))\cap \ov\phi_{j_3,j_1}Y(j_3,j_2)\!=\!\ov\phi_{j_1,j_2}Y(j_1,j_3,j_2).
$$  
By the uniqueness of $\ov{x}_1\!\in\! Y(j_1,j_3)$ and $\ov{x}_3\!\in\! Y(j_3,j_2)$, we conclude $\ov{x}_1\!\in\! Y(j_1,j_3,j_2)$ and $\ov\phi_{j_1j_3}(x_1)\!=\!\ov{x}_3$.
\eProof

\noindent
Let $\wp_i\colon\! Y(i)\!\!\lra\!\! Y(\mfL)$ be the obvious projection maps.
By (\ref{equ:sim-relation}) and Definition~\ref{level-system}.\ref{l:intersection-2}, the inverse footprint maps $\ov\psi_i^{-1}\colon F(i) \!\lra\! \ov{s}_i^{-1}(0)$ followed by $\wp_i$ can be glued on the overlaps into a global continuous one-to-one map
$$
\iota\colon M\lra Y(\mfL).
$$

\bDef{def:H-level-structure}
A dimensionally graded system $\mfL$ as in (\ref{level-tuple}) is called \textbf{Hausdorff}  if the quotient topology on the thickening
\bEqu{equ:thikening-space}
Y(\mfL)\equiv
\bigg(\coprod_{i\in \mfB} Y(i)\bigg)/\sim
\eEqu 
is Hausdorff.
\eDef

\bLem{lem:HtoP}
A dimensionally graded system  $\mfL$ as in (\ref{level-tuple}) is Hausdorff if and only if for every $i,j\!\in\!\mfB$, with $i\!>\!j$, the map
\bEqu{equ:HtoP}
 Y(i,j) \overset{\ov\iota_{ij}\times \ov\phi_{ij}}{\xrightarrow{\hspace*{1cm}}} Y(i)\times Y(j),
\eEqu
where $\ov\iota_{ij}$ is the inclusion map, is proper.
\eLem
\bProof
Assume $\ov\iota_{ij} \times \ov\phi_{ij}$ is not proper but $Y(\mfL)$ is Hausdorff. Then there exists a compact set $K\!\subset\! Y(i)\!\times\! Y(j)$ such that $\wt{K}\!\equiv\! (\ov\iota_{ij} \times \ov\phi_{ij})^{-1}(K)$ is not compact. 
Let $(\ov{x}_k)_{k\in\N}\!\subset\! \wt{K}$ be a sequence\footnote{Since $Y(i,j)$, $Y(i)$, and $Y(j)$ are metrizable, compactness is equal to sequential compactness.} with no subsequence limit in $\wt{K}$ and $\big((\ov{x}_k,\ov{y}_k)\big)_{k\in \N}\subset K$ be the image of this sequence under (\ref{equ:HtoP}). By assumption, and after possibly restricting to a subsequence, there exists $(\ov{x},\ov{y})\in K$ such that 
$$
\lim_{k\lra \infty} (\ov{x}_k,\ov{y}_k)\!=\!(\ov{x},\ov{y})\quad\tn{and}\quad \ov{x}\!\in\! Y(i)\!\setminus\! Y(i,j).
$$ 
By definition of quotient topology on $Y(\mfL)$ and since $\ov{x}_k\sim \ov{y}_k$, for all $k\!\in\!\N$, every two open sets around $[\ov{x}]\!=\!\wp_i(\ov{x})$ and $[\ov{y}]\!\in \!\wp_j(\ov{y})$ in $Y(\mfL)$ intersect non-trivially; thus, from Hausdorffness we deduce $[\ov{x}]\!=\![\ov{y}]$. This is a contradiction because 
$\ov{x}\!\notin\! Y(i,j)$. The other direction is similar.
\eProof

\subsection{Shrinking}\label{shrinking}

Given an DGS for $(M,\cK)$ as in Definition~\ref{level-system}, 
we often need to shrink the orbibundle pieces to accommodate certain perturbations and gluing-pasting arguments. For example, we use shrinking to achieve Hausdorffness and metrizability in the inductive construction of DGS. For the general case of good coordinate systems we refer to \cite{FOOO-shrink}.

\noindent
For the construction of VFC in Section~\ref{sec:VFC}, we need the thickening $Y(\mfL)$ of (\ref{equ:thikening-space}) to be metrizable; c.f. Definition~\ref{def:H-level-structure}. Even if the quotient topology on $Y(\mfL)$ is Hausdorff, it is often the case that this topology is not paracompact, thus not metrizable. An example of this issue is discussed in \cite[Example 6.1.11]{MW2}; we recall this example here. For simplicity, the underlying topological space $M$ is taken to be the non-compact manifold $\R$. Assume $\mfL$ is an DGS made of two vector bundles, $Y(1)\!=\!\R$ together with the trivial bundle, section, and footprint map, and 
$$
Y(2)=(0,\infty)\times \R,\quad E(2)=Y(2)\times \R, \quad s_2(x,y)=y, \quad \psi_2(x,y)=x.
$$
The embedding $\Phi_{12}$ is the obvious inclusion over $Y(1,2)\!=\!(0,\infty)\!\subset\! Y(1)$. Let $[0]$ be the image of $0\!\in\! Y(1)$ in $Y(\mfL)$. Then the quotient topology on $Y(\mfL)$ near $[0]$ strictly stronger than the subspace topology $Y(\mfL)\!\subset\! \R^2$. In fact, for every $\ep\!>\!0$ and every continuous function $f\colon (0,\ep)\lra (0,\infty)$, the set
$$
U_{f,\ep}\!=\!\big\{[x]\colon~x\!\in\! Y(1),~|x|\!<\!\ep\big\}\cup \big\{ [(x,y)]\colon ~(x,y)\! \in\! Y(2), ~|x|\!<\!\ep,~ y\!<\!f(x)\big\}
$$
is open in the quotient topology. Moreover, these sets form a basis around $[0]$ in the quotient topology.

\noindent
On the otherhand, by Proposition~\ref{pro:shrinking} below, if $Y(\mfL)$ is Hausdorff, then the induced topology on $Y(\wt\mfL)$, where $\wt\mfL$ is an arbitrary shrinking of $\mfL$, is metrizable. 

\bDef{def:shrinking}
With the notation as in Definition~\ref{level-system}, given a dimensionally graded system $\mfL$, a \textbf{shrinking} of $\mfL$ consist of relatively compact\footnote{i.e. the underlying topological inclusions are relatively compact.} open sub-orbifolds $\wt\cY(i)\subset \cY(i)$, for each $i\!\in\!\mfB$, such that 
\bEqu{e:covers}
M(\geq\!k)\!\subset\! \bigcup_{i\in \mfB(\geq k)}\hspace{-.1in} \wt{F}(i)\quad \forall k\!\geq\! 0,\quad\tn{with}~~~
\wt{F}(i)= \ov{\psi}_i(\ov{s}_i^{-1}(0)\cap \wt{Y}(i))\quad\forall i\!\in\! \mfB.
\eEqu
\eDef
\noindent
In this situation, we define $\wt\cE(i)$ to be the restriction to $\wt\cY(i)$ of $\cE(i)$. For each $i\!\in\!\mfB$, the section and footprint maps $s_i$ and $\psi_i$, respectively, restrict to section and footprint maps on $\wt\cY(i)$;  for simplicity, we will often denote them by the same notation.

\bLem{lem:shrinking}
Let $\mfL$ be an DGS as in Definition~\ref{level-system}. Given a shrinking $\{\wt\cY(i)\}_{i\in \mfB}$ of $\mfL$ as in Definition~\ref{def:shrinking}, with 
\bEqu{equ:tildeYij}
\wt\cY(i,j)\!=\! \wt\cY(i)\cap \cY(i,j)\cap \phi_{ij}^{-1}(\wt\cY(j)),\quad \wt\cE(i,j)\!=\!\cE(i,j)|_{\wt\cY(i,j)}\quad\forall i,j\!\in\! \mfB,~ i\!<\!j,
\eEqu
the shrunk collection 
\bEqu{shrunkL_e}
\wt\mfL=\big(\{\wt\cE(i),\wt\cY(i),s_i,\psi_i\}_{i\in \mfB},\{\wt\cE(i,j),\wt\cY(i,j),\Phi_{ij}|_{\wt\cY(i,j)}\}_{i<j,~i,j\in\mfB}\big)
\eEqu
is also an DGS for $(M,\cK)$. 
\eLem

\bProof
It immediately follows from the definition of $\wt\cY(i,j)$ and (\ref{e:covers}) that the new collection satisfies Conditions~\ref{l:dim},~\ref{l:footprint},~\ref{l:intersection-2},~\ref{l:tangent} of Definition~\ref{level-system}. For Condition~\ref{l:intersection-source}, we have
$$
\aligned
\ov\phi_{ij}^{-1}(\wt{Y}(j,k))\cap \wt{Y}(i)=&\ov\phi_{ij}^{-1}\bigg(\wt{Y}(j)\cap Y(j,k)\cap \ov\phi_{jk}^{-1}(\wt{Y}(k))\bigg)\cap \wt{Y}(i)=\\
& \ov\phi_{ij}^{-1}(\wt{Y}(j)) \cap  \ov\phi_{ij}^{-1}\bigg(Y(j,k) \cap \ov\phi_{jk}^{-1}(\wt{Y}(k))\bigg)\cap \wt{Y}(i)=\\
&\ov\phi_{ij}^{-1}(\wt{Y}(j)) \cap Y(i,j)\cap Y(i,k) \cap \ov\phi_{ik}^{-1}(\wt{Y}(k)) \cap \wt{Y}(i)= \\
&\wt{Y}(i,k)\cap \wt{Y}(i,j).
\endaligned
$$
For Condition~\ref{l:intersection-target}, if 
\bEqu{e:2-possibilities}
\aligned
\ov{x}\in \wt{Y}(i,j,k)=& \ov\phi_{ij}^{-1}(\wt{Y}(j)) \cap Y(i,j)\cap Y(i,k) \cap \ov\phi_{ik}^{-1}(\wt{Y}(k)) \cap \wt{Y}(i)=\\
&\ov\phi_{ij}^{-1}(\wt{Y}(j)\cap Y(j,k)) \cap \ov\phi_{ik}^{-1}(\wt{Y}(k)) \cap \wt{Y}(i),
\endaligned
\eEqu
then, by the first line of (\ref{e:2-possibilities}),
$$
\ov\phi_{ik}(\ov{x})\in \ov\phi_{ik}\bigg(\wt{Y}(i) \cap Y(i,j)\bigg)\cap \wt{Y}(k)=\ov\phi_{ik}(\wt{Y}(i,k)), 
$$
and, by the second line of (\ref{e:2-possibilities}),
$$
\ov\phi_{jk}(\ov{y})\in \ov\phi_{jk}\bigg(\wt{Y}(j) \cap Y(j,k)\bigg)\cap \wt{Y}(k)=\ov\phi_{jk}(\wt{Y}(j,k)), \quad \ov{y}=\ov\phi_{ij}(\ov{x}).
$$
Therefore, $\ov\phi_{ik}(\ov{x})\in \ov\phi_{ik}(\wt{Y}(i,k))\cap \ov\phi_{jk}(\wt{Y}(j,k))$.
Conversely, if 
$$
\ov{z}\in \ov\phi_{ik}(\wt{Y}(i,k)) \cap \ov\phi_{jk}(\wt{Y}(j,k))\subset \ov\phi_{ik}({Y}(i,k)) \cap \ov\phi_{jk}({Y}(j,k))=\ov\phi_{ik}(Y(i,j,k)) ,
$$
then 
$$
\ov{x}= \ov\phi_{ik}^{-1}(\ov{z})\in Y(i,j,k),~\wt{Y}(i),~\tn{and}~\phi_{ik}^{-1}(\wt{Y}(k)),
$$
and 
$$
\ov{x}=\ov\phi_{ij}^{-1}(\ov{y}), \quad \ov{y}= \ov\phi^{-1}_{jk}(\ov{z})\in \wt{Y}(j);
$$
i.e. $\ov{x}\in \ov{Y}(i,j,k)$. Condition~\ref{l:intersection-cycle} is just the restriction to the shrunk sub-orbifolds of the original cocycle condition. Finally, for Condition~\ref{l:compatible}, we can replace $\cU_{p,i}$ with a sub-chart centered at $\ov{p}$ of $\cU_{p,i}\cap \Phi_{pi}^{-1}(\wt\cY(i))$.
\eProof

\noindent
By (\ref{equ:tildeYij}) and Lemma~\ref{sim-thick},  the composite projection maps
$$
\wt{Y}(i)\hookrightarrow Y(i)\stackrel{\wp_i}{\lra} Y(\mfL)\quad \forall i\!\in\!\mfB
$$
descend to a continuous (both with quotient topologies) one-to-one map  
\bEqu{equ:tilde-to-Y}
Y(\wt\mfL)\!\hookrightarrow\!Y(\mfL).
\eEqu
In what follows, by $Y(\wt\mfL)\!\subset\!Y(\mfL)$ we mean the image of (\ref{equ:tilde-to-Y}), and 
we denote by $Y_{\mfL}(\wt\mfL)$ the induced topology given by this inclusion. The quotient topology on $Y(\wt\mfL)$ could be different than $Y_{\mfL}(\wt\mfL)$.

\bLem{lem:cl-im-vs-im-cl} 
Let  $\wt\mfL$ be a shrinking of $\mfL$ as in Lemma~\ref{lem:shrinking}. If $Y(\mfL)$ is Hausdorff, then $\tn{cl}(Y(\wt\mfL))\!\subset\!Y(\mfL)$ is compact and is equal to 
$$
\bigcup_{i\in\mfB} \wp_i(\tn{cl}(\wt{Y}(i))).
$$
\eLem
\bProof
Let $W$ be the complement of $\tn{cl}(Y(\wt\mfL))\!\subset\!Y(\mfL)$ and  
$$
W_i\!=\!\wp_i^{-1}(W)\quad \forall i\!\in\!\mfB.
$$ 
Every $W_i\!\subset\!Y(i)$ is open and is disjoint from $\wt{Y}(i)$; thus, it is disjoint from $\tn{cl}(\wt{Y}(i))$. We conclude that  
\bEqu{equ:ic-in-ci}
\bigcup_{i\in\mfB} \wp_i(\tn{cl}(\wt{Y}(i)))\!\subset\!\tn{cl}(Y(\wt\mfL)).
\eEqu
Since every $\tn{cl}(\wt{Y}(i))$ is compact and $\wp_i$ is continuous, we conclude that the left-hand side of (\ref{equ:ic-in-ci}) is a compact subset of $Y(\mfL)$. Therefore, if (\ref{equ:ic-in-ci}) is proper, there exists $[\ov{x}]\!\in\!Y(\wt\mfL)$ such that 
$$
[\ov{x}]\!\not\in\!\bigcup_{i\in\mfB} \wp_i(\tn{cl}(\wt{Y}(i))).
$$
By Haussdorfness of $Y(\mfL)$, there exists an open neighborhood $W$ of $[\ov{x}]$ in $Y(\mfL)$ which is disjoint from $\wp_i(\tn{cl}(\wt{Y}(i)))$, for all $i\!\in\!\mfB$. Therefore, $\wp_i^{-1}(W)$ is disjoint from $\wt{Y}(i)$, for all $i\!\in\!\mfB$. 
This is a contradiction because by definition of $Y(\wt\mfL)$, there exists some $i\!\in\!\mfB$ and $\ov{x}_i\!\in\!\wt{Y}(i)$ such that 
$[\ov{x}]\!=\!\wp_i(\ov{x}_i)$. 
\eProof

\bPro{pro:shrinking}
Let $\mfL$ be an DGS as in Definition~\ref{level-system} and $\wt\mfL$ be a shrinking of that as in (\ref{shrunkL_e}). If $Y(\mfL)$ with the quotient topology is Hausdorff, then $Y(\wt\mfL)$ with the quotient topology is Hausdorff and with the induced topology from $Y(\mfL)$ (i.e. $Y_\mfL(\wt\mfL)$) is metrizable.
\ePro

\bProof
By definition
$$
\wt{Y}(i,j)= Y(i,j)\cap \wt{Y}(i)\cap \ov{\phi}_{ij}^{-1} \wt{Y}(j)=(\ov\iota_{ij} \times \ov\phi_{ij})^{-1}(\wt{Y}(i)\times \wt{Y}(j)).
$$
Therefore, from Lemma~\ref{lem:HtoP} we conclude that if $Y(\mfL)$ with the quotient topology is Hausdorff, then so is $Y(\wt\mfL)$. This establishes the first claim. 

\noindent
In order to prove that $Y_\mfL(\wt\mfL)$ is metrizable, first, we show that $Y_\mfL(\wt\mfL)$ is regular. 
Let $\tn{pt}\!=\![\ov{x}]\!\in\! Y_\mfL(\wt\mfL)$ be a point and $C\!\subset\! Y_\mfL(\wt\mfL)$ be a non-empty closed set that does not include $\tn{pt}$. By Lemma~\ref{lem:cl-im-vs-im-cl},  as  subsets of $Y(\mfL)$,
$$
\tn{cl}(C)\subset \tn{cl}(Y(\wt\mfL))\subset Y(\mfL)
$$ 
is compact and $\tn{pt}\!\notin\! \tn{cl}(C)$. 
Since $Y(\mfL)$ is Hausdorff, for every other point $\tn{pt}'\!=\![\ov{y}]\in \tn{cl}(C)$, there exist open disjoint neighborhoods 
$$
\tn{pt}\in W_{\tn{pt}'}\quad\tn{and}\quad \tn{pt}'\in W'_{\tn{pt}'},
$$ 
separating $\tn{pt}$ and $\tn{pt}'$ from each other. 
Since $\tn{cl}(C)$ is compact, there exists a finite set of such points 
$$
J=\{\tn{pt}_1,\cdots,\tn{pt}_\ell\}\subset \tn{cl}(C)
$$ 
such that $W_C\!=\! \bigcup_{i=1}^\ell W'_{\tn{pt}_i}$ covers $\tn{cl}(C)$. In the induced topology, the open sets
$$
W=\bigcap_{i=1}^\ell W_{\tn{pt}_i}\quad \tn{and} \quad W_C
$$ 
separate $\tn{pt}$ and $C$. Next, we show that $Y_\mfL(\wt\mfL)$ is second countable. 
Since each $Y(i)$ is metrizable, $Y(i)$ is second-countable and we can take a countable basis $\{W_{i,\al}\}_{\al\in \Om_i}$ for it, including the empty set.
For every $i\!\in\! \mfB$, let 
\bEqu{equ:bad-stuff}
Y'(i)= Y(i)\setminus \bigcup_{j>i} \tn{cl}(\wt{Y}(i,j)),
\eEqu
and denote by $\Om_i'\!\subset\! \Om_i$ to be the set of those $\al\!\in \!\Om_i$ such that $W_{i,\al}\!\subset\! Y'(i)$.
For every $i\!\in\! \mfB$ and every $\al\!\in\! \Om_i'$, let 
\bEqu{equ:open-in-target}
[W_{i,\al}]\!=\!\wp_i(W_{i,\al});
\eEqu
this may not a priori be an open set in $Y(\mfL)$, but by (\ref{equ:bad-stuff}), its restriction to $Y(\wt\mfL)$ is open. Then, we claim that the set
$$
\{[W_{i,\al}]\cap Y(\wt\mfL)\}_{i\in \mfB,\al\in \Om_i'}
$$
is a basis for the induced topology. In fact, every open set of $Y_\mfL(\wt\mfL)$ is the restriction of some open set $W$ in $Y(\mfL)$. 
Let 
$$
W_i =\wp_i^{-1}(W)\subset Y(i).
$$ 
By (\ref{equ:bad-stuff}),
$$
Y(\wt\mfL)\cap W=Y(\wt\mfL)\cap \bigcup_{i\in \mfB} \wp_i(W_i) = Y(\wt\mfL)\cap \bigcup_{i\in \mfB}  \wp_i(W_i \cap Y'(i)).
$$ 
For every $i\!\in\! \mfB$, we have
$$
W_i \cap Y'(i)\!=\!\bigcup_{\al\in I_i} W_{i,\al},
$$ 
for some $I_i\!\subset\! \Om_i'$.
We conclude that 
$$
Y(\wt\mfL)\cap W= \bigcup_{i\in \mfB}\bigcup_{\al\in I_i} [W_{i,\al}].
$$
\noindent
Since $Y_\mfL(\wt\mfL)$ is already Hausdorff, it follows from Nagata-Smirnov's enhancement of Urysohn's theorem, see Remark~\ref{rem:top-remark} below, that $Y_\mfL(\wt\mfL)$ is separable and metrizable. 
\eProof

\bRem{rem:ulternative}
In Proposition~\ref{pro:shrinking}, a similar argument shows that the induced topology on $\tn{cl}(Y(\wt\mfL))\!\subset\!Y(\mfL)$ is also metrizable. For example, we can consider a shrinking $\mfL'$ of $\mfL$ such that $\wt\mfL$ is a shrinking of $\mfL'$. Then the conclusion follows from, first, applying Proposition~\ref{pro:shrinking} to $Y(\mfL')$ and then restricting to $\tn{cl}(Y(\wt\mfL))$. Moreover, in this case, the quotient topology and the induced topology on $\tn{cl}(Y(\wt{\mfL}))$ are the same. 
\eRem

\bRem{rem:top-remark}
There are some widely-recognized metrizability theorems. One of the first ones, due to Urysohn, states that every Hausdorff second-countable regular space is metrizable. The converse does not hold; there exist metric spaces that are not second countable. The Nagata-Smirnov theorem, provides a more specific theorem where the converse does hold. Urysohn's Theorem can be restated as: A topological space is separable and metrizable if and only if it is regular, Hausdorff, and second-countable. The Nagata-Smirnov metrizability theorem extends this to the non-separable case. It states that a topological space is metrizable if and only if it is regular, Hausdorff, and has a $\si$-locally finite base. A $\si$-locally finite base is a base which is a union of countably many locally finite collections of open sets.  A space is said to be locally metrizable if every point has a metrizable neighborhood. Smirnov proved that a locally metrizable space is metrizable if and only if it is Hausdorff and paracompact. Such spaces admit partitions of unity subordinate to any open cover. 
\eRem

\bRem{rem:Haus-later}
We have not yet shown that an arbitrary dimensionally graded system $\mfL$ admits a Haussdorf shrinking $\wt\mfL$ as in (\ref{shrunkL_e}).
We will do this in Theorem~\ref{exist-level-system}.
Proposition~\ref{pro:shrinking} only shows that once we build a Haussdorf DGS, the simple shrinking process of Lemma~\ref{lem:shrinking} preserves this property. 
\eRem

\subsection{Cobordism}\label{sec:cobordism}

Our construction of a VFC in Section~\ref{sec:VFC} depends on a choice of a dimensionally graded system $\mfL$, deformation of Kuranishi maps, etc. Before all these other auxiliary data, it depends on a choice of a Kuranishi structure $\cK$ on the topological space in question $M$. For example, in Section~\ref{sec:natural-KUR}, we build a class of natural Kuranishi structures on the moduli space of pseudoholomorphic maps; thus, it is important to show that the resulting VFC does not depend on the particular choice of a natural Kuranishi structure or the defining almost complex structure $J$. 
For this purpose we sketch a notion of a cobordism between two sets of choices. 

\noindent
The definition of an orbifold in Section~\ref{sec:orbifold} does not include charts with boundary and corners.
In order to define a cobordism, we need to at least include orbifold charts with boundary.  A bordered orbifold chart can be simply taken to be a product chart of the form
\bEqu{equ:B-charts}
\cV\!=\!(V\!\times\! [0,a),G,\psi);
\eEqu
where $(V,G,\psi)$ is an orbifold chart in the sense of Definition~\ref{def:local-orbi}, $a$ is some positive integer, and $G$ acts trivially on $[0,a)$. Then 
$$
\partial\cV\!=\!(V\!\times\! \{0\},G,\psi)\!\cong\! (V,G,\psi)\!\subset\! \cV\quad\tn{and}\quad \cV\!\setminus \!\partial\cV\!=\!(V\!\times\! (0,a),G,\psi)\!\subset\! \cV
$$
are orbifold charts in the sense of Definition~\ref{def:local-orbi}. A refinement of such $\cV$ at an interior point is a refinement of $\cV\!\setminus \!\partial\cV$ as in Definition~\ref{def:refinement}, and a refinement at a boundary point is another orbifold chart with boundary $\cV'\!=\!(V'\!\times\![0,a'),G',\psi')$ (of the same dimension ) together with a group homomorphism $h\colon\!G'\!\lra\! G$ and an $h$-equivariant smooth embedding 
\bEqu{equ:bordered-refinement} 
\phi \colon\!V'\!\times\![0,a') \!\lra\! V\!\times\![0,a)
\eEqu 
which takes $\partial\cV'$ to $\partial \cV$. Eventually, for a metrizable topological space $M$ with ``boundary'' $\partial M\!\subset\!M$, an $n$-dimensional smooth (effective) orbifold atlas on $(M,\partial M)$ is a countable collection 
\bEqu{equ:atlas-boundary}
\cA=\big\{\cV_\al=(V_\al,G_\al,\psi_\al)\big\}_{\al\in S}\cup \big\{\cV_\beta=(V_\beta\!\times[0,a_\beta),G_\beta,\psi_\beta)\big\}_{\beta\in S'}
\eEqu
of compatible $n$-dimensional orbifold charts covering $M$ such that 
$$
\psi_\al(V_\al)\cap \partial M\!=\!\emptyset\quad \forall \al\!\in\!S\quad\tn{and} \quad \psi_\beta^{-1}(\partial M)\!=\!V_\beta\!\times\! \{0\}\quad \forall\beta \!\in\! S'. 
$$
With this definition, the restriction 
$$
\cA|_{\partial M}\!=\!\big\{\cV_\beta=(V_\beta,G_\beta,\psi_\beta)\big\}_{\beta\in S'}
$$
defines an  orbifold structure $\partial \cM$ on $\partial M$ in the sense of Definition~\ref{def:orbi-atlas}. Oriented orbifold structures are defined similarly to Definition~\ref{def:orbi-oriented} and if an orbifold structure $\cA$ on $(M,\partial M)$ is oriented, by (\ref{equ:bordered-refinement}),  it induces an orientation of the induced orbifold structure $\partial \cM$. The latter depends on a sign convention; we consider the induced orientation on the boundary of (\ref{equ:B-charts}) such that
$$
T(V\!\times\! [0,a))|_{V\times \{0\}} = TV \oplus \R\cdot \frac{\partial}{\partial t},
$$
where $t$ is parametrizing $[0,a)$, is an isomorphism of oriented vector spaces. We can define other notions like orbifold morphisms, embeddings, orbibundles, multisections, etc. on bordered orbifolds, similarly. 

\noindent
Given an orbibundle $\cE$ over a bordered orbifold $(\cM,\partial\cM)$ and a section $s\colon \cM\lra \cE$ with compact support, the analogue of Proposition~\ref{pro:chain-cycle} gives us a rational homology class $e(\cE,s)$ in the relative singular homology group $H_{n}(M,\partial M,\Q)$, where $n\!=\!\dim(\cM)\!-\!\tn{rank}(\cE)$. The image of $e(\cE,s)$ under 
$$
H_{n}(M,\partial M,\Q)\lra H_{n-1}(\partial M,\Q)
$$
is the Euler class of the restriction $s|_{\partial \cM}$.

\noindent
The definition of a Kuranishi structure, an DGS for it, and other related notions analogously extend to the case of compact metrizable spaces with ``boundary". 

\bDef{def:cobo-Kur}
For $i\!=\!0,1$, let $\cK_i$ be a Kuranishi structure on the compact metrizable space $M_i$  and $\mfL_i$ be an DGS for $\cK_i$ as in Definition~\ref{level-system}. 
A \textbf{cobordism} between $(M_0,\cK_0,\mfL_0)$ and $(M_1,\cK_1,\mfL_1)$ consists of a compact metrizable space $\wt{M}$ with boundary 
$$ 
\wt{M}\!\supset\! \partial \wt{M}\cong M_0\sqcup M_1,
$$
a bordered Kuranishi structure $\wt\cK$ on $(\wt{M},\partial\wt{M})$, and an DGS $\wt\mfL$ for $\wt\cK$, such the restriction of $(\wt\cK,\wt\mfL)$ to $M_i$, with $i\!=\!0,1$, coincides with $(\cK_i,\mfL_i)$. In the oriented case, we require the orientation induced by the orientation of $\wt\cK$ on $\cK_0$ to be equal to its given orientation and the one induced on $\cK_1$ to be the opposite orientation.
\eDef

\noindent
In Definition~\ref{def:cobo-Kur}, we may forget about the DGS to define a cobordism between two Kuranishi structures. In particular, let $M$ be a fixed metrizable space as in the statement of Definition~\ref{def:Kur-structure}, and $\cK_0$ and $\cK_1$ be two Kuranishi structures on $M$. Then, by a \textbf{deformation equivalence} between $\cK_1$ and $\cK_2$, we mean a cobordism for which $\wt{M}\!=\!M\!\times\! [0,1]$.  
Similarly, let $\mfL_0$ and $\mfL_1$ be two dimensionally graded systems for a Kuranishi structure $(M,\cK)$. By $\cK\!\times\![0,1]$ we mean the product bordered Kuranishi structure on $M\!\times\![0,1]$ obtained by taking products of Kuranishi charts of $\cK$ with\footnote{More precisely, to comply with (\ref{equ:bordered-refinement}), by $[0,a)$ and $(b,1]\!\cong\![0,1-b)$, where $b\!<\!a$ are some fixed real numbers.} $[0,1]$. Then, by a deformation equivalence between $\mfL_1$ and $\mfL_2$, we mean a cobordism for which $\wt\cM\!=\!\cM\!\times\! [0,1]$ and $\wt\cK\!=\!\cK\!\times\![0,1]$.

\bRem{rmk:MW-cobordism}
In \cite[Section 6.2]{MW2}, a deformation equivalence between two Kuranishi structures is called a  ``concordance". McDuff-Wehrheim also use the notion of ``commensurability" which is an equivalence relation generated by calling $\cK_0$ and $\cK_1$ directly commensurate if they are both contained in the same Kuranishi structure $\cK_{01}$ on $M$. By \cite[Lemma 6.2.16]{MW2}, commensurability implies concordance (deformation equivalence), so it is a stronger notion. For a more detailed study of various notions of cobordism of Kuranishi structures we refer to \cite[Section 6.2]{MW2}.
\eRem

\noindent
We are now ready to state the existence result for dimensionally graded systems.

\bThm{exist-level-system}
Let $(M,\cK$) be a Kuranishi space\footnote{i.e. $\cK$ is a Kuranishi structure on $M$ with tangent bundle.} as in Definition~\ref{def:Kuranishi-tangent}. Then $(M,\cK)$ admits a plethora of ``natural" Hausdorff\footnote{With respect to the quotient topology. Unless we specifically mention, the topology considered on $Y(\mfL)$ is the quotient topology. } dimensionally graded systems. If $(M,\cK)$ is oriented, then the resulting natural DGS are also oriented.
Moreover, every two of such DGS are (orientably) ``deformation equivalent". \eThm
\noindent
We prove this theorem in the next section. As the proof shows, the same statements\footnote{i.e. existence of Hausdorff DGS and their deformation equivalence.} hold without the tangent bundle assumption. 
The word ``natural" in the statement means that the orbibundles 
$\{\cE(i)\!\lra\!\cY(i)\}_{i\in\mfB}$ are built by gluing Kuranishi charts of the same dimension along the overlaps. 
For the more general case of good coordinate systems, the construction of Hausdorff thickenings via the shrinking process is done in \cite{FOOO-shrink}.

\subsection{Existence of DGS (proof of Theorem~\ref{exist-level-system})}\label{sec:existence}
Our construction of natural dimensionally graded systems, out of a given Kuranishi structure or space $(M,\cK)$, is by reverse induction on $i\!\in\!\mfB$, where $\mfB\!\subset\!\Z$ is the index set of (\ref{equ:def-of-B}). With notation as in (\ref{L-starata}), starting from the maximal element $i_{\tn{max}}\!\in\! \mfB$, we first construct an orbibundle $\cE(i_{\tn{max}})\!\lra\!\cY(i_{\tn{max}})$ such that $M(i_{\tn{max}})\!\subset \!F(i_{\tn{max}})$ -recall that $M(i_{\tn{max}})\!\subset\! M$ is compact- and then proceed through lower values of $i\!\in\! \mfB$ inductively.  

\bLem{lem:compact-level}
Let $\cK$ be an $n$-dimensional Kuranishi structure on $M$.
For every $k\!\in\! \mfB$ and every compact subset $K\!\subset\! M(k)$, there exists an orbibundle $\pr\colon \cE\!\lra\! \cY$ together with a section $s\colon \cY\!\lra\! \cE$ and an open embedding $\ov\psi\colon\ov{s}^{-1}(0)\lra M$ with footprint $F$ such that $K\!\subset\! F$, $\tn{rank}~\cE\!=\!k$, $\dim \cY\!=\!n\!+\!k$, compatibility conditions\footnote{Between Conditions~\ref{l:compatible}(a) and~\ref{l:compatible}(b), only \ref{l:compatible}(a) is relevant in this case.}    \ref{l:compatible} of Definition~\ref{level-system} hold, and
$$
F\subset \bigcup_{\ov{p}\in K} F_p.
$$
If $(M,\cK)$ is a Kuranishi space then Condition~\ref{l:tangent} also holds.
\eLem
\noindent
We prove this lemma by induction, via the following lemma. In the light of  Example~\ref{exa:global-kur}, we call a tuple $(\pr\colon\!\cE\!\lra\!\cY,s,\psi)$ as in Lemma~\ref{lem:compact-level} an $n$-dimensional  \textbf{pure orbibundle structure} around $K$. In \cite[Definition 7.3]{FOOO-detail}, this is called a ``pure orbifold structure" around $K$. 

\bLem{lem:compact-level-2}
Let $\cK$ be an $n$-dimensional Kuranishi structure on $M$.
Give $k\!\in\! \mfB$ and compact subsets $K_1,K_2\!\subset\! M(k)$, for $a\!=\!1,2$, suppose $(\pr_a\colon\!\cE_a\!\lra\!\cY_a,s_a,\psi_a)$ is an $n$-dimensional pure orbibundle structure of rank $k$ with footprint $F_a$ around $K_a$ as described in Lemma~\ref{lem:compact-level}. Then $K\!=\!K_1\cup K_2$ admits a pure orbibundle structure around it.
\eLem

\newtheorem*{proofoflem:compact-level}{Proof of Lemma~\ref{lem:compact-level} via Lemma~\ref{lem:compact-level-2}}
\begin{proofoflem:compact-level}
We cover $K$ with finitely many Kuranishi charts of $\cK$ cantered at points of $K$,
$$
\lrc{ \cU_{p_i}\!=\!(\pr_{p_i}\colon\!U_{p_i}\!\lra\!V_{p_i},G_{p_i},s_{p_i},\psi_{p_i})}_{i\in[N]},\qquad \ov{p}_i\!\in\! K~~~\forall i\!\in\! [N].
$$
If $N\!=\!1$, we can simply take $\cE=\cU_{p_1}$.
If $N\!>\!1$, we fix a set of relatively compact open subsets $K_{p_i}^\circ\subset F_{p_i}$ such that 
$$
K\subset \bigcup_{i\in [N]} K^\circ _{p_i}.
$$
Let $K_{p_i}$ be the closure of $K_{p_i}^\circ$ in $F_{p_i}$. We then proceed by induction. By the $N\!=\!1$ case discussed above, $K\cap K_{p_1}$ admits such pure orbibundle structure around it. Assume $K\cap \big(\cup_{i=1}^{\ell-1} K_{p_i}\big)$ admits such pure orbibundle structure around it. Then by Lemma~\ref{lem:compact-level-2}, 
$$
K\cap \bigg( \bigcup_{i=1}^{\ell} K_{p_i}\bigg)=\bigg( K\cap \bigg(\bigcup_{i=1}^{\ell-1} K_{p_i}\bigg)\bigg) \cup \big( K\cap K_{p_\ell}\big)
$$
admits such orbibundle structure around it. 
\qed
\end{proofoflem:compact-level}

\newtheorem*{proofoflem:compact-level-2}{Proof of Lemma~\ref{lem:compact-level-2}}
\begin{proofoflem:compact-level-2}
We assume $\cK$ has tangent bundle, i.e. $(M,\cK)$ is a Kuranishi space; the general case is a simpler version of this. For every 
$\ov{p}\!\in \!F_1 \cap F_2$, by definition, there exist sub-charts $p\!\in\! \cU_{p,i}\!\subset \cU_p$, with $i\!=\!1,2$, such that the compatibility conditions \ref{l:compatible} and \ref{l:tangent} of Definition~\ref{level-system} for these sub-charts hold. Since both $V_{p,1}$ and $V_{p,2}$ are $G_p$-invariant neighborhoods of $p\!\in\! V_p$,  $V_{p,12}\!=\! V_{p,1}\cap V_{p,2}$ is $G_p$-invariant. If $V_{p,12}$ is disconnected, we replace it with the connected component containing the unique point $p$. Then $\cU_{p,12}=\cU_p|_{V_{p,12}}$ is a sub-chart of $\cU_p$ which can be embedded into both $\cY_1$ and $\cY_2$.

\noindent
Shrink each $\cU_{p,12}$ into some relatively compact sub-chart\footnote{The superscript ``m" used for these sub-charts denotes for ``main".} $\cU^m_{p,12}$ (i.e. $V^m_{p,12}\subset V_{p,12}$ is relatively compact), and take finitely many such points in 
$$
K_{12}= K_1\cap K_2\subset M(k),
$$ 
say $\lrc{\ov{p}_\al}_{\al\in I}\!\subset\! K_{12}$, such that 
\bEqu{equ:K12-alpha}
K_{12}\subset \bigcup_{\al\in I} F^m_{p_\al,12}.
\eEqu

\noindent
\textbf{Claim 1.} \textit{For every $\ov{q}\!\in\! K_{12}$, there exists a sub-chart\footnote{The superscript ``t" used for these sub-charts denotes for ``transition".} $\cU_{q,12}^t\!\subset \!\cU_{q,12}$ such that for $j\!=\!1,2$ and $\al\!\in\! I$,  if 
\bEqu{equ:main-transitive}
 \ov\phi_{qj}(\ov{V}^t_{q,12})\cap \ov\phi_{p_\al j}(\ov{V}^m_{p_\al,12})\neq \emptyset ,
\eEqu
 then $\ov{q}\in F_{p_\al}$ and $\ov{V}^t_{q,12}\subset \ov\phi_{qp_\al}^{-1}(\ov{V}_{p_\al,12})$.}

 \bProof
 Decompose $I$ into disjoint unions $I=I^{\tn{in}}\sqcup I^{\tn{out}}$ such that
$$
I^{\tn{in}}= \lrc{ \al \in I\colon ~\ov{q}\in \tn{cl}(F_{p_\al,12}^m))  }.
$$
Note that $\ov\psi_{p_\al,12}\colon\ov{s}_{p_\al,12}^{-1}(0)\!\lra\! F_{p_\al,12}$ is a homeomorphism, therefore we can think of $\ov{q}\in F_{p_\al,12}$ as a point in $$
\ov{s}_{p_\al,12}^{-1}(0)\subset \ov{V}_{p_\al,12}= V_{p_\al,12}/G_{p_\al,12}.
$$ 
Moreover, if $\ov{q}\!\in\! \tn{cl}(F_{p_\al,12}^m)$, since $\tn{cl}(F_{p_\al,12}^m)\!\subset\! F_{p_\al}$, the coordinate change map $\Phi_{qp_\al}$ of Kuranishi structure is defined.   
Then we take $\cU_{q,12}^t\subset \cU_{q,12}$ to be the sub-chart over the preimage\footnote{We would only consider the connected component containing $q$, if it is disconnected.} of 
$$
\ov{V}_{q,12}^t= 
\bigcap_{j=1,2} \bigg(\bigcap_{\al\in I^{\tn{in}}}    \bigg( \ov{\phi}_{qp_\al}^{-1}(\ov{V}_{p_\al,12} )\bigg)                \cap  
\bigcap_{\al\in I^{\tn{out}}}      \bigg(\ov{\phi}_{qj}^{-1}\big(Y_j\setminus \tn{cl}\lrp{ \ov{\phi}_{p_\al j}(\ov{V}^m_{p_\al,12})} \big)\bigg)            \bigg) .                                      
$$
It is easy to see that $\cU_{q,12}^t$ has the desired property.
\eProof

 \noindent
We take finitely many such points $\ov{q}$ in $K_{12}$, say $\lrc{\ov{q}_\beta}_{\beta\in J}\!\subset\! K_{12}$, such that 
\bEqu{equ:K12-beta}
K_{12}\subset \bigcup_{\beta\in J} F^t_{q_\beta,12}. 
\eEqu
Then, for $(i,j)\!=\!(1,2)$ and $(2,1)$ define
\bEqu{equ:Yij}
Y_{i,j}=\bigcup_{\al\in I}\bigcup_{\beta\in J}\lrp{
\ov\phi_{q_\beta i}(\ov{V}^t_{q_\beta,12}) \cap 
\ov\phi_{p_\al i}(\ov{V}^m_{p_\al,12})} \subset Y_i.
\eEqu
These are open subsets of $Y_1$ and $Y_2$. 
By restricting  the orbibundles $\cE_1\lra \cY_1$ and $\cE_2\!\lra\! \cY_2$ to $Y_{1,2}$ and $Y_{2,1}$, we obtain sub-orbibundles $\cE_{1,2} \!\lra\! \cY_{1,2}$ and $\cE_{2,1}\!\lra\! \cY_{2,1}$, respectively.
Let $s_{1,2}$, $s_{2,1}$,  and $F_{1,2}\!\cong\! F_{2,1}$ be the restriction of corresponding objects, respectively.
By (\ref{equ:K12-alpha}) and (\ref{equ:K12-beta}), $F_{1,2}\!\cong\! F_{2,1}$ is an open neighborhood of $K_{12}$.
 \vskip.1in
 \noindent
\textbf{Claim 2.} \textit{The tuples $(\cE_{1,2} \!\lra\!\cY_{1,2},s_{1,2},\psi_{1,2})$ and $(\cE_{2,1}\!\lra\! \cY_{2,1},s_{2,1},\psi_{2,1})$ are naturally isomorphic.}
\bProof
We show that the maps 
$$
    \xymatrix{ E_{1,2} \ar[rr]^{\ov{\mfD \varphi}}  \ar[d]^{\pr_{1,2}}&& E_{2,1} \ar[d]^{\pr_{2,1}} \\
               Y_{1,2} \ar[rr]^{\ov\varphi} && Y_{2,1} }
$$
which send
\bEqu{x-image}
 \ov{x}= \ov\phi_{q_\beta 1}(\ov{x}_\beta) = 
\ov\phi_{p_\al 1}(\ov{x}_\al),\quad \ov{x}_\beta \in \ov{V}^t_{q_\beta,12},~\ov{x}_\al\in \ov{V}^m_{p_\al,12},
\eEqu
to $\ov{y}=\ov\phi_{q_\beta 2}(\ov{x}_\beta) $ and
\bEqu{v-image}
\ov{v}= \ov{\mfD\phi}_{q_\beta 1}(\ov{v}_\beta) = 
\ov{\mfD\phi}_{p_\al 1}(\ov{v}_\al),\quad \ov{v}_\beta \in \ov{U}^t_{q_\beta,12},~\ov{v}_\al\in \ov{U}^m_{p_\al,12},
\eEqu
to $\ov{u}=\ov{\mfD\phi}_{q_\beta 2}(\ov{v}_\beta) $  are well-defined continuous maps and lift to an  orbibundle isomorphism $(\mfD \varphi,\varphi)$ commuting with the section and footprint maps.

\noindent
Given $\ov{x}\!\in \!Y_{1,2}$,  suppose there are two different pairs of indices $(\al,\beta),(\al',\beta')\!\in\! I\!\times\! J$ such that 
 $$
 \aligned
& \ov{x}= \ov\phi_{q_\beta 1}(\ov{x}_\beta) = 
\ov\phi_{p_\al 1}(\ov{x}_\al),\quad \ov{x}_\beta \in \ov{V}^t_{q_\beta,12},~\ov{x}_\al\in \ov{V}^m_{p_\al,12},\\
&\ov{x}= \ov\phi_{q_{\beta'} 1}(\ov{x}_{\beta'}) = 
\ov\phi_{p_{\al'} 1}(\ov{x}_{\al'}),\quad \ov{x}_{\beta'} \in \ov{V}^t_{q_{\beta'},12},~\ov{x}_{\al'}\in \ov{V}^m_{p_{\al'},12}.
\endaligned
$$
It follows from (\ref{equ:main-transitive}) that 
$$
\aligned
&\ov{x}\in \ov\phi_{q_\beta 1}(\ov{V}^t_{q_\beta,12})\cap  \ov\phi_{p_{\al'} 1}(\ov{V}^m_{p_{\al'},12})\neq \emptyset, \\
&\ov{x}\in \ov\phi_{q_{\beta'} 1}(\ov{V}^t_{q_{\beta'},12})\cap \ov\phi_{p_{\al} 1}(\ov{V}^m_{p_{\al},12})\neq \emptyset;
\endaligned
$$
therefore, 
$$
\ov{q}_\beta\in F_{p_{\al'}},~
\ov{q}_{\beta'}\in F_{p_{\al}},~
\ov{V}^t_{q_\beta,12}\subset \ov\phi_{q_\beta p_{\al'}}^{-1}(\ov{V}_{p_{\al'},12}),~~
\tn{and}~~
\ov{V}^t_{q_{\beta'},12}\subset \ov\phi_{q_{\beta'}p_\al}^{-1}(\ov{V}_{p_\al,12}).
$$
By Definition~\ref{level-system}.\ref{l:compatible}, for $i\!=\!1,2$, restricted to 
$$
\ov{V}^t_{q_\beta,12}\subset \ov{\phi}_{q_\beta p_{\al'}}^{-1}(\ov{V}_{p_{\al'},12})\cap \ov{V}_{q_\beta,12}\subset \ov{\phi}_{q_\beta p_{\al'}}^{-1}(\ov{V}_{p_{\al'},i})\cap \ov{V}_{q_\beta,i}
$$ 
we have 
\bEqu{equ:different-compositions}
\ov{\phi}_{p_{\al'}1}\circ \ov{\phi}_{q_\beta p_{\al'}}=\ov{\phi}_{q_\beta 1}\quad\tn{and}\quad \ov{\phi}_{p_{\al'}2}\circ \ov{\phi}_{q_\beta p_{\al'}}=\ov{\phi}_{q_\beta 2}.
\eEqu
From the first identity in (\ref{equ:different-compositions}) and the fact that $\ov{\phi}_{p_{\al'}1}$ is one-to-one we conclude  that
$$
\ov\phi_{q_\beta p_{\al'}}(\ov{x}_\beta)=\ov{x}_{\al'}.
$$
Then from the second identity in (\ref{equ:different-compositions}) we get
$$
\ov\phi_{q_\beta 2}(\ov{x}_\beta)=\ov\phi_{p_{\al'} 2}(\ov\phi_{q_\beta p_{\al'}}(\ov{x}_\beta))=
\ov\phi_{p_{\al'} 2}(\ov{x}_{\al'}).
$$
Therefore, (\ref{x-image}) is independent of the choice of $(\al,\beta)$ or $(\al',\beta')$. The case of  (\ref{v-image}) is similar.

\noindent
By definition, over 
$$
\ov\phi_{q_\beta 1}(\ov{V}^t_{q_\beta,12}) \cap 
\ov\phi_{p_\al 1}(\ov{V}^m_{p_\al,12}) \subset Y_1,
$$ 
$(\ov{\mfD\varphi},\ov\varphi)$ is equal to the restriction of $\ov\Phi_{p_\al 2}\circ \ov\Phi_{p_\al 1}^{-1}$; therefore, $(\ov{\mfD\varphi},\ov\varphi)$ is locally liftable to  smooth orbibundle isomorphisms $\Phi_{p_\al 2}\circ \Phi_{p_\al 1}^{-1}$ between equal rank orbibundle charts.
Compatibility of $(\mfD\varphi,\varphi)$ with the section and footprint maps follow from the corresponding set of assumptions in Definition~\ref{level-system}.\ref{l:compatible}.
\eProof

\noindent
In order to finish the proof of Lemma~\ref{lem:compact-level-2}, we identify $\pr_1\colon\cE_1\!\lra\! \cY_1$ and $\pr_2\colon\cE_2\!\lra\! \cY_2$ along the isomorphic sub-orbibundles  $\pr_{i,j}\colon\cE_{i,j}\lra \cY_{i,j}$, with $(i,j)=(1,2)$ and $(2,1)$, to obtain an  orbibundle 
$\pr\colon \cE\!\lra\! \cY$ together with a section $s\colon \cY\!\lra\! \cE$ and an open embedding $\ov\psi\colon\ov{s}^{-1}(0)\!\lra\! M$ with footprint $F\!=\!F_1\!\cup\! F_2$ such that $K\!=\!K_1\cup K_2\!\subset\! F$, $\tn{rank}~\cE\!=\!k$, $\dim \cY\!=\!n\!+\!k$, and the compatibility conditions \ref{l:compatible} and \ref{l:tangent} of Definition~\ref{level-system} hold.
However, without further restrictions on $Y_1$, $Y_2$, and $Y_{1,2}$, the resulting quotient space 
$$
Y= \bigg(Y_1\coprod Y_2\bigg)/Y_{1,2}\cup Y_{2,1}
$$
may not be a Hausdorff topological space, or the resulting continuous map $\ov\psi$ may not be injective. In order to obtain a pure orbibundle structure around $K$ as above, we further shrink $Y_1$, $Y_2$, and $Y_{1,2}\cong Y_{2,1}$ in several steps to get the right orbibundle structure.
\vskip.1in
\noindent
\textbf{First shrinking to make $\ov\psi$ injective.} For $i\!=\!1,2$, choose open sets $F_i^\one\!\subset\! F_i$ such that 
$$
K_i\subset F_i^\one~\tn{and}\quad F_{12}^\one= F_1^\one\cap F_2^\one \subset F_{1,2}=F_{2,1}.
$$ 
Take open sets $Y_i^\one\!\subset\! Y_i$  such that 
$$
Y_i^\one\cap \ov{s}_i^{-1}(0)\!= \!\ov\psi_i^{-1}(F_i^\one).
$$
Let 
$$
Y^\one_{1,2}=Y_{1,2}\cap Y_1^\one\cap \ov\varphi^{-1}(Y_{2,1}\cap Y_2^\one)~~\tn{and}~~Y^\one_{2,1}=Y_{2,1}\cap Y_2^\one\cap \ov\varphi(Y_{1,2}\cap Y_1^\one).
$$
Let $\cE_1^\one$, $\cE_2^\one$, and $\cE_{1,2}^\one\!\overset{\mfD\varphi}{\cong}\! \cE_{2,1}^\one$, together with the similarly denoted sections and footprint maps, be restrictions to $Y_1^\one$, $Y_2^\one$, and $Y_{1,2}^\one\overset{\ov\varphi}{\cong} Y_{2,1}^\one$ of the corresponding orbibundles, respectively. Since
$$
F_1^\one\cap F_2^\one= F_{1,2}^\one=F_{2,1}^\one,
$$
we conclude that the the induced map $\ov\psi$ on 
$$
\ov{s}^{-1}(0)\cap \lrp{Y^\one= \lrp{Y_1^0\cup_{Y_{1,2}^\one\cong Y_{2,1}^\one} Y_2^\one}}
$$ 
is injective.

\noindent
\textbf{Second shrinking to make the glued space Hausdorff.} \vskip.1in
\noindent
\textbf{Claim 3.} 
\textit{For $i\!=\!1,2$, there exist relatively compact open subsets $Y_i^\two\subset Y_i^\one$  and relatively compact open subsets $W_1\subset Y_{1,2}^\one$ and $W_2\subset Y_{2,1}^\one$ such that $K_i \subset F_i^\two$,
\bEqu{equ:sy}
\tn{cl}\!\lrp{Y_1^\two\cap \ov\varphi^{-1}(Y_{2,1}^\one\cap Y_2^\two)}\subset W_1,\quad 
\tn{cl}\!\lrp{Y_2^\two\cap \ov\varphi(Y_{1,2}^\one\cap Y_1^\two)}\subset W_2,
\eEqu
and $\ov\varphi(W_1)=W_2$.}

\bProof
We start from an arbitrary pair of relatively compact open subsets $\wt{Y}_i^\two\subset Y_i^\one$ and find $W_i$ such that the following pair of weaker conditions
\bEqu{equ:sy-weak}
\aligned
&\ov{s}_1^{-1}(0)\cap \tn{cl}\!\lrp{\wt{Y}_1^\two\cap \ov\varphi^{-1}(Y_{2,1}^\one\cap \wt{Y}_2^\two)}\subset W_1,\\
&\ov{s}_2^{-1}(0)\cap \tn{cl}\!\lrp{\wt{Y}_2^\two\cap \ov\varphi(Y_{1,2}^\one\cap \wt{Y}_1^\two)}\subset W_2,
\endaligned
\eEqu
hold. Then, we further shrink $\wt{Y}_i^\two$ such that the original inclusions of (\ref{equ:sy}) hold.
\noindent
Note that the left-hand side terms in (\ref{equ:sy-weak}) are compact subsets of $Y_1^\one$ and $Y_2^\one$, respectively. The existence of such $W_1$ and $W_2$ follows, if we  they are also subsets of $Y_{1,2}^\one$ and $Y_{2,1}^\one$, respectively.
Since $\wt{Y}_1^\two\subset Y_1^\one$ is relatively compact and $\ov{s}_1^{-1}(0)\cap Y_1^\one$ is closed
$$
\ov{s}_1^{-1}(0)\cap \tn{cl}\!\lrp{\wt{Y}_1^\two\cap \ov\varphi^{-1}(Y_{2,1}^\one\cap \wt{Y}_2^\two)}\subset \tn{cl}\!\lrp{\wt{Y}_1^\two\cap Y_{1,2}^\one}\subset Y_1^\one.
$$
For every 
$$
\ov{x}\in \ov{s}_1^{-1}(0)\cap \tn{cl}\!\lrp{\wt{Y}_1^\two\cap \ov\varphi^{-1}(Y_{2,1}^\one\cap \wt{Y}_2^\two)}
$$ 
there exists a sequence 
\bEqu{page-ref}
(\ov{x}_k)_{k=1}^\infty \subset \ov{s}_1^{-1}(0)\cap \ov\varphi^{-1}(Y_{2,1}^\one\cap \wt{Y}_2^\two) \subset Y_{1,2}^\one
\eEqu
such that $\lim_{k\lra \infty} \ov{x}_k\! =\!\ov{x}$. Since $\wt{Y}_2^\two\!\subset\! Y_2$ is relatively compact, there exist $\ov{y}\!\in\! Y_2^\one$ such that the sequence 
$$
(\ov\varphi(\ov{x}_k))_{k=1}^\infty \subset Y_{2,1}^\one\cap \wt{Y}_2^\two
$$ 
converges to $\ov{y}$; moreover, $s_2(\ov{y})\!=\!0$. Since 
$F_1^\one\cap F_2^\one=F_{1,2}^\one$ and $\ov\psi_1(\ov{x})=\ov\psi_2(\ov{y})$,
we conclude that $\ov{x}\in Y_{1,2}^\one$. The argument for the inclusion of 
$$
\ov{s}_2^{-1}(0)\cap \tn{cl}\!\lrp{\wt{Y}_2^\two\cap \ov\varphi(Y_{1,2}^\one\cap \wt{Y}_1^\two)}
$$ 
in $Y_{2,1}^\one$ is similar. Fix such $W_1$, $W_2$.

\noindent
We now shrink $\wt{Y}_1^\two$ and $\wt{Y}_2^\two$ into relatively compact open subsets $Y_1^\two$ and $Y_2^\two$, still covering $K_1$ and $K_2$, respectively, such that (\ref{equ:sy}) hold. 
\vskip.1in
\noindent
For every $\ov{x}\!\in\! \ov{s}_1^{-1}(0)\cap \wt{Y}_1^\two$ such that $\ov\psi_1(\ov{x})\!\in\! K_1$, we can find a sufficiently small open neighborhood $W_{\ov{x}}\subset \wt{Y}^\two_1$ such that
\bEnum
\item 
$
\ov{x}\notin \tn{cl}\!\lrp{\wt{Y}_1^\two\cap \ov\varphi^{-1}(Y_{2,1}^\one\cap \wt{Y}_2^\two)}
\Longrightarrow 
\tn{cl}\!\lrp{W_{\ov{x}}} \cap  \tn{cl}\!\lrp{\wt{Y}_1^\two\cap \ov\varphi^{-1}(Y_{2,1}^\one\cap \wt{Y}_2^\two)} =\emptyset,
$
\item 
$
\ov{x}\in \tn{cl}\!\lrp{\wt{Y}_1^\two\cap \ov\varphi^{-1}(Y_{2,1}^\one\cap \wt{Y}_2^\two)}
\Longrightarrow
\tn{cl}\!\lrp{W_{\ov{x}}} \subset W_1.
$ 
\eEnum
Cover $K_1$ by finitely many such neighborhoods and let $Y_1^\two$ be their union. We construct $Y_2^\two$ similarly. If 
$$
\ov{y}\! \in\! \tn{cl}\!\lrp{Y_1^\two\cap \ov\varphi^{-1}(Y_{2,1}^\one\cap Y_2^\two)},
$$ 
then there is one of those finitely many $\ov{x}$ such that 
$$
\ov{y} \in \tn{cl}\!\lrp{W_{\ov{x}}\cap \ov\varphi^{-1}(Y_{2,1}^\one\cap Y_2^\two)}\subset \tn{cl}\!\lrp{W_{\ov{x}}} \cap  \tn{cl}\!\lrp{\wt{Y}_1^\two\cap \ov\varphi^{-1}(Y_{2,1}^\one\cap \wt{Y}_2^\two)} .
$$ 
In this case, $\ov{x}$ should be of type (2), hence, $\ov{y}\in W_1$.
\eProof

\noindent
Let
$$
\aligned
Y^\two=
\bigg(
&Y_1^{\tn{out}}= \lrp{Y_1^\two\setminus  \tn{cl}\!\lrp{Y_1^\two\cap \ov\varphi^{-1}(Y_{2,1}^\one\cap Y_2^\two)}} \coprod\\
& 
Y_2^{\tn{out}}= \lrp {Y_2^\two \setminus \tn{cl}\!\lrp{Y_2^\two\cap \ov\varphi(Y_{1,2}^\one\cap Y_1^\two)} }
\coprod \lrp{ (W_1\overset{\ov\varphi}{\cong} W_2)}
\bigg)/\sim~,
\endaligned
$$ 
where the equivalence relation is given by the intersections 
\bEqu{equ:q-rel}
Y_1^{\tn{out}}\supset  W_{11}= Y_1^{\tn{out}} \cap W_1 \subset W_1, \quad 
Y_2^{\tn{out}}\supset W_{22}= Y_2^{\tn{out}} \cap W_2 \subset W_2.
\eEqu
Note that  $\ov\varphi (W_{11})\cap W_{22}=\emptyset$. In fact, 
$$
\ov\varphi (W_{11})=\ov\varphi\lrp{
(Y_1^\two\cap W_1)\setminus  \tn{cl}\!\lrp{Y_1^\two\cap \ov\varphi^{-1}(Y_{2,1}^\one\cap Y_2^\two)}
}
 \subset \ov\varphi\lrp{ Y_{1,2}^\one \cap Y_1^\two};
 $$
therefore, $\ov\varphi (W_{11})\cap Y_2^{\tn{out}} \!=\!\emptyset$.
Since $W_1$ and $Y_1^{\tn{out}}$ are open subsets of the same Hausdorff space $Y_1$, the quotient topology with respect to only the first relation in (\ref{equ:q-rel}) is Hausdorff. Similarly, since $W_2$ and $Y_2^{\tn{out}}$ are open subsets of the same Hausdorff space $Y_2$, the quotient topology with respect to only the second relation in (\ref{equ:q-rel}) is Hausdorff. Finally, since the intersections $W_{11}$ and $W_{22}$ in $W_1\overset{\ov\varphi}{\cong} W_2$ are disjoint, the quotient topology on $Y^\two$ with respect to both relations in (\ref{equ:q-rel}) is Hausdorff. Moreover, since both $Y_1$ and $Y_2$ are metrizable, $Y^\two$ is metrizable; see Remark~\ref{rem:top-remark}.
\vskip.1in
\noindent
We identify the sub-orbibundles $\cE_1|_{Y_1^{\tn{out}}}$, $\cE_2|_{Y_2^{\tn{out}}}$, and $\cE_1|_{W_1}\cong \cE_2|_{W_2}$ to obtain an open orbibundle $\pr\colon \cE^\two\lra \cY^\two$ together with a section $s^\two\colon \cY^\two\lra \cE^\two$ and an open embedding $\ov\psi^\two\colon{\ov{s}^\two}^{-1}(0)\lra M$ with footprint $F^\two$ that contains $K\!=\!K_1\cup K_2$, $\tn{rank}~\cE^\two\!=\!k$, $\dim \cY^\two\!=\!n\!+\!k$, and the compatibility conditions of Definition~\ref{level-system}.(8) hold. This finishes the proof of Lemma~\ref{lem:compact-level-2}.

\noindent
\end{proofoflem:compact-level-2}

\newtheorem*{proofoflevel-system}{Proof of Theorem~\ref{exist-level-system}}
\begin{proofoflevel-system}
Proof of existence is by reverse induction on the index set $\mfB$.
By the argument after Definition~\ref{def:Kur-structure}, for the largest index $i_{\tn{max}}\in \mfB$, $M(i_{\tn{max}})\subset M$ is compact. Therefore, by Lemma~\ref{lem:compact-level}, there exists a pure orbibundle structure around $K\!=\!M(i_{\tn{max}})$.

\noindent
For every $i\!\in\! \mfB$, let 
$$
\mfB(>\!i)\!=\!\{j\!\in\! \mfB\colon j\!>\!i\}.
$$
For some fixed $i\!\in\! \mfB$ and every $j \!\in\! \mfB(>\! i)$, suppose we have built an orbibundle $\pr_{j}\colon \cE(j) \lra \cY(j)$ together with a section $s_{j}\colon \cY(j)\lra \cE(j)$ and an open embedding $\ov\psi_{j}\colon\ov{s}_{j}^{-1}(0)\lra M$ with footprint $F(j)$ such that the conditions of Definition~\ref{level-system} and Definition~\ref{def:H-level-structure}, with $\mfB(>\! i)$ instead of $\mfB$, are satisfied; i.e.  
$$
\mfL(>\! i)\equiv \bigg(\{\cE(j),\cY(j)\}_{ j\in \mfB(> i)},\{\cE(j_1,j_2),\cY(j_1,j_2),\Phi_{j_1j_2}\}_{\substack{j_1,j_2\in\mfB(> i)\\j_1<j_2}}\bigg)
$$
is an DGS on
$$
M(>\! i)\subset F(>\!i)= \bigcup_{j \in \mfB(> i)} F(j) \subset M
$$ 
such that corresponding thickening $Y(\mfL(>\! i))$ as in (\ref{equ:thikening-space}) is Hausdorff.

\noindent
Let $K_{i}$ be a compact subset\footnote{Possibly empty, in which case we proceed to lower values of $\mfB$.} of $M(i)$ such that 
\bEqu{equ:inclusion}
\tn{cl}\!\lrp{M(\geq\! i)\!\setminus\! K_{i}} \subset F(>\! i).
\eEqu 
In fact, we can start from a relatively compact open subset $F'\!\subset\! F(>\! i)$ such that $M(>\! i)\!\subset\! F'$  and set $K_{i}\!=\!M(i)\!\setminus \!F'$.
Let 
$$
\xymatrix{
 \cE(i)   \ar@/^/[rrr]|{\pr_{i}}
 &&& \cY(i)  \ar@/^/[lll]|{s_{i}} }
 ,\quad \ov\psi_{i}\colon\ov{s}_{i}^{-1}(0)\lra F(i)\subset M,
$$ 
be a pure orbibundle structure as in Lemma~\ref{lem:compact-level} around $K_{i}\!\subset\! M(i)$. 
After possibly shrinking $F(i)$, we may assume 
$$
\tn{cl}\!\lrp{F(i)}\cap M(>\!i)=\emptyset;
$$ 
cf.~Definition~\ref{level-system}.(2). 
In what follows, we shrink each each $Y(j)$, for $ j \!\in\! \mfB(>\! i)$, and build the required transition sub-charts $\cE(i,j)\lra \cY(i,j)$, such that the conditions of Definition~\ref{level-system} and Definition~\ref{def:H-level-structure} for the extended DGS over $\mfB(\geq\!i)$ are satisfied. 

\noindent
We first shrink $\{ Y(j)\}_{j\in \mfB(> i)}$, respectively $Y(i)$, as in Definition~\ref{def:shrinking}, into relatively compact open subsets $\{ Y^\one(j)\}_{i\in \mfB(> i)}$, respectively $Y^\one_0(i)$, such that the shrunk system 
$$
\mfL^\one(>\! i)\equiv \big(\{\cE^\one(j),\cY^\one(j)\}_{ j\in \mfB(> i)},\{\cE^\one(j_1,j_2),\cY^\one(j_1,j_2),\Phi_{j_1j_2}\}_{\substack{j_1,j_2\in\mfB(> i)\\ j_1<j_2}}\big)
$$
as in Lemma~\ref{lem:shrinking}
still covers $M(>\!\!i)$, $K_{i}\!\subset\! F^\one_0(i)$, all the conditions of Definition~\ref{level-system} and Definition~\ref{def:H-level-structure}, with $\mfB(>\!i)$ instead of $\mfB$, are still satisfied, and the inclusion condition (\ref{equ:inclusion}) is still valid; i.e. the shrunk system still has all the properties of last three paragraphs. 
\noindent
Let 
$$
O=M(i)\cap F^\one_0(i) \cap F^\one(>\!i)\subset M.
$$ 
For every $\ov{p}\! \in\! \tn{cl}\!\lrp{O}\!\subset\! F(i)\cap F(>\!i)$, let $\cU_{p,O}\!\subset\! \cU_{p}$ be a sufficiently small sub-chart centered at $\ov{p}$ such that the following conditions hold. 
\bEnum
\item\label{l:inclusion} For every $j\!\in\! \mfB(>\!i)$, if $\ov{p}\!\in\! F(j)$, then $\cU_{p,O}\!\subset\! \cU_{p,j}$, where $\cU_{p,j}$ are the compatibility sub-charts with $\mfL(>\!i)$ of  Definition~\ref{level-system}.\ref{l:compatible} for $\cE(j)\!\lra \!\cY(j)$. 

\item\label{l:no-inclusion} For every $j\!\in \!\mfB(>\!i)$, if $\ov{p}\!\notin\!  \tn{cl}\!\lrp{F^\one(j)}$, then $F_{p,O}\cap F^\one(j) \!=\!\emptyset$; moreover, if $\ov{p}\!\in\! F(j)\!\setminus\!\tn{cl}\!\lrp{F^\one(j)}$, then 
\bEqu{equ:no-inclusion-1}
\ov\phi_{pj}(\ov{V}_{p,O})\cap Y^\one(j) =\emptyset.
\eEqu

\item\label{l:inclusion2} For every $j_1,j_2\!\in\! \mfB(>\!i)$, if $j_1\!<\!j_2$ and $\ov{p}\!\in\!  F(j_1)\cap  F(j_2)$, then 
\bEqu{equ:for(4)to-hold}
\ov\phi_{pj_1}(\ov{V}_{p,O})\subset Y(j_1,j_2);
\eEqu
\item\label{l:no-inclusion2} For every $j,k\!\in\! \mfB(>\!i)$, if $\ov{p}\!\in\!  F(j)$ but $\ov{p}\!\notin\! \tn{cl}\lrp{F^\one(k)}$, then 
\bEqu{equ:for(5)to-hold}
\begin{cases} 
\ov\phi_{pj}(\ov{V}_{p,O})\cap \ov\phi_{kj}(Y^\one(j,k))=\emptyset  &\mbox{if } j\!>\!k, \\ 
\ov\phi_{pj}(\ov{V}_{p,O})\cap Y^\one(j,k)=\emptyset & \mbox{if } j\!<\!k. 
 \end{cases} 
\eEqu
\item\label{l:inclusion-i-1} Finally, $\cU_{p,O}\subset \cU_{p,i}$. 
\eEnum

\noindent 
For every $\ov{p} \!\in\! \tn{cl}\!\lrp{O}\!\subset\! M$, shrink $\cU_{p,O}$ into a relatively compact sub-chart $\cU^m_{p,O}\subset \cU_{p,O}$.
For every finite set of points $S\!=\!\{\ov{p}_\al\}_{\al\in I}\! \subset\! \tn{cl}\!\lrp{O}$ and every $j\!\in\! \mfB(>\!i)$ let 
$$
I_j=\{\al\in I\colon \ov{p}_\al \in S\cap \tn{cl}\!\lrp{F^\one(j)}\}.
$$
Since $\tn{cl}\!\lrp{O}$ is compact, we can choose a finite set $S$ as above such that
\bEqu{equ:cover1}
\tn{cl}\!\lrp{O}\cap \tn{cl}\!\lrp{F^\one(j)} \subset \bigcup_{\al\in I_j} F^m_{p_\al,O}\quad \forall j\!\in\!\mfB(>\!i).
\eEqu
\vskip-.2in
\noindent
Similar to Claim~1 in the proof of Lemma~\ref{lem:compact-level-2}, for every $\ov{q}\in \tn{cl}\!\lrp{O}$ there also exists a sub-chart $\cU^t_{q,O}\!\subset\! \cU_{q,O}$ such that for every $j\!\in\! \mfB(>\!i)$, every $\al\! \in\! I_j$, and $a\!=\!i~\tn{or}~j$, 
\bEqu{equ:transition-condition}
\ov{\phi}_{q a}(\ov{V}^t_{q,O})\cap \ov\phi_{p_\al a}(\ov{V}_{p_\al,O}^m)\neq \emptyset \Longrightarrow 
\ov{q}\in F_{p_\al}~\tn{and}~ \ov{V}^t_{q,O}\subset \ov\phi^{-1}_{qp_\al}(\ov{V}_{p_\al,O}).
\eEqu

\noindent
Choose a finite set of points $T\!=\! \{\ov{q}_\beta\}_{\beta\in J}\! \subset\! \tn{cl}\!\lrp{O}$ such that 
\bEqu{equ:cover2}
\tn{cl}\!\lrp{O}\cap \tn{cl}\!\lrp{F^\one(j)} \subset \bigcup_{\beta\in J_j} F^t_{q_\beta,O}\quad\forall  j\!\in\!\mfB(>\!i).
\eEqu
Then, similar to (\ref{equ:Yij}), for every $j\!\in\!\mfB(>\!i)$, let
\bEqu{equ:overlap-charts}
Y(i,j)= \bigcup_{\beta \in J_j,\al\in I_j}\lrp{Y_{\al,\beta}= \lrp{ \ov\phi_{q_\beta i}(\ov{V}^t_{q_\beta,O})\cap \ov\phi_{p_\al i}(\ov{V}^m_{p_\al,O}) }}.
\eEqu
For every $j_1,j_2\! \in \mfB(>\! i)$, by (\ref{equ:overlap-charts}) and Condition~\ref{l:no-inclusion} in Page~\pageref{l:no-inclusion},
\bEqu{equ:need-for-cond3}
F(i,j_1)\cap F^\one(j_2)=\hspace{-0.1in} \bigcup_{\substack{\beta \in J_{j_1}\cap J_{j_2}\\ \al\in I_{j_1}\cap I_{j_2}}}\hspace{-0.15in}( F_{q_\beta,O}^{t}\cap F_{p_\al,O}^{m}\cap  F^\one(j_2))\subset F(i,j_2)\cap F^\one(j_2).
\eEqu
Let $\cE(i,j)\lra \cY(i,j)$ (with similarly denoted section and footprint maps) be sub-orbibundles of $\cE(i)\lra \cY(i)$ obtained via restriction to $Y(i,j)\!\subset\!  Y(i)$. By (\ref{equ:cover1}), (\ref{equ:cover2}), and (\ref{equ:overlap-charts}), we have 
$$
\tn{cl}\!\lrp{O}\cap \tn{cl}\!\lrp{F^\one(j)} \subset F(i,j) \quad \forall j\!\in\! \mfB(>\! i).
$$

\noindent
\textbf{Claim 1.} \textit{For every $j\!\in\!\mfB(>\!i)$, there exists a natural orbibundle embedding 
$$
\Phi_{i j}\equiv (\mfD\phi_{i j},\phi_{ij})\colon\lrp{\cE(i,j),\cY(i,j)}\lra \lrp{\cE(j),\cY(j)}
$$
that commutes with the section and footprint maps. Moreover, for $j_1,j_2\!\in\!\mfB(>\!i)$, with $j_1\!<\!j_2$, restricted to $\ov\phi_{ij_1}^{-1}(Y(j_1,j_2))\cap Y(i,j_2)$ we have 
$$
\ov\Phi_{ij_2}= \ov\Phi_{j_1j_2}\circ \ov\Phi_{i j_1}
$$ 
and Definition~\ref{level-system}.\ref{l:compatible} holds.}
\bProof
The proof is similar to the proof of Claim~2 in the proof of Lemma~\ref{lem:compact-level-2}.
For $j\!\in\!\mfB(>\!i)$, $\al\!\in\! I_j$, and $\beta\!\in\! J_j$, we show that the maps
$$
(\ov{\mfD\phi}^{\al,\beta}_{ij},\ov\phi^{\al,\beta}_{ij})\colon \lrp{E(i,j),Y(i,j)}|_{Y_{\al,\beta}}\lra \lrp{E(j), Y(j)}
$$
which send
\bEqu{x-image-2}
 \ov{x}= \ov\phi_{q_\beta i}(\ov{x}_\beta) = 
\ov\phi_{p_\al i}(\ov{x}_\al),\quad \ov{x}_\beta \in \ov{V}^t_{q_\beta,O},~\ov{x}_\al\in \ov{V}^m_{p_\al,O},
\eEqu
to $\ov{y}=\ov\phi_{q_\beta j}(\ov{x}_\beta) $, and
\bEqu{v-image-2}
\ov{v}= \ov{\mfD\phi}_{q_\beta i}(\ov{v}_\beta) = 
\ov{\mfD\phi}_{p_\al i}(\ov{v}_\al),\quad \ov{v}_\beta \in \ov{U}^t_{q_\beta,O},~\ov{v}_\al\in \ov{U}^m_{p_\al,O},
\eEqu
to $\ov{u}=\ov{\mfD\phi}_{q_\beta j}(\ov{v}_\beta) $  are compatible on the overlaps  and define a well-defined natural continuous map $\ov\Phi_{ij}$ as stated which lifts to a bundle isomorphism 
$$
\Phi_{ij}=(\mfD \phi_{ij},\phi_{ij})
$$ 
with the required properties.

\noindent
Given $\ov{x}\!\in\! Y(i,j)$,  suppose there are two different pairs of indices $(\al,\beta),(\al',\beta')\!\in\! I_j\!\times\! J_j$ such that 
 $$
 \aligned
& \ov{x}= \ov\phi_{q_\beta i}(\ov{x}_\beta) = 
\ov\phi_{p_\al i}(\ov{x}_\al),\quad \ov{x}_\beta \in \ov{V}^t_{q_\beta,O},~\ov{x}_\al\in \ov{V}^m_{p_\al,O},\\
&\ov{x}= \ov\phi_{q_{\beta'} i}(\ov{x}_{\beta'}) = 
\ov\phi_{p_{\al'} i}(\ov{x}_{\al'}),\quad \ov{x}_{\beta'} \in \ov{V}^t_{q_{\beta'},O},~\ov{x}_{\al'}\in \ov{V}^m_{p_{\al'},O}.
\endaligned
$$
By (\ref{equ:transition-condition}), it follows from 
$$
\aligned
&\ov{x}\in \ov\phi_{q_\beta i}(\ov{V}^t_{q_\beta,O})\cap \tn{cl}\!\lrp{ \ov\phi_{p_{\al'} i}(\ov{V}^m_{p_{\al'},O})}\neq \emptyset, \\
&\ov{x}\in \ov\phi_{q_{\beta'} i}(\ov{V}^t_{q_{\beta'},O})\cap \tn{cl}\!\lrp{ \ov\phi_{p_{\al} i}(\ov{V}^m_{p_{\al},O})}\neq \emptyset;
\endaligned
$$
that
$$
\ov{q}_\beta\in F_{p_{\al'}},~~\ov{q}_{\beta'}\in F_{p_{\al}},~~\ov{V}^t_{q_\beta,O}\subset \ov\phi_{q_\beta p_{\al'}}^{-1}(\ov{V}_{p_{\al'},O}),~~\tn{and}~ \ov{V}^t_{q_{\beta'},O}\subset \ov\phi_{q_{\beta'}p_\al}^{-1}(\ov{V}_{p_\al,O}).
$$
By Definition~\ref{level-system}.\ref{l:compatible}.(a), with $a\!=\!i$ or $j$, restricted to 
$$
\ov{V}^t_{q_\beta,O}\subset\ov{\phi}_{q_\beta p_{\al'}}^{-1}(\ov{V}_{p_{\al'},O})\cap \ov{V}_{q_\beta,O}
\subset \ov{\phi}_{q_\beta p_{\al'}}^{-1}(\ov{V}_{p_{\al'},a})\cap \ov{V}_{q_\beta,a},
$$
we have 
$$
\ov{\Phi}_{p_{\al'}a}\circ \ov{\Phi}_{q_\beta p_{\al'}}=\ov{\Phi}_{q_\beta a}.
$$ 
Therefore, with $a\!=\!i$ and since $\ov{\Phi}_{p_{\al'}i}$ is one-to-one, we get
$$
\ov{x}_{\al'}= \ov{\Phi}_{q_\beta p_{\al'}}(\ov{x}_\beta);
$$
and then, with $a\!=\!j$, we get 
$$
\ov{\phi}_{q_\beta j}(\ov{x}_\beta)=\ov{\phi}_{p_{\al'} j}(\ov{x}_{\al'}).
$$ 
We conclude that (\ref{x-image-2}) is independent of the choice of $(\al,\beta)$ or $(\al',\beta')$. The case of  (\ref{v-image-2}) is similar. 

\noindent
By definition, restricted to $Y_{\al,\beta}$, 
\bEqu{e:local-equation}
(\ov{\mfD\phi}_{ij},\ov\phi_{ij})=\ov\Phi_{p_\al j}\circ \ov\Phi_{p_\al i}^{-1};
\eEqu
therefore, it is locally liftable to a smooth orbibundle embedding between orbibundle charts. 
The cocycle condition $\ov\Phi_{ij_2}\!=\! \ov\Phi_{j_1j_2}\circ \ov\Phi_{ij_1}$ is a consequence of (\ref{e:local-equation}).
Compatibility with the section and footprint maps and  Definition~\ref{level-system}.\ref{l:compatible} follow from the corresponding compatibility conditions of Definition~\ref{level-system} for $\mfL(>\! i)$.
\eProof

\noindent
The extended system of pure orbibundle structures
$$
\mfL(\geq\! i)\equiv 
\bigg(\{\cE(j),\cY(j),s_j,\psi_j\}_{ j\in \mfB(\geq i)},\{\cE(j_1,j_2),\cY(j_1,j_2),\Phi_{j_1j_2}\}_{\substack{j_1,j_2\in\mfB(\geq i)\\ j_1<j_2}}\bigg)
$$
constructed  so far does not satisfy all the requirements of Definition~\ref{level-system}. Also, the corresponding thickening may not be Hausdorff. In the remaining part of the proof, we shrink these orbibundles in a way that the reduced system satisfies all the conditions of Definition~\ref{level-system} and Definition~\ref{def:H-level-structure}, as well. 

\noindent
First, for $j\!>\!i$, the intersection condition (\ref{equ:Fij=FicapFj}) does not necessarily  hold. In fact,
$$
\tn{cl}\lrp{M(i)\cap F^\one_0(i)\cap F^\one(j)}\subset F(i,j)\subset F(i)\cap F(j)
$$
and the second inclusion could be proper. Let 
$$
K=K_{i}\setminus \hspace{-.1in}  \bigcup_{j\in \mfB(> i)}\hspace{-.1in} F(i,j)\subset M(i);
$$
this is a compact set which has empty intersection with $\tn{cl}(F^\one(>\! i))$.
For $j\!>\!i$, we replace $Y(j)$ with $Y^\one(j)$. For $i\!<\!j_1\!<\!j_2$, we replace $Y(j_1,j_2)$ with $Y^\one(j_1,j_2)$. We replace $Y(i)$ with 
\bEqu{e:Y1}
{Y}^\one(i)= W\cup \bigcup_{j\in \mfB(> i)} Y(i,j)
\eEqu
\vskip-.2in
\noindent
such that $W$ has the following properties:
\bEqu{e:W-condition}
W\subset Y(i), \quad K\subset \ov\psi_{i}(W), ~~ \tn{and}~~W \cap \bigcup_{j\in \mfB(> i)}\ov\phi_{ij}^{-1}(Y^\one(j))=\emptyset.
\eEqu
For every $j\!\in\! \mfB(>\! i)$, we replace $Y(i,j)$ with
\bEqu{equ:overlap-charts2}
Y^\one(i,j)=  \ov\phi_{ij}^{-1}(Y^\one(j))\subset Y(i,j)\subset Y^\one(i).
\eEqu
It is immediate that the new collection
$$
\mfL^\one(\geq\! i)\equiv \bigg(\{\cE^\one(j),\cY^\one(j),s_j,\psi_j\}_{ j\in \mfB(\geq i)},\{\cE^\one(j_1,j_2),\cY^\one(j_1,j_2),\Phi_{j_1j_2}\}_{\substack{j_1,j_2\in\mfB(\geq i)\\ j_1<j_2}}\bigg)
$$
satisfies Conditions~\ref{l:dim}, \ref{l:footprint}, \ref{l:compatible}, and \ref{l:tangent}  of Definition~\ref{level-system}. It is a bit tedious, but we show that $\mfL^\one(\geq\! i)$ also satisfies Conditions~\ref{l:intersection-2}, \ref{l:intersection-source}, and \ref{l:intersection-target}. Then it follows that it also satisfies \ref{l:intersection-cycle}. Eventually, we achieve Hausdorffness by another appropriate shrinking.

\noindent
\textbf{Condition~\ref{l:intersection-2} of Definition~\ref{level-system}.}
For every $j\!\in\! \mfB(>\! i)$, by (\ref{equ:need-for-cond3}), (\ref{e:Y1}), and (\ref{e:W-condition}), we have
\bEqu{equ:new-intersection}
\aligned
 &F^\one(i)\cap F^\one(j)=\big( \ov\psi_{i}(W)\cup \bigcup_{k\in \mfB(> i)} F(i,k)\big) \cap F^\one(j)=\\
& \bigcup_{k\in \mfB(> i)} (F(i,k) \cap F^\one(j))=  F(i,j)\cap F^\one(j)=F^\one(i,j).
\endaligned
\eEqu
Therefore, $\mfL^\one(\geq\! i)$ satisfies (\ref{equ:Fij=FicapFj}), i.e. Definition~\ref{level-system}\ref{l:intersection-2} holds. 

\noindent
\textbf{Claim 2.} \textit{For $j_1,j_2\!\in\!\mfB(>\!i)$, 
$$
Y^\one(i,j_1)\cap Y^\one(i,j_2)=\ov\phi_{ij_1}^{-1}(Y^\one(j_1))\cap \ov\phi_{ij_2}^{-1}(Y^\one(j_2)) \cap  \hspace{-.1in}
\bigcup_{\substack{\beta \in J_{j_1}\cap J_{j_2}\\\al \in I_{j_1}\cap I_{j_2}}} \hspace{-.1in} Y_{\al,\beta} 
$$}
\vspace{-.1in}
\bProof 
The inclusion of right-hand side in the left-hand side is clear from (\ref{equ:overlap-charts}) and (\ref{equ:overlap-charts2});
moreover, 
$$
Y^\one(i,j_1)\cap Y^\one(i,j_2)=\phi_{ij_1}^{-1}(Y^\one(j_1))\cap \phi_{ij_2}^{-1}(Y^\one(j_2)) \cap
 \bigcup_{\substack{\beta_1 \in J_{j_1}, \beta_2\in J_{j_2}\\\al_1\in I_{j_1},\al_2\in I_{j_2}}} \hspace{-.1in}\bigg(Y_{\al_1,\beta_1} \cap Y_{\al_2,\beta_2}\bigg).
$$
We may assume $j_1\!<\!j_2$. For some $\beta_1\! \in\! J_{j_1}$, $\beta_2\!\in\! J_{j_2}$, $\al_1\!\in\! I_{j_1}$, and $\al_2\!\in\! I_{j_2}$, if
\bEqu{equ:one-term}
\phi_{ij_1}^{-1}(Y^\one(j_1))\cap \phi_{ij_2}^{-1}(Y^\one(j_2)) \cap Y_{\al_1,\beta_1} \cap Y_{\al_2,\beta_2}\neq \emptyset,
\eEqu
it follows from (\ref{equ:transition-condition}) and (\ref{equ:for(4)to-hold}) that 
$$
\ov{q}_{\beta_1}\in F^\one(j_1)\cap F(j_2),\quad \ov{q}_{\beta_2}\in F(j_1)\cap F^\one(j_2),
$$ 
and 
\bEqu{equ:composition-defined}
\ov\phi_{ij_1}(Y_{\al_1,\beta_1}),\ov\phi_{ij_1}(Y_{\al_2,\beta_2}) \subset Y(j_1,j_2).
\eEqu
If $\ov{q}_{\beta_1}\!\notin\! F^\one(j_2)$, by (\ref{equ:no-inclusion-1}), 
$$
\ov\phi_{ij_2}(Y_{\alpha_1,\beta_1})\cap Y^\one(j_2)\subset \ov\phi_{q_{\beta_1} j_2}(\ov{V}^t_{q_{\beta_1},O})\cap Y^\one(j_2)=\emptyset;
$$
which contradicts (\ref{equ:one-term}). Thus, $\ov{q}_{\beta_1}\!\in\! F^\one(j_2)$. Similarly, we conclude $\ov{q}_{\beta_2}\!\in\! F^\one(j_1)$. 
If $\ov{p}_{\alpha_1}\!\notin\! F^\one(j_2)$, by (\ref{equ:for(5)to-hold}), 
$$
\ov\phi_{ij_1}(Y_{\al_1,\beta_1})\cap Y^\one(j_1,j_2)\subset \ov\phi_{p_{\alpha_1} j_1}(\ov{V}^m_{p_{\alpha_1},O})\cap Y^\one(j_1,j_2)=\emptyset.
$$
Therefore, by (\ref{equ:composition-defined}),
$$
Y_{\al_1,\beta_1}\cap \phi_{ij_1}^{-1}(Y^\one(j_1))\cap \phi_{ij_2}^{-1}(Y^\one(j_2))=\emptyset;
$$
which contradicts (\ref{equ:one-term}). Thus, $\ov{p}_{\alpha_1}\!\in\! F^\one(j_2)$. Similarly, we conclude $\ov{p}_{\alpha_2}\!\in\! F^\one(j_1)$. This finishes the proof of Claim 2.
\eProof

\noindent
\textbf{Condition \ref{l:intersection-source} of Definition~\ref{level-system}.} 
For $j_1,j_2\!\in\!\mfB(>\!i)$, with $j_1\!<\!j_2$, by Claim~2 and (\ref{equ:for(4)to-hold}), 
\bEqu{equ:left-to-right}
Y^\one(i,j_1)\cap Y^\one(i,j_2)\subset \ov\phi_{ij_1}^{-1}(Y^\one(j_1,j_2)).
\eEqu
On the other hand, for every $\alpha \!\in\! I_{j_1}$ and $\beta\!\in\! J_{j_1}$, if 
$$
\ov{p}_\al\!\notin\! \tn{cl}(F^\one(j_2))\quad \tn{or}\quad \ov{q}_\beta\!\notin \!\tn{cl}(F^\one(j_2)),
$$ 
by the second case of (\ref{equ:for(5)to-hold}) we get
\bEqu{equ:contributing-parts}
\ov\phi_{ij_1}(Y_{\al,\beta})\cap Y^\one(j_1,j_2)=\emptyset; 
\eEqu
i.e. 
\bEqu{equ:what-matters}
\ov\phi_{ij_1}^{-1}(Y^\one(j_1,j_2))\cap  Y_{\al,\beta}=\emptyset.
\eEqu
Therefore, by (\ref{equ:what-matters}) and the definition of $Y^\one(j_1,j_2)$ in the proof Lemma~\ref{lem:shrinking}, 
$$
\aligned
 &\ov\phi_{ij_1}^{-1}(Y^\one(j_1,j_2))=\ov\phi_{ij_1}^{-1}\big(Y^\one(j_1)\cap \ov\phi_{j_1j_2}^{-1}(Y^\one(j_2)) \cap Y(j_1,j_2)\big)\subset \\
 & \phi_{ij_1}^{-1}(Y^\one(j_1))\cap \phi_{ij_2}^{-1}(Y^\one(j_2)) \cap \bigcup_{\substack{\beta \in J_{j_1}\cap J_{j_2}\\ \al \in I_{j_1}\cap I_{j_2}}} \hspace{-.2in} Y_{\al,\beta}= Y^\one(i,j_1)\cap Y^\one(i,j_2).
 \endaligned
$$
Together with (\ref{equ:left-to-right}), we conclude that the equality in Definition~\ref{level-system}.\ref{l:intersection-source} is satisfied.

\noindent
Let
$$
Y(i,j_1,j_2)=\ov\phi_{ij_1}^{-1}(Y^\one(j_1,j_2))= Y^\one(i,j_1)\cap Y^\one(i,j_2)
$$
as in Definition~\ref{level-system}.\ref{l:intersection-source}.

\noindent
\textbf{Condition~\ref{l:intersection-target} of Definition~\ref{level-system}.}
For $j_1,j_2\!\in\!\mfB(>\!i)$, with $j_1\!<\!j_2$, by (\ref{equ:overlap-charts}),
\bEqu{e:image-1}
\aligned
&\ov\phi_{ij_2}(Y^\one(i,j_2))= Y^\one(j_2)\cap \bigcup_{\beta \in J_{j_2},\al\in I_{j_2}}\lrp{ \ov\phi_{q_\beta j_2}(\ov{V}^t_{q_\beta,O})\cap \ov\phi_{p_\al j_2}(\ov{V}^m_{p_\al,O}) },\\
&\ov\phi_{j_1j_2}(Y^\one(i,j_2))= Y^\one(j_2)\cap \ov\phi_{j_1,j_2}(Y^\one(j_1)\cap Y(j_1,j_2)),
\endaligned
\eEqu
By Claim~2 above,
\bEqu{e:image-2}
\aligned
\ov\phi_{ij_2}(Y^\one(i,j_1,j_2))= &Y^\one(j_2)\cap \ov\phi_{j_1,j_2}(Y^\one(j_1)\cap Y(j_1,j_2))\cap \\
&\bigcup_{\substack{\beta \in J_{j_1}\cap J_{j_2}\\\al\in I_{j_1}\cap I_{j_2}}}\lrp{ \ov\phi_{q_\beta j_2}(\ov{V}^t_{q_\beta,O})\cap \ov\phi_{p_\al j_2}(\ov{V}^m_{p_\al,O}) }
\endaligned
\eEqu
For every $\alpha\! \in\! I_{j_2}$ and $\beta\!\in\! J_{j_2}$, if 
\bEqu{equ:non-empty-over-lap}
Y^\one(j_2)\cap  \ov\phi_{q_\beta j_2}(\ov{V}^t_{q_\beta,O})\cap \ov\phi_{p_\al j_2}(\ov{V}^m_{p_\al,O}) 
\cap \ov\phi_{j_1,j_2}(Y^\one(j_1)\cap Y(j_1,j_2))\neq \emptyset,
\eEqu
by (\ref{equ:for(5)to-hold}),  $\ov{q}_\beta,\ov{p}_\al\!\in\! \tn{cl}\lrp{F^\one(j_1)}$. It immediately follows from (\ref{equ:non-empty-over-lap}), (\ref{e:image-1}), and (\ref{e:image-2}), that
$$
\ov\phi_{ij_2}(Y^\one(i,j_1,j_2))= \ov\phi_{ij_2}(Y^\one(i,j_2))\cap \ov\phi_{j_1j_2}(Y^\one(j_1,j_2)).
$$
Hence, Definition~\ref{level-system}.\ref{l:intersection-target} is also satisfied.

\noindent
\textbf{Achieving of Hausdorffness.}
The thickening $Y(\mfL^\one(\geq\! i))$ has the form
$$
Y(\mfL^\one(\geq\! i))=\bigg(Y^\one(i) \coprod Y(\mfL^\one(>\! i))\bigg)/\sim,
$$
where, by induction, $Y(\mfL^\one(>\! i))$ is Haussdorf and 
$$
Y^\one(i) \ni \ov{x} \sim [\ov{y}] \in Y(\mfL^\one(>\! i))
$$
whenever there exists $j\!\in\!\mfB(>\! i)$ such that $\ov{y}\!\in\! Y^\one(j)$, $[\ov{y}]\!=\!\wp_j (\ov{y})$, and $\ov{y}\!=\!\ov\phi_{ij}(\ov{x})$.
In order to achieve Hausdorffness, we shrink  $Y^\one(i)$ and $Y(\mfL^\one(>\! i))$ into relatively compact open sets 
$Y^\two(i)$ and $\wt{Y}(\mfL^\one(>\! i))$, respectively, similarly to the argument after the proof of Claim~3 in page~\pageref{page-ref}. 
For $j\!\in\!\mfB(>\!i)$, let 
$$
Y^\two(j)= \wp_j^{-1}(\wt{Y}(\mfL^\one(>\! i))), 
$$
where $\wp_i\colon Y^\one(j)\lra Y(\mfL^\one(\geq i))$ is the inclusion map as in the paragraph before Definition~\ref{def:H-level-structure}. The shrunk DGS  
$$
\aligned
&\mfL^\two(\geq\! i)\equiv \\
&\bigg(\{\cE^\two(j),\cY^\two(j),s_j,\psi_j\}_{ j\in \mfB(\geq i)},\{\cE^\two(j_1,j_2),\cY^\two(j_1,j_2),\Phi_{j_1j_2}\}_{\substack{j_1,j_2\in\mfB(\geq i)\\ j_1<j_2}}\bigg)
\endaligned
$$
covers $M(\geq\!i)$, is Hausdorff , and satisfies all the conditions of Definition~\ref{level-system}.
This completes the induction step and the proof of the first statement in Theorem~\ref{exist-level-system}.

\noindent
For every $i\!\in\!\mfB$, the orbibundle $\cE(i)\!\lra\!\cY(i)$ is made of Kuranishi charts of the same rank $i$. Therefore, if $(M,\cK)$ is oriented, by the last condition in Definition~\ref{def:Kuranishi-orientable}, every $\cE(i)\!\lra\!\cY(i)$ is relatively oriented. Compatibility of (\ref{equ:index-bundle-i-0}) and (\ref{equ:iso-det-bundles}) in Definition~\ref{def:Kuranishi-orientable-2} follow from the corresponding condition in Definition~\ref{def:Kuranishi-orientable}.

\noindent
It is just left to show that every two such dimensionally graded systems are deformation equivalent. In fact, it follows from the construction of the proof that every two of such DGS admit a common refinement\footnote{This is analogue of commensurability in \cite{MW2}.} in the following sense: for every two dimensionally graded systems 
$\mfL_0$ and $\mfL_1$ constructed as above, there exist shrinkings $\wt\mfL_{0}$ and  $\wt\mfL_{1}$ of $\mfL_0$ and $\mfL_1$, respectively, which are isomorphic in the following sense. 
There are orbibundle isomorphisms
$$
\Phi_{01;i}\colon \colon \wt\cE_0(i) \lra \wt\cE_1(i)\quad \forall i\!\in\!\mfB
$$ 
such that 
\bItem
\item each $\Phi_{i}$ commutes with the section and footprint maps,
\item for every $\ov{p}\!\in\!F_0(i)\!=\!F_1(i)$, the compatibility embeddings of Definition~\ref{level-system}.\ref{l:compatible} (taken over a sub-chart in the intersection of their domains)
$$
\Phi_{0;pi}, \Phi_{1;pi}\colon ([U_{p,i}/ G_{p,i}],[V_{p,i}/ G_{p,i}])\lra (\cE_0(i),\cY_0(i)),(\cE_1(i),\cY_1(i))
$$
satisfy $\Phi_{1;pi}\!=\!\Phi_{01;i}\!\circ\!\Phi_{0;pi}$.
\eItem

\noindent
Such a common refinement gives us a deformation equivalence\footnote{This is analogue of \cite[Lemma 6.2.16]{MW2}.} dimensionally graded system $\wt\mfL$ for $(M\!\times\![0,1],\cK\!\times\![0,1])$ (see the argument after Definition~\ref{def:cobo-Kur}) with orbibundle pieces of the form
$$
\aligned
(\wt\cE(i),\wt\cY(i))=~&(\cE_0(i),\cY_0(i))\!\times\![0,1/3)\coprod (\wt\cE_0(i),\wt\cY_0(i)) \!\times\! (0,1) \coprod \\
& (\wt\cE_1(i),\wt\cY_1(i)) \!\times \!(0,1)\coprod (\cE_1(i),\cY_1(i))\!\times\!(2/3,1] /\sim \quad \forall i\!\in\!\mfB,
\endaligned
$$
where the gluing relation $\sim$ between $(\wt\cE_0(i),\wt\cY_0(i)) \!\times\! (0,1) $ and $(\wt\cE_1(i),\wt\cY_1(i)) \!\times \!(0,1)$ is by $\Phi_{01;i}\!\times\!\tn{id}$ and the gluing relations in the first and second lines are given by obvious inclusions restricted to $(0,1/3)$ and $(2/3,1)$, respectively. The footprint map and section on each piece are the obvious product footprint map and section. This finishes the proof of Theorem~\ref{exist-level-system}.

\qed
\end{proofoflevel-system}

\subsection{Deformations of Kuranishi maps}\label{sec:def}
Let $\mfL$ be an oriented Hausdorff DGS as in Definition~\ref{level-system} for an $n$-dimensional oriented Kuranishi  structure $(\cK,M)$ as in Section~\ref{sec:level}. In Section~\ref{sec:VFC}, with an approach similar to the orbibundle case of Section~\ref{sec:euler}, we will build a virtual fundamental class for $M$ with respect to $\mfL$. More precisely, after shrinking $\mfL$ finitely many times into a sufficiently small $\wt\mfL$,  we build a rational singular homology $n$-cycle 
$$
\vfc(\wt\mfL)\!\in\! H_n(Y(\wt\mfL),\Q)
$$ 
which will play the role of fundamental class of $M$ in calculations.
If $Y(\wt\mfL)$ admits a continuous map into some topological manifold, say $f\colon Y(\wt\mfL)\lra X$, we can push-forward this homology class to a homology class in $X$. This is for instance the case in the construction of Gromov-Witten VFC. In this particular application, the evaluation and forgetful maps, (\ref{equ:eval}) and (\ref{equ:st}), naturally extend to the thickening $Y(\wt\mfL)$ and the Gromov-Witten VFC of $\ov\cM_{g,k}(X,A)$, as we define it in this article, is the image under $\ev\!\times\! \st$ of 
$\vfc(\wt\mfL)$ in $H_*(X^k\times \ov\cM_{g,k},\Q)$. In the case of abstract Kuranishi structures, up to cobordism, VFC is independent of the particular choice of natural $\mfL$ in Theorem~\ref{exist-level-system}; thus we may denote it by $\vfc(\cK)$. In the particular case of GW theory, $\vfc(\cK)$ does not depend on the particular choice of natural $\cK$ in Theorem~\ref{thm:VFC}, either; thus we denote it by $\vfc(\ov{\cM}_{g,k}(X,A,J))$ or $[\ov{\cM}_{g,k}(X,A,J)]^\vfc$.

\noindent
With notation as in (\ref{level-tuple}), in order to construct $\vfc(\mfL)$, we start from the smallest  value $i_1\!\in\! \mfB$, shrink $\mfL$ into $\mfL^\one$, perturb the kuranishi map $s_{i_1}$ into a sequence of transversal lifted multisections $(\mfs_{i_1;a}\!=\! s_{i_1}\!+\mft_{i_1;a})_{a\in \N}$ on $\tn{cl}(Y^\one(i_1))$, where $\lim_{a\lra \infty} \mft_{i_1;a}\equiv 0$ as in Proposition~\ref{pro:estimate}, and inductively extend the multisections  $\mft_{i_1;a}$ to $\cE(i)\lra \cY(i)$ for higher values of $i\!\in\! \mfB$, again via Proposition~\ref{pro:estimate}. 
The result of this process is a compatible set of sequences of deformations of Kuranishi maps, 
$$
(\mfs_{i;a}= s_i+\mft_{i;a})_{i\in \mfB,a\in \N}, $$ 
described in Definition~\ref{def:pert-system-sections} below. 

\noindent
For this induction process, we need to ``compatibly" extend every sub-orbibundle 
$$\mfD\phi_{ij}(\cE(i,j))\subset \cE(j)
$$ 
to a sub-orbibundle 
$$
\cE_{\tn{ext}}(i,j)\subset \cE(j)
$$ 
over a tubular neighborhood $\cN(i,j)$ of $\phi_{ij}(\cY(i,j))\subset \cY(j)$. 
We need some compatibility conditions between these extensions in order for the resulting deformations of Kuranishi maps to be equal one the overlaps. Following definitions abstractly characterize such extensions. 

\noindent
For a sub-orbifold $\cM_1\!\subset\!\cM_2$ as in Definition~\ref{def:tub-nbhd}, let 
$$
\pr\colon\cN_{\cM_1}\cM_2\lra \cM_1
$$ 
be the orbifold normal bundle of that as in Example~\ref{eax:OrbiNormal}. Similar to the case of manifolds, the restriction of orbifold tangent bundles $T\cM_2$ and $T\cN_{\cM_1}\cM_2$ to $\cM_1$ are naturally isomorphic to each other. We make use of this identification in (\ref{equ:f-map}). In what follows, by a \textbf{tubular neighborhood} of $\cM_1$ in $\cN_{\cM_1}\cM_2$ we mean the restriction $\cN'_{\cM_1}\cM_2$ of $\cN_{\cM_1}\cM_2$  to some open set $ N'_{M_1}M_2$ including $M_1$. Here, by convention of Notation~\ref{not:calconv}, we are thinking about $\cN_{\cM_1}\cM_2$ as an orbifold with underlying topological space $N_{M_1}M_2$.

\bRem{rmk:unusual}
Our definition of tubular neighborhood is weaker than its usual meaning in the literature. In our definition, $N'_{M_1}M_2$ does not even have to be connected; it is just an open set in $N_{M_1}M_2$ containing $M_1$. 
\eRem

\bDef{def:tub-nbhd}
Given a sub-orbifold $\cM_1\!\subset\! \cM_2$, a \textbf{tubular neighborhood} of $\cM_1$ in $\cM_2$ is an open sub-orbifold $\cW\!\subset\! \cM_2$ containing $\cM_1$ with an orbifold smooth map $\pi\colon\cW\lra \cM_1$ such that there exists a tubular neighborhood $\cN'_{\cM_1}\cM_2$ of the zero section $\cM_1\!\subset\! \cN_{\cM_1}\cM_2$ and an orbifold diffeomorphism 
$$
f\colon \cN'_{\cM_1}\cM_2\lra \cW
$$
satisfying 
\bEqu{equ:f-map}
\pr=\pi\circ f~~\tn{and}~~\nd f|_{T\cM_2|_{\cM_1}}=\tn{id}.
\eEqu
\eDef

\noindent
In the context of Definition~\ref{def:tub-nbhd}, we call such $f$ a \textbf{regularization} of the tubular neighborhood $\pi\colon\cW\lra \cM_1$.

\bDef{def:sub-nbhd}
Given a tubular neighborhood $\pi\colon \cW\lra \cM_1$ for $\cM_1\!\subset\! \cM_2$, a \textbf{sub-tubular neighborhood} of $\cW$ is some open\footnote{Together with the restriction of the corresponding projection map.} sub-orbifold $\cW'\subset \cW$ containing $\cM_1$.
\eDef
\noindent
For a tubular neighborhood $\cW$ as in Definition~\ref{def:tub-nbhd} and any $\cM_1'\!\subset\!\cM_1$, let
$$
\cW|_{\cM_1'}= \pi^{-1}(\cM_1')\subset \cW;
$$
if $\cM_1'$ is an open sub-orbifold of $\cM_1$, then $\cW|_{\cM_1'}$ is a tubular neighborhood of $\cM_1'\!\subset\!\cM_2$.
\bDef{def:ext-obs}
Let $\mfL$ be an DGS as in Definition~\ref{level-system} and $\ell\!\in\!\mathbb{N}$. An \textbf{extension of the obstruction bundles (\textbf{EOB})} $\cN$ of \textbf{width $\ell$} for $\mfL$ consists of a set of tubular neighborhoods 
$$
\pi_{ji}\colon\cN(i,j)\lra \phi_{ij}(\cY(i,j))
$$ 
for $\phi_{ij}(\cY(i,j))\!\subset\!\cY(j)$ and orbibundle embeddings
\bEqu{equ:embedd-bundles}
\iota_{ij}\colon \!\pi_{ji}^*(\mfD\phi_{ij}\cE(i,j))\!\lra\! \cE(j),\quad \iota_{ij}|_{\phi_{ij}(\cY(i,j))}=\tn{id},
\eEqu
for all $i,j\!\in\!\mfB$, with $0\!<\!j\!-\!i\!\leq \!\ell$, such that the following \textbf{compatibility} conditions, for every $i,j,k\in \mfB$, with $i\!<\!j\!<\!k$ and $k\!-\!i\leq \ell$, hold.
\bEnum
\item\label{l:nested1} With 
$\cN(i,j,k)\!=\! \cN(i,k)\cap \cN(j,k)$, 
we have 
\bEqu{l:like-5}
\phi_{jk}^{-1}(\pi_{kj}(\cN(i,j,k)))\subset \cN(i,j).
\eEqu
\item\label{l:nested2} Restricted to $\cN(i,j,k)$, we have 
\bEqu{equ:nested4}
(\phi_{ik}^{-1}\circ \pi_{ki})=(\phi_{ij}^{-1}\circ \pi_{ji})\circ(\phi_{jk}^{-1}\circ \pi_{kj}).
\eEqu

\item\label{l:nested3} For every $x\!\in\! \cN(i,j,k)$, with $y\!=\!\phi_{jk}^{-1}(\pi_{kj}(x))$ and $z\!=\!\phi_{ki}^{-1}(\pi_{ki}(x))$,
$$
\iota_{ik}\lrp{x,\mfD\phi_{ik}(v)}=\iota_{jk}\lrp{x,\mfD\phi_{jk}\lrp{\iota_{ij}\lrp{y,\mfD\phi_{ij}(v)}}}\quad \forall v\!\in\! \cE(i,k)_z.
$$
Here $\lrp{x,\mfD\phi_{ik}(v)}$, and the similar notation on the right-hand side of the equation, denotes the pull-back of vector 
$\mfD\phi_{ik}(v)\!\in\! \mfD\phi_{ik}(\cE(i,k))$ to the fiber of $\pi_{ki}^*\mfD\phi_{ik}(\cE(i,k))$ at $x\!\in\!\pi_{ki}^{-1}(\phi_{ik}(z))$. 
\eEnum
\eDef

\noindent
\bRem{rmk:clarification}
With $\cY(i,j,k)$ as in Definition~\ref{level-system}\ref{l:intersection-source}, (\ref{l:like-5}) and (\ref{equ:nested4}) imply that 
$\cN(i,j,k)$ is a sub-tubular neighborhood of $
\cN(i,k)|_{\phi_{ik}(\cY(i,j,k))}
$; i.e. the following diagram is well-defined and commutes:
$$
    \xymatrix{ \cN(i,j,k)  \ar[rrd]^{\phi_{ik}^{-1}\circ\pi_{ki}}\ar[rr]^{\phi_{jk}^{-1}\circ \pi_{kj}} &&  \cN(i,j)|_{\phi_{ij}(\cY(i,j,k))} \ar[d]^{\phi_{ij}^{-1}\circ\pi_{ji}}\\
                                && \cY(i,j,k).  }
$$
Figure~\ref{fig:tubular} illustrates a sequence of compatible neighborhoods and their overlaps.
\eRem

\noindent
In the context of Definition~\ref{def:ext-obs}, for every $i,j\!\in\!\mfB$, with $0\!<\!j\!-\!i\!\leq \!\ell$, let 
$$
\cE_{\tn{ext}}(i,j)= \iota_{ij}(\pi_{ji}^*(\mfD\phi_{ij}\cE(i,j)))\subset \cE(j)|_{\cN(i,j)}.
$$ 

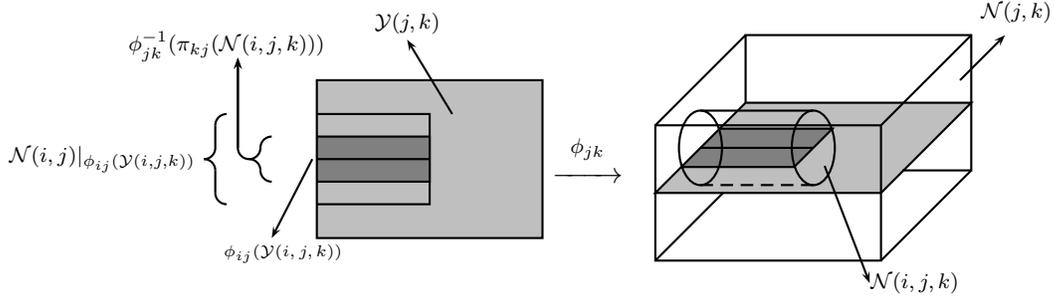
\begin{figure}
\begin{pspicture}(-6,-5)(10,-2)
\psset{unit=.3cm} 
\pspolygon*[linecolor=lightgray](18,-13)(8,-13)(12,-9)(22,-9)
\pspolygon*[linecolor=gray](10.85,-10.15)(15.85,-10.15)(14.15,-11.85)(9.15,-11.85)

\pspolygon(18,-13)(8,-13)(12,-9)(22,-9)
\psline(10.85,-10.15)(15.85,-10.15)(14.15,-11.85)(9.15,-11.85)

\pspolygon(18,-10)(18,-16)(8,-16)(8,-10)
\psline(18.95,-12)(22,-12)(22,-6)(12,-6)(12,-9)
\psline(18,-10)(22,-6)
\psline(18,-16)(22,-12)
\psline(8,-16)(11,-13)
\psline(8,-10)(12,-6)

\psellipse(15,-11)(1,1.7)
\psellipse(10,-11)(1,1.7)
\psline(10,-9.35)(15,-9.35)
\psline[linestyle=dashed](10,-12.65)(15,-12.65)
\psline(10,-11)(15,-11)
\pspolygon*[linecolor=lightgray](3,-15)(-7,-15)(-7,-8)(3,-8)
\pspolygon*[linecolor=gray](-2,-12.5)(-7,-12.5)(-7,-10.5)(-2,-10.5)

\pspolygon(3,-15)(-7,-15)(-7,-8)(3,-8)
\psline(-7,-12.5)(-2,-12.5)(-2,-10.5)(-7,-10.5)
\psline(-2,-11.5)(-7,-11.5)

\psline(-7,-13.5)(-2,-13.5)(-2,-9.5)(-7,-9.5)

\psline{->}(-1,-9.5)(-3,-6.2)
\rput(-3,-5.5){\scriptsize{$\cY(j,k)$}}

\psarc(-12,-11){.5}{270}{360}\psline(-11.5,-11)(-11.5,-10)\psarc(-11,-10){.5}{90}{180}
\psarc(-12,-12){.5}{0}{90}\psline(-11.5,-12)(-11.5,-13)\psarc(-11,-13){.5}{180}{270}
\rput(-16.5,-11.5){\scriptsize{$\cN(i,j)|_{\phi_{ij}(\cY(i,j,k))}$}}

\psarc(-10,-11){.5}{180}{360}\psarc(-9,-11){.5}{90}{180}
\psarc(-10,-12){.5}{0}{90}\psarc(-9,-12){.5}{180}{270}
\psline{->}(-10.5,-11)(-10.5,-7)
\rput(-11,-6.5){\scriptsize{$\phi_{jk}^{-1}(\pi_{kj}(\cN(i,j,k)))$}}

\psline{->}(-7.2,-11.5)(-9,-15)
\rput(-8.5,-15.7){\tiny{$\phi_{ij}(\cY(i,j,k))$}}

\rput(5,-11.5){$\stackrel{\phi_{jk}}{\xrightarrow{\hspace*{.7cm}}}$}
\psline{->}(-10.5,-11)(-10.5,-7)

\rput(19.5,-17){\scriptsize{$\cN(i,j,k)$}}
\psline{->}(15.5,-11.8)(17.5,-17)

\rput(24,-5){\scriptsize{$\cN(j,k)$}}
\psline{->}(21.5,-8)(23.5,-6)

\end{pspicture}
\caption{Overlap of tubular neighborhoods.}
\label{fig:tubular}
\end{figure}

\noindent
For simplicity, by an EOB, i.e. whenever we do not specify its width, we mean an EOB of the width infinity (or any number larger than the width of $\mfB$). The case of EOB with bounded width appears in the inductive proof Theorem~\ref{thm:exists-extension}; otherwise, we only consider EOB of the width infinity. The case of $\ell\!=\!1$ corresponds to a set of tubular neighborhoods and extension of the corresponding orbibundles without any compatibility condition on the overlaps.

\bDef{def:shrink-pair}
Let $\mfL$ be an DGS as in Definition~\ref{level-system}. Let $\cN$ be an EOB of width $\ell$ for $\mfL$ as in Definition~\ref{def:ext-obs}. A shrinking of $(\mfL,\cN)$ consist of a shrinking $\wt\mfL$ of $\mfL$ as in Lemma~\ref{lem:shrinking} and an EOB $\wt\cN$ of width $\ell$  for $\wt\mfL$, such that for all $i,j\!\in\!\mfB$, with $0\!<\!j\!-\!i\!\leq \!\ell$,
$\wt\cN(i,j)$ is a sub-tubular neighborhood of $\cN(i,j)|_{\phi_{ij}(\wt\cY(i,j))}$. 
\eDef

\noindent
The extensions $\wt\cE_{\tn{ext}}(i,j)$ in this situation are simply the restrictions of the corresponding bundles to $\wt\cN(i,j)$.

\bLem{lem:shrunk-EOB}
Let $(\mfL,\cN)$ be a pair as in Definition~\ref{def:shrink-pair}. For every shrinking $\wt{\mfL}$ of $\mfL$ as in Lemma~\ref{lem:shrinking}, the sub-tubular neighborhoods 
\bEqu{equ:wtN}
\wt\cN(i,j)= \cN(i,j)|_{\phi_{ij}(\wt\cY(i,j))}\cap \wt\cY(j)\quad \forall i,j\!\in\!\mfB,~\tn{with }i\!<\!j,
\eEqu
give rise to a shrinking $(\wt\mfL,\wt\cN)$ of $(\mfL,\cN)$ as in Definition~\ref{def:shrink-pair}. 
\eLem
\bProof
We just need to verify (\ref{l:like-5}), the rest follows by restriction. By definition we have 
\bEqu{equ:wtcNijk}
\wt\cN(i,j,k)=\cN(i,k)|_{\phi_{ik}(\wt\cY(i,k))}\cap \cN(j,k)|_{\phi_{jk}(\wt\cY(j,k))}\cap \wt\cY(k).
\eEqu
Suppose there exists $x\!\in\!\wt\cN(i,j,k)$ such that 
$$
y= \phi_{jk}^{-1}(\pi_{kj}(x))\!\notin\! \wt\cN(i,j)\!=\! \cN(i,j)|_{\phi_{ij}(\wt\cY(i,j))}\cap \wt\cY(j).
$$
By (\ref{equ:wtcNijk}), 
\bEqu{equ:yin}
y\!\in\! \wt\cY(j,k)\!\subset \!\wt\cY(j).
\eEqu
Therefore, 
$$
y\notin  \cN(i,j)|_{\phi_{ij}(\wt\cY(i,j))}.
$$
Since $y\!\in\!\cN(i,j)$ (because (\ref{l:like-5}) holds for $\cN$), we conclude that 
$$
z=\phi_{ij}^{-1}(\pi_{ji}(y))\notin \wt\cY(i,j)= \phi_{jk}^{-1}(\wt\cY(j))\cap \cY(i,j)\cap \wt\cY(i).
$$
It follows from (\ref{equ:yin}) that $z\!\notin\! \cY(i,j)\!\cap \!\wt\cY(i)$. On the other hand, it follows from (\ref{equ:wtcNijk}), (\ref{equ:nested4}), and Remark~\ref{rmk:clarification} that 
$$
z\in \cY(i,j,k)\cap \wt\cY(i)\cap \phi_{ik}^{-1}(\wt\cY(k)) \subset  \cY(i,j)\cap \wt\cY(i);
$$
this is a contradiction.
\eProof

\noindent
The following lemma is one of the reasons we consider EOB. More importantly, we will make use of EOB to consistently deform Kuranishi maps into traversal orbifolds multisections.

\bLem{lem:ds-ext}
Let $\mfL$ be an DGS as in Definition~\ref{level-system} and $\cN$ be an EOB for $\mfL$ as in Definition~\ref{def:ext-obs}. Via $\cN$, for every $i,j\!\in\!\mfB$, with $i\!<\!j$, the normal direction derivative map $\nd s_{j/i}$ in  (\ref{comm-tangent}), originally defined along the zero locus, extends to a similarly denoted map
\bEqu{extended-der}
\nd s_{j/i}\colon \cN_{\phi_{ij}(\cY(i,j))}\cY(j) \lra \cE(j)/\mfD\phi_{ij}\cE(i,j).
\eEqu
\eLem

\bProof
The argument is local. For simplicity we forget about the orbifold structure and consider manifolds and vector bundles instead. The argument readily extends to local charts with group actions. Thus, let $M_1\subset M_2$ be a submanifold, $\pi\colon E_2\lra M_2$ be a vector bundle, $\pi\colon E_1 \lra M_1$ be a sub-bundle of $E_2|_{M_1}$, and $E'_1$ be an extension of $E_1$ to a sub-bundle of $E_2$ over a neighborhood $M_{12}$ of $M_1\subset M_2$. Let $s_2$ be a section of $E_2$ such that restricted to $M_1$, $s_1\!=\! s_2|_{M_{1}}$, is a section of $E_1$. Fix some $x\in M_1$. Let $U_{12}\!\subset\! M_{12}$ be a sufficiently small chart around $x$ such that 
\bEqu{equ:another}
E'_1|_{U_{12}}\cong E_1|_x \times U_{12},\quad E_2|_{U_{12}}\cong E_2|_x \times U_{12}.
\eEqu
Given a choice of such trivialization, it determines a normal direction derivative map
\bEqu{equ:ds12-before}
\nd s_{2/1}\colon N_{M_1}M_2|_x \lra (E_2/E_1)|_{x}.
\eEqu
A different choice of trivialization in (\ref{equ:another}) corresponds to some 
$$
g\colon M_{12}\lra \tn{End}(E_2|_x,E_2|_x)
$$ 
such that 
$$
g(y)(E_1|_x)= E_1|_x \quad \forall~y\in M_{12}.
$$
For the new choice of trivialization, the map (\ref{equ:ds12-before}) changes to 
$$
\nd^{\tn{new}} s_{2/1} ([v])= [g(x)][\nd s_{2/1}(v)]+[\nd_x g(v) s_2(x) ].
$$
Since $g$ preserves $E_1'$ and $s_2(x)\in E_1$, the second term on the left-hand side is zero;
therefore, $\nd s_{2/1}$ is independent of the choice of trivialization.
\eProof

\noindent
By definition, if $\mfL$ has tangent bundle, $\nd s_{j/i}$ restricted to the zero set of Kuranishi map is an isomorphism. Therefore, after possibly shrinking 
$(\mfL,\cN)$, we may assume  that (\ref{extended-der}) is an isomorphism along entire $\phi_{ij}(\cY(i,j))$.

\noindent
For $i,j\!\in\!\mfB$, with $i\!<\!j$, let 
$$
f\colon \cN'_{\phi_{ij}(\cY(i,j))}\cY(j)\lra \cN(i,j)
$$
be a regularization of $\cN(i,j)$ as in Definition~\ref{def:tub-nbhd}. Fix metrics $g$ on $\cN_{\phi_{ij}(\cY(i,j))}\cY(j)$ and $h$ on $\cE(j)$. If $\nd s_{j/i}$ is an isomorphism on the entire $\phi_{ij}(\cY(i,j))$, by the second assumption in (\ref{equ:embedd-bundles}) and after restricting to a sufficiently small sub-neighborhood $\cN''_{\phi_{ij}(\cY(i,j))}\cY(j)$, there exists a positive continuous function 
$$
c \colon \phi_{ij}(\cY(i,j))\lra \R^+
$$ 
such that
\bEqu{equ:ineq1}
|[s_j(x)]| \geq c |v| \quad \forall~x\!=\!f(v)\in\cN(i,j),~~v\in \cN''_{\phi_{ij}(\cY(i,j))}\cY(j).
\eEqu
Here, $|v|$ is the norm of $v$ with repeat to $g$, $[s_j(x)]$ is the class of $s_j(x)$ in $\cE(j)/\cE_{\tn{ext}}(i,j)$, and $|[s_j(x)]|$ is the norm of $[s_j(x)]$ with respect to the metric induced by $h$. Changing the regularization $f$ or the metrics $g$ and $h$ changes the function $c$ but the existence of a lower bound (\ref{equ:ineq1}) is a property  of the tubular neighborhood independent  of such choices.

\bDef{def:goodN}
Let $\mfL$ be an DGS as in Definition~\ref{level-system} with tangent bundle.
A \textbf{perfect} EOB $\cN$ for $\mfL$ is an EOB such that for every $i,j\!\in\!\mfB$, with $i\!<\!j$, the normal direction derivative map $\nd s_{j/i}$  along entire $\phi_{ij}(\cY(i,j))$ is an isomorphism and (\ref{equ:ineq1}) holds. 
\eDef

\bRem{rmk:orient-vs-perfect}
Via a perfect EOB, and after possibly some shrinking so that $s_{i,j}^{-1}(0)$ over each connected component of $\cY(i,j)$ is non-empty, we can extend the isomorphism (\ref{equ:iso-det-bundles}) to an isomorphism over entire $\phi_{ij}(\cY(i,j))$.
\eRem

\bThm{thm:exists-extension}
Let $\mfL$ be an DGS as in Definition~\ref{level-system} with tangent bundle. There exists a shrinking of $\mfL$ that admits a perfect EOB. 
\eThm
\noindent
We postpone the proof to Section~\ref{sec:exists-extension}. Given such an EOB, in the rest of this section, we introduce the notion of deformation of Kuranishi map. Note that further shrinkings of an DGS with a perfect EOB as in Lemma~\ref{lem:shrunk-EOB} remain perfect.

\noindent
Similar to Section~\ref{sec:multisection}, in order to deform the Kuranishi maps into something transversal, we substitute the orbifold sections with transversal orbifold multisections. We first need to extend the notion of normal direction derivative map in (\ref{global-dsij}) to the case of multisections. Since the argument is local, suppose 
$$
\cU_i=(\pr_i\colon\!U_i\lra V_i,G_i,s_i,\psi_i),\quad i=1,2,
$$
are two Kuranishi charts of dimension $n\!=\!\dim V_i\!-\!\tn{rank}~ U_i$ centered at $\ov{x}$, and 
$$
\Phi=(\mfD\phi,\phi,h)\colon \cU_1\lra \cU_2
$$ 
is a coordinate change map as in Definition~\ref{def:intersection}. Since $\cU_1$ and $\cU_2$ are local charts centered at $\ov{x}$, we have 
$$
(G_1)_x\!\cong \!G_1\!\cong\! h(G_1)\!\cong\! G_2\! \cong\! (G_2)_x
$$ 
and we simply denote these isomorphic groups by $G$. We also identify $\phi(V_1)$ and $\mfD\phi(U_1)$ with their images and think of $V_1$ as an $G$-invariant sub-manifold of $V_2$. Assume $U_1^{\tn{ext}}$ is an extension of $U_1$ to a $G$-invariant sub-bundle of $U_2$. Suppose $\mfs_2\!=\![s^1_2,\ldots,s^k_2]$ is a smooth lifted orbifold $k$-section of $\cU_2$ such that 
$$
s^i_1= s^i_2|_{V_1}\in U_1\quad \forall i\!=\!1,\ldots,k.
$$
Similar to Lemma~\ref{lem:ds-ext}, the extension $U_1^{\tn{ext}}$ results in a set of well-defined normal direction derivative maps
\bEqu{equ:branch-wise}
\nd s^i_{2/1}\colon N_{V_1}V_2 \lra U_2/U_1 \quad \forall i\!=\!1,\ldots,k.
\eEqu
While every individual map in (\ref{equ:branch-wise}) may not be $G$-invariant, the set $\{\nd s^i_{2/1}\}_{i\in[k]}$ is $G$-invariant. Putting together, we obtain a well-defined map
\bEqu{equ:multi-ds}
\nd \mfs_{2/1}\colon S^k(N_{V_1}V_2)\lra S^k (U_2/U_1),\quad \mfs_1=\mfs_2|_{V_1},
\eEqu
which we call the \textbf{normal direction derivative map} for lifted multisections. Then $\nd \mfs_{2/1}$ is an isomorphism if and only if all branches in (\ref{equ:branch-wise}) are isomorphisms. 

\bDef{def:pert-system-sections}
Let $\cN$ be an EOB for an DGS $\mfL$ as in Definition~\ref{def:ext-obs} and $m\!\in\!\mathbb{N}$. A set of \textbf{compatible sequences of deformations} of Kuranishi maps $\{s_{i}\}_{i\in \mfB}$, is a set of sequences of transversal (thus lifted) orbifold multisections,
$$
\mfs_{i;a}\in C^m_{\tn{multi}}(\cE(i)) \quad \forall i\!\in\! \mfB,~a\!\in\! \N,
$$
satisfying the following conditions.
\bEnum
\item For every $i,j\!\in\! \mfB$, with $i\!<\!j$, and every $a\!\in\! \N$, 
$$
\mfD\phi_{ij}(\mfs_{i;a}|_{\cY(i,j)})\cong \mfs_{j;a}|_{\phi_{ij}(\cY(i,j))}
$$
as lifted multisections, and the normal derivative map $\nd \mfs_{j;a/i;a}$ (\ref{equ:multi-ds}) along $\phi_{ij}(\cY(i,j))$, defined via $\cN$, is an isomorphism. 
\item\label{l:zerozero} For every $i,j\!\in\! \mfB$, with $i\!<\!j$, and every $a\!\in\! \N$, 
$$
\mfs_{j;a}^{-1}(0)\cap \cN(i,j)=\phi_{ij}(\mfs_{i;a}^{-1}(0)\cap \cY(i,j)).
$$
\item\label{l:valval}  With the function $\tn{val}_{\mfs}$ as in (\ref{equ:val-function})
$$
\tn{val}_{\mfs_{i;a}}(\ov{x})\!=\!\tn{val}_{\mfs_{j;a}}(\ov\phi_{ij}(\ov{x}))\quad \forall i,j\!\in\! \mfB, i\!<\!j, a\!\in\! \N, \ov{x}\!\in\!Y(i,j).
$$
\item For every $i\!\in\! \mfB$, as $a\lra \infty$, $\mfs_{i;a}$ $C^m$-converges  to $s_i$. 
\eEnum
\eDef

\bThm{thm:pert-system-sections}
Let $\mfL$ be an  DGS with tangent bundle as in Definition~\ref{level-system} and $m\!\in\!\mathbb{N}$ is sufficiently large. 
Then there exists a shrinking $\mfL'$ of $\mfL$, as in Lemma~\ref{lem:shrinking}, and a prefect EOB $\cN$ for $\mfL'$ such that $(\mfL',\cN)$ admits a compatible set of sequences of deformations as in Definition~\ref{def:pert-system-sections}.
\eThm

\bProof
After possibly shrinking $\mfL$, we equip $\mfL$ with a perfect EOB $\cN$ as in Theorem~\ref{thm:exists-extension}. Let $i_1$ be the smallest value in $\mfB$. Let $(\mfL^\one,\cN^\one)$ be a shrinking of $(\mfL,\cN)$ as in Lemma~\ref{lem:shrunk-EOB}; thus, $\tn{cl}\!\lrp{Y^\one(i_1)} \subset Y(i_1)$ is compact. 
Then, by Proposition~\ref{pro:estimate}, with $K_1\!=\!W_1\!=\!\emptyset$ and $K_2\!=\!\tn{cl}(Y^\one(i_1))$ in $M\!=\!Y(i_1)$, there exists a sequence of lifted multisections 
$$
\big(\mft_{i_1;a}\in C^{m}_{\tn{multi}}(\cE^\one(i_1))\big)_{a\in \N}$$ 
such that 
$$
\lrp{\mfs_{i_1;a}= s_{i_1}|_{\cY^\one(i_1)}+\mft_{i_1;a}}_{a\in \N}
$$ 
is a sequence of transversal orbifold multisections $C^m$-converging to $s_{i_1}|_{\cY^\one(i_1)}$.

\noindent
Next, for every $j\!\in\! \mfB(>\!i_1)$ and $a\!\in\! \N$, restricted to $\cN^\one(i_1,j)$,
define
\bEqu{e:extend-deformation}
\aligned
&\mft_{j,i_1;a}= \iota_{i_1j}( \pi_{ji_1}^*(\mfD\phi_{i_1j}(\mft_{i_1;a}))),\quad \mfs_{j,i_1;a}= s_j|_{\cN^\one(i_1,j)}+\mft_{j,i_1;a}, \\
& \mft_{j,i_1;a},\mfs_{j,i_1;a} \in C^{m}_{\tn{multi}}(\cE^\one(j)|_{\cN^\one(i_1,j)}).
\endaligned
\eEqu
By the second condition in (\ref{equ:embedd-bundles}), restricted to $\phi_{i_1j}(\cY^\one(i_1,j))$ we have 
$$
\mfs_{j,i_1;a}= \mfD\phi_{i_1j} \mfs_{i_1;a}\quad \forall j\!\in\! \mfB(>\!i_1).
$$ 
Moreover, the pull-back of every lifted multisection is lifted; therefore, $\mft_{j,i_1;a}$ and  $\mfs_{j,i_1;a}$ are lifted as well.
It is also clear from (\ref{e:extend-deformation}) that
$$
\tn{val}_{\mfs_{i_1;a}}(\ov{x})\!=\!\tn{val}_{\mfs_{j,i_1;a}}(\ov\phi_{i_1j}(\ov{x}))\quad \forall j\!\in\! \mfB(>\!i_1), a\!\in\! \N, \ov{x}\!\in\!Y^\one(i_1,j).
$$

\noindent
\textbf{Claim 1.}
\textit{For every $j\!\in\! \mfB(>\!\!i_1)$, $(\mfs_{j,i_1;a})_{a\in\N}$ is a of sequence transversal multisections over $\cN^\one(i_1,j)$ converging in $C^m$-norm to $s_j|_{\cN^\one(i_1,j)}$. Moreover, 
\bEqu{equ:equal-zero-set}
\mfs_{j,i_1;a}^{-1}(0)=\phi_{i_1j}(\mfs_{i_1;a}^{-1}(0)\cap \cY^\one(i_1,j)).
\eEqu
}

\bProof
By (\ref{e:extend-deformation}), the $\cE^\one_{\tn{ext}}(i_1,j)^\perp$ component of $\mfs_{j,i_1;a}$ is equal to the $\cE^\one_{\tn{ext}}(i_1,j)^\perp$ component of $s_j$; therefore, (\ref{equ:equal-zero-set}) follows from (\ref{equ:ineq1}).\

\noindent
Consider an arbitrary point $x\!\in\! \mfs_{j,i_1;a}^{-1}(0)$. The $\cE^\one_{\tn{ext}}(i_1,j)^\perp$ component 
$$
\nd^\perp \mfs_{j,i_1;a/i_1;a}(x)
$$ 
of $\nd \mfs_{j,i_1;a/i_1;a}(x)$ is equal to $\nd s_{j/i_1}(x)$, which surjects onto $\cE^\one_{\tn{ext}}(i_1,j)^\perp_x$ by the tangent bundle condition. The 
$\cE^\one_{\tn{ext}}(i_1,j)$ component 
$$
\nd^{i_1} \mfs_{j,i_1;a/i_1;a}(x)
$$
of $\nd \mfs_{j,i_1;a/i_1;a}(x)$ 
is equal to $\nd (s_{i_1}\!+\mft_{i_1;a})(x)$, which surjects onto $\cE^\one_{\tn{ext}}(i_1,j)_x$, because $(s_{i_1}\!+\mft_{i_1;a})_{a\in \N}$ is a sequence of transversal orbifold multisections for $\cE^\one(i_1)$.
\eProof
\noindent
We proceed by induction.  For $k\!\geq\! 1$, let $\mfB_k$ be the set of first $k$ values in $\mfB$. Suppose $(\mfL^\sk,\cN^\sk)$ is a shrinking of $(\mfL,\cN)$,  as in Lemma~\ref{lem:shrunk-EOB},  admitting a compatible set of sequences of deformations 
$$
\mfs_{i;a}=s_i|_{\cY^{(k)}}\!+\! \mft_{i;a}\in C^m_{\tn{multi}}(\cE^\sk(i)) \quad\forall i\!\in\! \mfB_k,~a \!\in\! \N,
$$
for a family of lifted multisections 
$$
\mft_{i;a}\in C^m_{\tn{multi}}(\cE^\sk(i))\quad \forall i\!\in\! \mfB_k,~a\!\in\! \N,
$$ 
such that for every $i,i'\!\in \!\mfB_k$, with $i\!<\!i'$ and $\mft_{i',i;a}$ defined over $\cN^\sk(i,i')$ as in (\ref{e:extend-deformation}), 
$$
\mft_{i';a}|_{\cN^\sk(i,i')}=\mft_{i',i;a}\quad \forall a\!\in\! \N.
$$
Let $i_{k+1}$ be the smallest number in $\mfB\! \setminus\! \mfB_k$ and $(\mfL^\skk,\cN^\skk)$ be a shrinking of $(\mfL^\sk,\cN^\sk)$ as in Lemma~\ref{lem:shrunk-EOB}. Set 
$$
K^\circ_1 = \bigcup_{i\in \mfB_k} N^\skk(i,i_{k+1}) \subset Y^\sk(i_{k+1}), \quad 
W_1= \bigcup_{i\in \mfB_k} N^\sk(i,i_{k+1}) \subset Y^\sk(i_{k+1}),
$$
and $K_1\!=\!\tn{cl}(K_1^\circ)\!\subset\! W_1$. For every $a\!\in\!\N$, by Definition~\ref{def:ext-obs}.\ref{l:nested3}, the multisection
$$
\mft'_{i_{k+1};a}\in  C^m_{\tn{multi}}(\cE^\sk(i_{k+1})|_{W_1}), \quad \mft'_{i_{k+1};a}|_{\cN^\sk(i,i_{k+1})}= \mft_{i_{k+1},i;a}\quad \forall i\in \mfB_k,
$$
is well-defined and lifted. By an argument similar to Claim~1 above and the paragraph before that, for every $a\!\in\!\N$,  the lifted multisection 
$$
\mfs'_{i_{k+1};a}\in  C^m_{\tn{multi}}(\cE^\sk(i_{k+1})|_{W_1}),\quad \mfs'_{i_{k+1};a}= s_{i_{k+1}}\!+\!\mft'_{i_{k+1};a}, 
$$
is transverse to the zero-section, 
$$
(\mfs'_{i_{k+1};a})^{-1}(0) \cap \cN^\sk(i,i_{k+1})=\phi_{ii_{k+1}}(\mfs_{i;a}^{-1}(0)\cap \cY^\sk(i,i_{k+1}))\quad \forall i\!\in\! \mfB_k,
$$
and 
$$
\tn{val}_{\mfs_{i;a}}(\ov{x})\!=\!\tn{val}_{\mfs'_{i_{k+1},i;a}}(\ov\phi_{ii_{k+1}}(\ov{x}))\quad \forall i\!\in\! \mfB_k, a\!\in\! \N, \ov{x}\!\in\!Y^\sk(i,i_{k+1}).
$$
Therefore, by Proposition~\ref{pro:estimate}, with $K_1$ and $W_1$ above and 
$$K_2=\tn{cl}(Y^\sk(i_{k+1}))\subset M\!=\!Y(i_{k+1}),$$ 
there exists a sequence of lifted multisections 
$$
\big(\mft_{i_{k+1};a}\in C^{m}_{\tn{multi}}(\cE^\sk(i_{k+1}))\big)_{a\in \N}
$$ 
such that 
$$
(\mfs_{i_{k+1};a}= s_{i_{k+1}}|_{\cY^\sk(i_{k+1})}+\mft_{i_{k+1};a})_{a\in \N}
$$ 
is a sequence of transversal orbifold multisections $C^m$-converging to $s_{i_{k+1}}|_{\cY^\sk(i_{k+1})}$ and 
$$
\mft_{i_{k+1};a}|_{\cN^\skk(i,i_{k+1})}=\mft'_{i_{k+1};a}|_{\cN^\skk(i,i_{k+1})}=\mft_{i_{k+1},i;a}\quad \forall a\!\in\!\N,~i\!\in\! \mfB_k.
$$
After the restriction to $(\mfL^\skk,\cN^\skk)$,  this establishes the induction step and thus completes the proof of Theorem~\ref{thm:pert-system-sections}.
\eProof

\subsection{Construction of perfect EOB (proof of Theorem~\ref{thm:exists-extension})}\label{sec:exists-extension}
Let $\mfL$ be an DGS with tangent bundle as in Definition~\ref{level-system}. 
For $\ell\!=\!1$, an EOB of width $1$ is simply a set of tubular neighborhoods $\cN(i,i\!+\!1)$ and extensions of obstruction bundles $\cE(i,i\!+\!1)$, for every $i\!\in\! \mfB$ with $i\!+\!1\!\in\!\mfB$, without any compatibility condition.
Construction of such extensions, via exponential map and parallel extension of sub-orbibundle to a neighborhood can be done in the following way.

\noindent
We equip every $T\cY(i)$ with an orbifold Riemannian metric $g_i$. Let $\nabla^{\cY(i)}$ the Levi-Civita connection of $g_i$. Also, we equip every obstruction bundle $\cE(i)$ with a metric $h_i$ and a compatible orbibundle connection $\nabla^{\cE(i)}$. 
For every $i,j\!\in\! \mfB$,  with $j\!=\!i\!+\!1$, identify  $\cN_{\phi_{i}(\cY(i,j))}\cY(j)$ with the orthogonal complement $T\cY(i,j)^\perp$ of the sub-orbibundle 
$$
\nd\phi_{ij}(T\cY(i,j))\subset T\cY(j)|_{\phi_{ij}(\cY(i,j))}.
$$ 
Then, via the exponential map of $g_j$, restricted to $T\cY(i,j)^\perp$, we obtain an isomorphism 
\bEqu{exa-map}
\tn{exp}^{ij}\colon \cN^{\de_{ij}}_{\phi_{ij}(\cY(i,j))}\cY(j)\lra  \cN(i,j)
\eEqu
from a tubular neighborhood of the zero section in $\cN_{\phi_{ij}(\cY(i,j))}\cY(j)$ to a neighborhood of $\phi_{ij}(\cY(i,j))$ in $\cY(j)$. In (\ref{exa-map}), $\de_{ij}\colon Y(i,j)\lra \R^+$ is a continuous positive function and $
\cN^{\de_{ij}}_{\phi_{ij}(\cY(i,j))}\cY(j)$ is the set of vectors of the length less than $\de_{ij}$ with respect to the metric induced by $g_j$. 
Let 
$$
 \cN^{\de_{ij}}_{\phi_{ij}(\cY(i,j))}\cY(j)\stackrel{\pi_{ji}}{\longrightarrow}  \phi_{ij}(\cY(i,j)),
$$ 
be the projection map. Via the orbifold diffeomorphism $\tn{exp}^{ij}$, $\pi_{ji}$ gives rise to a similarly denoted projection map 
$$
\pi_{ji}\colon\cN(i,j)\lra \phi_{ij}(\cY(i,j)).
$$

\noindent
Then the orbibundle embedding 
$$
\iota_{ij}\colon \pi_{ji}^*(\mfD\phi_{ij}\cE(i,j))\xlongrightarrow{\cong}  \cE_{\tn{ext}}(i,j) \subset \cE(j)|_{\cN(i,j)}
$$
at 
$$
\pi_{ji}^*(\mfD\phi_{ij}\cE(i,j))|_y,\quad y=\tn{exp}^{ij}_{x}(v), ~~ x\in \phi_{ij}(\cY(i,j)),~~v\in \cN^{\de_{ij}}_{\phi_{ij}(\cY(i,j))}\cY(j)|_{x},
$$ 
is given by parallel-wise extension of $\mfD\phi_{ij}\cE(i,j)|_{x}$, with respect to $\nabla^{\cE(j)}$,  along the geodesic ray $\{\tn{exp}^{ij}_{x}(tv)\}_{t\in [0,1]}$. \\

\noindent
Let $\cN$ be an EOB of width $\ell\!-\!1$ for $\mfL$ for some $\ell\!>\!1$, $\wt{\mfL}$ be a shrinking of $\mfL$ as in Lemma~\ref{lem:shrinking}, and $\wt\cN$ be a the shrinking of $\cN$ associated to $\wt\mfL$ as in (\ref{equ:wtN}). 
Let $i,k\!\in\!\mfB$, with $k\!-\!i\!=\!\ell$. For every $j\!\in\!\mfB$, with $i\!<\!j\!<\!k$, the tubular neighborhoods $\cN(i,j)$ and $\cN(j,k)$ and the extension of orbibundles $\cE^{\tn{ext}}(i,j)$ and $\cE^{\tn{ext}}(j,k)$  are given by assumption. Choose 
$$
\cW_{ik;j}\subset \cN(i,j)|_{\phi_{ij}(\cY(i,j,k))}
$$ 
to be sub-tubular neighborhood of 
$$
\phi_{ij}(\cY(i,j,k))\subset \cN(i,j)\cap\cY(j,k)
$$
and set
$$
\cN(i,k)_j= \pi_{kj}^{-1}\big(\phi_{jk}(\cW_{ik;j})\big)\subset \cN(j,k).
$$
Let
$$
\pi_{ki;j}\colon \cN(i,k)_j \lra  \phi_{ik}(\cY(i,j,k)), \quad \pi_{ki;j}=\phi_{jk}\circ\pi_{ji}\circ\phi_{jk}^{-1}\circ\pi_{kj}.
$$
Then, for $j,j'\!\in\!\mfB$, with $i\!<\!j\!<\!j'\!<\!k$, by (\ref{equ:nested4}) applied to $j,j',k$, we get
\bEqu{equ:matching-retraction}
\pi_{ki;j}(x)=\pi_{ki;j'}(x)\quad \forall x\!\in\! \cN(i,k)_j\cap \cN(i,k)_{j'}.
\eEqu

\noindent
By (\ref{equ:matching-retraction}), the projection maps $\{\pi_{ki;j}\}_{i<j<k}$ give rise to a projection map
$$
\pi_{ki}'\colon \cN'(i,k)= \bigcup_{\substack{j\in\mfB\\~i<j<k}} \cN(i,k)_j\lra \phi_{ik}(\cY'(i,k)),\quad
\cY'(i,k)= \bigcup_{\substack{j\in\mfB\\~i<j<k}} \cY(i,j,k).
$$
Similarly, by Definition~\ref{def:ext-obs}.\ref{l:nested3}, with 
$$
\cE(i,k)_j= \pi_{kj}^*\mfD\phi_{jk}(\iota_{ij}(\pi_{ji}^*(\mfD\phi_{ij}(\cE(i,j)|_{\cY(i,j,k)}))))\subset \pi_{kj}^*(\mfD\phi_{jk}(\cE(j,k)))|_{\cN(i,k)_j},
$$ 
the orbibundle embedding 
$$
\iota_{ik}'\colon (\pi_{ki}')^*\mfD\phi_{ik}(\cE(i,k)|_{\cY'(i,k)})\lra \cE(k)|_{\cN'(i,k)},\quad
\iota_{ij}'|_{\cE(i,k)_j}=\iota_{jk}|_{\cE(i,k)_j},
$$
is well-defined.

\bLem{lem:merging-bundles}
There exists a tubular neighborhood 
$$\pi_{ki}\colon \cN(i,k)\lra \phi_{ik}(\cY(i,k))$$
with an orbibundle embedding 
$$
\iota_{ik}\colon \!\pi_{ki}^*(\mfD\phi_{ik}\cE(i,k))\!\lra\! \cE(k),\quad \iota_{ik}|_{\phi_{ik}(\cY(i,k))}=\tn{id},
$$
such that $\big(\pi_{ki}, \cN(i,k)|_{\phi_{ik}(\wt\cY(i,j,k))}\big)$ is a sub-tubular neighborhood of  
$$
\big(\pi'_{ki}, \cN'(i,k)|_{\phi_{ik}(\wt\cY(i,j,k))}\big),
$$
for all $i\!<\! j\!<\! k$, and 
$$
(\iota_{ik}=\iota_{ik}')|_{\pi_{ki}^*(\mfD\phi_{ik}\wt\cE(i,j,k))}.
$$
\eLem
\bProof
For every $j\!\in\!\mfB$, with $i\!<\!j\!<\!k$, $\phi_{ik}(\wt\cY(i,j,k))\!\subset\!\cY(k)$ is relatively compact. 
With $$
\wt\cM_1\!= \!\bigcup_{\substack{j\in\mfB\\~i<j<k}}\hspace{-.1in} \phi_{ik}(\wt\cY(i,j,k)),\quad \cM_1'\!=\!\phi_{ik}(\cY'(i,k)),\quad \cW'\!=\!\cN'(i,k), \quad \cM_1\!=\!\phi_{ik}(\cY(i,k)),
$$
and $\cM_2\!=\!\!\cY(k)$, the existence of $\cN(i,k)$ and $\iota_{ik}$ follows from Proposition~\ref{pro:ext-nbhd} below.
\eProof
\noindent
After replacing $\cN(i,k)$ with $\cN(i,k)\cap \wt\cY(k)$ and restricting everything to $\wt\mfL$ and $\wt\cN$, we obtain a set of tubular neighborhoods $\wt\cN(i,j)$ for orbifold embeddings $\phi_{ij}(\wt\cY(i,j))\!\subset\!\wt\cY(j)$ and orbibundle embeddings $\iota_{ij}$ as in (\ref{equ:embedd-bundles}),  for all $i,j\!\in\!\mfB$ with $j\!-\!i\!\leq\! \ell$. 
For every $i,j,k\in \mfB$, with $i\!<\!j\!<\!k$, Definition~\ref{def:ext-obs}.\ref{l:nested1}-\ref{l:nested3} hold by induction, if $k\!-\!i\leq \ell-1$, and they hold by Lemma~\ref{lem:merging-bundles}, if $k\!-\!i\!=\!\ell$.
This finishes the inductive step and thus the proof of Theorem~\ref{thm:exists-extension}. 
\qed

\bPro{pro:ext-nbhd}
Let $\cM_1\!\subset\!\cM_2$ be an orbifold embedding. For $i\!=\!1,2$, let $\tn{pr}_i\colon \cE_i\!\to\! \cM_i$ be orbibundles such that $\cE_1$ is a sub-orbibundle of $\cE_2|_{\cM_1}$.
Let $\cM'_1$ be an open sub-orbifold of $\cM_1$ and $\wt\cM_1$ be a relatively compact open sub-orbifold of $\cM'_1$.
Let $\pi'\colon\cW'\lra \cM'_1$ be a tubular neighborhood of the orbifold embedding $\cM'_1\!\subset\!\cM_2$ and 
$$
\iota'\colon (\pi')^*(\cE_1|_{\cM'_1})\lra \cE_2|_{\cW'}, \quad \iota'|_{\cM'_1}=\tn{id},
$$ 
be an embedding of the orbibundles.
Then there exists a tubular neighborhood $\pi\colon\cW\lra \cM_1$  and an orbibundle embedding 
$$
\iota \colon \pi^*\cE_1\lra \cE_2|_{\cW}, \quad \iota|_{\cM_1}=\tn{id},
$$ 
such that $(\pi|_{\wt\cM_1}, \cW|_{\wt\cM_1})$ is a sub-tubular neighborhood of $(\pi'|_{\wt\cM_1}, \cW'|_{\wt\cM_1})$ and 
$$
(\iota=\iota')|_{\pi^*(\cE_1|_{\wt\cM_1})}.
$$
\ePro

\bProof.
Let 
$$
\pi_0\colon \cW_0\lra \cM_1
$$
be a sufficiently small tubular neighborhood of $\cM_1\!\subset\!\cM_2$ and 
$$
\iota_0\colon \pi_0^*\cE_1\lra \cE_2|_{\cW_0}, \quad \iota_0|_{\cM_1}=\tn{id},
$$
be an embedding of orbibundles  such that 
\bEqu{equ:WWW}
\quad \cW_0|_{\tn{cl}(\wt\cM_1)}\subset \cW'.
\eEqu
Such tubular neighborhood and orbibundle embedding exist by the exponential map construction of the beginning of the section.
The inclusions of open orbifolds in (\ref{equ:WWW}) is not necessarily compatible with $\pi_0$ and $\pi'$. 
%
%
%
%

\noindent
\bLem{lem:diff-change}
There exists a sufficiently small sub-neighborhood $\cW\!\subset \!\cW_0$ of $\cM_1$ in $\cM_2$ with an orbifold embedding
$$
\mc{F}\colon \cW\lra\cW_0, \quad \tn{d}\mc{F}|_{T\cW|_{\cM_1}}=\tn{id},
$$
such that 
\bEqu{equ:comapre-pi}
(\pi_0\circ \mc{F}=\pi')|_{ \pi'^{-1}(\wt\cM_1)\cap\cW}.
\eEqu
\eLem

\noindent
For such $\cW$ and $\mc{F}$, let
\bEqu{equ:lastpi}
 \pi= \pi_0\circ \mc{F}\colon \cW\lra \cM_1, \quad \iota = \mc{F}^* \iota_0\colon \pi^*\cE_1\lra \cE_2|_{\cW}.
\eEqu
The projection map and orbibundle embedding (\ref{equ:lastpi}) have the required properties of Proposition~\ref{pro:ext-nbhd}.
\eProof

\noindent
We finish this section by proving Lemma~\ref{lem:diff-change} which follows from an orbifold version of the isotropy extension theorem. 
\newtheorem*{proofofLem:diff-change}{Proof of Lemma~\ref{lem:diff-change}}
\begin{proofofLem:diff-change}
Let
$$
f'\colon \cN'_{\cM_1'}\cM_2\lra \cW'\quad \tn{and}\quad  f_0\colon \cN^0_{\cM_1}\cM_2\lra \cW_0
$$
be arbitrary regularizations of $\cW'$ and $\cW_0$ as in Definition~\ref{def:tub-nbhd}, respectively. Let
$$
\wh\cW\!\subset\!\cW'_0= (\cW'\cap \cW_0)
$$
be a sufficiently small relatively compact open sub-orbifold such that $\tn{cl}(\wt\cM_1)\!\subset\!\wh\cW$ and
\bEqu{equ:cF}
\wh{\mc{F}}= f'\circ f_0^{-1} \colon \tn{cl}(\wh\cW)\lra \cW'_0
\eEqu
is defined. 

\noindent
\textbf{Claim.} \textit{For sufficiently small $\wh\cW$, there exist a diffeotopy (i.e. a one-parameter family of orbifold diffeomorphisms)
\bEqu{equ:Palais-embedding}
\mc{F}_t\colon \cW'_0\lra \cW'_0, \quad \tn{d}\mc{F}_t|_{T\cM_2|_{\cW'_0\cap\cM_1}}=\tn{id}, \quad \forall t\!\in\![0,1],
\eEqu
supported in $\cW'_0$ such that
\bEqu{equ:restriction-of-cF1}
\mc{F}_0\!=\!\tn{id}\quad\tn{and}\quad\mc{F}_1|_{\wt\cW}\!=\!\wh{\mc{F}}|_{\wt\cW},
\eEqu
where $\wt\cW\!\subset \wh\cW$ is some smaller open neighborhood of  $\tn{cl}(\wt\cM)$.}
\vspace{-.1in}
\bProof
In the case of manifolds, this follows from the isotopy extension theorem \cite[Theorem 1.4, Ch 8]{Hirsch}. More precisely, by definition of regularization and (\ref{equ:cF}) we have
$$
\tn{d}\wh{\mc{F}}|_{T\cM_2|_{\tn{cl}(\wh\cW\cap\cM_1)}}=\tn{id};
$$ 
in particular, $\wh{\mc{F}}|_{\tn{cl}(\wh\cW\cap\cM_1)}=\tn{id}$. Therefore, given a Riemannian metric $g$ on $\cM_2$, for sufficiently small $\wh\cW$, there is a unique vector field
\bEqu{equ:1st-2nd-jets}
\xi\in \Gamma(\cM_2|_{\tn{cl}(\wh\cW)},T\cM_2|_{\tn{cl}(\wh\cW)}),\quad \xi|_{\tn{cl}(\wh\cW\cap\cM_1)}\equiv 0,\quad  \nabla\xi|_{\tn{cl}(\wh\cW\cap\cM_1)}\equiv 0,
\eEqu
such that 
\bEqu{equ:restricted-embedding}
\wh{\mc{F}}(x)=\tn{exp}_x(\xi(x))\quad \tn{and}\quad \mc{F}_t(x)= \tn{exp}_x(t\xi(x))\!\in\! \cW'_0\quad \forall t\!\in\![0,1],~x\!\in\! \tn{cl}(\wh\cW).
\eEqu
For $\wh\cW$ sufficiently small, by the last two conditions of (\ref{equ:1st-2nd-jets}) and the implicit function theorem, $(\mc{F}_t)_{t\in[0,1]}$ is one-parameter family of embeddings of $\wh\cW$ in $\cW'_0$.
Then by \cite[Theorem 1.4, Ch 8]{Hirsch}, we can extend $(\mc{F}_t)_{t\in[0,1]}$ to a diffeotopy supported in $\cW'_0$ such that (\ref{equ:Palais-embedding}) over some smaller open neighborhood $\wt\cW\!\subset \wh\cW$ of  $\tn{cl}(\wt\cM)$ holds. The proof is by cutting off the time dependent vector field of the isotopy (\ref{equ:restricted-embedding}) away from a neighborhood of $\tn{cl}(\wt\cM)$ to obtain a time dependent vector field on entire $\cW'_0$ (such that the last two conditions of (\ref{equ:1st-2nd-jets}) remain valid).
In the orbifold case, the same argument applies. We start with an orbifold Riemannian metric and define (\ref{equ:restricted-embedding}) similarly.  We perform the same operation on the resulting time dependent vector field and consider the orbifold ODE flow of that to get (\ref{equ:restriction-of-cF1}).
\eProof

%
%
%
%
%
\noindent
Going back to the proof of Lemma~\ref{lem:diff-change}, choose a sufficiently small sub-neighborhood $\cW\!\subset\!\cW_0$ such that 
$$
\cW|_{\wt\cM_1}\!\subset\! \wt\cW.
$$
Define
\bEqu{equ:finalcF}
\mc{F}\colon \cW\lra\cW_0, \quad \mc{F}(x)\!=\!\mc{F}_{1}(x)\quad \forall x\!\in\! \cW\cap\cW'_0, \quad \mc{F}(x)\!=\!x\quad \forall x\!\in\! \cW\setminus \cW'_0. 
\eEqu
Since $\mc{F}_1$ is supported in $\cW'_0$, (\ref{equ:finalcF}) is a well-defined orbifold embedding.
Finally, it follows from (\ref{equ:restriction-of-cF1}) and (\ref{equ:cF}) that (\ref{equ:comapre-pi}) holds.
This finishes the proof of Lemma~\ref{lem:diff-change}.
\end{proofofLem:diff-change}

\subsection{Kuranishi vector bundles}\label{sec:bundle-kuranishi}
A Kuranishi structure, as we defined in Section~\ref{sec:basics}, is locally given by a pair of an orbibundle and a section, such that the difference of the dimension of the base and the rank of orbibundle is some fixed constant. In this section, we define the notation of Kuranishi vector bundle. 
It consists of an additional orbibundle of some fixed rank (unlike obstruction bundle) on the base orbifold of every Kuranishi chart. As we show below, we can incorporate this extra bundle into the obstruction bundle and obtain a new augmented Kuranishi structure with a bigger obstruction bundle for each chart. Then, by the results of Section~\ref{sec:VFC}, we define the Euler class of such Kuranishi vector bundle to be  the VFC of corresponding augmented Kuranishi structure. In Section~\ref{sec:GW-VFC}, we use the notion of Kuranishi vector bundle to define GW invariants involving $\psi$-classes.

\vskip.1in
\noindent
\bDef{def:Kur-bundle}
For a Kuranishi structure $\cK$ of real dimension $n$ on a compact metrizable topological space $M$ as in Definition~\ref{def:Kur-structure}, with notation as in Definition~\ref{def:kur-chart}, a \textbf{Kuranishi vector bundle} $\mf{E}$ of rank $m$ on $\cK$ consists of an additional rank $m$ $G_p$-equivariant vector bundle $\pr'_p\colon U_p'\!\lra\! V_p$, for every Kuranishi chart $\cU_p\!=\!(\pr_p\colon\! U_p\!\lra\! V_p,G_p,s_p,\psi_p)$, such the following conditions hold.
\bEnum
\item For every $\ov{q}\!\in\! F_p$, the coordinate change map $\phi_{qp}$ lifts to an equivariant vector bundle isomorphism 
\bEqu{equ:iso-bundle}
\mfD\phi'_{qp}\colon U'_{q,p}\equiv U'_q|_{V_{q,p}} \lra U'_p|_{\phi_{qp}(V_{q,p})}.
\eEqu
\item Whenever $\ov{r}\in F_p\cap F_q$ and $\ov{q}\in F_p$, then $\mfD\phi'_{qp}$, $\mfD\phi'_{qp}$, and $\mfD\phi'_{rp}$ satisfy the cocycle condition of Definition~\ref{def:cocycle}.
\eEnum
\eDef
\noindent
Similarly to Section~\ref{sec:existence}, the local defining bundles of a  Kuranishi vector bundle $\mf{E}$ can be glued 
into a \textbf{dimensionally graded bundle} $\mf{E}_\mfL$,  or \textbf{DGB} for short, over different orbifolds $\cY(i)$, with $i\!\in\!\mfB$, of an associated natural DGS $\mfL$. 
In other words, $\mf{E}_\mfL$ consist of an additional rank $m$ orbibundle $\pr'_i\colon\mf{E}(i)\!\lra\! \cY(i)$, for every orbifold $\cY(i)$ in $\mfL$, such that compatibility conditions similar to Conditions~\ref{l:intersection-2}, \ref{l:intersection-cycle}, and \ref{l:compatible} of Definition~\ref{level-system} hold.

\noindent
Given an $n$-dimensional Kuranishi structure $\cK$ on $M$, and a rank $m$ Kuranishi vector bundle $\mf{E}$ on $\cK$, we build a ``new" Kuranishi structure $\cK^{\mf{E}}$ of dimension $n\!-\!m$ on $M$ in the following way. For every $\ov{p}\!\in\! M$, we define the Kuranishi chart $\cU^{\mf{E}}_p$ of $\cK^{\mf{E}}$ centered at $\ov{p}$ to be
\bEqu{equ:Augmented-chart}
\cU^{\mf{E}}_p=\{ \pr^{\mf{E}}_p\equiv \pr_p\oplus \pr'_p\colon\! U^{\mf{E}}_p\equiv U_p\oplus U'_p \!\lra\! V_p, {s}^{\mf{E}}_p\equiv s_p\oplus 0,\psi_p\}.
\eEqu
By definition, $({s}^{\mf{E}}_p)^{-1}(0)=s^{-1}_p(0)$; therefore, $\psi_p$ still makes sense and has the same footprint. We define the change of coordinate maps similarly. Moreover, by definition, for the resulting coordinate change maps $\Phi^{\mf{E}}_{qp}\colon \cU^{\mf{E}}_{q,p}\lra \cU^{\mf{E}}_p$, since 
$$
({U}^{\mf{E}}_p/\mfD\phi^{\mf{E}}_{qp}({U}^{\mf{E}}_{q,p}))|_{\phi^{\mf{E}}_{qp}(({s}^{\mf{E}})^{-1}_{q,p}(0))}= 
(U_p/\mfD\phi_{qp}(U_{q,p}))|_{\phi_{qp}(s^{-1}_{q,p}(0))},
$$
we conclude that if $\cK$ has tangent bundle, then $\cK^{\mf{E}}$ has tangent bundle as well. Similarly, if $\cK$ and $\mf{E}$ are oriented, $\cK^{\mf{E}}$ is oriented as well. It is clear from (\ref{equ:Augmented-chart}) and the construction of Section~\ref{sec:existence} that every Kuranishi vector bundle gives rise to an DGB over every associated natural DGS. 

\noindent
We then define the Euler class of $\mf{E}$ to be the VFC of the augmented Kuranishi space $(M,\cK^{\mf{E}})$, as defined in Section~\ref{sec:VFC}. In Section~\ref{sec:GW-VFC}, we use  vector bundles over Kuranishi structures and their Euler class to define GW invariants involving $\psi$-classes and other similarly defined insertions.

\section{VFC for abstract Kuranishi spaces}\label{sec:VFC}
\noindent
The notion of a virtual fundamental class, or simply a VFC,  is a generalization of the fundamental class of topological manifolds to more general compact topological spaces with singularities and with components of various dimensions.

\noindent
For a suitable kind of topological space $M$, with some extra structure $\cM$, a VFC is ideally a homology class 
$$[\cM]^\vir \!\in\! H_{\dim_\vir(\cM)}(M),
$$ 
where $\dim_\vir(\cM)$ is the ``expected dimension" of $\cM$. Depending on the application, the homology theory used above can be the de Rham, singular, or the \v{C}ech homology/cohomology.
In the case of extracting invariants from moduli spaces, evaluating certain natural cocycles against the VFC class should be thought of as the integration over that moduli space as in (\ref{equ:GW}).

\noindent
For an oriented Kuranishi space $(M,\cK)$ as in Definition~\ref{def:Kuranishi-tangent} and \ref{def:Kuranishi-orientable}, we first build $[\cK]^\vir$ as a singular homology class inside a thickening $Y(\mfL)$ as defined in (\ref{equ:thikening-space}), where $\mfL$ is a dimensionally graded system for $\cK$ as in Definition~\ref{level-system}. Therefore, a priori  $[\cK]^\vir$ is not supported in $M$ and depends on the auxiliary data of the thickening.
After this step, which is common between various constructions in the literature, there are different  methods for getting a VFC independent of the choice of the thickening (and other auxiliary data). We explain two different methods: one by looking at the image cycle under certain natural evaluation maps as in (\ref{equ:evst}) that extend to Kuranishi structure, and the other, by taking the inverse limit, due to McDuff-Wehrheim  \cite[Theorem B]{MW2}, which gives us a \v{C}ech homology class over $M$, itself.

\subsection{The construction of a VFC in a thickening}\label{sec:VFC1}
\noindent
Let $(M,\cK)$ be an oriented $n$-dimensional Kuranishi space, $\mfL$ be an oriented Hausdorff DGS for $\cK$ as in 
(\ref{level-tuple}),  and $\cN$ be a perfect EOB for $\mfL$ as in Definition~\ref{def:goodN}. By Theorem~\ref{thm:pert-system-sections}, such an EOB gives rise to a compatible set of sequences of deformations $\mfs_a\!\equiv\! \{\mfs_{i,a}\}_{i\in \mfB, a\in \N}$ of Kuranishi maps $s\!\equiv\! \{s_{i}\}_{i\in \mfB}$ of $\mfL$.

\noindent
Similar to Section~\ref{sec:euler}, for sufficiently large $a\!\in\! \N$, after restricting to some shrinking $\mfL'$ of $\mfL$, we build a singular $n$-cycle supported in the zero set of $\mfs_a$ inside the thickening $Y(\mfL')$. 

\noindent
With notation similar to (\ref{equ:sero-set}), let 
$$
\begin{aligned}
Z_{i}\!=\! Z(s_{i})\!\subset\! Y(i)\qquad \forall i\!\in\! \mfB,& \qquad M\!=\!Z(s)\!=\!\bigcup_{i\in\mfB} \wp_i(Z_{i})\!\subset\! Y(\mfL),\\
Z_{i,a}\!=\! Z(\mfs_{i,a})\!\subset\! Y(i)\qquad \forall i\!\in\! \mfB,~a\!\in\! \N,& \qquad Z(\mfs_a)\!=\!\bigcup_{i\in\mfB} \wp_i(Z_{i,a})\!\subset\! Y(\mfL)\quad\forall a\!\in\! \N.
\end{aligned}
$$
By Definition~\ref{def:pert-system-sections}.\ref{l:zerozero},
\bEqu{equ:deformed-zeros}
\ov\phi_{ij}(Z_{i,a}\cap Y(i,j))\!=\!Z_{j,a}\cap N(i,j) \qquad \forall i,j\!\in\! \mfB,~i\!<\!j,~a\!\in\! \N.
\eEqu
Let 
\bEqu{equ:shrink123}
\mfL^\two\!\subset\!\mfL^\one\!\subset\!\mfL^\zero\equiv \mfL
\eEqu
be a sequence of shrinkings of $\mfL$.
The natural inclusion maps
\bEqu{equ:bc}
Y(\mfL^\tb)\!\subset\!Y(\mfL^\tc)\quad \forall~0\!\leq \!c\!<\!b\!\leq\! 2
\eEqu
are relatively compact embeddings, but they are not open\footnote{Otherwise, the quotient topology would have been the same as induced topology, which by the example before Definition~\ref{def:shrinking} we know it is not the case.}. 
For $b\!=\!1,2$, let  $s^\tb$, $\mfs^\tb_{i,a}$, etc, be the restriction of corresponding objects to the corresponding sets for $\mfL^\tb$. By Remark~\ref{rem:ulternative}, the induced and the quotient topologies on 
$$
\tn{cl}(Y(\mfL^\tb))\subset Y(\mfL),\quad \tn{cl}(Y(\mfL^\tb))=\bigcup_{i\in\mfB}\wp_i(\tn{cl}\big(Y^\tb(i))\big),\quad b\!=\!1,2,
$$
are the same and are metrizable.

\bLem{lem:comp-support}
For sufficiently large $a\!\in\!\N$, 
$$
Z(\mfs^\two_{a})= Z(\mfs^\one_{a}); 
$$
in particular $Z(\mfs^\one_{a})$ is compact.
\eLem

\bProof
Since 
$$
M\!=\!Z(s^\tb)\!\subset\! Y(\mfL^\tb)\quad \forall b\!\in\!\{0,1,2\},
$$
for every $p\!\in\!M$, there exists a maximal index $i_{b,p}\!\in\!\mfB$ such that 
\bEqu{equ:dfnxibp}
p\!\in\!\wp_{i_{b,p}}(\ov{x}_{i_{b,p}})\quad\tn{for some (unique)}\quad \ov{x}_{i_{b,p}}\!\in\! Z^\tb_{i_{b,p}}\!\subset\!Y^\tb(i_{b,p}).
\eEqu
For $b\!<\!c$, by (\ref{equ:shrink123}) and (\ref{equ:sim-relation}), 
\bEqu{equ:propertiesxibp}
i_{b,p}\geq i_{c,p}, \quad \ov{x}_{i_{c,p}}\!\in\!Y(i_{c,p},i_{b,p}), \quad \ov\phi_{i_{c,p}i_{b,p}}(\ov{x}_{i_{c,p}})=\ov{x}_{i_{b,p}}.
\eEqu
For $b\!=\!0$, we will write $i_{0,p}$ as $i_p$.
For every $p\!\in\!M$, let 
$$
W_{\ov{x}_{i_{p}}}\subset Y(i_{p})
$$
be an open set such that 
\bEqu{equ:Uxp}
W_{\ov{x}_{i_{p}}}\cap \tn{cl}(Y^\one(i_{p},j))\!=\!\emptyset \quad \forall j>i_{p},\qquad W_{\ov{x}_{i_{p}}}\subset \cN(j,i_{p})\quad \forall j<i_{p}.
\eEqu
Although 
$$
W_p:= \wp_{i_p}(W_{\ov{x}_{i_{p}}})\subset Y(\mfL)
$$
is not necessarily open, but similarly to the case of (\ref{equ:open-in-target}) in the proof of Proposition~\ref{pro:shrinking}, for $b\!=\!1,2$, the restriction 
$W_p\!\cap\! \tn{cl}(Y(\mfL^\tb))$ is open in $\tn{cl}(Y(\mfL^\tb))$.
For a sequence of integers $(a_k)_{k=1}^\infty$ with $\lim_{k\to\infty}a_k\!=\!\infty$, suppose there exists 
$$
p_k\!\in\! Z(\mfs_{a_k}^\one) \quad\tn{s.t}\quad p_k\!\not\in\! Z(\mfs_{a_k}^\two)\quad \forall k\in\N;
$$
i.e. $p_k\!\not\in\! Y(\mfL^\two)$, for all $k\!\in\!\N$.
By compactness of $\tn{cl}(Y(\mfL^\one))$ there exists 
$$
p\!\in\! \tn{cl}(Y(\mfL^\one))\subset Y(\mfL)
$$
such that $\lim_{k\to\infty}p_k\!=\!p$ with respect to a metric for the metrizable topology of $\tn{cl}(Y(\mfL^\one))$. Since $\mfs_{i,a_k}$ is $C^1$-converging to $s_i$, for all $i\!\in\!\mfB$, we conclude that $p\!\in\! M\!=\!Z(s)$; in particular, $p\!\in\!Y(\mfL^\two)$. 

\noindent
For $k$ sufficiently large, $p_k\!\in\!W_p$ and by the left-hand side of (\ref{equ:Uxp}) there exists (a unique) $\ov{x}_k\!\in\! W_{\ov{x}_{i_p}}$ such that $\wp_{i_p}(\ov{x}_k)\!=\!p_k$. By (\ref{equ:propertiesxibp}) and (\ref{equ:dfnxibp}),
$$
i_{2,p}\leq i_{0,p}\equiv i_p, \quad \ov{x}_{i_{2,p}}\!\in\!Y(i_{2,p},i_{p})\cap Y^\two(i_{2,p}), \quad \ov\phi_{i_{2,p}i_{p}}(\ov{x}_{i_{2,p}})=\ov{x}_{i_{p}}.
$$
By the right-hand side of (\ref{equ:Uxp}),
$$
\ov{x}_k\!\in\! N(i_{2,p},i_{p});
$$
therefore, by (\ref{equ:deformed-zeros}),
$$
\ov{x}_k\!\in\!\ov\phi_{i_{2,p}i_{p}} Y(i_{2,p},i_{p}).
$$ 
Since $\ov{x}_k$ is converging to $\ov{x}_{i_{p}}$, we conclude that $\ov{x}_k'\!\equiv\! \ov\phi_{i_{2,p}i_{p}}^{-1}(\ov{x}_k)$ is converging to $\ov{x}_{i_{2,p}}\!\in\!Y^\two(i_{2,p})$. Therefore, $\ov{x}_k'\!\in\!Y^\two(i_{2,p})$, for sufficiently large $k$. This implies that 
$$
p_k\!=\!\wp_{i_{2,p}}(\ov{x}_k')\!\in\!Y(\mfL^\two)
$$  
for sufficiently large $k$, which is a contradiction.

\eProof

\vskip.1in
\noindent
In the rest of this section, for every $a\!\in\!\N$ such that the conclusion of  Lemma~\ref{lem:comp-support} holds, we build a singular $n$-cycle\footnote{Where $n$ is the dimension of $\mfL$.} with $\Q$-coefficients supported in $Z(\mfs^\one_a)$ whose homology class in $Y(\mfL)$ is independent of the choice of the perturbation. This gives us a VFC as a singular homology class inside the thickening $Y(\mfL)$.  To this end, we use a generalization of the concept of the resolution from Section~\ref{sec:resolution} and construct $\vfc(\mfL)\!\in\! H_n(Y(\mfL),\Q)$ similarly to Section~\ref{sec:euler}.  In the case of pure orbibundle Kuranishi structures, the resulting homology class is simply the Euler class of the orbibundle as in (\ref{equ:Euler-class}). The thickening, and hence the resulting VFC homology class depends on the choice of $\mfL$. We remove this indeterminacy  in Sections~\ref{sec:VFC2} and \ref{sec:VFC3}, in two different ways.

\vskip.1in
\noindent
Starting with $\cK$, $\mfL$, and $\cN$ as in the beginning of this section, we replace $\mfL$ with $\mfL^\one$ of 
Lemma~\ref{lem:comp-support}. Similarly, we replace $\cN$ with the corresponding shrinking of Lemma~\ref{lem:shrunk-EOB}. We will also denote $\mfL^\two$ by $\mfL'$. For the following discussion, let $a\!\in\!\N$ be some fixed sufficiently large integer such that the conclusion of Lemma~\ref{lem:comp-support} holds. For simplicity, we write $\mfs$ for $\mfs_a$. For every $p\!\in\!Y(\mfL)$, let 
$$
\mfB_p=\{i\in \mfB\colon p\!\in\!\wp_i(Y(i))\}.
$$
For every $p\!\in\!Y(\mfL)$ and $i,j\!\in\!\mfB$, with $i\!<\!j$, if 
$$
p=\wp_i(\ov{x}_i),\quad \ov{x}_i\in Y(i),\quad\tn{and}\quad  p=\wp_j(\ov{x}_j),\quad \ov{x}_j\in Y(j),
$$
then 
$$
\ov{x}_i\!\in\!Y(i,j), \quad \phi_{ij}(\ov{x}_{i})=\ov{x}_j,
$$ 
and by Definition~\ref{def:pert-system-sections}.\ref{l:valval}
\bEqu{equ:valivalj}
\tn{val}_{\mfs_{i}}(\ov{x}_i)=\tn{val}_{\mfs_{j}}(\ov{x}_j),\quad\tn{where}~~~ \mfs_i=\mfs|_{\cY(i)}\quad \forall i\!\in\!\mfB;
\eEqu
i.e. the positive integer 
$$
\tn{val}_{\mfs}(p)\equiv \tn{val}_{\mfs_{i}}(\ov{x}_i),\quad \forall i\!\in\!\mfB_p, \ov{x}_i\!\in\!Y(i)\quad \tn{s.t}\quad \wp_i(\ov{x}_i)\!=\!p,
$$
is well-defined. Similar to Section~\ref{sec:resolution}, the restricted function $\tn{val}_{\mfs}\colon\! Y(\mfL')\!\lra\! \N$ is upper 
semi-continuous. Therefore, $\tn{val}_{\mfs}$ restricted to the compact set $M\!=\!Z(\mfs)\!\subset \!Y(\mfL')$ takes only finitely many values 
$$
c_1\!<\!\cdots\!<\!c_N.
$$
Let
$$
Z(\mfs)_\ell=\{p\!\in\! Z(\mfs)\colon \tn{val}_{\mfs}(p)\!=\!c_\ell\}\quad \forall \ell\!\in\![N];
$$
in particular, $Z(\mfs)_N$ is compact. 

\noindent
Let 
$$
\mfB\!=\!\{i_1<\ldots<i_S\}, \quad Z(\mfs)_{i;\ell}\equiv Z(\mfs)_\ell\cap \wp_i(\tn{cl}(Y'(i)))\!\subset\!Y(\mfL')\quad \forall i\!\in\!\mfB, \ell\!\in\![N].
$$
By definition, $Z(\mfs)_{i;\ell} \!=\!\wp_i( Z(\mfs_{i})_\ell )$, where 
\bEqu{equ:Zmfsiell}
Z(\mfs_{i})_\ell\equiv\{\ov{x}\in Z(\mfs_{i}) \cap  \tn{cl}(Y'(i))\colon \tn{val}_{\mfs_{i}}(\ov{x})\!=\!c_\ell\}\!\subset\!Y(i)\quad \forall i\!\in\!\mfB, \ell\!\in\![N].
\eEqu
Since $Z(\mfs)_N$ is compact,  $Z(\mfs_{i})_{N}\!\subset\! Y(i)$, for all $i\!\in\! \mfB$, 
is compact. 

\noindent
Starting with maximal elements $i_S\!\in\!\mfB$ and $N\!\in\![N]$, by Lemma~\ref{lem:cover}, there exist a neighborhood $W_{i_S;N}$ of 
$$
K(\mfs_{i_S})_N\equiv Z(\mfs_{i_S})_N
$$ 
in $Y(i_S)$, and an orbifold $\wt\cW_{i_S;N}$ with a $c_N$-covering map 
$$
q_{i_S;N}\colon \wt\cW_{i_S;N}\lra \cW_{i_S;N}= \cY(i_S)|_{W_{i_S;N}},
$$ 
such that the multisection $\mfs_{i_S}$, restricted to $W_{i_S;N}$, is equivalent to the weighted push forward, with weight $w_{i_S;N}$, of some orbifold section $s_{i_S;N}$ of 
$$
\wt\cE_{i_S;N}= q_{i_S;N}^*(\cE|_{W_{i_S;N}}).
$$
Moving down in $\mfB$ to $i_{S-1}$, let
$$
K(i_{S-1})_N:= Z(\mfs_{i_{S-1}})_{N} \setminus \ov\phi^{-1}_{i_{S-1}i_S}(W_{i_S;N})\subset Y(i_{S-1}).
$$
Once again, by Lemma~\ref{lem:cover}, there exist a neighborhood $W_{i_{S-1};N}$ of 
$
K(\mfs_{i_{S-1}})_N
$ 
in $Y(i_{S-1})$, and an orbifold $\wt\cW_{i_{S-1};N}$ with a $c_N$-covering map 
$$
q_{i_{S-1};N}\colon \wt\cW_{i_{S-1};N}\lra \cW_{i_{S-1};N}= \cY(i_{S-1})|_{W_{i_{S-1};N}},
$$ 
such that the multisection $\mfs_{i_{S-1}}$, restricted to $W_{i_{S-1};N}$, is equivalent to the weighted push forward, with weight $w_{i_{S-1};N}$, of some orbifold section $s_{i_{S-1};N}$ of 
$$
\wt\cE_{i_{S-1};N}= q_{i_{S-1};N}^*(\cE|_{W_{i_{S-1};N}}).
$$
By (\ref{equ:valivalj}) and (\ref{equ:Zmfsiell}), 
$$
K(i_{S-1})_N\cap \tn{cl}(Y'(i_{S-1},i_S))\!=\!\emptyset;
$$
therefore, after possibly replacing $W_{i_{S-1};N}$ with a smaller neighborhood, we may assume that
\bEqu{equ:WS-1S}
W_{i_{S-1};N}\cap \tn{cl}(Y'(i_{S-1},i_S))\!=\!\emptyset.
\eEqu
Continuing downward inductively in $\mfB$, we obtain a set of $c_N$-covering maps 
$$
q_{i;N}\colon \wt\cW_{i;N}\lra \cW_{i;N}= \cY(i)|_{W_{i;N}}\quad \forall i\!\in\!\mfB,
$$ 
such that the multisection $\mfs_{i}$, restricted to $W_{i;N}$, is equivalent to the weighted push forward, with weight $w_{i;N}$, of a natural orbifold section $s_{i;N}$ of 
$$
\wt\cE_{i;N}= q_{i;N}^*(\cE|_{W_{i;N}}), 
$$
and 
$$
W_{i;N}\cap \tn{cl}(Y'(i,j))\!=\!\emptyset\quad \forall i\!<\!j.
$$
Moreover, by the explicit construction of proof of Lemma~\ref{lem:cover}, the orbibundle embeddings 
$\Phi_{ij}\!=\!(\mfD\phi_{ij},\phi_{ij})$ in Definition~\ref{level-system} lift to orbibundle embeddings $\wt\Phi_{ij;N}=(\mfD\wt\phi_{ij;N},\wt\phi_{ij;N})$ between the covering spaces such that the following diagram commutes,
$$
\xymatrix{
\wt\cE_{i;N}|_{\ov{q}_{i;N}^{-1}(\ov\phi_{ij}^{-1}(W_{j;N})\cap W_{i;N})}\ar[rr]^{\qquad\mfD\wt\phi_{ij;N}}	&&\wt\cE_{j;N}\\
\wt\cW_{i;N}|_{\ov{q}_{i;N}^{-1}(\ov\phi_{ij}^{-1}(W_{j;N})\cap W_{i;N})}\ar[rr]^{\qquad\wt\phi_{ij;N}}	 \ar[u]^{s_{i;N}}	&&\wt\cW_{j;N}\ar[u]^{s_{j;N}}\,.
}
$$

\noindent
We next move down to $c_{N\!-\!1}$; i.e. we move down inductively with respect to the lexicographic order 
$$
(i_b,c_\ell)<(i_{b'},c_\ell') \quad\Leftrightarrow \quad \ell\!<\!\ell' \quad \tn{or}\quad (\ell\!=\!\ell'\quad \tn{and} \quad b\!<\!b').
$$
Let 
\bEqu{equ:WN}
W_N\!=\!\bigcup_{i\in\mfB} \wp_i(W_{i;N})\subset Y(\mfL);
\eEqu
this is not necessarily an open set in $Y(\mfL)$, but similarly to the case of (\ref{equ:open-in-target}) in the proof of Proposition~\ref{pro:shrinking}, by (\ref{equ:WS-1S}), the restriction $W_N\!\cap\!Y(\mfL')$ is an open set in $Y(\mfL')$. Let
$$
K_{N-1}= Z(\mfs)_{N-1}\!\setminus\! W_N.
$$
Since $Z(\mfs)\!\subset\!Y(\mfL')$, $W_N\!\cap\!Y(\mfL')$ is open, and $\tn{val}_\mfs$ is upper semi-continuous, $K_{N-1}$
is compact. With $N\!-\!1$ instead of $N$ and $K_{N-1}$ instead of $Z(\mfs)_N$ in the previous argument, we inductively build 
$$
q_{i;N-1}\colon \wt\cW_{i;N-1}\lra \cW_{i;N-1}= \cY(i)|_{W_{i;N-1}}\quad \forall i\!\in\!\mfB,
$$ 
similarly. Continuing inductively with respect to the lexicographic order on $\mfB\!\times\![N]$, we obtain a double-indexed resolution  of $\mfs$ similar to Definition~\ref{def:resolution-conditions}, described in Definition~\ref{def:resolution-KUR} below. 

\bDef{def:resolution-KUR} 
With $\mfL$, $\mfL'$, $\mfs$, and $N$ as above,
a \textbf{resolution} of $\mfs$ consists of a set of $c_\ell$-covering maps
\bEqu{equ:seq-covering2}
q_{i;\ell}\colon\wt\cW_{i;\ell}\lra \cW_{i;\ell},\quad \forall i\!\in\!\mfB,~\ell\!\in\![N],
\eEqu
over open subsets $W_{i;\ell}\!\subset\! Y(i)$ such that the following properties hold.
\bEnum
\item The image of open sets $W_{i;\ell}$ cover $M$ in the sense that 
$$
M\!=\!Z(\mfs)\subset \bigcup_{i=\ell}^N W_{\ell},\quad W_\ell\!=\!\bigcup_{i\in\mfB} \wp_i(W_{i,\ell}).
$$ 
Moreover, for every $\ell\!\in\![N]$, $W_\ell\!\cap\!Y(\mfL')$ is open in $Y(\mfL')$, 
$$
W_\ell\cap \big(Z(\mfs)_{>\ell}\equiv \{p\!\in\!Z(\mfs)\colon \tn{val}_\mfs(p)\!>\!c_\ell\}\big)\!=\!\emptyset,~~ \tn{and}\quad  Z(\mfs)_{\geq \ell}\subset \bigcup_{b=\ell}^N W_b.
$$
\item For every $\ell\!\in\! [N]$ and $i\!\in\!\mfB$, the multisection $\mfs_i$, restricted to $W_{i;\ell}$, is equivalent to the weighted push forward with weight $w_{i;\ell}$ of some orbifold section $s_{i;\ell}$ of $\wt\cE_{i;\ell}\!\equiv\!q_{i;\ell}^*(\cE(i)|_{\cW_{i;\ell}})$.
\item For every $i\!\in\!\mfB$ and $\ell \!>\!\ell'$, on the overlap 
\bEqu{equ:Wiljl}
W_{i;\ell,\ell'}= W_{i;\ell',\ell}= W_{i;\ell'} \cap W_{i;\ell},
\eEqu
there exists a component-wise covering map
$$
q_{i;\ell,\ell'}\colon \wt\cW_{i;\ell}|_{\ov{q}_{i;\ell}^{-1}(W_{i;\ell,\ell'})}\lra \wt\cW_{i;\ell'}|_{\ov{q}_{i;\ell'}^{-1}(W_{i;\ell',\ell})},
$$
such that the analogue of Definition~\ref{def:resolution-conditions}.\ref{l:Wij-qij} holds. 
\item\label{l:ijfixell} For every $\ell\!\in\![N]$ and  $i,j\!\in\!\mfB$ with $i\!<\!j$, on the overlap 
\bEqu{equ:Wiljl2}
W_{i,j;\ell}= \ov\phi_{ij}^{-1}(W_{j;\ell})\cap W_{i;\ell},
\eEqu
the orbifold embedding $\phi_{ij}$ lifts to an orbifold embedding $\wt\phi_{ij;\ell}$ 
$$
\wt\phi_{ij;\ell}\colon \wt\cW_{i;\ell}|_{\ov{q}_{i;\ell}^{-1}(W_{i,j;\ell})}\lra \wt\cW_{j;\ell}.
$$
\item For $i\!< \!j$ and $\ell_i\!<\!\ell_j$, there exists a $c_{\ell_j}$-covering 
$$
q_{j;\ell_j\mid i;\ell_i}\colon \wt\cW_{i;\ell_i\mid j;\ell_j}\lra \cW_{i;\ell_i\mid j;\ell_j}, \quad W_{i;\ell_i\mid j;\ell_j}=\ov\phi_{ij}^{-1}(W_{j;\ell_j})\cap W_{i;\ell_i},$$
such that the orbifold embedding $\phi_{ij}$ lifts to an orbifold embedding 
$$
\wt\phi_{i;\ell_i\mid j;\ell_j}\colon \wt\cW_{i;\ell_i\mid j;\ell_j}\lra \wt\cW_{j;\ell_j},
$$
commuting with $q_{j;\ell_j\mid i;\ell_i}$ on the source and $q_{j;\ell_j}$ on the target, and $q_{j;\ell_j\mid i;\ell_i}$ factors through $q_{i;\ell_i}$:
\bEqu{equ:comm1}
\xymatrix{
 \wt\cW_{i;\ell_i\mid j;\ell_j}\ar[d]^{q'_{j;\ell_j\mid i;\ell_i}}\ar[rr]^{\wt\phi_{i;\ell_i\mid j;\ell_j}} \ar@/^-6pc/[dd]|-{q_{j;\ell_j\mid i;\ell_i}}&& 
 \wt\cW_{j;\ell_j}|_{\ov{q}_{j;\ell_j}^{-1}(W_{j;\ell_j\mid i;\ell_i})}\ar[dd]|-{q_{j;\ell_j}}\\
\wt\cW_{i;\ell_i}|_{\ov{q}_{i;\ell_i}^{-1}(W_{i;\ell_i\mid j;\ell_j})}\ar[d]^{q_{i;\ell_i}}&& \\
\cW_{i;\ell_i\mid j;\ell_j}\ar[rr]^{\phi_{ij}}		&&		 \cW_{j;\ell_j\mid i;\ell_i}\,.
}
\eEqu
\item\label{l:i<j-ell_i>ell_j} For $i\!< \!j$ and $\ell_i\!> \!\ell_j$, there exists a $c_{\ell_j}$-covering 
$$
q_{j;\ell_j\mid i;\ell_i}\colon \wt\cW_{i;\ell_i\mid j;\ell_j}\lra \cW_{i;\ell_i\mid j;\ell_j}, \quad W_{i;\ell_i\mid j;\ell_j}=\ov\phi_{ij}^{-1}(W_{j;\ell_j})\cap W_{i;\ell_i},$$
such that the orbifold embedding $\phi_{ij}$ lifts to an orbifold embedding 
$$
\wt\phi_{i;\ell_i\mid j;\ell_j}\colon \wt\cW_{i;\ell_i\mid j;\ell_j}\lra \wt\cW_{j;\ell_j},
$$
commuting with $q_{j;\ell_j\mid i;\ell_i}$ on the source and $q_{j;\ell_j}$ on the target, and $q_{i;\ell_i}$ factors through $q_{j;\ell_j\mid i;\ell_i}$:
\bEqu{equ:comm2}
\xymatrix{
\wt\cW_{i;\ell_i}|_{\ov{q}_{i;\ell_i}^{-1}(W_{i;\ell_i\mid j;\ell_j})}\ar[d]^{q'_{j;\ell_j\mid i;\ell_i}}\ar@/^-4pc/[dd]|-{q_{i;\ell_i}}&& \\
  \wt\cW_{i;\ell_i\mid j;\ell_j}\ar[rr]^{\wt\phi_{i;\ell_i\mid j;\ell_j}} \ar[d]^{q_{j;\ell_j\mid i;\ell_i}}
&& \wt\cW_{j;\ell_j}|_{\ov{q}_{j;\ell_j}^{-1}(W_{j;\ell_j\mid i;\ell_i})}\ar[d]^{q_{j;\ell_j}} \\
\cW_{i;\ell_i\mid j;\ell_j}\ar[rr]^{\phi_{ij}}		&&		 \cW_{j;\ell_j\mid i;\ell_i}\,.
}
\eEqu
\item The overlap maps $q_{i;\ell,\ell'}$, $q_{i,j;\ell}$, and $q_{i;\ell_i\mid j;\ell_j}$ are compatible on triple intersections similar\footnote{The precise statement involves some complicated commutative diagrams  similar to (\ref{equ:comm1}) and (\ref{equ:comm2}). } to (\ref{equ:third-condition}).
\item In \ref{l:ijfixell}-\ref{l:i<j-ell_i>ell_j}, the lifted maps commute with the sections $s_{i;\ell}$ on the source and  the target.
\item The twisted weight-ratio functions $(\ov{w}_{i;\ell}/\mf{I}_{\cW_{i;\ell}})$ are equal on the overlaps as in (\ref{equ:wIiwIj}).
\eEnum
\eDef

\noindent
We will not further go into the details of existence of such a resolution and its properties. The argument is similar to that of Section~\ref{sec:resolution}, but as the rather long statement of this definition indicates, it is more tedious to write down the proof in details. 

\noindent
Finally, similarly to Section~\ref{sec:euler}, in order to construct $\vfc(\mfL)$, we compatibly triangulate the zero set of orbifold sections of a given resolution in the following way. 
The key point is that at the level of zero sets, the double-indexed 
data of Definition~\ref{def:resolution-KUR} simplifies into a set of orbifolds simply indexed by $\ell\!\in\![N]$.
For every $\ell\!\in\![N]$, although the orbifolds $\wt\cW_{i,\ell}$ and $\wt\cW_{j,\ell}$, for different $i,j\!\in\!\mfB$, have different dimensions,
the zero sets $\wt{Z}(s_{i;\ell})$ and $\wt{Z}(s_{j;\ell})$ inherit, possibly non-effective, oriented orbifold structures $\wt{\cZ}(s_{i;\ell})$ and $\wt{\cZ}(s_{j;\ell})$ of the same dimension\footnote{This is by the transversality assumption, Definition~\ref{level-system}.\ref{l:dim}, and the fact that each orbibundle is relatively oriented.} $n$ from $\wt\cE_{i,\ell}\!\lra\!\wt\cW_{i,\ell}$ and $\wt\cW_{j,\ell}\!\lra\!\wt\cW_{j,\ell}$, respectively.
Therefore, via the restriction of orbifold embedding maps $\wt\phi_{ij;\ell}$ in Definition~\ref{def:resolution-KUR}.\ref{l:ijfixell} to $\wt{\cZ}(s_{i;\ell})$ and $\wt{\cZ}(s_{j;\ell})$, we can glue these orbifolds and obtain an oriented, and possibly non-effective, orbifold $\wt\cZ_\ell$ as in (\ref{equ:zero-orbifold}).

\noindent
For $\ell\!>\!\ell'$, we can glue the overlap maps $q_{i;\ell,\ell'}$ and $q_{i;\ell \mid j;\ell'}$, on the intersection of overlaps with zero sets of the orbifold sections, to obtain covering maps
$$
q_{\ell,\ell'}\colon \wt\cZ_{\ell,\ell'} \lra \wt\cZ_{\ell',\ell}
$$ 
similar to the covering maps used in Definition~\ref{def:triangulation}. Therefore, Definition~\ref{def:triangulation} still applies to the induced resolution of zero sets,
$$
\mf{R}_0\equiv \lrp{
\{\ov{q}_\ell\colon \wt{Z}_\ell\lra Z(\mfs)_\ell\cap W_\ell\}_{\ell\in [N]}, \{q_{\ell,\ell'}\colon \wt\cZ_{\ell,\ell'}\lra \wt\cZ_{\ell',\ell} \}_{1\leq \ell' < \ell \leq N}
},
$$
and gives us an ``admissible'' triangulation that via Proposition~\ref{pro:chain-cycle} gives us a singular homology class 
$$
\vfc(\mfL)\in H_n(Y(\mfL),\Q).
$$

\bRem{rem:double-admissible}
We could also define a double-indexed version of admissible triangulation for every resolution in Definition~\ref{def:resolution-KUR}, and generalize Definition~\ref{def:triangulation}, similarly.
However, by the trick of previous paragraph, we avoided such complicated definition by first gluing the zero sets $\wt\cZ(\mfs_{i;\ell})$, for every fixed $\ell$, into a single orbifold $\wt\cZ_\ell$ to which Definition~\ref{def:triangulation} and Proposition~\ref{pro:chain-cycle} apply.
\eRem

\noindent
We summarize the results of this section in the following  proposition. 
\bPro{pro:VFC-cL}
Let $\mfL$ be an oriented Hausdorff dimensionally graded system. Let $\mfL'$ be a sufficiently small shrinking of $\mfL$ admitting a  perfect EOB $\cN$, and $\mfs$ be a transversal multisection deformation of the Kuranishi map $s\!=\!\{s_i\}_{i\in\mfB}$ provided by Theorem~\ref{thm:pert-system-sections} with $Z(\mfs)\!\subset\!Y(\mfL')$. Then the homology class 
$$
\vfc(\mfL)\!\in\! H_n(Y(\mfL),\Q)
$$ 
constructed above is independent of the choice of $\mfL'$, $\cN$, $\mfs$, a resolution as in Definition~\ref{def:resolution-KUR},  and the choice of admissible triangulation.
\ePro

\bRem{rem:restriction-of-val}
The proof of independence is by considering the trivial deformation equivalence $\mfL\!\times\![0,1]$ and extending two given sets of auxiliary data over the boundaries to the entire space. There is however one delicate issue that we need to be careful about.
The $\tn{val}_\mfs$ function used in the construction of a resolution does not generally behave well under the restriction to boundary, i.e. if $\mfs$ is a multisection of an orbibundle $\cE$ over an orbifold $\cM$ with boundary $\partial\cM$ and $\ov{p}\!\in\! \partial M$, then $\tn{val}_\mfs(\ov{p})$ can be larger than $\tn{val}_{\mfs|_{\partial M}}(\ov{p})$; i.e. a set of  branches of $\mfs$ which are different in every neighborhood of $\ov{p}$ in $M$ may have equal restrictions over some neighborhood of $\ov{p}$ in $\partial M$. In this situation, the restriction of a resolution corresponding to a given tuple 
$(\cE\!\lra\!\cM,\mfs)$ to the boundary $(\cE|_{\partial M}\!\lra\!\partial \cM,\mfs|_{\partial M})$ is not a resolution of the latter. In order to avoid this issue, we may only consider multisections that over a sufficient small bordered chart $V\!\times\![0,a)$ around each point of $\partial M$ are trivial extensions of a multisection on $V$ to the normal direction.
\eRem

\subsection{VFC via evaluation maps}\label{sec:VFC2}

In this section, we consider Kuranishi structures for a pair consisting of a compact metrizable topological space and a map into another space.
We use this specialization to define VFC of a given pair as a singular homology class inside the target space. For example, we use this approach to define the GW virtual fundamental class in Theorem~\ref{thm:VFC}  for the pair of moduli space of $J$-holomorphic maps and evaluation and forgetful maps.

\bDef{def:rho-Kur-structure}
Let $M$ be a compact metrizable topological space, $W$ be a topological space, and 
$
\ov\vr\colon M\!\lra\! W
$
be a continuous map.
A Kuranishi structure $\cK$ on $M$ compatible with $\ov\vr$, called a \textbf{$\ov\vr$-Kuranishi structure}, is a Kuranishi structure on $M$ as in Definition~\ref{def:Kur-structure} such the the map $\ov\vr$ consistently lifts to all Kuranishi charts in the following way.
For every Kuranishi chart $\cU_{p}\!=\!(\pr_p\colon\!U_p\!\lra\!V_p,G_p, s_p,\psi_p)$, the map 
$$
\ov\vr_p=\ov\vr|_{F_p}\colon F_p=\psi_p(s_p^{-1}(0))\lra W 
$$
lifts to a continuous map $\ov\vr_p\colon (V_p/G_p)\lra W$ in a way that 
$$
\ov\vr_p\circ \ov\phi_{qp}=\ov\vr_q,
$$ 
for every coordinate change map $\Phi_{qp}\!=\!(\mfD\phi_{qp},\phi_{qp})\colon\cU_{q,p}\lra \cU_{p}$.
If $W$ is the underlying topological space of an orbifold $\cW$, we say $\cK$ is an \textbf{$\ov\vr$-smooth Kuranishi structure}, if the restrictions $\ov\vr_p$ consistently lift to orbifold smooth maps $\vr_p$ between $[V_p/ G_p]$ and $\cW$.
\eDef

\noindent
We define $\ov\vr$-DGS and $\ov\vr$-smooth DGS similarly. It is clear from the proof of Theorem~\ref{exist-level-system} that if $\cK$ is a 
$\ov\vr$-Kuranishi structure (resp. $\ov\vr$-smooth Kuranishi structure), then every associated natural DGS $\mfL$ is canonically an $\ov\vr$-DGS (resp. $\ov\vr$-smooth DGS). It is also clear from the definition of the thickening $Y(\mfL)$ in (\ref{equ:thikening-space}) that if $\mfL$ is an $\ov\vr$-DGS, then $\ov\vr$ extends to a continuous map 
\bEqu{equ:vrcLW}
\ov\vr\colon Y(\mfL)\lra W.
\eEqu

\noindent
\bPro{pro:VFC-map}
Let $M$ be a compact metrizable topological space, $W$ be a topological space, and 
$
\ov\vr\colon M\!\lra\! W
$
be a continuous map. Assume $(M,\cK)$ is an oriented $n$-dimensional $\ov\vr$-Kuranishi space and $\mfL$ is a natural oriented Hausdorff $\ov\vr$-DGS for $\cK$ as in Theorem~\ref{exist-level-system}.
Let $\vfc(\mfL)$ be the $n$-dimensional homology class of Proposition~\ref{pro:VFC-cL}.
Then the push forward\footnote{Via (\ref{equ:vrcLW}).} homology class 
$$
\ov\vr_*(\vfc(\mfL))\in H_n(W,\Q)
$$ 
is independent of the choice of a natural $\ov\vr$-DGS $\mfL$ for $\cK$ and the cobordism class of $\ov\vr$-Kuranishi structure  $\cK$; thus, we denote it by $\vfc(\cK,\ov\vr)$. 
\ePro

\noindent
The invariance part of Proposition~\ref{pro:VFC-map} is by extending two given sets of auxiliary data on the boundary of a cobordism between two given Kuranishi structures to the entire cobordism Kuranishi structure. We will not go into the details of such cobordism argument in this article as it is a relative version of what we have done so far.

\subsection{\v{C}ech homology VFC}\label{sec:VFC3}
The VFC of Proposition~\ref{pro:VFC-cL} is not a priori supported in the topological space $M$ and depends on the auxiliary choice of $\mfL$. Another way of resolving this issue is by taking the inverse limit over an infinite sequence of dimensionally graded systems such that $Y(\mfL)$ converges to $M$. This gives us a VFC supported in $M$, at the expense of replacing singular homology with \v{C}ech homology; this is how the VFC is defined in \cite{MW2,MW3,Pa}. We outline the process here.

\noindent
Let $(M,\cK)$ be an oriented $n$-dimensional Kuranishi space, $\mfL$ be an oriented Hausdorff DGS for $\cK$, $\cN$ be a perfect EOB for $\mfL$, and $\mfs_a\!\equiv\! \{\mfs_{i,a}\}_{i\in \mfB, a\in \N}$ be a compatible set of sequences of deformations of Kuranishi maps $s\!\equiv\! \{s_{i}\}_{i\in \mfB}$ of $\mfL$ as in the beginning of Section~\ref{sec:VFC1}.

\noindent
Let 
\bEqu{equ:shrink123infty}
 \mfL\!\equiv\!\mfL^\zero \!\supset\!\mfL^\one\! \supset\!\mfL^\two\!\supset\! \cdots 
\eEqu
be an infinite sequence of shrinkings of $\mfL$ such that
$$
\bigcap_{b=1}^\infty Y(\mfL^\tb)= M.
$$
By Lemma~\ref{lem:comp-support}, for every $b\in\N$, there exists $m_b\!\in\!\N$ such that 
$$
Z(\mfs^\tb_a)=  Z(\mfs^{\tn{\tiny{(b+1)}}}_a)\quad \forall a\!>\!m_b.
$$
Thus, by the construction of Section~\ref{sec:VFC1}, we obtain a sequence of homology classes 
$$
\vfc(\mfL^\tb) \in H_n(Y(\mfL^\tb),\Q)\quad \forall b\in \N,
$$
such that 
\bEqu{equ:btoc}
(\iota_{bc})_*(\vfc(\mfL^\tb))=\vfc(\mfL^\tc) \qquad \forall b\!>\!c,
\eEqu
where 
$$
\iota_{bc}\colon Y(\mfL^\tb)\lra Y(\mfL^\tc)
$$
is the inclusion map. The system of inclusion maps
$$
H_n(M,\Q)\lra H_n(Y(\mfL^\tb),\Q) \lra H_n(Y(\mfL^\tc),\Q)\quad \forall b\!>\!c,
$$
gives rise to an inclusion map 
\bEqu{equ:sing-limit}
H_n(M,\Q)\lra \lim_{\longleftarrow} H_n(Y(\mfL^\tb),\Q),
\eEqu
which is not necessary an isomorphism. 
Therefore, with this approach, we do not necessarily obtain a singular homology VFC in $H_n(M,\Q)$ by applying the inverse limit to (\ref{equ:btoc}).
For every $\tn{CW}$-complex $Y$, the \v{C}ech and singular homology groups $\check{H}_n(Y,\Q)$ and $H_n(Y,\Q)$, respectively, are isomorphic.
For every $b\!\in\!\N$, $\tn{cl}(Y(\mfL^\tb))$ has the structure\footnote{In \cite{MW2,MW3}, McDuff-Wehrheim use Poincare duality over the cycle realizing $\tn{VFC}(\cL)$ which has the structure of a CW-complex and then push forward the resulting \v{C}ech (co-)homology class.} of a CW-complex.  Therefore, the analogue of (\ref{equ:sing-limit}) with  \v{C}ech homology groups instead is an isomorphism. Thus, after replacing singular homology with \v{C}ech homology and taking inverse limit, we obtain a \v{C}ech homology virtual fundamental class
\bEqu{checkVFC_e}
\check{\vfc}\big((\mfL^{\tiny{(b)}})_{b=0}^\infty)\in \check{H}_n(M,\Q).
\eEqu
In conclusion, the following theorem is the \v{C}ech homology analogue of Proposition~\ref{pro:VFC-cL}; see \cite[Remark 8.2.4]{MW2}.
\bPro{pro:VFC-cech}
Let $(M,\cK)$ be an oriented $n$-dimensional Kuranishi space, $\mfL$ be a natural oriented Hausdorff DGS for $\cK$ as in Theorem~\ref{exist-level-system}, $(\mfs_a)_{a=1}^\infty$ be a compatible set of sequences of deformations of Kuranishi maps of $\mfL$ as in Theorem~\ref{thm:pert-system-sections}, and $(\mfL^{\tiny{(b)}})_{b\in \N}$ be a sequence of shrinkings of $\mfL$ as in (\ref{equ:shrink123infty}) .
Then the resulting homology class (\ref{checkVFC_e})
is an invariant of oriented cobordism class of $\cK$.
\ePro

\newpage
 \section{Moduli spaces of stable maps}\label{sec:moduli}
Stable compactification (or Gromov compactification) $\ov\cM_{g,k}(M,A)$ of $\cM_{g,k}(X,A)$ includes pseudoholomorphic maps with nodal domains.  
In order to cover the entire $\ov\cM_{g,k}(X,A)$ with Kuranishi charts, we need to first understand the local structure of the Deligne-Mumford space at an arbitrary nodal curve.

\noindent
In Section~\ref{sec:stable}, we set up some notation and review the definition of the moduli space of stable maps.
In Section~\ref{sec:DM}, following the analytic approach of \cite{RS}, we review the orbifold structure of $\ov\cM_{g,k}$. 
We recall the notions of universal families and of universal curve, and collect some well-known results about these objects. 
In Section~\ref{sec:topology}, we fix some more notation and state the Gromov's Convergence Theorem. 
In section~\ref{sec:main}, we will discuss a different realization of the Gromov topology via the notion of $\ep$-close maps.

\subsection{Stable curves and stable maps}\label{sec:stable}

Let $\Gamma$ be a graph, possibly with open edges, called \textbf{flags}, i.e.  edges that have a vertex at one end and are open at the other end. Let $V_\Gamma$, $E_\Gamma$, and $E^\circ_\Gamma$ be the set of vertices, edges, and flags of $\Gamma$, respectively.  A \textbf{labeling} of $\Gamma$ consist of two functions
$$
g\colon V_\Gamma\!\lra\! \Z^{\geq 0},\quad v\lra g_v,\quad \tn{and}\quad o\colon E^\circ_\Gamma\lra \{1,\ldots, |E^\circ_\Gamma|\}.
$$
For such a labeled graph $\Gamma$ let
\bEqu{equ:genus}
g_\Gamma= \sum_{v\in V_\Gamma} g_v\! +\! \tn{rank} ~H_1(\Gamma,\Z)
\eEqu
to be the (arithmetic) genus of $\Gamma$, where $H_1(\Gamma,\Z)$ is the first homology group of the underlying topological space of $\Gamma$. Figure~\ref{fig:labled-graph}-left illustrates a labeled graph with $2$ flags.

\begin{figure}[h]
\begin{pspicture}(8,1.3)(11,3.8)
\psset{unit=.3cm}

\pscircle*(35,6){.25}\pscircle*(39,6){.25}
\pscircle*(35,9){.25}\pscircle*(41.5,9){.25}\pscircle*(39,9){.25}
\psline[linewidth=.05](39,6)(35,9)
\psline[linewidth=.05](39,6)(39,9)\psline[linewidth=.05](39,6)(41.5,9)
\psarc[linewidth=.05](33,7.5){2.5}{-36.9}{36.9}\psarc[linewidth=.05](37,7.5){2.5}{143.1}{216.0}
\psline[linewidth=.05](41.5,9)(40.5,10.5)\psline[linewidth=.05](41.5,9)(42.5,10.5)
\rput(40.5,11.2){1}\rput(42.5,11.2){2}
\rput(35,5){$g_4$}\rput(39.5,5){$g_5$}
\rput(34.5,10){$g_1$}\rput(43.7,8.5){$g_3$}
\rput(39,10){$g_2$}

\pscircle*(55,6){.25}\pscircle*(59,6){.25}
\pscircle*(55,9){.25}\pscircle*(61.5,9){.25}\pscircle*(59,9){.25}
\psline[linewidth=.05](59,6)(55,9)
\psline[linewidth=.05](59,6)(59,9)\psline[linewidth=.05](59,6)(61.5,9)
\psarc[linewidth=.05](53,7.5){2.5}{-36.9}{36.9}\psarc[linewidth=.05](57,7.5){2.5}{143.1}{216.0}
\psline[linewidth=.05](61.5,9)(60.5,10.5)\psline[linewidth=.05](61.5,9)(62.5,10.5)
\rput(60.5,11.2){1}\rput(62.5,11.2){2}
\rput(55,5){$(g_4,A_4)$}\rput(59.5,5){$(g_5,A_5)$}
\rput(54,10){$(g_1,A_1)$}\rput(63.8,8.5){$(g_3,A_3)$}
\rput(58.3,10){$(g_2,A_2)$}
\end{pspicture}
\caption{On left, a labeled graph $\Ga$ representing elements of $\ov\cM_{g,2}$. On right, a labeled graph $\Ga$ representing elements of $\ov\cM_{g,2}(X,A)$.}
\label{fig:labled-graph}
\end{figure}
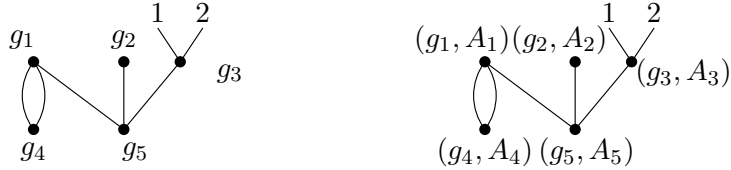

\noindent 
Such labeled graphs $\Gamma$ characterize different topological types of nodal marked surfaces $(\Si,\vec{z}\!=\!(z^1,\ldots,z^k))$ in the following way.
Each vertex $v\!\in\!V_\Gamma$ corresponds to a smooth\footnote{We mean a smooth closed oriented surface.} component $\Si_v$ of $\Si$ with genus $g_v$. Each edge $e\!\in\!E_\Gamma$ corresponds to a node $q_e$ obtained by connecting $\Sigma_v$ and $\Sigma_{v'}$ at the points $q_{e,v}$ and $q_{e,v'}$, where $e\!=\!\ll v,v'\rr$; note that $e$ can be a loop connecting $v$ to itself (i.e $v'\!=\!v$). Finally, each flag $e^\circ$ connected to the vertex $v(e^\circ)$ corresponds to a marked point $z^{o(e^\circ)}\!\in \Sigma_{v(e^\circ)}$ disjoint from the connecting nodes. Then $g\!=\!g_\Gamma$ is the arithmetic genus of $\Si$. Thus we get
\bEqu{equ:nodalcurve}
(\Si,\vec{z})\! =\! \coprod_{v\in V_\Gamma}(\Si_v,\vec{z}_v,{q}_v)/ \{ q_{e,v_1}\!\sim\! q_{e,v_2},~~\forall e\!=\!\ll v_1,v_2\rr \in E_\Gamma \},
\eEqu
where 
$$
\vec{z}_v\!=\!\vec{z}\cap \Sigma_v\quad\tn{and}\quad \quad {q}_v=\{q_{e,v}\colon \forall e\!\in\!E_\Gamma, ~v\!\in\! e\},
$$
counting loops twice.
In this situation we say $\Gamma$ is the \textbf{dual graph} of $(\Si,\vec{z})$.
A complex structure $\mfj$ on $\Sigma$ is a set of complex structures $(\mfj_v)_{v\in V_\Gamma}$ on its components. By a (complex) marked nodal curve, we mean a  marked nodal real surface together with a complex structure $(\Si,\mfj,\vec{z})$. 
Figure~\ref{fig:labled-curve} illustrates a nodal curve with $(g_1,g_2,g_3,g_4,g_5)=(0,2,0,1,0)$ corresponding to Figure~\ref{fig:labled-graph}-left. Note that the nodal points ${q}_v$ are not a priori ordered; if we order them someway, we denote the ordered set by $\vec{q}_v$.

\vskip.1in
\noindent 
Similarly, for nodal marked surfaces with maps into $X$, we consider similar labeled graphs where the vertices carry an additional labeling $A_v\!\in\! H_2(X,\Z)$, recording the homology class of the image of the corresponding component. Figure~\ref{fig:labled-graph}-right illustrates a dual graph associated to a marked nodal map over the graph on the left. 
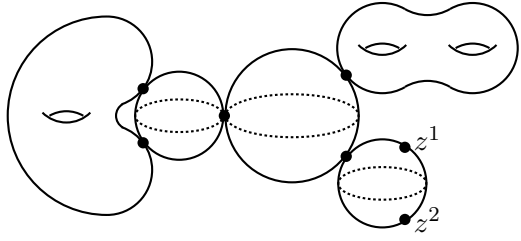
\begin{figure}

\begin{pspicture}(8,1.7)(11,4.5)
\psset{unit=.3cm}

\pscircle*(52.4,8.2){.25}\pscircle*(52.4,11.8){.25}
\pscircle(50,10){3}\psellipse[linestyle=dashed,dash=1pt](50,10)(3,1)
\pscircle(54,7){2}\psellipse[linestyle=dashed,dash=1pt](54,7)(2,.77)
\pscircle*(55,8.6){.25}\rput(56,9){$z^1$}
\pscircle*(55,5.4){.25}\rput(56,5.4){$z^2$}

\pscircle(45,10){2}\psellipse[linestyle=dashed,dash=1pt](45,10)(2,.77)
\pscircle*(47,10){.25}\pscircle*(43.4,8.8){.25}\pscircle*(43.4,11.2){.25}
\psarc(41.8,7.6){2}{270}{71}
\psarc(41.8,12.4){2}{-71}{90}
\psarc(41.8,10){4.4}{90}{270}\psarc(42.8,10){.6}{120}{240}
\psarc(40,8.7){1.5}{60}{120}
\psarc(40,11.2){1.5}{225}{315}

\psarc(54,13){2}{60}{300}
\psarc(58,13){2}{240}{120}
\psarc(56,16.464){2}{240}{300}
\psarc(56,9.546){2}{60}{120}
\psarc(54,11.7){1.5}{60}{120}
\psarc(54,14.2){1.5}{225}{315}
\psarc(58,11.7){1.5}{60}{120}
\psarc(58,14.2){1.5}{225}{315}

\end{pspicture}
\caption{A nodal curve in $\ov\cM_{4,2}$.}
\label{fig:labled-curve}
\end{figure}
For every $v\!\in\!V_\Gamma$, let $\Gamma_v$ be the labled graph with $V_{\Gamma_v}\!=\!\{v\}$, $E_{\Gamma_v}\!=\!\emptyset$, and 
\bEqu{equ:v-flags}
E^\circ_{\Gamma_v}= \{ e^\circ\!\in\! E^\circ_\Gamma\colon v(e^\circ)\!=\!v\} \cup \{ e\!\in\! E_\Gamma\colon v\!\in\! e\},
\eEqu
counting loops twice. The order function on $E^\circ_{\Gamma_v}$ depends on a choice of ordering on $q_v$.  Each sub-graph $\Gamma_v$ describes the topological type of the component corresponding to $v$ in the normalization of the corresponding nodal marked  surface.

\noindent
With notation as above, for every nodal marked curve $C\!\equiv\! \big(\Si,\mfj, \vec{z}\big)$ with dual graph $\Gamma$, let 
\bEqu{equ:auto-nodal}
\aut(C)\!=\!
\lrc{h\!\equiv\! \{h_v\colon \Si_v \!\to\!\Si_{h(v)}\}_{\{v\in V_\Gamma\}}\colon h(z^i)\!=\!z^i~~\forall z^i\!\in\!\vec{z},~ h^*\mfj_{h(v)}\!=\!\mfj_v~~\forall v\!\in\!V_\Gamma}
\eEqu
be the group of biholomorphic automorphisms of $C$.
A marked nodal curve $C$ is called \textbf{stable} if $\aut(C)$ is finite; this is the case if and only if all $(\Si_v,\mfj_v,\vec{z}_v\cup\vec{q}_v)$ are stable. 

\bRem{rem:permuting-auto}
For every $v\!\in\!V_\Gamma$ and $h\!\in\!\aut(C)$ such that $h(v)\!=\!v$, the automorphism $h$ preserves ${q}_v$ (as a set) but it may permute the points non-trivially; thus, if we fix an ordering $\vec{q}_v$, for all $v\!\in\!V_\Gamma$, $\aut(C)$ contains the product 
$$
\prod_{v\in V_\Gamma} \aut(\Si_v,\mfj_v,\vec{z}_v\cup\vec{q}_v)
$$ 
but it may have more elements permuting different components or some of $q_v$. 
\eRem

\noindent 
A nodal map $u$ from $\big(\Si,\mfj, \vec{z}\big)$ into $X$ is a collection of (sufficiently smooth) maps 
$$
u_v \colon \Si_v\lra X~~\forall v\!\in\!V_\Gamma , \quad u_{v_1}(q_{e,v_1})\!=\!u_{v_2}(q_{e,v_2})\quad \forall e\!=\!\ll v_1,v_2\rr \!\in\! E_\Gamma.
$$
If $(X,J)$ is an almost complex manifold, then a nodal map is $J$-holomorphic if every $u_v$ is $(J,\mfj_v)$-holomorphic as in (\ref{equ:J-holo}).
Similarly, for a nodal $(J,\mfj)$-holomorphic map $f\!=\!\big(u,\Si,\mfj, \vec{z}\big)$, let
\bEqu{equ:auto-nodal-map}
\aut(f)\!=\!\{h\!\in\! \aut(\Si,\mfj, \vec{z})\colon~u\circ h\!=\!u\};
\eEqu
such $f$ is called stable if $\aut(f)$ is finite. We can think of a stable curve, as a stable map with the trivial map to a point.

\noindent
Let $(X,\om)$ be a closed symplectic manifold. Given $g,k \!\in\! \Z^{\geq 0}$, $J\!\in\! \cJ(X,\om)$, and $A\!\in\! H_2(X,\Z)$, the moduli space $\ov\cM_{g,k}(X,A,J)$ is the set of  equivalence classes of stable $J$-holomorphic maps over nodal $k$-marked genus $g$ (not necessarily stable) curves. 
If $X$ is a point $\tn{pt}$, 
$$\ov\cM_{g,k}\equiv \ov\cM_{g,k}(\tn{pt},0,\tn{null})$$
is the Deligne-Mumford space of stable $k$-marked genus $g$  curves.  For simplicity, we often fix a choice of $J$ and drop it from the notation; i.e. we denote the moduli space by $\ov\cM_{g,k}(X,A)$. 

\vskip.1in
\noindent 
If $2g+k\!<\!3$ (resp. $A\!=\!0$ and $2g\!+\!k\!<\!3$), there is no stable curve (resp. map) of type $(g,k)$. 
In these cases, for our notion to be consistent, we define $\ov\cM_{g,k}$ (resp. $\ov\cM_{g,k}(X,0)$) to be a single point (resp. $X$). 
The moduli space $\ov\cM_{g,k}(X,A)$ can be decomposed into pieces $\cM_{g,k}(X,A)_\Gamma$, called \textbf{strata}, depending on the dual graph of the map.
For $J\!\in\! \cJ(X,\om)$, the symplectic area (energy) of a $J$-holomorphic map $u\colon (\Si,\mfj) \lra X$ is given by 
\bEqu{equ:area}
\om(u)\!=\!\sum_{v\in \Gamma_V}\int_{\Si_v} u_v^*\om = \frac{1}{2}  \sum_{v\in \Gamma_V} \int_{\Si_v} |du_v|_J^2 ~\nd \tn{vol}_{\Si_v}.
\eEqu
Since the symplectic area of a non-trivial $J$-holomorphic maps has a positive lower bound, see \cite[Proposition 4.1.4]{MS2004},  it follows from the Gromov's Compactness Theorem in Section~\ref{sec:topology} that the stratification 
\bEqu{equ:Gamma-M}
\cM_{g,k}(X,A)\!=\!\coprod_{\Gamma} \cM_{g,k}(X,A)_{\Gamma}
\eEqu
is finite. The forgetful map (\ref{equ:st}) takes each stratum $\cM_{g,k}(X,A)_\Gamma$ to the stratum $\cM_{\st(\Gamma)}\!\subset\! \ov\cM_{g,k}$, where $\st(\Gamma)$ is the dual graph uniquely obtained via the stabilization process of shrinking unstable components and may have fewer vertices and edges but equal number of flags. Note that each stratum $\cM_\Gamma$ (and similarly $\cM_{g,k}(X,A)_\Gamma$) is of the form
$$
\big(\prod_{v\in V_\Gamma} \cM_{g_v,k_v}\big)/\aut(\Gamma)
$$
where $k_v\!=\!|E_{\Gamma_v}^\circ|$ (see (\ref{equ:v-flags})) and $\aut(\Gamma)$ is the symmetry group of $\Gamma$.
\subsection{Orbifold structure of the Deligne-Mumford space}\label{sec:DM}
In this section, following the construction of Robbin-Salamon \cite{RS}, we review the complex (K\"ahler) orbifold structure of the Deligne-Mumford space. Except the two cases of $(g,k)\!=\!(1,1),(2,0)$,  $\ov\cM_{g,k}$ has the structure of a complex effective orbifold. In those two special cases, the orbifold structure is not effective but it has an effective reduction; see Example~\ref{exa:M11}.
From the algebraic geometry point of view, $\ov\cM_{g,k}$ has the structure of a complex projective variety representing families of stable marked nodal curves modulo isomorphism; see \cite{HM}. 

\noindent
A \textbf{holomorphic family} of complex curves consists of connected open complex manifolds $\cB$ and $\cC$, with $\dim \cC\!=\!\dim_\C \cB\!+\!1$, and a proper  holomorphic projection map $\pi\colon \cC\! \lra\! \cB$ with  complex one dimensional fibers. By the holomorphic Implicit Function Theorem, for every regular point $p\!\in\! \cC$, there are local holomorphic coordinates $w\!=\!(w_0,\ldots,w_{n})$ around $p$, with $n\!=\!\dim \cB$, such that $\pi(w)\!=\!(w_1,\ldots,w_n)$.  A critical point $p$ of $\pi$ is called a \textbf{nodal point} if there are local holomorphic coordinates $w\!=\!(w_0,\ldots,w_{n})$ around $p$ such that $\pi(w)\!=\!(w_0w_1,w_2,\ldots,w_n)$. A holomorphic family $\cC\!\lra\!\cB$ is called \textbf{nodal} if all the critical points are nodal. 

\noindent
Given $k\!\in\! \Z^{\geq 0}$, a \textbf{$k$-marked nodal family} of complex curves consists of a nodal family $\pi\colon\cC\lra \cB$ and an ordered set of holomorphic sections 
$$
\vec{\mfz}=\{\mfz^i\colon \mc{B} \lra \cC\quad \forall i\!\in\! [k]\}
$$ 
away from the nodal points, such that the graph of every two $\mfz^i$ are mutually disjoint. 
For every $a\!\in\! \cB$, let $\Si_a\!=\!\pi^{-1}(a)$, $\mfj_a$ be the restriction of complex structure of $\cC$ to $\Si_a$, and 
$$
\vec{z}_a=(\mfz^1(a),\ldots,\mfz^k(a))\!\subset\! \Si_a;
$$ 
then, $C_a\!=\!(\Si_a,\mfj_a,\vec{z}_a)$ is a nodal curve of some non-negative genus $g$ with $k$ marked points. The arithmetic genus of fibers of a nodal family is constant, c.f. \cite[Lemma 4.6]{RS}; therefore, the genus of a nodal family is well-defined. A marked family is called \textbf{stable} if each of its fibers is stable. The stability is an open condition. 

\noindent
Given two $k$-marked nodal families 
$$
\cV_i=\big(\pi_i\colon\cC_i\lra \cB_i, \quad \vec{\mfz}_i\colon \cB_i\lra \cC_i\big),\quad  \tn{with } i\!=\!1,2,
$$ 
a \textbf{morphism} between $\cV_1$ and $\cV_2$ is a commutative diagram of holomorphic maps
$$
\xymatrix{
& \cC_1  \ar[rr]^{\varphi^\cC} \ar[d]_{\pi_1} && \cC_2 \ar[d]_{\pi_2}  \\
&  \cB_1 \ar@/_/[u]_{(\mfz_1^i)_{i\in [k]}} \ar[rr]^{ \varphi^\cB}                 && \cB_2 \ar@/_/[u]_{(\mfz_2^i)_{i\in [k]}}}
$$
such that the restriction of $\varphi^\cC$ to each fiber of $\cC_1$ is a holomorphic isomorphism. 

\bDef{def:versal-family}
Given a marked nodal curve $C\!=\!(\Si,\mfj,\vec{z})$, an \textbf{unfolding} of $C$ consists of a marked nodal family 
\bEqu{equ:unfolding}
\lrp{ \pi\colon\cC\lra \cB, \vec{\mfz}, b\!\in\! \cB}
\eEqu 
and an identification of $C$ with the base curve $C_b\!=\!(\Si_b,\mfj_b,\vec{z}_b)$. An unfolding is called a \textbf{universal family} around $C$ if for every other unfolding of $C$,
$$
\lrp{ \pi'\colon\cC'\lra \cB', \vec{\mfz}', b'\!\in\! \cB'},
$$ 
and any holomorphic identification 
\bEqu{equ:base-varphi}
\varphi\colon C'_{b'}\!\lra\! C_b,
\eEqu
there is some open neighborhood  $b'\!\in\!\cB''\!\subset\! \cB'$ and a unique base point preserving morphism (a germ of a morphism)
\bEqu{equ:universal-morphism}
(\varphi^\cC,\varphi^\cB)\colon\lrp{ \pi'\colon\cC'|_{\cB''} \lra \cB'', \vec{\mfz}'|_{\cB''}, b\!\in\! \cB''} \lra \lrp{ \pi\colon\cC\lra \cB, \vec{\mfz}, b\!\in\! \cB}
\eEqu
with $\varphi^\cC|_{C'_{b'}}\!=\!\varphi$. 
\eDef
\noindent
By \cite[Theorems 5.5,5.6]{RS}, every stable marked curve admits a unique (up to isomorphism) germ of universal family. The base $\cB$ of such universal families cover $\ov\cM_{g,k}$ and form an (possibly non-effective) orbifold atlas of complex dimension $3g-3+k$ for $\ov\cM_{g,k}$. In fact, for every $h\!\in\!\aut(C_b)$, with $\cC'\!=\!\cC$, $b'\!=\!b$, and $\varphi\!=\!h$ in (\ref{equ:base-varphi}), 
let\footnote{By abuse of notation, and after possibly replacing $\cB$ with a smaller neighborhood of $b$, we write $\cB$ instead of $\cB''$. } 
\bEqu{equ:induced-h}
(\varphi_{h}^\cC,\varphi_{h}^\cB)\colon (\cC,\cB)\lra (\cC,\cB) 
\eEqu 
be the resulting (germ of) isomorphism of (\ref{equ:universal-morphism}). This gives us a holomorphic action of $\aut(C_b)$ on (the germ of) 
$\cB$ (resp. $\cC$) such that 
\bEqu{equ:DM-chart}
[\cB/ \aut(C_b)]
\eEqu
is an orbifold chart centered at $[C_b]\!\in \!\ov\cM_{g,k}$. By \cite[Theorems 5.3,5.4]{RS}, if (\ref{equ:unfolding}) is a universal family around $C_b$, then for every $b'\!\in\! \cB$ sufficiently close to $b$, (\ref{equ:unfolding})  is also a universal family around $C_{b'}$. Moreover, the equality
\bEqu{bbprime-aut_e}
\aut(C_{b'})=\{\varphi_h^\cC|_{C_{b'}}\colon \varphi_h^\cB(b')=b',~~h\in \aut(C_b)\} 
\eEqu
identifies $\aut(C_{b'})$ with a subgroup of $\aut(C_b)$; the size of automorphism group is an upper semi-continuous function.

\noindent
If $C$ is smooth, there exists a natural identification
\bEqu{equ:Dm-tangent}
T_b\cB\!=\!H^1(\cT\Si_{\mfj}(-\vec{z}))\cong \tn{Ext}^1(\Om_{\Sigma,\mfj}^{1,0}(\vec{z}),\cO_{\Sigma,\mfj})
\eEqu
where $\cT\Si_{\mfj}(-\vec{z})$ is the logarithmic\footnote{$\cT\Si_{\mfj}(-\vec{z})$ is isomorphic to the tensor product of $\cT\Si_{\mfj}$ and the line bundle $\cO_{\Sigma,\mfj}(-\vec{z})$. It corresponds to the sheaf of holomorphic tangent vector fields vanishing at $\vec{z}$.} holomorphic tangent bundle of $(\Si,\mfj,\vec{z})$, $\Om^{1,0}_{\Sigma,\mfj}(\vec{z})$ is the sheaf of logarithmic holomorphic $1$-forms with at most simple poles at $\vec{z}$, $\cO_{\Sigma,\mfj}$ is the structure sheaf of $(\Si,\mfj)$, and  $\tn{Ext}^*$ are hyper-cohomology functors. We will avoid using these facts as much as possible in our calculations. The second identification has an analogue for nodal curves as well. The conclusion is that $\ov\cM_{g,k}$ is an orbifold of complex dimension 
$$
3(g\!-\!1)\!+\!k=h^1(\cT\Si_{\mfj}(-\vec{z})).
$$

\noindent
We refer to Definition 6.2 and 6.4, Proposition 6.3, and Theorem 6.5 and 6.6 in \cite{RS} for further details on the (Lie groupoid) orbifold structure of $\ov\cM_{g,k}$ obtained via the base of universal families above.

\noindent
Considering the total space of universal families $\cC$ above, instead of just the bases $\cB$, and identifying them along their intersections, we obtain a so called \textbf{universal curve} $\ov\cC_{g,k}$ over $\ov\cM_{g,k}$ which is an effective complex orbifold of one dimension higher. The universal curve admits a projection map $\pi\colon\ov\cC_{g,k}\lra \ov\cM_{g,k}$ such that for every marked curve $C_a=(\Si_a,\mfj_a,\vec{z}_a)$ representing $[C_a]\!\in\! \ov\cM_{g,k}$, 
$$
\pi^{-1}([C_a])\cong [\Si_a/ \aut(C_a)].
$$ 
In fact, $\ov\cC_{g,k}\cong \ov\cM_{g,k+1}$ with $\pi\!=\!\pi_{k+1}$ as in (\ref{equ:pi-I}).

\bRem{rem:smoothly-trivial-family}
Let $C\!=\!(\Si,\mfj,\vec{z})$ be a smooth stable marked curve; i.e. $[C]\!\in\!\cM_{g,k}$.
In this case, there exists a sufficiently small universal family $\lrp{\pi\colon \cC\lra \cB, \vec{\mfz},b}$ around $C$ such that $\cC$ is smoothly trivial, i.e. there exists an $\aut(C)$-equivariant\footnote{The first factor of the map $\varphi$ could be defined to send every $x\!\in\!\cC$ to the closest point in  $\Si$ with respect to some properly defined metric $\mc{G}$ on $\cC$. If we choose $\mc{G}$ to be $\aut(C)$-invariant, then $\varphi$ becomes equivariant.} diffeomorphism
$$
\varphi\colon  \cC\lra\Si\times \cB
$$ 
such that $\pi\circ \varphi^{-1}$ is the projection to the second component, 
$$
(\varphi\circ \mfz^i)(x)\equiv (z^i\!\times\! x,x)\quad \forall x\in \cB,~i\!\in\![k],
$$
(i.e. every section $\varphi_*( \mfz^i)$ is constant), $\varphi_b\!=\!\tn{id}_{\Si}$, and the action of $\aut(C)$ on $\Si\!\times\!\cB$ is the product action.
Here, for every $a\!\in\!\cB$, $\varphi_a$ is the restriction of $\varphi$ to $\Si_a\!=\!\pi^{-1}(a)\!\subset\!\cC$.
In this sense, we can think of $\cB$ as a smooth family of complex structures $\{\mfj_{\varphi,a}=(\varphi_a)_*(\mfj_{a})\}_{a\in \cB}$ on the fixed marked surface $(\Si,\vec{z})$ with $\mfj_{\varphi,b}\!=\!\mfj_b\!=\!\mfj$. Note that the $\aut(C)$-equivariance implies that for every $a\!\in\!\cB$ and $h\!\in\!\aut(C)$, with $a'\!=\!\varphi_h^\cB(a)$ as in (\ref{equ:induced-h}), the map 
\bEqu{hjphia_e}
h\colon (\Si, \mfj_{\varphi,a})\!\lra\! (\Si, \mfj_{\varphi,a'})
\eEqu 
is holomorphic.
\eRem

\noindent
\bRem{rem:family-with-unmarked}
In what follows, we also need to consider marked curves $C\!=\!(\Si,\mfj,\vec{z},q)$ where $q\!=\!\{q_1,\ldots,q_\ell\}$ is an un-ordered set of disjoint points away from $\vec{z}$; e.g. they appear  in decomposition of nodal curves, in which case $q$ is the pre-image of nodal points in a smooth component. Let
\bEqu{equ:un-order-group}
\aut(\Si,\mfj,\vec{z},q)\!=\!
\big\{h\!\in\! \aut(\Si,\mfj,\vec{z}): h(\{q_1,\ldots,q_\ell\})\!=\!\{q_1,\ldots,q_\ell\}\big\}.
\eEqu
After fixing an ordering $\vec{q}$ of $q$; the automorphism group $\aut(\Si,\mfj,\vec{z}\cup\vec{q})$ may become a proper subgroup\footnote{This subgroup is independent of choice of ordering.} of (\ref{equ:un-order-group}); see Remark~\ref{rem:permuting-auto}. Assume $C'\!=\!(\Si,\mfj,\vec{z}\cup\vec{q})$ is stable and let 
$$
\lrp{ \pi\colon\cC\lra \cB, \vec{\mfz}\cup\vec{\mfq}, b\!\in\! \cB}
$$ 
be a universal family around $C'$. After forgetting the ordering on $\vec{\mfq}$, the action of $\aut(C')$ on $\cC$ discussed prior to (\ref{equ:DM-chart}) extends to an action of (\ref{equ:un-order-group}) on 
$$
\lrp{ \pi\colon\cC\lra \cB,~ \vec{\mfz}, \mfq\colon \cB\lra \cC, ~b\!\in\! \cB}.
$$
Therefore, by a Universal family around $C$ we mean a universal family around $C'$, in which we forget the added data of ordering on $\vec{\mfq}$ and extend the group action to (\ref{equ:un-order-group}). 
\eRem

\noindent
Let $C\!=\!(\Si,\mfj,\vec{z})$ be a stable nodal curve with dual graph $\Gamma$. For every $v\!\in \!V_\Gamma$, let 
\bEqu{equ:v-families}
\lrp{\pi_v\colon \cC_v\!\to\! \cB_v, \vec\mfz_v,\mfq_v,b_v}
\eEqu
be a sufficiently small universal family around $C_v\!\equiv\!(\Si_v,\mfj_v,\vec{z}_v, q_v)$ as above with no nodal point, see Remark~\ref{rem:smoothly-trivial-family}.
Then there exists a \textbf{standard} extension of (\ref{equ:v-families}) to a universal family $\lrp{\pi\colon \cC\lra \cB, \vec{\mfz}, b}$ around $C\!=\!C_b$ with 
\bEqu{equ:B-smmothing}
\cB\!=\!\prod_{v\in V_\Gamma}\cB_v \times \De^{E_\Gamma}\quad \tn{and}\quad b=(b_v)_{v\in V_\Gamma}\times 0^{E_\Gamma}\in \cB,
\eEqu
where $\De$ is a disk of sufficiently small radius $\de$ in $\C$ and $\De^{E_\Gamma}\!=\!\prod_{e\in E_\Gamma}\De$.
More precisely, $\cC$ is obtained from $\{\cC_v\}_{v\in V_\Gamma}$ in the following way.

\noindent
For every pair\footnote{If $e$ is a loop  at $v$, then by definition there are two pairs $(e',v)$ and $(e'',v)$ corresponding to $e$ and $v$.} $(e,v)$ of an edge and a vertex on one end of it, the section
\bEqu{equ:divisor}
\mc{D}_{e,v}=\mfq_v^e(\cB_v)\subset \cC_v
\eEqu
is a (Cartier) divisor. In (\ref{equ:divisor}), via (\ref{equ:v-flags}), we think of $e$ as a flag in $E^\circ_{\Gamma_v}$ and $\mfq_v^e$ is the corresponding section of $\cC_v$. Therefore, for $\{\cB_v\}_{v\in V_\Gamma}$ sufficiently small, there exist $\ve\!>\!0$, a set of mutually disjoint open neighborhoods 
\bEqu{equ:Uev}
\mc{D}_{e,v}\!\subset\!U_{e,v}\!\subset\! \cC_v\quad \forall e\!\in\!E_\Gamma,~v\!\in\!e,\qquad U_{e,v}\cap \mfq_v^{e^\circ}(\cB_v)\!=\!\emptyset \quad \forall e^\circ\!\in\!E^\circ_\Gamma, v\!=\!v(e^\circ),
\eEqu
and holomorphic functions 
\bEqu{equ:holo-functions}
w_{e,v}\colon U_{e,v}\lra \C
\eEqu
such that $\mc{D}_{e,v}\!=\!w_{e,v}^{-1}(0)$ and
$$
\tn{cl}(U^\ve_{e,v})\!\subset\!U_{e,v},\quad \tn{where}\quad U_{e,v}^\ve=\{y\!\in\!U_{e,v}\colon |w_{e,v}(y)|\!<\!\ve\}.
$$

\bLem{lem:symmetric-gluing}
There exists a choice of universal families $\{\cC_v\}$ and holomorphic functions $\{w_{e,v}\}$ which are compatible with the action of $\aut(C)$ in the following sense.
\bEnum
\item\label{l:v} For $v\!\in\!V_\Gamma$ and $h\!\in\!\aut(C)$, we have $\cC_v\!\cong\!\cC_{h(v)}$.
\item\label{l:h} For every $h\!\in\!\aut(C)$ 
$$
U_{h(e),h(v)}\!=\! \varphi^\cC_h(U_{e,v}),
$$
where $\varphi^\cC_h\colon \cC_v \lra \cC_{h(v)}\cong \cC_v$ is the isomorphism corresponding to $h$ in (\ref{equ:induced-h}).
\item\label{l:symmetric-w} For every $e\!=\!\ll v_1,v_2\rr\!\in E_\Gamma$ (or similarly for the pairs $(e',v)$ and $(e'',v)$ that correspond to a loop $e$) and every $h\!\in\!\aut(C)$ of order $m_{h}$, there exists an $m_{h}$-th root of unity $\mu_{h,e}\in \C^*$  such that
\bEqu{equ:whw}
(w_{h(e),h(v_1)}\circ \varphi^\cC_h) \cdot (w_{h(e),h(v_1)}\circ \varphi^\cC_h)\!=\! \mu_{h,e} w_{e,v_1}\cdot w_{e,v_2}.
\eEqu
\eEnum
\eLem

\bProof
Only the third condition is somewhat non-trivial. Assuming Conditions~\ref{l:v} and \ref{l:h} hold, we construct $\{w_{e,v}\}$ in a way that Condition~\ref{l:symmetric-w} holds as well. Fix some $e\!\in\!E_\Gamma$ and let $\aut(C)\!\cdot\! e$ be the orbit space of the action of $\aut(C)$ on $e$. For every $e'\!=\!\ll v'_1,v'_2\rr\!\in\!\aut(C)\!\cdot\! e$, fix an $h'\!\in\!\aut(C)$ such that $e'\!=\!h'(e)$; given $w_{e,v_1}$ and $w_{e,v_2}$ define  $w_{e',v'_1}$ and $w_{e',v'_2}$ via (\ref{equ:whw}) with $h'$ in place  of $h$.  
If $h''$ is another automorphism with the same property, there exists $h$ such that $h(e)\!=\!e$ and $h''\!=\!h'\!\circ h$. 
Thus, we reduce to the subgroup $\aut(C)_e\!\subset\!\aut(C)$ of the authomorphisms fixing $e$. 
If $e$ is a loop at $v$ and there exists $h\!\in\!\aut(C)_e$ exchanging the two pairs $(e',v)$ and $(e'',v)$ corresponding to that, given $w_{e',v}$ and a fixed choice of such $h$, we define $w_{e'',v}\!=\!w_{e',v}\!\circ\!\varphi^\cC_h$. Thus, we reduce to the subgroup $\aut(C)_{(e,v)}\!\subset\!\aut(C)_e$ of those automorphisms fixing a pair $(e,v)$. For every $h\!\in\! \aut(C)_{(e,v)}$, $\varphi^\cC_h(\mc{D}_{e,v})\!=\!\mc{D}_{e,v}$  and $U_{e,v}$ is $\varphi^\cC_h$-invariant. Choose an arbitrary defining equation $w'_{e,v}$ for $\mc{D}_{e,v}$ and define
$$
\wt{w}_{e,v}=\prod_{h\in\aut(C)_{(e,v)}} w'_{e,v}\circ \varphi^\cC_h.
$$
If $m\!=\!|\aut(C)_{(e,v)}|$, any branch $w_{e,v}\!=\!\sqrt[\leftroot{-2}\uproot{2}m]{\wt{w}_{e,v}}$ has the desired property.\eProof

\noindent
Going back to the construction of $\cC$,
for every $ e\!=\!\ll v_1,v_2\rr\!\in\! E_\Gamma$, let
$$
\aligned
&\cC_{e}\!= \bigg\{
\big(y_{1},y_{2},\la_{e}\big)\!\in\! 
 U_{e,v_1}\!\times \!  U_{e,v_2}\!\times \! \De^e \colon w_{e,v_1}(y_{1})\cdot w_{e,v_2}(y_{2})\!=\!\la_{e} 
\bigg\}\!\times\!\prod_{\substack{v\in V_\Gamma\\ v\not\in e} }\cB_v\! \times\! \De^{E_\Gamma\setminus\{e\}},\\
&\cC_{e,v_1}\!= \cC_e \cap \bigg( \big( U_{e,v_1}\setminus \tn{cl}(U^\ve_{e,v_1})\big)\!\times \!  U_{e,v_2}\!\times \! \De^e\!\times\!\prod_{\substack{v\in V_\Gamma\\ v\not\in e} }\cB_v\! \times\! \De^{E_\Gamma\setminus\{e\}}\bigg).
\endaligned
$$
For every $v\!\in\! V_\Gamma$, let
\bEqu{equ:cCcircv}
\cC_v^\circ\!=\!\big(\cC_v\setminus \bigcup_{e\in E_\Gamma, v\in e}  \tn{cl}(U^\ve_{e,v})\big)\! \times\! \prod_{v'\in V_\Gamma\setminus \{v\}} \cB_{v'} \! \times\!\De^{|E_\Gamma|}.
\eEqu
Then, for sufficiently small $\de$, 
\bEqu{equ:smoothing-type-family}
\cC=\bigg(\coprod_{v\in V_\Gamma} \cC^\circ_v \coprod_{e\in E_\Gamma} \cC_e\bigg)/\sim\quad \tn{and}\quad \vec\mfz\!=\!(\mfz^{e^\circ})_{e^\circ\in E^\circ_\Gamma},
\eEqu
with the natural projection map, where 
$$
\aligned
&\mfz^{e^\circ}\big((x_v)_{v\in V_\Gamma}, (\la_e)_{e\in E_\Gamma}\big)\!=\![\mfz_{v(e^\circ)}^{e^\circ}(x_{v(e^\circ)}),(x_{v'})_{v'\in V_\Gamma\setminus\{v\}}, (\la_e)_{e\in E_\Gamma}]\!\in\![\cC^\circ_v]\quad \forall e^\circ\!\in\!E^\circ_\Gamma\\
&\tn{and}\quad\cC_{e,v_1}\!\ni\! \bigg( \big(y_{1},y_{2},\la_{e}\big), (x_v)_{\substack{v\in V_\Gamma\\ v\not\in e} },(\la_{e'})_{e'\in E_\Gamma\setminus\{e\}} \bigg)\sim
 \\
&\quad\qquad \qquad\qquad\qquad\bigg(y_{1}, \pi_{v_2}(y_2),(x_v)_{\substack{v\in V_\Gamma\\ v\not\in e} },(\la_{e})_{e\in E_\Gamma}\bigg) \!\in\!\cC^\circ_{v_1}\quad \forall e\!=\!\ll v_1,v_2\rr\!\in\!E_\Gamma,
\endaligned
$$
is a universal family around $C$.
In particular, the curve over $((x_v)_{v\in V_\Gamma},(\la_e)_{e\in E_\Gamma})$ 
is obtained from the curves $\{\pi_v^{-1}(x_v)\}_{v\in V_\Gamma}$ via the gluing identifications
$$
 w_{e,v_1}\cdot w_{e,v_2}=\la_e \quad \forall e\!=\!\ll v_1,v_2\rr \!\in\! E_\Gamma.
$$
If $w_{e,v}$ is compatible with $\aut(C)$ as in Lemma~\ref{lem:symmetric-gluing}, the action of $\aut(C)$ on $\prod_{v\in V_\Gamma}\cB_v$ (as well as $\prod_{v\in V_\Gamma}\cC_v$)  described in the argument after Definition~\ref{def:versal-family} and Remark~\ref{rem:family-with-unmarked} naturally extends to (\ref{equ:B-smmothing}) and (\ref{equ:smoothing-type-family}); for every $h\!\in\!\aut(C)$ and $e\!\in \!E_\Gamma$, 
\bEqu{equ:gluing-action}
h(\la_e)\!=\!\mu_{h,e}\la_{h(e)}.
\eEqu

\noindent
In Section~\ref{sec:main}, we will use universal families of standard type as in (\ref{equ:smoothing-type-family}) to build Kuranishi charts around maps with nodal domain. 

\subsection{Gromov Topology}\label{sec:topology}
The underlying topological space of the orbifold structure, constructed in \cite{RS}, on the Deligne-Mumford space $\ov\cM_{g,k}$, has a basis for the open sets the collection of all  sets 
\bEqu{quotient-open_e}
(\cB/\aut(C))\!\subset\! \ov\cM_{g,k}
\eEqu
where each $\cB$ is the base of a universal family as in Definition~\ref{def:versal-family}. With respect to this topology $\ov\cM_{g,k}$ is compact and Hausdorff; the proof uses the notion of DM-convergence of a sequence of marked nodal curves, see \cite[Section~14]{RS}. 
In this section, after a quick review of  the convergence problem in the Deligne-Mumford space, we recall the more general notion of Gromov convergence for a sequence of marked nodal stable maps.  
This sequential convergence topology is compact and metrizable as well.
\vskip.1in
\noindent

\bDef{def:cuttingset}
Given a $k$-marked genus $g$ (possibly not stable) nodal surface $C\equiv(\Si,\vec{z})$ with dual graph $\Gamma$, a 
\textbf{cutting configuration} with dual graph $\Gamma'$ is a set of disjoint embedded circles 
$$
\gamma\equiv\{\gamma_e\}_{e\in E_{\Gamma',\Gamma}}\!\subset\! \Si,
$$ 
away from the nodes and marked points, such that  the nodal marked surface $(\Si',\vec{z}')$ obtained by pinching every $\gamma_e$ into a node $q_e$ has dual graph $\Gamma'$.
\eDef

\noindent
Thus, a cutting configuration corresponds to a continuous map
$$
\pi_\gamma\colon C\lra C',
$$ 
denoted by a $\gamma$-\textbf{degeneration}\footnote{It is called deformation in \cite{RS}.} in what follows, onto a $k$-marked genus $g$ nodal surface $C'$ with dual graph $\Gamma'$,
such that $\vec{z}'\!=\! \pi_\gamma(\vec{z})$, the preimage of every node of $\Si$ is either a node in $\Si'$ or a circle in $\gamma$, and the restriction 
$$
\pi_\gamma\colon \Si\!\setminus\! \gamma \lra \Si'\! \setminus \!(\pi_\gamma (\gamma )\!\equiv \! \{q_e\}_{e\in E_{\Gamma',\Gamma}})
$$
is a diffeomorphism. Let
\bEqu{equ:deg-graphs}
\pi_\gamma^*\colon \Gamma'\lra \Gamma
\eEqu
be the map corresponding to $\pi_\gamma$ between the dual graphs. 
We have
$$
E_{\Gamma'}\!\cong \!E_\Gamma\cup E_{\Gamma',\Gamma}, \qquad E^\circ_{\Gamma'}\!\cong\!E^\circ_\Gamma,
$$
such that $\pi_\gamma^*|_{E_\Gamma\subset E_{\Gamma'}}$ and $\pi_\gamma^*|_{E^\circ_{\Gamma'}}$ are isomorphisms and 
$$
\pi_\gamma^*\colon E_{\Gamma',\Gamma}\lra V_{\Gamma}
$$
sends the edge $e$ corresponding to $\gamma_e$ to $v$, if $\gamma_e\!\subset\!\Sigma_v$.
For every ${v}'\!\in\!V_{\Gamma'}$ there exists a unique $v\!\in\!V_{\Gamma}$ and a connected component $U_{v'}$ of $\Sigma_v\!\setminus\! \{\gamma_e\}_{e\in E_{\Gamma',\Gamma}}$ such that $\Sigma'_{v'}\!\subset\!\Sigma'$ is obtained by collapsing the boundaries of $\tn{cl}({U}_{v'})$. This determines the surjective map 
$$
\pi_\gamma^*\colon V_{\Gamma'}\lra V_{\Gamma},\quad {v}'\lra v.
$$
From this perspective, a cutting configuration corresponds to expanding each vertex $v\!\in\!\Gamma$ into a sub-graph $\Gamma'_v\!\subset\!\Gamma'$ with vertices $(\pi_\gamma^*)^{-1}(v)$ and edges $(\pi_\gamma^*)^{-1}(v)\!\cap \!E_{\Gamma',\Gamma} $; moreover, $g_v\!=\!g_{\Gamma'_v}$, the ordering of marked points are as before, and 
\bEqu{equ:a-division}
A_v\!=\!\sum_{v'\in (\pi_\gamma^*)^{-1}(v)} A_{v'}.
\eEqu
Figure~\ref{fig:cut-conf} illustrates a cutting configuration over a $1$-nodal curve of genus $3$ and the corresponding dual graphs. 

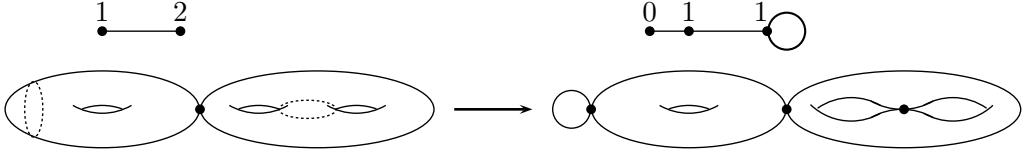
\begin{figure}
\begin{pspicture}(7.8,1)(11,2.5)
\psset{unit=.26cm}

\pscircle*(35,9){.25}\pscircle*(39,9){.25}
\psline[linewidth=.07](35,9)(39,9)
\rput(35,10){$1$}\rput(39,10){$2$}
\psellipse[linewidth=.07](35,5)(5,2)
\psarc[linewidth=.07](35,7.8){3}{240}{300}
\psarc[linewidth=.07](35,2.2){3}{70}{110}
\psellipse[linestyle=dashed,dash=1pt,linewidth=.07](31.5,5)(0.5,1.46)
\pscircle*(40,5){.25}
\psellipse[linewidth=.07](46,5)(6,2)
\psarc[linewidth=.07](43,7.8){3}{240}{300}
\psarc[linewidth=.07](43,2.2){3}{70}{110}
\psarc[linewidth=.07](48,7.8){3}{240}{300}
\psarc[linewidth=.07](48,2.2){3}{70}{110}
\psellipse[linestyle=dashed,dash=1pt,linewidth=.07](45.5,5)(1.5,.5)

\psline[linewidth=1pt]{->}(53,5)(57,5)

\pscircle*(63,9){.25}\pscircle*(65,9){.25}\pscircle*(69,9){.25}
\psline[linewidth=.07](65,9)(69,9)\psline[linewidth=.07](63,9)(65,9)\pscircle(70,9){1}
\rput(63,10){$0$}\rput(65,10){$1$}\rput(68.7,10){$1$}
\psellipse[linewidth=.07](65,5)(5,2)
\psarc[linewidth=.07](65,7.8){3}{240}{300}
\psarc[linewidth=.07](65,2.2){3}{70}{110}
\pscircle[linewidth=.07](59,5){1}\pscircle*(60,5){.25}
\pscircle*(70,5){.25}
\psellipse[linewidth=.07](76,5)(6,2)\pscircle*(76,5){.25}
\psarc[linewidth=.07](73.3,7.4){3}{226}{300}
\psarc[linewidth=.07](75.9,1.8){3.2}{60}{120}
\psarc[linewidth=.07](78.5,7.4){3}{240}{313}
\psarc[linewidth=.07](73.3,2.6){3}{60}{125}
\psarc[linewidth=.07](75.9,8.2){3.2}{240}{300}
\psarc[linewidth=.07](78.5,2.6){3}{55}{120}

\end{pspicture}
\caption{On left, a $1$-nodal curve of genus $3$ and a cutting set made of two circles. On right, the resulting pinched curve.}
\label{fig:cut-conf}
\end{figure}

\noindent 
A set $\{\pi_{\gamma_a}\colon C_a\lra C'\}_{a\in S}$ of degenerations of marked nodal curves is called \textbf{monotonic} if $\Gamma(C_a)\!\cong\!\Gamma$ for some fixed $\Gamma$ and the induced maps $\pi_{\gamma_a}^*\colon \Gamma\!\lra\! \Gamma'$ are all the same. In this situation, the underlying marked nodal surfaces are diffeomorphic; i.e.,
\bEqu{equ:indentify-domians}
(C_a,\gamma_a)\cong \big((\Sigma,\mfj_a,\vec{z}),\gamma\big)\quad \forall a\in \N,
\eEqu
for some fixed marked surface $(\Si,\vec{z})$ with dual graph $\Gamma$ and cutting configuration $\gamma$. In addition, with this identification, the set of nodes 
$$
 \{q_e\}_{e\in E_{\Gamma',\Gamma}}\subset\!\Sigma'
$$ 
is independent of $a\!\in\!\N$ and we denote the complement by $\Si'_*$.

\begin{definition}[{\cite[Definition 13.3]{RS}}]\label{def:DM-convergence}
We say a  sequence $\lrc{ C_a\!\equiv\!(\Si_a,\mfj_a,\vec{z}_a) }_{a=1}^\infty$ of genus $g$ $k$-marked nodal curves  \textbf{monotonically} converges to $C'\!\equiv \!(\Si', \mfj',\vec{z}' )$, if there exist a sequence of cutting configurations $\gamma_a$ on $C_a$ of type $\Gamma'$ and a monotonic sequence $\pi_{\gamma_a}\colon C_a\!\lra\! \!C'$  of $\gamma_a$-degenerations such that the sequence $(\pi_{\gamma_a}|_{\Si_a\setminus \gamma_a})_*\mfj_a $ converges to 
$\mfj'|_{\Si'_*}$ in the $C^\infty$-topology\footnote{Uniform convergence on compact sets with all derivatives.}. 
\eDef

\noindent
By \cite[Section 13]{RS}, the topology underlying the holomorphic orbifold structure of Section~\ref{sec:DM} on $\ov\cM_{g,k}$ is equivalent to the sequential $\tn{DM}$-convergence topology: a sequence $\lrc{ C_a}_{a=1}^\infty$ of genus $g$ $k$-marked stable nodal curves $\tn{DM}$-converges to $C'$ if a subsequence of that monotonically converges to $C'$.

\bExa{exa:deformation}
For $a\!\in\!\N$ sufficiently large, let $(\Si_a,\mfj_a)\!=\! \P^1$ with $\vec{z}_a=(0,\frac{1}{a},1,\infty)$, 
$$
(\Si,\mfj)=\P^1_0~{}_{\infty}\!\!\cup_{0} \P^1_{\infty}
$$ 
be a $1$-nodal genus $0$ curve made of two copies $\P^1_0$ and $\P^1_{\infty}$ of $\P^1$ and 
$$
\vec{z}=\big((0,1)\subset (\P^1_0)^2\big)\cup \big((1,\infty)\subset (\P^1_{\infty})^2\big),
$$ 
and $\gamma_a$ be the cutting configuration given by the circle 
$$
\{z\!\in\! \P^1\cong \C\cup \{\infty\} \colon |z|\!=\!\frac{1}{\sqrt{a}}\}.
$$ 
Let $\eta\colon\R^{\geq 0}\!\lra\! \R^{\geq 0}$ be an increasing cutoff function such that $\eta\!\equiv\! 0$ on $[0,1]$ and $\eta\!\equiv \!1$ on $[2,\infty)$. Define $\pi_{\gamma_a}\colon\Si_a \lra \Si$ by
$$
\pi_{\gamma_a}(z)=
\begin{cases} 
\eta(\sqrt{a}|z|)z \in \P^1_\infty&\mbox{if } |z| > \frac{1}{\sqrt{a}}, \\ 
\eta((\sqrt{a}|z|)^{-1}) a z \in \P^1_0 & \mbox{if } |z| < \frac{1}{\sqrt{a}}, \\
  \P^1_0 \ni \infty \sim 0\in \P^1_\infty& \mbox{if } |z| = \frac{1}{\sqrt{a}} . 
 \end{cases} 
$$
Then $\pi_{\gamma_a}(\vec{z}_a)\!=\!\vec{z}$ and $\pi_{\gamma_a}$ is holomorphic away from the neck 
$$
A_a\!=\!\{z\colon \frac{1}{2\sqrt{a}}\leq |z| \leq \frac{2}{\sqrt{a}}\}.
$$ 
Since the image of $\P^1\!-\!A_a$ is exhausting, the sequence $(\pi_{\gamma_a}|_{\Si_a\setminus \gamma_a})_*\mfj_a $ converges to 
$\mfj|_{(P^1_0-\infty) \cup (\P^1_\infty-0)}$ in the $C^\infty$-topology. Therefore, $\{C_a\}_{a=1}^\infty$ is $\tn{DM}$-converging to $C$. 
\eExa

\noindent
The following result, known as Gromov's Compactness Theorem \cite[Theorem 1.5.B]{Gromov}, describes a convergence topology on $\ov\cM_{g,k}(X,A,J)$  which is compact and metrizable;  see also \cite[Theorem 1.2]{Hu}, \cite[Theorem 0.1]{Ye}, \cite[Theorem 2.2.1]{AL}, and \cite[Theorem 5.3.1]{MS2004} for the genus $0$ case.
In the special case of Deligne-Mumford space, Gromov convergence is equal to DM-convergence discussed above.

\bThm{thm:gromov}
Let $(X,\om)$ be a compact symplectic manifold and $J\!\in\! \cJ(X,\om)$. 
Let 
$$
\lrc{f_i \equiv \lrp{ u_i,C_i \equiv \lrp{\Si_i,\mfj_i,\vec{z}_i }}}_{i=1}^\infty
$$ 
be a sequence of $k$-marked genus $g$ stable $J$-holomorphic maps of bounded symplectic area. After passing to a subsequence, still denoted by $\{f_i\}_{i=1}^\infty$, there exists a unique (up to automorphism) $k$-marked genus $g$ stable $J$-holomorphic map 
$$
f'\!\equiv\!(u',C'\equiv(\Si',\mfj',\vec{z}'))
$$ 
such that $\{C_i\}_{i=1}^\infty$ monotonically converges to $C'$,
and the corresponding monotonic subsequence, still denoted by $\{C_i\}_{i=1}^\infty$, has the following properties. 
\bEnum
\item\label{l:smooth-conv} Via $\gamma_i$-degeneration maps $\pi_{\gamma_i}\colon C_i \!\lra\! C'$, the restriction $u_i|_{\Sigma_i\setminus \gamma_i}\!\circ\! \pi_{\gamma_i}^{-1}|_{\Si'_*}$ uniformly converges to $u|_{\Si'_*}$ over compact sets.
\item With $\Gamma\cong \Gamma(C_i)$ as in definition of monotonic sequences, for every $e\!\in\!E_{\Gamma',\Gamma}$, 
$$
\lim_{i\lra \infty} u_i(\gamma_{i,e})\!=\!u'(q_e),\qquad q_e\!\equiv \!\pi_{\gamma_i}(\gamma_{i,e}).
$$
\item\label{l:symp-area} Symplectic area of $f'$ coincides with the symplectic area of $f_i$, for all $i\!\in\!\N$.
\eEnum
\eThm
\noindent
It follows from the properties \ref{l:smooth-conv} and \ref{l:symp-area} that for every $v'\!\in\! \Gamma'\!=\!\Gamma(C')$, with $U_{i,v'}\!\subset\! \Sigma_i$ as in the definition of $\pi_{\gamma_i}^*$,
$$
\lim_{i\lra\infty} \int_{\tn{cl}(U_{i;v'})} {u_i}^*\om = \int_{\Si'_{v'}} (u')^*\om.
$$
Moreover, the stronger identity (\ref{equ:a-division}) holds. 
With respect to the identification of the domains and degeneration maps 
$$
(\pi_{\gamma_i}\colon \Si_i\lra \Si')\cong (\pi_{\gamma}\colon\Si\lra \Si')
$$ 
as in (\ref{equ:indentify-domians}), the second property implies that the sequence $(u_i\!\colon\! \Si\lra X)_{i=1}^\infty$
$C^0$-converge to $u\circ \pi_{\gamma}$. The construction of primary Kuranishi charts in  Section~\ref{sec:main} gives a rather different description of the Gromov topology, generalizing open sets of the form (\ref{quotient-open_e}); see Remarks~\ref{rem:varphi-dependence} and~\ref{rmk:Gromov2}, and Section~\ref{sec:natural}.

\section{Kuranishi structure over moduli space of stable maps}\label{sec:main}
In this section, for every $\ov\cM_{g,k}(X,A,J)$, we sketch the construction of a class of natural Kuranishi structures that together with the results of Section~\ref{sec:VFC}, gives us the VFC of Theorem~\ref{thm:VFC}.

\noindent
To this end, in Sections~\ref{sec:canonical-smooth}, \ref{sec:canonical-nodal}, and \ref{sec:natural}, we construct a class of Kuranishi charts, called \textbf{primary charts}, around  maps with smooth stable domain, nodal stable domain, and unstable domain, respectively.
The construction of every primary chart depends (at most) on the choice of a so called obstruction space, stabilizing data, slicing divisors, smooth trivializations of the universal families, and holomorphic coordinates at the nodes in (\ref{equ:holo-functions}); see the list of all auxiliary data in Page~\pageref{small-enough_c}.

\noindent
In Section~\ref{sec:natural-KUR}, we construct a class of natural Kuranishi structures on  $\ov\cM_{g,k}(X,A,J)$ in the sense that every Kuranishi chart is obtained from the auxiliary data  of a finite set of primary charts. We call these \textbf{induced charts}. In Section~\ref{sec:transition}, we show that the induced charts admit coordinate change maps.

\noindent
The construction of Kuranishi charts around maps with nodal domain involves a gluing theorem with respect to the gluing parameters of (\ref{equ:B-smmothing}). In order for the Kuranishi and change of coordinate maps to be smooth, we need a strong gluing theorem and few other analytical results that, due to the lack of space, we will only recall from \cite{FOOO-detail3} in Sections~\ref{sec:canonical-nodal} and~\ref{sec:natural-KUR}.

\subsection{Analytics preliminaries}\label{sec:analysis}
Let $(X,\om)$ be a closed symplectic manifold. Fix $J\!\in\! \cJ(X,\om)$ and let $g_J$ be the associated metric\footnote{Not to be confused with the genus.}, $\nabla$  be the Levi-Civita connection of $g_J$,  and 
$$
\wt\nabla_v w := \nabla_v w\! -\! \frac{1}{2}J (\nabla_vJ) w\quad \forall v,w\!\in \Gamma(X,TX)
$$ 
 be the associated complex linear \textbf{Hermitian} connection. The Hermitian connection $\wt\nabla$ coincides with $\nabla$ if and only if $(X,\om,J)$ is K\"ahler, i.e. $\nabla J\equiv 0$. The torsion $\wt{T}$ of $\wt\nabla$ is related to Nijenhueis tensor (\ref{equ:Nij}) by
 $$
 \wt{T}(w,v)=-\frac{1}{4}N_J(w,v)\quad \forall w,v\!\in\! TX;
 $$
 see \cite[Section 2.1]{MS2004}.
 
\noindent 
Fix $\ell\!\in\!\N$ and $p\!>\!1$ such that 
\bEqu{equ:ell-p}
\ell p\!>\!2;
\eEqu
for the necessity of this assumptions see \cite[Section 3.1]{MS2004}.
For the proof of gluing theorem in \cite{FOOO-detail3}, the authors consider $p\!=\!2$ and $\ell$ sufficiently large.  
For a fixed domain $\Si$, let $W^{\ell,p}(\Si,X)$ be the separable Banach manifold of $(\ell,p)$-smooth maps\footnote{For $(\ell,p)$ as in (\ref{equ:ell-p}), with $\ell\!-\!\frac{2}{p}\!>\!m$, by Sobolev inequality, every such $W^{\ell,p}$-smooth map is at least $C^{m}$-smooth; see \cite[Theorem B.1.11]{MS2004}.} from $\Si$ into $X$. The local charts\footnote{With smooth transition maps.} of the Banach manifold structure $W^{\ell,p}(\Si,X)$ are constructed in the following way; see the proof of \cite[Proposition 3.2.1]{MS2004} for more details.

\bRem{rmk:notation}
If $\Si$ is the underlying smooth domain of a marked curve $C$, we will often  write $W^{\ell,p}(\Si,X)$ as $W^{\ell,p}(C,X)$. The Banach manifold structure does not depend on the choice of the marked complex structure $(\mfj,\vec{z})$, but later, we will consider certain additional structures and evaluation maps on this space which do depend on the extra data.
For instance, the group $\aut(C)$ of the automorphisms of $C$ acts on $W^{\ell,p}(C,X)$ by 
$$
u \lra u\circ h\qquad \forall ~h\!\in\!\aut(C),~u\!\in\!W^{\ell,p}(C,X).
$$

\eRem

\noindent
Given a map $u \colon\! \Si\!\lra\! X$, for every $x\!\in\! \Si$, let 
$$
\tn{exp}_{u(x)}\colon T_{u(x)}X \lra X
$$
be the exponentiation map of $\wt\nabla$. For $\ep\!>\!0$ sufficiently small, depending on the injectivity radius of $g_J$ across $\Si$, $\tn{exp}_{u(x)}$ is a diffeomorphism from the $\ep$-ball $B_\ep(x)$ around $0\!\in\! T_{u(x)}X$ onto its image, for all $x\!\in\! \Si$.
For $\ep\!>\!0$ sufficiently small, given an $(\ell,p)$-smooth vector field 
\bEqu{zeta_e}
 \ze \!\in\! \Gamma^{\ell,p}(\Si,u^*TX),\quad\tn{with}\quad |\ze(x)|\!<\!\ep,
\eEqu
let 
\bEqu{equ:exp_uze}
u_\ze\!\equiv\! \exp_{u}\ze
\eEqu 
be the associated close-by $(\ell,p)$-smooth map\footnote{Note that two different maps in $W^{\ell,p}(\Si,X)$ can have identical images in $X$.} in $W^{\ell,p}(\Si,X)$. Therefore, for $\ep$ sufficiently small, the exponential map (\ref{equ:exp_uze}) gives an identification
\bEqu{ep-chart_e}
\tn{exp}\colon \big( B_\ep(0) \subset \Gamma^{\ell,p}(\Si,u^*TX)\big) \lra \big(B_\ep(u)\subset W^{\ell,p}(\Si,X)\big)
\eEqu
of an $\ep$-neighborhood of $0$ in the Banach space $\Gamma^{\ell,p}(\Si,u^*TX)$ with a neighborhood, denoted by $B_\ep(u)$, of a fixed map $u$ in $W^{\ell,p}(\Si,X)$. 

\noindent
For every differentiable map $u\colon \Si \!\lra\! X$ and complex structure $\mfj$ on $\Si$,
$$
\nd u \!\in \!\Gamma(\Si,\Om_\Si^1\otimes u^*TX)
$$ 
decomposes into complex linear and anti-complex linear parts,
\bEqu{dbar_e}
\aligned
&\partial u\equiv \partial_{J,\mfj} u\! =\!\frac{1}{2}( \nd u \!-\! J \circ \nd u \circ \mfj) \in \Gamma(\Sigma, \Om^{1,0}_{\Si,\mfj}\otimes_\C u^*TX)\quad \tn{and}\\
&\dbar u\equiv \dbar_{J,\mfj} u =\frac{1}{2}( \nd u \!+\! J \circ \nd u \circ \mfj)\in \Gamma(\Si,\Om^{0,1}_{\Si,\mfj}\otimes_\C u^*TX),
\endaligned
\eEqu
respectively, where $\Om^{1,0}_{\Si,\mfj}$ and $\Om^{0,1}_{\Si,\mfj}$ are the sheaves/spaces of $(1,0)$ and $(0,1)$-forms on $\Si$ with respect to $\mfj$, respectively; $u$ is $(J,\mfj)$-holomorphic whenever $\dbar u\!\equiv\!0$. 

\bRem{lack-of-Symm_rmk}
If $f\!=\!(u,C)$ is a $J$-holomorphic map, the automorphism group $\aut(f)$ of $f$, defined in (\ref{equ:auto-nodal-map}), acts on $\Gamma^{\ell,p}(\Sigma, u^*TX)$ by 
\bEqu{autActiononGamma_e}
\ze\lra h^*\ze \!\in\! \Gamma^{\ell,p}(\Sigma, h^*u^*TX)\!=\!\Gamma^{\ell,p}(\Sigma, (u\circ h)^*TX))\!=\!\Gamma^{\ell,p}(\Sigma, u^*TX)\quad \forall h\!\in\!\aut(f).
\eEqu
The exponential map identification of (\ref{ep-chart_e}) is $\aut(f)$-equivariant with respect to (\ref{autActiononGamma_e}) and the action of $\aut(f)\!\subset\!\aut(C)$ on $W^{\ell,p}(C,X)$ as in Remark~\ref{rmk:notation}, i.e.
$$
u_{h^*\ze}=u_\ze\circ h\qquad \forall h\!\in\!\aut(f),~ \ze\!\in\! \Gamma(\Sigma, u^*TX).
$$
\eRem

\noindent
Motivated by (\ref{dbar_e}), for a fixed curve $C\!=\!(\Si,\mfj,\vec{z})$, let
\bEqu{E-bundle_e}
E^{\ell-1,p}(C,X)\lra W^{\ell,p}(C,X)
\eEqu
be the infinite dimensional smooth Banach bundle whose fiber over every map $u\colon \!C\!\lra\! X$ is the Banach space of $(\ell\!-\!1,p)$-smooth $u^*TX$-valued $(0,1)$-forms
$$
\Gamma^{\ell-1,p}(\Sigma,\Om^{0,1}_{\Sigma,\mfj}\otimes_\C u^*TX).
$$ 
Similarly, the tangent bundle of $W^{\ell,p}(\Si,X)$,
\bEqu{TW-bundle_e}
TW^{\ell,p}(\Si,X)\lra W^{\ell,p}(\Si,X),
\eEqu
is the infinite dimensional smooth Banach bundle whose fiber over every map $u\colon\! \Si\!\lra\! X$ is the Banach space of $(\ell,p)$-smooth vector fields  
$
\Gamma^{\ell,p}(\Sigma, u^*TX).
$
Similarly to Remark~\ref{rmk:notation} and (\ref{autActiononGamma_e}), the action of $\aut(C)$ on $C$ extends to $E^{\ell-1,p}(C,X)$ and $TW^{\ell,p}(C,X)$ such the projection maps to $W^{\ell,p}(C,X)$ are $\aut(C)$-equivariant. For instance, the action of $\aut(C)$ on  $E^{\ell-1,p}(C,X)$ is given by
\bEqu{autActiononGamma2_e}
(\eta,u)\lra(h^*\eta,u\circ h)\qquad \forall h\!\in\! \aut(C).
\eEqu
The local defining charts for  $E^{\ell-1,p}(C,X)$ (and similarly for $TW^{\ell,p}(\Si,X)$) around a fixed $J$-holomorphic map $f\!=\!(u,C)$ are constructed in the following way.
For a fixed map $u\colon\! C\!\lra\! X$, $\ze$ as in (\ref{zeta_e}), and $u_\ze$ as in (\ref{equ:exp_uze}), let
$$
\tn{P}_{u,\ze}\colon u^*TX\lra u_\ze^*TX
$$
denote the complex line bundle isomorphism given by the parallel transport along the geodesic
$$
\big(s\!\lra\! \tn{exp}_{u(x)}(s\ze(x)\big)_{s\in [0,1]}\quad \forall x\!\in\! \Si
$$ 
with respect to the complex linear connection $\wt\nabla$.
This induces a similarly denoted complex linear isomorphism
\bEqu{equ:Pzeu}
\tn{P}_{u,\ze}\colon E^{\ell-1,p}(C,X)|_{u}\lra E^{\ell-1,p}(C,X)|_{u_\ze}.
\eEqu
These parallel translation maps, together with the local charts (\ref{ep-chart_e}) for $W^{\ell,p}(C,X)$, give us local product charts 
\bEqu{E-localchart_e}
\aligned
&\Gamma^{\ell-1,p}(\Sigma,\Om^{0,1}_{\Sigma,\mfj}\otimes_\C u^*TX)\times  B_\ep(0)\stackrel{P_u}{\lra}E^{\ell-1,p}(C,X)|_{B_\ep(u)},\\
&(\eta,\zeta)\lra (\eta',u'),\quad  u'\!=\!u_\ze,\quad \eta'\!=\!\tn{P}_{u,\ze}(\eta),
\endaligned
\eEqu
for the Banach bundle (\ref{E-bundle_e}). Since 
$
\tn{P}_{u,h^*\ze}(h^*\eta)\!=\!h^*\tn{P}_{u,\ze}(\eta),
$
similarly to Remark~\ref{lack-of-Symm_rmk}, the map $P_u$ is $\aut(f)$-equivariant.

\noindent
We now discuss\footnote{Under the assumption that the domain is stable.} the case where the marked complex structure $(\mfj,\vec{z})$ on the domain $\Si$ can vary as well.
Assume $C\equiv (\Si,\mfj, \vec{z}\big)$ is a smooth stable marked curve and let
\bEqu{equ:a-universal-family}
\lrp{\pi\colon\cC\lra \cB, \vec{\mfz}, b\!\in\! \cB}
\eEqu 
be a sufficiently small universal family around $\![C]\!\in\! \cM_{g,k}$ as in Section~\ref{sec:DM}. 
Let 
\bEqu{equ:Banach-smooth-domain}
W^{\ell,p}(\cC,X):= \big\{ (u,a)\colon ~~a\!\in\! \cB,~ u\! \in\! W^{\ell,p}(C_a,X)\big\},
\eEqu
be the set of pairs $(u,a)$ such that 
$a$ is a point of $\cB$,  $\Si_a\!=\!\pi^{-1}(a)\!\subset\!\cC$, and $u$ is an $(\ell,p)$-smooth  map from 
$\Si_a$ to $X$.
Recall from (\ref{equ:induced-h}) that (for sufficiently small $\cC$) the action of $\aut(C)$ on $C$ (uniquely) extends to $\cC$. This action further extends to $W^{\ell,p}(\cC,X)$ by 
$$
(u,a) \lra (u\circ \varphi_h^\cC, (\varphi_h^\cB)^{-1}(a))\qquad \forall  h\in\aut(C),
$$
where $(\varphi_h^\cC,\varphi_h^\cB)\colon (\cC,\cB)\lra (\cC,\cB)$ is the automorphism corresponding to $h$ in (\ref{equ:induced-h}).

\noindent
By definition, each fiber of the projection map $W^{\ell,p}(\cC,X)\!\lra\!\cB$ has a Banach manifold structure 
but the total space does not a priori come with any natural smooth structure. 
In order to define a Banach manifold  structure on (\ref{equ:Banach-smooth-domain}),
we need to fix a smooth trivialization of $\cC$. Such a trivialization gives us a trivialization of $W^{\ell,p}(\cC,X)$ and thus a product Banach manifold structure on that.
More precisely, by Remark~\ref{rem:smoothly-trivial-family}, for $\cB$ sufficiently small, $\cC$ is smoothly trivial, i.e. there exists an $\aut(C)$-equivariant  diffeomorphism 
\bEqu{equ:smooth-trivial-phi}
\varphi\colon \cC\!\lra\!\Si\!\times\! \cB
\eEqu
such that $\pi\circ\varphi^{-1}$ is the projection onto the second factor, each section $\varphi\!\circ\!\mfz^i$ is constant, and $\varphi_b\!=\!\tn{id}_{\Si}$.
In this sense, we can think of $\cB$ as a smooth family of complex structures $\{\mfj_{\varphi,a}\!=\!(\varphi_a)_*(\mfj_{a})\}_{a\in \cB}$ on the fixed marked surface $(\Si,\vec{z})$ with $\mfj_{\varphi,b}\!=\!\mfj_{b}\!=\!\mfj$. A choice of such $\varphi$ induces an $\aut(C)$-equivariant trivialization
\bEqu{equ:Bspace-triv}
W^{\ell,p}(C,X)\!\times\!\cB\!\stackrel{T_\varphi}{\lra}\! W^{\ell,p}(\cC,X), \quad  (u,a)\lra (u\circ \varphi_a,a),
\eEqu
and thus a Banach manifold structure which we denote by $W^{\ell,p}_\varphi(\cC,X)$. 
In particular, for a fixed $J$-holomorphic map $u\colon \!C\!\lra\!X$, the exponentiation map of (\ref{ep-chart_e}) and the product structure (\ref{equ:Bspace-triv}) give us
 $\aut(u,C)$-equivariant local charts 
\bEqu{LocalChart_e0}
(B_\ep(0)\subset \Gamma^{\ell,p}(\Sigma, u^*TX)) \times\cB \stackrel{\tn{exp}_{\varphi}}{\xrightarrow{\hspace*{.8cm}}} W^{\ell,p}_\varphi(\cC,X)   ,\quad
(\ze,a)\lra (u_\ze\circ \varphi_a,a).
\eEqu
If $\varphi'$ is another such  smooth trivialization, then $\varphi\!=\rho\!\circ\!\varphi'$, where
\bEqu{equ:CTM}
\rho\colon \Si\times \cB\lra \Si\times \cB, \quad \rho(x,a)=(\rho_a(x),a), 
\eEqu
is a fiber-preserving diffeomorphism fixing the constant sections $\vec{z}$. Unless $\rho_a$ is constant in $a\!\in\!\cB$, the change of trivialization map 
\bEqu{equ:Frho}
\xymatrix{
W^{\ell,p}(C,X)\!\times\!\cB \ar[d]^{T_{\rho}} \ar[rr]^{T_\varphi} &&  W^{\ell,p}(\cC,X)  \\
 W^{\ell,p}(C,X)\!\times\!\cB \ar[rru]^{T_{\varphi'}}&&
}
, \quad T_{\rho}(u,a)=(u\circ\rho_a,a),
\eEqu
is continuous but not smooth. Therefore, the Banach manifold structure (\ref{equ:Bspace-triv}) depends on a choice of a smooth trivialization. Whenever the choice of $\varphi$ is fixed or irrelevant in a discussion, for simplicity, we may omit $\varphi$ from $W^{\ell,p}_{\varphi}(\cC,X)$ and simply write it as $W^{\ell,p}(\cC,X)$. For more on this  differentiability issue see \cite[Section 3.1]{MW2}.

\noindent
Similarly to (\ref{equ:Banach-smooth-domain}), for a universal family $\cC\!\lra\!\cB$ around $C$ as in (\ref{equ:a-universal-family}), let 
\bEqu{equ:Banach-Bundle-cC}
\aligned
&E^{\ell-1,p}(\cC,X):= \big\{ (\eta,u,a)\colon ~(u,a)\!\in\!W^{\ell,p}(\cC,X),~ \eta\! \in\! E^{\ell-1,p}(C_a,X)|_{u}\big\}\quad \tn{and}\\
&TW^{\ell,p}(\cC,X):= \big\{ (\ze,u,a)\colon ~(u,a)\!\in\!W^{\ell,p}(\cC,X),~ \ze\! \in\! TW^{\ell,p}(C_a,X)|_{u}\big\}.
\endaligned
\eEqu
For $\cB$ sufficiently small, a choice of smooth trivialization $\varphi$ as in (\ref{equ:smooth-trivial-phi}) also gives us trivializations\footnote{For simplicity, we use $T_\varphi$ for both the trivialization map of $W^{\ell,p}_\varphi(\cC,X)$ and $E^{\ell-1,p}_\varphi(\cC,X)$. }  
\bEqu{LT-ETW_e}
\aligned
&E^{\ell-1,p}(C,X)\times \cB\stackrel{T_\varphi}{\lra} E^{\ell-1,p}_\varphi(\cC,X), \quad (\eta,u,a)\!\lra\!\big(\varphi_a^*[\eta]^{0,1}_{\mfj_{\varphi,a}},u\circ \varphi_a,a\big),\quad \tn{and}\\
&TW^{\ell,p}(C,X)\times \cB\stackrel{T_\varphi}{\lra} TW^{\ell,p}_\varphi(\cC,X), \quad (\ze,u,a)\!\lra\!\big(\varphi_a^*\ze,u\circ \varphi_a,a\big),
\endaligned
\eEqu
where 
\bEqu{01part_e}
[\eta]^{0,1}_{\mfj_{\varphi,a}}= \frac{1}{2} \big(\eta + J\circ  \eta \circ \mfj_{\varphi,a}\big)
\eEqu
is the $(0,1)$-part of $\eta$ with respect to the complex structure $\mfj_{\varphi,a}$ on $\Si$. For $\mfj_{\varphi,a}$ sufficiently close to $\mfj$, the map $\eta\!\lra\![\eta]^{0,1}_{\mfj_{\varphi,a}}$ is an isomorphism between the Banach spaces
$$
\Gamma^{\ell-1,p}(\Sigma,\Om^{0,1}_{\Sigma,\mfj}\otimes_\C u^*TX)\quad\tn{and}\quad \Gamma^{\ell-1,p}(\Sigma,\Om^{0,1}_{\Sigma,\mfj_{\varphi,a}}\otimes_\C u^*TX).
$$ 
Also note that $\varphi_a^*[\eta]^{0,1}_{\mfj_{\varphi,a}}\!=\![\varphi_a^*\eta]^{0,1}_{\mfj_{a}}$.
For every $h\!\in\! \aut(C)$ and $a\!\in\!\cB$, with $a\!=\!\varphi_h^\cB(a')$, by (\ref{hjphia_e}) and the $\aut(C)$-equivariance of $\varphi$ we have 
$$
[h^*\eta]^{0,1}_{\mfj_{\varphi,a'}}= h^*[\eta]^{0,1}_{\mfj_{\varphi,a}}\quad \tn{and}\quad \varphi_a\circ \varphi_h^\cC=h\circ \varphi_{a'};
$$
i.e. the trivialization maps in (\ref{LT-ETW_e}) are $\aut(C)$-equivariant. 
With respect to these smooth structures
\bEqu{equ:infinite-Banach-bundle}
\aligned
&E^{\ell-1,p}_\varphi(\cC,X)\lra W^{\ell,p}_\varphi(\cC,X), \quad (\eta,u,a)\lra (u,a),\quad\tn{and}\\
&TW^{\ell,p}_\varphi(\cC,X)\lra W^{\ell,p}_\varphi(\cC,X), \quad (\ze,u,a)\lra (u,a),
\endaligned
\eEqu
are smooth $\aut(C)$-equivariant Banach bundles. In particular, for a fixed map $u\colon \!C\!\lra\!X$, with notation as in (\ref{LocalChart_e0}), local trivializations  of $E^{\ell-1,p}(C,X)$ as in (\ref{E-localchart_e}) and the product structure (\ref{LT-ETW_e}) give us $\aut(f)$-equivariant local charts 
\bEqu{ELocalChart_e}
\aligned
& \Gamma^{\ell-1,p}(\Sigma,\Om^{0,1}_{\Sigma,\mfj}\otimes_\C u^*TX) \times B_\ep(0)\times\cB   
\stackrel{P_{\varphi} }{\lra}  E^{\ell-1,p}_\varphi(\cC,X)|_{\tn{exp}_\varphi(B_\ep(0)\times \cB)} ,\\
&(\eta,\ze,a)\lra (\eta',u',a), \quad u'= u_\ze\circ \varphi_a, \quad \eta'= \varphi_a^*P_{u,\ze,a}(\eta):= \varphi_a^*[P_{u,\ze}(\eta)]^{0,1}_{\mfj_{\varphi,a}}.
\endaligned
\eEqu
\bRem{rmk:PPchange}
Similarly to (\ref{equ:Frho}), for two different trivializations $\varphi'$ and $\varphi$, with the change of parametrization map $\rho$ as in (\ref{equ:CTM}), the resulting change of trivialization map $T_\rho$ in the commutative diagram below
\bEqu{equ:Frho22}
\xymatrix{
E^{\ell-1,p}(C,X)\!\times\!\cB \ar[d]^{T_{\rho}} \ar[rr]^{T_\varphi} && E^{\ell-1,p}(\cC,X)  \\
E^{\ell-1,p}(C,X)\!\times\!\cB \ar[rru]^{T_{\varphi'}}&&
}
\eEqu
is simply the pull-back map $(\eta,u,a)\!\lra\!(\rho_a^*\eta,u\circ \rho_a,a)$. However, it is rather impossible to  write down explicit equations for the change of coordinate maps $\tn{exp}^{-1}_{\varphi'}\!\circ\tn{exp}_\varphi$ corresponding to (\ref{LocalChart_e0}) and $P^{-1}_{\varphi'}\!\circ\!P_\varphi$ corresponding to (\ref{equ:Frho22}).
\eRem
\subsection{Case of smooth stable domain}\label{sec:canonical-smooth}

In this section, we outline the construction of primary Kuranishi charts around maps with smooth stable domain.
Assume $C\equiv (\Si,\mfj, \vec{z}\big)$ is a smooth stable marked curve and let
$$
\lrp{\pi\colon\cC\lra \cB, \vec{\mfz}, b\!\in\! \cB}
$$ 
be a sufficiently small universal family around $\![C]\!\in\! \cM_{g,k}$ as in (\ref{equ:a-universal-family}). The non-linear Cauchy-Riemann operator 
\bEqu{equ:CR-section}
\dbar\colon W^{\ell,p}(\cC,X)\!\lra\! E^{\ell-1,p}(\cC,X), \quad u\lra \dbar u,
\eEqu
given by (\ref{dbar_e}), is an $\aut(C)$-equivariant section of (\ref{equ:infinite-Banach-bundle}) whose zero locus is the set of $J$-holomorphic maps on the fibers of $\cC\!\lra\!\cB$.
For any choice of smooth trivialization $\varphi$ of (\ref{equ:a-universal-family}) as in (\ref{equ:smooth-trivial-phi}), the Cauchy-Riemann section  (\ref{equ:CR-section}) is smooth with respect to the smooth structures $W^{\ell,p}_\varphi(\cC,X)$ and $E^{\ell-1,p}_\varphi(\cC,X)$; see \cite[Section 3]{MS2004} for more details.

\noindent
Deformation theory of $J$-holomorphic maps close to a particular $J$-holomorphic map $u$ in $W^{\ell,p}(C,X)$ (i.e. if the domain is fixed) is described by the first order $\aut(u,C)$-equivariant\footnote{Recall from Section~\ref{sec:analysis} that the finite group $\aut(u,C)$ acts on both the domain and target of $\tn{D}_u\dbar$.} approximation or linearization of (\ref{equ:CR-section})\footnote{Along the zero set of the section $\dbar$, we have a canonical decomposition $TE^{\ell-1,p}(C,X)\cong E^{\ell-1,p}(C,X)\oplus TW^{\ell,p}(C,X)$.}, which is of the form
\bEqu{equ:linearization}
\aligned
 &\tn{D}_u\dbar\colon  TW^{\ell,p}(C,X)|_{u}\lra E^{\ell-1,p}(C,X)|_{u},\\
 &\tn{D}_u\dbar (\ze)\!=\! \frac{1}{2}(\wt\nabla \ze + J \wt\nabla \ze \circ \mfj )+\frac{1}{4}N_J(\ze,\partial u);
 \endaligned
 \eEqu
see \cite[Chapter 3.1]{MS2004}. The connection 
$$
\wt\nabla^{0,1}\equiv \frac{1}{2}(\wt\nabla  + J \wt\nabla  \circ \mfj )
$$
is the $(0,1)$-part of $\wt\nabla$ and defines a $\dbar$-operator, thus a holomorphic structure, on $u^*TX$; see \cite[Remark C.1.1]{MS2004}.
The Nijenhuis term in $\tn{D}_u\dbar$ is a zeroth order complex anti-linear compact operator. Therefore, $\tn{D}_u\dbar$ is a 
deformation of some $\dbar$-operator by a compact operator; see \cite[Definition C.1.5]{MS2004}. It follows that $\tn{D}_u\dbar$ is Fredholm, c.f. \cite[Theorem A.1.5]{MS2004}, i.e. it has a finite dimensional kernel and cokernel.

\noindent
Define
$$
\Def(u)\!=\!\ker(\tn{D}_u\dbar)\quad \tn{and} \quad \Obs(u)\!=\!\coker(\tn{D}_u\dbar).
$$
Both $\Def(u)$ and $\Obs(u)$ are finite dimensional $\aut(u,C)$-invariant real vector spaces. 
The first space corresponds to infinitesimal deformations of $u$ (over the fixed domain $C$) and the second one is the obstruction space for integrating elements of $\Def(u)$ to actual deformations. 
By Riemann-Roch Theorem \cite[Theorem C.1.10]{MS2004},
\bEqu{equ:RR}
\dim_\R \Def(u)\!-\!\dim_\R \Obs(u)\!=\!2\big( c_1(u^*TX) \!+\! \dim_\C\! X (1\!-\!g)\big).
\eEqu
If $\Obs(u)\!\equiv\! 0$, by Implicit Function Theorem \cite[Theorem A.3.3]{MS2004}, in a small neighborhood $B_\ep(u)$ of $u$ in $W^{\ell,p}(C,X)$ as in (\ref{ep-chart_e}), the set of $(J,\mfj)$-holomorphic maps
\bEqu{equ:Vu}
V_u \equiv \dbar^{-1}(0)\cap B_\ep(u)
\eEqu
is a smooth $\aut(u,C)$-invariant manifold of real dimension $\dim_\R \Def(u)$, all the elements of $\Def(u)$ are smooth (i.e. independent of choice of $\ell$ and $p$), and 
$
T_u V_u\cong \Def(u);
$
c.f. \cite[Theorem 3.1.5]{MS2004}. The manifold $V_u$ carries a natural orientation: starting with the $\dbar$-part of $\tn{D}_u\dbar $, both the kernel and cokernel are complex linear and thus naturally oriented. By going from $t\!=\!0$ to $t\!=\!1$ in the one-parameter family of Fredholm operators
$$
\big(\wt\nabla^{0,1}+t \frac{1}{4}N_J(\cdot ,\partial u)\big)_{t\in [0,1]},
$$
\cite[Proposition A.2.4]{MS2004} gives us a natural orientation on $\Def(u)$.

\bRem{rem:actual-Obs}
Let $f\!=\!(u,C)$. More generally, deformation theory of $f$, i.e. if we allow deformations of both $u$ and $C$, is given by a long exact sequence of the form
\bEqu{equ:long-def}
\aligned
0\lra &\Def(u)\lra \Def(f) \lra \Def(C) \stackrel{\de}{\lra} \\
	&\Obs(u)  \lra  \Obs(f) \lra 0,
\endaligned
\eEqu
where 
$$
\Def(C)\cong T_b\cB\cong H^1(\cT\Si_\mfj(-\vec{z}))
$$ 
is as in (\ref{equ:Dm-tangent}). The long exact sequence (\ref{equ:long-def}) is the hypercohomology of a short exact sequence of complexes of fine sheaves constructed in the following way. 

\noindent
For every holomorphic or complex vector bundle $E\!\lra\! (\Si,\mfj)$, let $\Om^0(E)$ and  $\Om^{0,1}_{\Sigma,\mfj}(E)$ denote the associated fine sheaves of smooth sections of $E$ and of smooth $E$-valued $(0,1)$-forms, respectively. 

\noindent
The map $\nd u\colon T\Si\!\lra\!TX$ gives rise to a similarly denoted map 
$$
\nd u\colon \Om^0(\cT\Si_\mfj(-\vec{z}))\lra \Om^0(u^*TX).
$$
Away from the marked points $\vec{z}$, it maps a local generating section $ \frac{\partial}{\partial w}$ to $\nd u( \frac{\partial}{\partial w})$\footnote{With $\partial w\!=\!\frac{1}{2}(\!\frac{\partial}{\partial x}\!-\!\mfi \frac{\partial}{\partial y})$, $\nd u( \frac{\partial}{\partial w})$ means $\frac{1}{2}(\nd u(\frac{\partial}{\partial x})-J \nd u(\frac{\partial}{\partial y}))\!=\!\nd u(\frac{\partial}{\partial x})$.}, and in a holomorphic chart $w$ around a marked point $z^i\!\equiv\! (w\!=\!0)$, it maps the local generating section $ \frac{\partial}{\partial w}$ to $w \nd u( \frac{\partial}{\partial w})$. The following commutative diagram has exact rows:
$$
\xymatrix{
 0\ar[r]\ar[d]& \Om^0(\cT\Si_\mfj(-\vec{z}))\ar[r]  \ar[d]^{\nd u\oplus \dbar} & \Om^0(\cT\Si_\mfj(-\vec{z})) \ar[d]^{\dbar}\\
 \Om^0(u^*TX) \ar[r]\ar[d]^{\tn{D}_u\dbar}  &\Om^0(u^*TX)\oplus \Omega^{0,1}_{\Sigma,\mfj}(\cT\Si_\mfj(-\vec{z}))\ar[r] \ar[d]^{\tn{D}_u\dbar - \nd u}& \Omega^{0,1}_{\Sigma,\mfj}(\cT\Si_\mfj(-\vec{z}))\ar[d]\\
 \Omega^{0,1}_{\Sigma,\mfj}(u^*TX)\ar[r] & \Omega^{0,1}_{\Sigma,\mfj}(u^*TX)\ar[r] & 0 \; ;
}
$$
i.e. it is an exact sequence of chain complexes given by the columns. Then (\ref{equ:long-def}) is the  hypercohomology of this diagram.
We refer to the article of Siebert-Tian \cite[Section 3.2]{TS} for an explicit  and more detailed description of (\ref{equ:long-def}); also see \cite[Section 24.4]{mirror} for more details in the algebraic case.  
\eRem

\noindent
In the light of Remark~\ref{rem:actual-Obs}, similarly to (\ref{equ:Vu}), if $\Obs(u)\!=\!0$, then 
\bEqu{equ:defu+defC}
\Def(f) \cong \Def(u)\oplus \Def(C)
\eEqu
and every tangent direction in  $\Def(f)$ is integrable; i.e. in a small neighborhood $B_\ep(f)$ of $f$ in\footnote{For any choice of smooth structure $W^{\ell,p}_\varphi(\cC,X)$ on $W^{\ell,p}(\cC,X)$. } $W^{\ell,p}(\cC,X)$, the set of $J$-holomorphic maps
$
V_f \equiv \dbar^{-1}(0)\cap B_\ep(f)
$
is a smooth $\aut(f)$-invariant manifold of real dimension $\dim_\R \Def(f)$, all the elements of $\Def(f)$ are smooth (i.e. independent of choice of $\ell$ and $p$), and 
$
T_f V_f\cong \Def(f).
$
A similar conclusion holds if $\Obs(f)\!=\!0$, i.e. if the map $\de$ in (\ref{equ:long-def}) is surjective; however, we will work with $\Obs(u)$ for simplicity\footnote{As Remark~\ref{rem:actual-Obs} indicates, a generalization of $\tn{D}_u\dbar$ in (\ref{equ:linearization}) to some Fredholm map $\tn{D}_f\dbar$ that includes deformations of complex structure of the domain is rather difficult to define and work with.  }. 

\noindent
The main idea of this article is that if $\Obs(u)\!\neq\!0$, we can still reduce the right-hand side of the first equation in (\ref{equ:linearization}) to a finite dimensional sub-space modulo which $\tn{D}_u \dbar$ is surjective and then we define $V_u$ and $V_f$ similarly.
The auxiliary  finite dimensional space then appears as an orbibundle over $V_f$ and gives rise to a natural Kuranishi chart on the moduli space in the sense of Definition~\ref{def:kur-chart} .

\bDef{def:obs-space}
Let $C\!=\!(\Si,\mfj,\vec{z})$ be a smooth marked curve\footnote{Not necessarily stable.} and $u\colon C\!\lra\!X$ be a $J$-holomorphic map. An \textbf{obstruction space} $E_f$ for $f\!=\!(u,C)$ is a finite dimensional complex linear (with respect to the induced action of $J$) sub-vector space
$$
E_f\subset \Gamma(\Si,\Om^{0,1}_{\Si,\mfj}\otimes_\C u^*TX)
$$
satisfying the  following conditions.
\bEnum
\item Every $\eta\!\in\! E_f$ is smooth;
\item\label{l:surj-E} with respect to the projection map
$$
\pi_{E_f}\colon \Gamma^{\ell-1,p}(\Si,\Om^{0,1}_{\Si,\mfj}\otimes_\C u^*TX) \lra \Gamma^{\ell-1,p}(\Si,\Om^{0,1}_{\Si,\mfj}\otimes_\C u^*TX)/E_f,
$$
the composition $ \pi_{E_f}\!\circ\! \tn{D}_u\dbar$ is surjective;
\item\label{l:inv-E} and $E_f$ is $\aut(u,C)$-invariant.
\eEnum
\eDef

\bRem{rmk:exists-Obs}
If $\dim X\!>\!0$, every $J$-holomorphic map admits a plethora of obstruction spaces. 
In fact, for any open set $U\!\subset\!\Si$, there exists an obstruction space supported in $U$. For a proof see \cite[Lemma 17.11]{FOOO-detail}; it can also be deduced from the proof of \cite[Proposition 3.2.1]{MS2004}.
\eRem

\bRem{rmk:minimal-obs}
In the light of Remark~\ref{rem:actual-Obs}, it is often more practical for calculations to replace Definition~\ref{def:obs-space}.\ref{l:surj-E} with the weaker assumption that
$$
[E_f]\lra \tn{Obs}(f),\quad [E_f]\!=\!E_f/\tn{Image}(D_u{\dbar}) \subset \tn{Obs}(u),
$$
is surjective. For example, it sometimes happens that the moduli space $\ov\cM_{g,k}(X,A)$ has the structure of an oriented orbifold of larger than expected dimension, $\tn{Obs}(f)$ has constant rank throughout $\ov\cM_{g,k}(X,A)$, and they are fibers of an orbibundle on the moduli space. In these examples, we can take the obstruction space to be precisely (isomorphic to) $\tn{Obs}(f)$\footnote{i.e. we take the minimal possible choice.} and we obtain a pure orbibundle Kuranishi structure on the moduli space (with Kuranishi map $s\!\equiv\!0$) that significantly simplifies the calculations; see the examples in Sections~\ref{subsec:degree0} and \ref{subsec:elliptic}.

\eRem

\noindent
Assume $f\!=\!(u,C)$ is a fixed $J$-holomorphic map with smooth stable domain, $E_f$ is an obstruction space for $f$,  $(\pi\colon\!\cC\!\lra\!\cB, b\!\in\! \cB)$ is a universal family around $C$, and  $\varphi$ is a smooth trivialization for $\cC$ as in (\ref{equ:smooth-trivial-phi}). 
With respect to the local charts (\ref{LocalChart_e0}) for $W^{\ell,p}(\cC,X)$ and (\ref{LT-ETW_e}) for $E^{\ell-1,p}(\cC,X)$, let
\bEqu{equ:Vf}
\aligned
&\wt{V}_f=\wt{V}_{f}(\varphi)=\{ (\ze,a) \!\in\! B_\ep(0)\!\times\! \cB\colon~ \tn{P}_{u,\ze,a}^{-1}(\dbar_{J,\mfj_{\varphi,a}} u_\ze) \!\in E_{f} \} \quad \tn{and}\\
&\wt{s}_f\colon \wt{V}_f \lra E_f,\quad (\ze,a)\lra \tn{P}_{u,\ze,a}^{-1}(\dbar_{J,\mfj_{\varphi,a}} u_\ze).
\endaligned
\eEqu
In other form, with notation as in (\ref{ELocalChart_e}) and
\bEqu{Bfep_e}
B_\ep(f)\!=\!\tn{exp}_\varphi(B_\ep(0)\times\cB),
\eEqu 
by Definition~\ref{def:obs-space}.\ref{l:inv-E} and $\aut(C)$-equivariance of $P_\varphi$, we obtain a finite rank $\aut(f)$-equivariant \textit{obstruction bundle} $E_{f,\varphi}$,
\bEqu{EVphitilde_e}
\xymatrix{
E_{f,\varphi}=P_\varphi( E_f \times B_\ep(0)\times\cB) \ar@{^{(}->}[r]\ar[rd]& E^{\ell-1,p}_\varphi(\cC,X)|_{B_\ep(f)}\ar[d]\\
& B_\ep(f)\subset W^{\ell,p}(\cC,X).
}
\eEqu
Let
\bEqu{EVtilde_e}
\aligned
&V_f=\{  f' \!\in\! B_\ep(f)\colon~ \dbar f'\!\in E_{f,\varphi} \},\quad\quad U_f=E_{f,\varphi}|_{{V}_f},\\
& \tn{and} \quad {s}_f \colon {V}_f \lra U_f, \quad f' \lra \dbar f'.
\endaligned
\eEqu
Then the tuple\footnote{By abuse of notation, we use $\wt{s}_f$ to also denote for the section $x\!\in\! \wt{V}_f \lra (x,\wt{s}_f(x))\!\in\! \wt{V}_f\!\times\! \wt{E}_f$.} $(\wt{U}_f\!=\!E_f\!\times\! \wt{V}_f\! \lra\! \wt{V}_f, \wt{s}_f)$ in (\ref{equ:Vf}) is the $\aut(f)$-equivariant trivialization of $({U}_f\!\lra\! {V}_f, {s}_f=\dbar)$ with respect to the local charts (\ref{LocalChart_e0}) and (\ref{ELocalChart_e}), i.e.
$$
{V}_f=\tn{exp}_\varphi(\wt{V}_f),\quad {U}_f=P_\varphi(E_f\times \wt{V}_f), \quad\tn{and}\quad s_f=P_\varphi\circ \wt{s}_f.
$$

\noindent
For every choice of obstruction space as in Definition~\ref{def:obs-space}, the analogue of Riemann-Roch formula (\ref{equ:RR}) asserts that
\bEqu{equ:RR2}
\dim_\R \ker  (\pi_{E_f}\circ \tn{D}_u\dbar)\!=\!2\big( c_1(u^*TX) \!+\! \dim_\C\! X (1\!-\!g) \!+\!\dim_\C E_f\big).
\eEqu
Then, similarly to the case of $\Obs(u)\!=\!0$ before, for sufficiently small $\ep$ and $\cB$, 
it follows from Definition~\ref {def:obs-space}.\ref{l:surj-E} and Implicit Function Theorem that $V_f$
 is a smooth $\aut(f)$-invariant manifold of real dimension 
\bEqu{equ:dim-with-E}
2\big( c_1(u^*TX) \!+\! (\dim_\C\! X\!-\!3) (1\!-\!g) \!+\!\dim_\C E_f\big)
\eEqu
independent of choice $(\ell,p)$ and 
$$
T_{f} {V}_f \cong T_{0\times \{b\}} \wt{V}_f \cong \ker (\pi_{E_f}\!\circ\!\tn{D}_u\dbar)\oplus T_b\cB.
$$
The zero set of $s_f$ is a set of actual $J$-holomorphic maps close to $f$; see Remark~\ref{rem:varphi-dependence}. 
Putting together, by (\ref{equ:dim-with-E}), we obtain a Kuranishi chart  
\bEqu{equ:Can-Kur-chart}
\cU_f\equiv (\pr_f\colon U_f\lra V_f , G_f\!=\!\aut(f), s_f,\psi_f), \qquad \psi_f(f')=[f']\quad \forall f'\!\in\! s_f^{-1}(0),
\eEqu
of the expected dimension (\ref{equ:exp-dim}) centered at $[f]\!\in\! \cM_{g,k}(X,A)$.  We will denote the trivialization (\ref{equ:Vf}) of $\cU_f$ by 
\bEqu{equ:Can-Kur-chart-triv}
\aligned
&\wt\cU_f\equiv (\wt\pr_f\colon \wt{U}_f=E_f\times \wt{V}_f\lra \wt{V}_f , G_f\!=\!\aut(f),  \wt{s}_f,\wt{\psi}_f), \\
&\wt\psi_f(\ze,a)=[u_\ze,C_a]\quad \forall (\ze,a)\!\in\! \wt{s}_f^{-1}(0);
\endaligned
\eEqu
Depending on the situation, it is some times more convenient to work with $\wt\cU_f$ and vice versa.
The manifold $V_f$ admits an $\aut(f)$-equivariant smooth fibration 
$$
\tn{st}_f\colon V_f\! \lra\! \cB,\quad\tn{st}_f((u',C_a))\!=\!a,
$$  
where the fiber over each point is as in (\ref{equ:Vu}). 
Therefore, since $E_f$ and $T_b\cB$ are complex linear, $V_f$ carries a natural orientation similarly to the case of (\ref{equ:Vu}).

\bRem{aut-limit_rmk}
Similarly to (\ref{bbprime-aut_e}), for sufficiently small ${V}_f$, it follows from Theorem~\ref{thm:gromov} and (\ref{bbprime-aut_e}) that for every $J$-holomorphic map $f'\!=\!(u',C_a)\!\in\! ({s}_f^{-1}(0)\!\cap\! {V}_f)$ we have an identification 
$$
\aut(f')\cong\{h\colon \colon \varphi_h^\cB(a)=a,~~u'\circ \varphi_h^\cC|_{C_a}=u',~~h\!\in \!\aut(f)\!\subset\! \aut(C)\}
$$
between $\aut(f')$ and a subgroup of $\aut(f)$. Therefore, for a sufficiently small open neighborhood ${V}_{f;f'}\!\subset\!{V}_f$, the restriction 
$$
\cU_{f;f'}\equiv (\pr_f\colon U_f|_{V_{f;f'}}\lra V_{f;f'}, \aut(f'),s_f|_{V_{f;f'}}, \psi_f|_{V_{f;f'}})
$$
defines a sub-chart of $\cU_f$ around $[f']\!\in\!\cM_{g,k}(X,A)$ in the sense of Definition~\ref{Sub-Chart_dfn}, with the correct isotropy group.
The same argument holds for the primary Kuranishi charts of the next two sections.
\eRem

\bRem{rem:chart-varphi}
If $\varphi'$ is another trivialization of the universal curve, the obstruction bundle $E_{f,\varphi'}$ would be different from $E_{f,\varphi}$ and thus there is no obvious identification of $V_f(\varphi')$ and $V_f(\varphi)$ away from the subset of $J$-holomorphic maps.
\eRem

\bRem{rem:lift-ev}
The evaluation and forgetful maps in (\ref{equ:evst}) naturally lift to $W^{\ell,p}(\cC,X)$ and thus to $V_f$. 
The lifts are smooth. Therefore, the resulting Kuranishi charts satisfy the requirements of Definition~\ref{def:rho-Kur-structure}. Since the domain is stable in this case, the lift of forgetful map is the composition of smooth projection map $\st_f$ above and the quotient map $\cB\!\lra\!\cM_{g,k}$. If $\dim X\!>\!0$, for sufficiently large $E_f$, the lifted evaluation map
$$
\tn{ev}\colon V_f \lra X^k,\quad \tn{ev}((u',\Si_a,\mfj_a,\vec{z}_a))\!=\!u'(\vec{z}_a),
$$
becomes a submersion. 
\eRem

\bRem{rem:effectivise}
If $\dim X\!>\!0$, for appropriate choices of $E_f$, the action of $\aut(f)$ on $V_f$ would be effective. This way, we can fortify the possible non-effectiveness of the action of $\aut(f)\!\subset\!\aut(C)$ on $\cB$. However, for simplicity, in the calculations of Section~\ref{sec:examples}, we avoid such enlargements and allow non-effective actions in the natural way.
\eRem

\noindent
\bRem{rem:embedding}
For every fixed $J$-holomorphic map, inclusion defines a partial order on the set of admissible obstruction spaces in Definition~\ref{def:obs-space}. If $E'_f\!\subset\!E_f$, then the germ of Kuranishi chart  $\cU'_f$ corresponding to $E'_f$ as in (\ref{equ:Can-Kur-chart}) embeds in the Kuranishi chart  $\cU_f$ corresponding to $E_f$ (with respect to the same fixed $\varphi$) as in Definition~\ref{def:intersection}.  This embedding satisfies the tangent bundle condition of Definition~\ref{def:Kuranishi-tangent}. In fact, if $\wt{s}_f$ and $\wt{s}'_f$ are the Kuranishi maps of $\wt\cU_f$ and $\wt\cU'_f$ as in (\ref{equ:Vf}), respectively, then the corresponding normal direction derivative map in (\ref{dspq}),
$$
\nd \wt{s}_{f}/\wt{s}'_{f}\colon\cN_{\wt{V}'_{f}}{\wt{V}_f}|_{f}\lra E_f/E'_f,
$$
is a restriction of $\tn{D}_u\dbar$. Since $\tn{D}_u\dbar$ is surjective modulo $E'_f$, $\nd \wt{s}_{f}/\wt{s}'_{f}$ is surjective and thus an isomorphism (by dimensional reason).
\eRem

\noindent
The natural Kuranishi structure constructed in Section~\ref{sec:natural-KUR} is made of charts \textbf{induced} by the defining equations of the charts $V_f$ constructed in this section and Sections~\ref{sec:canonical-nodal} and \ref{sec:natural}.
For this reason, we call the Kuranishi charts of Sections~\ref{sec:canonical-smooth}, \ref{sec:canonical-nodal}, and \ref{sec:natural}, \textbf{primary charts}.
The following observation explains the meaning and necessity of the induced charts.

\noindent
Assume $f'\!=\!(u',C_{a})$ is a $J$-holomorphic map in ${V}_f$ sufficiently close to $f$. Then the subspace $E_{f'}\!=\! E_{f,\varphi}|_{f'}$ is an obstruction space for $f'$; see Remark~\ref{aut-limit_rmk} and \cite[Proposition 17.22]{FOOO-detail}.
By the argument after Definition~\ref{def:versal-family}, for $a\!\in\!\cB$ sufficiently close to $b$,  $\cC\!\lra\! \cB$ is also a universal family around $C_{a}$. 
However, unless $\wt\nabla$ is a flat connection, the manifold 
$$
{V}_{f'}\!\subset W^{\ell,p}(\cC,X)
$$ 
constructed via $E_{f'}$ is not naturally isomorphic to a neighborhood of $f'$ in ${V}_f$.
In this simple case, an ``induced'' chart at $f'$ is the restriction to some neighborhood ${V}_{f;f'}\!\subset\! {V}_f$ of  $\cU_f$; see  Remark~\ref{aut-limit_rmk}.
In Section~\ref{sec:natural-KUR}, we will generalize this trick to the cases where $f'$ belongs to several primary charts $\{{V}_{f_i}\}_{i=1}^N$.
The defining equation of an induced chart around $f'$ will be a sum of Cauchy-Riemann equations centered at $f_i$, with $i\!\in \![N]$.

\bRem{rem:family-with-unmarked2}
Similarly to Remark~\ref{rem:family-with-unmarked}, we also need to consider marked maps $(u,C\!=\!(\Si,\mfj,\vec{z},q))$ where $q\!=\!\{q_1,\ldots,q_\ell\}$ is an un-ordered set of disjoint marked points away from $\vec{z}$. In certain cases, we will require the obstruction space to be supported away from $q$. Other than that,   the construction of primary charts above readily extends to this situation. In this case, we use a universal family around $C\!=\!(\Si,\mfj,\vec{z},q)$ as in Remark~\ref{rem:family-with-unmarked} and construct a Kuranishi chart via such a universal family, similarly. 
The only difference between a chart for $(u,(\Si,\mfj,\vec{z},\vec{q}))$ and the resulting chart for $(u,(\Si,\mfj,\vec{z},q))$ is that the latter may have a larger group of automorphisms.
\eRem

\bRem{rem:varphi-dependence}
Although the exponentiation map identification
\bEqu{equ:Bepf}
B_\ep(f)\!\stackrel{\tn{exp}_{\varphi}}{\cong} \!\{ (\ze,a) \in  \Gamma^{\ell,p}(\Si,u^*TX)\!\times\!\cB\colon~~|\ze|\!<\! \ep\}
\eEqu
of a neighborhood of $f$ in $W^{\ell,p}(\cC,X)$ with a neighborhood of 
$$
0\times \{b\}\in\Gamma^{\ell,p}(\Si,u^*TX)\!\times \!\cB
$$ 
depends on the choice of a smooth trivialization $\varphi$ in (\ref{equ:smooth-trivial-phi}), as we pointed out earlier, the topology defined by the neighborhoods $B_\ep(f)$ is independent of the choice of $\varphi$.
Therefore, the following notion of convergence to $f$ is independent of the choice of $\varphi$. We say a sequence $(f_i=(u_i,C_i=(\Si_i,\mfj_{i},\vec{z}_i))_{i=1}^\infty$ of $J$-holomorphic maps $\ep$-\textbf{converges} to $f$ (or they are $\ep$-\textbf{close} to $f$) if 
for every $i\!>\!\!>\!1$ there exists $a_i\!\in\! \cB$ and $\ze_i\!\in\!u^*TX$ such that $C_i\!\cong\!C_{a_i}$, 
$u_i\!=\!\tn{exp}_{u}\ze_i\circ \varphi_{a_i}$, and 
$$
\lim_{i\lra \infty} a_i\!=\!b, \quad \lim_{i\lra \infty} |\ze_i|_{\ell,p}\!=\!0.
$$
This notion of convergence coincides with the Gromov's convergence of Theorem~\ref{thm:gromov}; see \cite{Ye}. 
Therefore, with notation as in (\ref{equ:Can-Kur-chart}), the topology on $\cM_{g,k}(X,A)$ determined by\footnote{If one $f$ in $\cM_{g,k}(X,A)$ is domain-stable, then every element of $\cM_{g,k}(X,A)$ is domain-stable.} ``open sets'' $\psi_f (s_f^{-1}(0))$ coincides with the Gromov topology.
In Sections~\ref{sec:canonical-nodal} and \ref{sec:natural}, we will extend this notion of convergence to nodal maps with possibly unstable domain; see Remark~\ref{rmk:Gromov2}. 
This alternative description of the Gromov convergence is the description of the choice in \cite{FO,FOOO} and directly shows that our footprint maps are homeomorphisms onto open sets of $\ov\cM_{g,k}(X,A)$. The identification (\ref{equ:Bepf}) also shows that the Gromov topology is locally (and thus by Remark~\ref{rem:top-remark} globally) metrizable.

\eRem

\subsection{Case of stable nodal domain}\label{sec:canonical-nodal}
In this section we extend the construction of previous section to primary charts for nodal maps with stable domain. The construction depends on the choice of an obstruction space $E_f$, smooth trivializations of universal families of components of the nodal domain, and the holomorphic functions in (\ref{equ:holo-functions}).

\noindent
With notation as in Sections~\ref{sec:stable} and \ref{sec:DM}, let $C\!=\!(\Si,\mfj,\vec{z})$ be a stable marked nodal curve with dual graph $\Gamma$ and $u\colon C\!\lra\!X$ be a $J$-holomorphic map. The analogue of linearization map (\ref{equ:linearization}) for $f\!=\!(u,C)$, still denoted by $\tn{D}_u\dbar$, is the restriction of $\bigoplus_{v\in V_\Gamma} \tn{D}_{u_v}\dbar$ to 
$$
\aligned
&\tn{Domain}(\tn{D}_u\dbar )=TW^{\ell,p}(C,X)|_{u}\equiv \Gamma^{\ell,p}(\Sigma, u^*TX):=\\
& \bigg\{\ze\!\equiv\!(\ze_v)_{v\in V_\Gamma}\!\in\! \bigoplus_{v\in V_\Gamma}\Gamma^{\ell,p}(\Sigma, u_v^*TX)\colon \ze_{v}(q_{e,v})\!=\!\ze_{v'}(q_{e,v'})\!\in\! T_{u(q_e)}X\quad \forall e\!=\!\ll v,v'\rr\!\in\! E_\Gamma \bigg\}
\endaligned
$$
with image in the direct sum
\bEqu{target_e}
\tn{Target}(\tn{D}_u\dbar )\!=\! \Gamma^{\ell-1,p}(\Si,\Om^{0,1}_{\Si,\mfj}\otimes_\C u^*TX):= \bigoplus_{v\in V_\Gamma} \Gamma^{\ell-1,p}(\Si_v,\Om^{0,1}_{\Si_v,\mfj_v}\otimes_\C u_v^*TX).
\eEqu
Recall from the notation of (\ref{equ:nodalcurve}) that $q_{e,v}$ is the special point on $\Si_v$ corresponding to the node $q_e\!\in\!\Si$ and the set of points $q_v\!=\!\{q_{e,v}\}_{e=\ll v,v'\rr}$ form a set of un-ordered marked points on $\Si_v$ away from $\vec{z}_v$.
Note that by \cite[Theorem B.1.11]{MS2004} and (\ref{equ:ell-p}), every $\ze_{v}$ is continuous and thus the evaluation $\ze_{v}(q_{e,v})$ is well-defined. Similarly to (\ref{autActiononGamma_e}) and (\ref{autActiononGamma2_e}), the automorphism group $\aut(f)$, defined in (\ref{equ:auto-nodal-map}), acts on the domain and target of $D_u\dbar$,  and $D_u\dbar$ is $\aut(f)$-equivariant.

\bDef{def:obs-space-2}
With $f\!=\!(u,C)$ as above, an \textbf{obstruction space} for $f$ is a direct sum
$$
E_f\!=\!\bigoplus_{v\in V_\Gamma} E_{f_v} \subset \Gamma^{\ell-1,p}(\Si,\Om^{0,1}_{\Si,\mfj}\otimes_\C u^*TX),
$$ 
satisfying the following conditions.
\bEnum
\item For every $v\!\in\! V_\Gamma$, 
$$
E_{f_v}\!\subset\!\Gamma^{\ell-1,p}(\Si_v,\Om^{0,1}_{\Si_v,\mfj_v}\otimes_\C u_v^*TX)
$$ 
is finite dimensional, made of smooth sections, and is supported away from $q_v$ (see Remark~\ref{rem:family-with-unmarked2}).
\item\label{l:surj-E2} With respect to the projection map
$$
\pi_{E_f}\colon \Gamma^{\ell-1,p}(\Si,\Om^{0,1}_{\Si,\mfj}\otimes_\C u^*TX) \lra \Gamma^{\ell-1,p}(\Si,\Om^{0,1}_{\Si,\mfj}\otimes_\C u^*TX)/E_f,
$$
the composition $ \pi_{E_f}\!\circ\! \tn{D}_u\dbar$ is surjective.
\item\label{l:inv-E2} $E_f$ is $\aut(f)$-invariant.
\eEnum
\eDef

\noindent
Fix an obstruction space $E_f$. With notation as in Section~\ref{sec:DM}, let 
\bEqu{equ:used-family}
\big(\pi\colon \cC \lra \cB,\vec{\mfz}=(\mfz^{e^\circ})_{e^\circ\in E^\circ_\Gamma}, b\big)
\eEqu
be a sufficiently small standard universal family as in (\ref{equ:smoothing-type-family}) around $C_b\!\equiv\!C$ such that 
\bEqu{equ:support-off-neck}
\tn{Support}(E_{f_v})\!\subset\! \Si_v\!\setminus \hspace{-.1in} \!\bigcup_{e\in E_\Gamma, v\in e} \hspace{-.1in}(\Si_v\cap U_{e,v}) \qquad \forall v\!\in\!V_\Gamma.
\eEqu
In this equation, $U_{e,v}$ are the domains\footnote{Open sets around nodal points.} of holomorphic functions $w_{e,v}$ in (\ref{equ:holo-functions}), which we used to smooth out the nodes. Recall from the construction following (\ref{equ:v-families}) that $\cC$ is explicitly made from a set of universal families 
$$
\lrp{\pi_v\colon \cC_v\!\to\! \cB_v, \vec\mfz_v,\mfq_v,b_v}\qquad \forall v\in V_\Gamma,
$$
by gluing them at the nodal points (sections) $\mfq_v$, $\cB\!=\!\prod_{v\in V_\Gamma} \cB_v\times\De^{E_\gamma}$, and the action of $\aut(C)$ on $\prod_{v\in V_\Gamma} \cB_v$ and $\prod_{v\in V_\Gamma} \cC_v$ extends to $(\cC,\cB)$ according to (\ref{equ:gluing-action}).

\noindent
For every $v\!\in\!V_\Gamma$, let 
$$
\cU_{f_v}\equiv (\pr_v\colon\! U_{v}\!\lra\! V_{v} , G_{v}, s_v,\psi_{v})
$$
be the primary Kuranishi chart\footnote{For simplicity, we write $(\pr_v\colon\! U_{v}\!\lra\! V_{v} , G_{v}, \psi_{v}, s_{v})$ instead of $(\pr_{f_v}\colon\! U_{f_v}\!\lra \!V_{f_v} , G_{f_v}, \psi_{f_v},  s_{v})$ in (\ref{equ:Can-Kur-chart}). } of (\ref{equ:Can-Kur-chart}) centered at $[f_v]$ with respect to the obstruction space $E_{f_v}$, the universal family $\cC_v$ around $C_v$, and a fixed choice of smooth trivializations $\varphi_v$ of $\cC_v$ as in (\ref{equ:smooth-trivial-phi}). 
Recall from (\ref{EVtilde_e}) that each $V_v$ is a finite dimensional submanifold of $W^{\ell,p}_{\varphi_v}(\cC_v,X)$ and we denote the trivialization of this manifold in $\Gamma^{\ell,p}(\Sigma_v, u_v^*TX)\!\times\!\cB_v$ under the exponential map $\tn{exp}_{\varphi_v}$ by $\wt{V}_v$; see (\ref{equ:Vf}). In the following discussion, depending on the situation and for the sake of convenience, we will switch from $\cU_v$ to its trivialization $\wt{\cU}_v$ as in (\ref{equ:Can-Kur-chart-triv}) and vice versa. 

\noindent
Let 
$$
\ev_v\colon {V}_v\lra X^{Q_v}\!=\!\hspace{-.05in}\prod_{(e,v)\in Q_v}\hspace{-.05in} X,\quad Q_v\!=\!\{(e,v)\colon v\!\in\! e, ~e\!\in\! E_\Gamma\},
$$
(and similarly denoted $\ev_v\colon V_v\lra X^{Q_v}$) be the evaluation map at the nodal points $\mf{q}_v$ in $\cC_v$, i.e.
$$
\ev_v^e(f'_v)=u'(\mf{q}_v^e(a_v))\quad \forall f'_v\!=\!\big(u'_v,C_{a_v}\big)\!\in\! {V}_{v},~v\!\in\! V_\Gamma. 
$$
Let $\De_X\!\subset\! X\!\times\!X$ be the diagonal submanifold, 
$$
\tn{Diag}_\Gamma(X)\!=\!(\De_X)^{E_\Gamma}\!\subset \! ( X\!\times\!X)^{E_\Gamma},
$$
and 
\bEqu{equ:Fiber-Product}
{V}_{f;0}\!=\!\bigg\{ (f'_v)_{v\in V_\Gamma}\in \prod_{v\in {V}_\Gamma} {V}_{v}: \ev^e_{v_1}(f'_{v_2})\!=\!\ev^e_{v_2}(f'_{v_2})\!\in\!X\quad \forall e\!=\!\ll v_1,v_2\rr\!\in\!E_\Gamma\bigg\}
\eEqu
be the fiber product of $\{{V}_v\}_{v\in V_\Gamma}$ over $\tn{Diag}_\Gamma(X)$; i.e.
\bEqu{equ:fp-intersection}
{V}_{f;0}= \big(\prod_{v\in {V}_\Gamma} {V}_{v}\big) \cap \big(\hspace{-.2in}\prod_{e=\ll v_1,v_2\rr\in E_\Gamma} \hspace{-.2in}\ev_{v_1}^e\times \ev_{v_2}^e\big)^{-1}(\tn{Diag}_\Gamma(X)).
\eEqu
In other form, by (\ref{equ:Vf}), via the exponential map identification of (\ref{LocalChart_e0}) we have
\bEqu{tildeVvsVf0_e}
\aligned
{V}_{f;0}\stackrel{(\tn{exp}_{\varphi_v})_{v\in V_\Gamma}}{\cong} &\wt{V}_{f;0}= \bigg\{ (\ze_v,a_v)_{v\in V_\Gamma}\! \in\! \tn{Domain}(\tn{D}_u\dbar)\!\times\!\prod_{v\in V_\Gamma} \cB_v\colon\\
 & \tn{P}_{u_v,\ze_v,a_v}^{-1}(\dbar_{J,\mfj_{\varphi_v,a_v}} (u_v)_{\ze_v})\!\in\!E_{f_v},~~ |\ze_v|\!<\! \ep,~~ \forall v\!\in\! V_\Gamma\bigg\}.
 \endaligned
\eEqu
Similarly to (\ref{equ:Vf}), in what follows, we will write the defining equations on the right-hand side of (\ref{tildeVvsVf0_e}) as 
\bEqu{family-P_e}
\tn{P}_{u,\ze,a}^{-1}(\dbar_{J,\mfj_{\varphi,a}} u_{\ze})\!\in\!E_{f}.
\eEqu

\bRem{rem:transversality-vs}
Condition~\ref{l:surj-E2} in Definition~\ref{def:obs-space-2} is equivalent to the intersection (\ref{equ:fp-intersection}) to be transverse to the diagonal. 
In \cite[Definition 10.1.1]{MS2004} and \cite[Assumption 10.3]{FOOO-detail}, the transversally condition is assumed instead of Definition~\ref{def:obs-space-2}.\ref{l:surj-E2}. That  the transversality condition implies Definition~\ref{def:obs-space-2}.\ref{l:surj-E2} is discussed after \cite[Proposition 10.5.1]{MS2004}; the other direction follows from a simple dimension counting as well. 
\eRem

\noindent
By  Definition~\ref{def:obs-space-2} and (\ref{equ:fp-intersection}), for sufficiently small $\{\cB_v\}$ and $\ep$, $V_{f;0}$ (isomorphically $\wt{V}_{f;0}$) is a smooth $\aut(f)$-invariant oriented submanifold of the fiber product 
\bEqu{Bepf0_e}
B_\ep(f;0)= \big(\prod_{v\in {V}_\Gamma} B_\ep(f_v)\big) \cap \big(\hspace{-.2in}\prod_{e=\ll v_1,v_2\rr\in E_\Gamma} \hspace{-.2in}\ev_{v_1}^e\times \ev_{v_2}^e\big)^{-1}(\tn{Diag}_\Gamma(X)).
\eEqu
(isomorphically $\tn{Domain}(\tn{D}_u\dbar)\!\times\!\prod_{v\in V_\Gamma} \cB_v$) with
$$
T_fV_{f;0}\cong T_{0\times \{b_v\}_{v\in V_\Gamma}}\wt{V}_{f;0}\cong \ker( \pi_{E_f}\!\circ\!\tn{D}_u\dbar) \oplus \bigoplus_{v\in V_\Gamma} T_{b_v}\cB_v.
$$
Let  
\bEqu{equ:bundle-0regular}
\pr_{f;0}\colon U_{f;0}\!=\! \boxtimes_{v\in V_\Gamma} U_v \lra V_{f;0}\quad \tn{and}\quad
s_{f;0}\!=\!\oplus_{v\in V_\Gamma} s_v\colon V_{f;0}\lra U_{f;0}
\eEqu
be the box-product\footnote{With respect to the projection maps $\rho_v\colon \!V_{f;0}\! \lra\! V_v$, the box-product bundle is defined by $\boxtimes_{v\in V_\Gamma} U_v\!=\!\bigoplus_{v\in V_\Gamma} \rho_v^* U_v$.  } Kuranishi bundle and direct sum Kuranishi map, respectively. With respect to the trivialization (\ref{tildeVvsVf0_e}) and (\ref{ELocalChart_e}), (\ref{equ:bundle-0regular}) is isomorphic to the  product bundle and direct sum Kuranishi map
\bEqu{equ:bundle-0}
\wt\pr_{f;0}\colon \wt{U}_{f;0}\!=\! \wt{V}_{f;0}\times E_f\lra \wt{V}_{f;0}\quad \tn{and}\quad
\wt{s}_{f;0}\!=\!\oplus_{v\in V_\Gamma} \wt{s}_v\colon \wt{V}_{f;0}\lra \wt{U}_{f;0}.
\eEqu
The forgetful map 
$$
\st_{f;0}\colon V_{f;0}\lra \cB
$$
has image the complex codimension $|E_\Gamma|$ submanifold 
$$
\prod_{v\in V_\Gamma}\hspace{-.03in}\cB_v\times 0^{E_\Gamma}\!\subset\!\cB
$$
of marked nodal curves of topological type $\Gamma$ sufficiently close to $C$. Similarly, the footprint map 
\bEqu{equ:psif0}
\psi_{f;0}\!=\!(\psi_v)_{v\in V_\Gamma}\colon \big(V_{f;0}\cap s_{f;0}^{-1}(0)\big)\lra \ov\cM_{g,k}(X,A)
\eEqu
has image the set of $J$-holomorphic maps in the stratum $\ov\cM_{g,k}(X,A)_\Gamma$ (see (\ref{equ:Gamma-M})) sufficiently close to $f$.
By Theorem~\ref{thm:gluing-thm-main} below, we $\aut(f)$-equivariantly extend the Kuranishi chart $\cU_{f;0}\!=\!(\pr_{f;0}\colon \!U_{f;0}\!\lra\! V_{f;0}, \aut(f), s_{f;0},\psi_{f;0})$ on $\cM_{g,k}(X,A)_{\Gamma}$ to a Kuranishi chart on  $\cM_{g,k}(X,A)$. In order to state the gluing theorem and describe the extended chart, we need to fix some notation.

\vskip.1in
\noindent
In Section~\ref{sec:DM}, for every $(a\!=\!(a_v),\la\!\equiv\!(\la_e))\!\in\!\cB$
the marked curve 
$$
C_{a,\la}\!=\!\big(\Si_{a,\la}\!=\!\pi^{-1}(a,\la), \vec{z}_{a,\la}=\mfz(a,\la)\big)\!\subset\! (\cC, \vec\mfz(\cB)) 
$$ 
is obtained from the curves $\{C_{a_v}\}_{v\in V_\Gamma}$ via the gluing identifications
\bEqu{equ:smoothing-nodes}
 w_{e,v_1}\cdot w_{e,v_2}=\la_e \quad \forall e\!=\!\ll v_1,v_2\rr \!\in\! E_\Gamma.
\eEqu
For sufficiently small $\de$ in the gluing construction of $\cC$, the defining equation (\ref{equ:smoothing-type-family}) gives us a holomorphic identification between the \textbf{core} regions
$$
C_{a,\la}^\circ\!\equiv\!C_{a,\la}\cap \bigcup_{v\in V_\Gamma} \cC_v^\circ \quad \tn{and}\quad C_{a}^\circ\!=\!    C_{a,0}^\circ;
$$
Therefore, in the following equations we will often write $C_a^\circ$ instead of $C_{a,\la}^\circ$.
We refer to the complement
$$
C_{a,\la}^\tn{neck}\!=\!C_{a,\la}\setminus C^\circ_{a,\la} \!\subset \!C_{a,\la}\cap \bigcup_{e\in E_\Gamma} \cC_e
$$
as the \textbf{neck} region. 
Let 
$$
C_{a}^\circ\!=\!\bigcup_{v\in V_\Gamma}C_{a_v}^\circ, \quad C_{a_v}^\circ\!=\!C_{a}\cap \cC_v^\circ\quad \forall v\!\in\! V_\Gamma, a\!\in\!\cB,
$$
be the decomposition of the core region $C_{a}^\circ\!\cong\!C_{a,\la}^\circ$ into connected components labeled by $v\!\in\!V_\Gamma$; see Figure~\ref{fig:neck-decomp}.
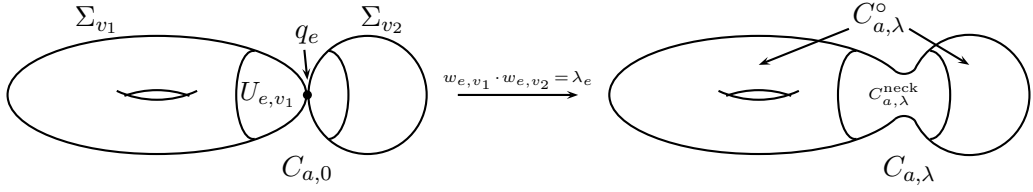
\begin{figure}
\begin{pspicture}(8,-4)(10,-5)
\psset{unit=.4cm}

\psellipse(25,-10)(5,2)\pscircle(32,-10){2}  
\pscircle*(30,-10){.15}
\rput(30,-12.5){$C_{a,0}$}
\rput(23,-7.5){$\Si_{v_1}$}\rput(32.5,-7.5){$\Si_{v_2}$}
\rput(28.7,-10){\small{$U_{e,v_1}$}}
\rput(30,-8){$q_e$}\psline{->}(29.9,-8.5)(30,-9.5)
\psellipticarc(28.2,-10)(.6,1.5){88}{272}
\psellipticarc(30.7,-10)(.7,1.5){-90}{90}
\psarc(25,-13){3.16}{70}{110}
\psarc(25,-7){3.16}{245}{295}

\psline{->}(35,-10)(39,-10)\rput(37,-9.5){\tiny{$ w_{e,v_1}\!\cdot\! w_{e,v_2}\!=\!\la_e$}}

\psellipticarc(45,-10)(5,2){10}{350}
\rput(50,-12.5){$C_{a,\la}$}
\rput(49,-7.5){$C_{a,\la}^\circ$}\psline{->}(48,-7.8)(45,-9)\psline{->}(50,-7.8)(52,-9)
\rput(49.5,-10){\tiny{$C_{a,\la}^\tn{neck}$}}
\psarc(52,-10){2}{-155}{155}  
\psarc(49.83,-8.8){.5}{-127}{-44}
\psarc(49.83,-11.2){.5}{44}{127}
\psellipticarc(48.2,-10)(.6,1.5){88}{272}
\psellipticarc(50.7,-10)(.7,1.5){-90}{90}
\psarc(45,-13){3.16}{70}{110}
\psarc(45,-7){3.16}{245}{295}

\end{pspicture}
\caption{Neck-Core decomposition of a smoothing.}
\label{fig:neck-decomp}
\end{figure}
\noindent
On the other hand, by (\ref{equ:cCcircv}), for suitable choices of $U_{e,v}$ in (\ref{equ:Uev}), the smooth trivializations $\{\varphi_v\}_{v\in V_\Gamma}$ of $\{\cC_v\}_{v\in V_\Gamma}$ give a smooth trivialization 
$$
\varphi=(\varphi_v)_{v\in V_\Gamma}\colon \big(\cC^\circ \lra \cB,\vec{\mfz}, b\big) \lra \big(\Si^\circ \!\times\!\cB \lra \cB,\vec{z}, b\big)
$$
of the core region 
$$
\cC^\circ\!=\!\bigcup_{v\in V_\Gamma} \cC_v^\circ
$$ 
of the universal family $\cC$. Thus
$$
(\Sigma_{a_v}^\circ,\vec{\mfz}_v(a_v))\!\cong \!(\Sigma^\circ_v,\vec{z}_v)
$$
as marked smooth surfaces; i.e. the smooth structure of $C_{a,\la}^\circ$ does not depend on the choice of $\la$ or $a$ and thus we may identify the underlying smooth domains with $\Si^\circ$ and write $(\Si_{a,\la}^\circ,\mfj_{a,\la},\vec{\mfz}(a,\la))$ as\footnote{By the first holomorphic identification above, restricted to the core region, the complex structure is determined only by $a$.} $(\Si^\circ,\mfj_{\varphi,a,\la},\vec{z})$, with $(a,\la)\!\in\!\cB$. The following construction depends on the choice of smooth trivialization $\varphi$.

\noindent
As a generalization of (\ref{Bfep_e}), and in order to extend (\ref{Bepf0_e}) over $\cC$, define $B_\ep(f)$ to be the set of $(\ell,p)$-smooth maps $u'\colon C_{a,\la}\!\lra\!X$, with $(a,\la)\!\in\! \cB$, such that 
\bEnum
\item\label{l:smooth-part} $u'|_{C_{a,\la}^\circ}\!=\!\tn{exp}_{u}(\xi)\circ\varphi$ for  some
$\xi\!\in\! \Gamma^{\ell,p}(\Si^\circ,u^* TX)$ with $|\xi|\!<\!\ep$,
\item\label{l:ball-inequlity}
and 
\bEqu{equ:ball-inequlity}
\tn{dist}(u(q_e),u'|_{C^e_{a,\la}})\!<\!\ep,
\eEqu 
where $C^e_{a,\la}\!=\!C^\neck_{a,\la}\cap \cC_e$
and (\ref{equ:ball-inequlity}) means that the image of $u'$ restricted to the neck region $C^e_{a,\la}$ lies in a ball of radius $\ep$ around $u(q_e)$.
\eEnum
Let $B_\ep^{\tn{hol}}(f)\!\subset\!B_\ep(f)$ be the subset where in addition
\bEnum
\setcounter{enumi}{2}
\item $\dbar u'\! \equiv \! 0$ on $C_{a,\la}^\neck$.
\eEnum 
The automorphism group $\aut(f)$ acts on both $B_\ep(f)$ and $B_\ep^{\tn{hol}}(f)$. Since the obstruction space $E_f$ is supported away from the nodes, $B_\ep^{\tn{hol}}(f)$ includes $B_\ep(f;0)$ (see (\ref{Bepf0_e})).

\noindent
\bRem{rmk:Gromov2}
For sufficiently small $\ep$ (etc.\footnote{There are various other parameters in the construction of $\cC$ that 
we want them to be sufficiently small as well.}), given $(a,\la)\!\in\! \cB$ and $\xi\!\in\! \Gamma^{\ell,p}(\Si^\circ,u^* TX)$, there is at most one extension $u'\colon C_{a,\la}\!\lra\!X$ such that $(u',C_{a,\la})\!\in\!B_\ep^{\tn{hol}}(f)$ and $u'|_{C_{a,\la}^\circ}\!=\!\tn{exp}_{u}(\xi)\circ\varphi$.
In other words, there is at most one pseudoholomorphic extension of a given map on the core region to the neck region.
This is a consequence of Aronszajn's Unique Continuation Theorem \cite[Theorem 2.3.4]{MS2004} applied to each neck, as every neck is isomorphic to an annulus in $\C$.
Therefore, we get an $\aut(f)$-equivariant embedding of $B_\ep^{\tn{hol}}(f)$ into the Banach manifold\footnote{Since $\Si^\circ$ is open, by abuse of notation in (\ref{equ:BB}), $W^{\ell,p}(\Si^\circ,X)$ means the space of maps on the closure $\tn{cl}(\Si^\circ)$. }
\bEqu{equ:BB}
W^{\ell,p}(\Si^\circ,X)\times \cB\cong W^{\ell,p}_{\varphi}(\cC^\circ,X).
\eEqu
While this embedding can be very wild, similarly to Remark~\ref{rem:varphi-dependence}, it defines a topology on $B_\ep^{\tn{hol}}(f)$ that restricted to the subset of pseudoholomorphic maps gives an open neighborhood of $[f]$ in $\ov\cM_{g,k}(X,A)$.  
The resulting topology coincides with Gromov topology.
This extends the notion of $\ep$-close maps and $\ep$-convergence in Remark~\ref{rem:varphi-dependence} to a neighborhood of maps with stable nodal domain in the moduli space. From this characterization of the Gromov topology we get the last claim before Theorem~\ref{thm:gluing-thm-main} below. 
\eRem

\noindent
Similarly to (\ref{EVphitilde_e}), via trivialization (\ref{tildeVvsVf0_e}), we obtain an $\aut(f)$-equivariant obstruction bundle 
\bEqu{ObsEvarphif_e}
E_{\varphi,f}\lra B_\ep(f)
\eEqu
whose fiber over each point is obtained from parallel translations of $E_v$ over the corresponding component of the core region. In other words, for $f'\!=\!(u',C_{a,\la})\!\in\!B_\ep(f)$ with
$$
u'|_{C_{a_v}^\circ}\!=\!\tn{exp}_{u_v}(\xi_v) |_{\Si_v^\circ} \circ \varphi|_{C^\circ_{a_v}},
$$
and $P_\varphi\!=\!(P_{\varphi_v})_{v\in V_\Gamma}$ as in (\ref{ELocalChart_e}), let
\bEqu{equ:translated-space}
E_{\varphi,f}|_{f'}\!=\!P_\varphi\big(E_f\big):=\bigoplus_{v\in V_\gamma} P_{\varphi_v}\big(E_v\big)\!\subset \! \Gamma^{\ell-1,p}\big(\Si_{a,\la},\Om^{0,1}_{\Si_{a,\la},\mfj_{a,\la}}\otimes_\C (u')^*TX\big),
\eEqu
where the right-hand side is defined as in (\ref{target_e}) and the direct sum extends by zero to the neck region. The latter is well-defined by the first condition in definiton~\ref{def:obs-space-2}. Restricted to $V_{f;0}\!\subset\!B_\ep(f)$, $E_{\varphi,f}$ coincides with the obstruction budnle $U_{f;0}$ in (\ref{equ:bundle-0regular}).
Similarly to the paragraph before Remark~\ref{rem:family-with-unmarked2}, by \cite[Proposition 17.22]{FOOO-detail}, if $\cB$ and $\ep$ are sufficiently small, for every $f'\!\in\! B_\ep(f)$ with $\dbar u'\!=\!0$, the subspace (\ref{equ:translated-space}) is an obstruction space for $f'$ supported away from the neck region. 
Define 
\bEqu{OverAllChart_e}
\cU_f=(\pr_f\colon U_f\lra {V}_f, \aut(f), s_f,\psi_f)
\eEqu
as in (\ref{EVtilde_e}), i.e.
\bEqu{nodaltildechart_e}
{V}_f= \dbar^{-1}(E_{\varphi,f})\subset B_\ep^{\tn{hol}}(f), \quad {U}_f=E_{\varphi,f}|_{{V}_f},\quad \quad {s}_f=\dbar|_{{V}_f},
\eEqu
the action of $\aut(f)$ is the  restriction of the action of $\aut(f)$ to the zero set of $\aut(f)$-equivariant section $\dbar$, and $\psi_f$ sends $f'$ to $[f']$.
The restriction of $\cU_f$ to  $\la\!=\!0\!:=\!(0)_{e\in E_\Gamma}\! \in\! \De^{E_\gamma}$ coincides with the smooth fiber product chart $\cU_{f;0}$.  
The inclusion $\dbar^{-1}(E_{\varphi,f})\!\subset\! B_\ep^{\tn{hol}}(f)$ follows from (\ref{equ:translated-space}) and the third assumption in the defintion of $B_\ep(f)$ above. 
In the following, for $f'\!=\!(u',C_{a,\la})\in B_\ep(f)$, we write 
\bEqu{equ:nodal-dbar}
\dbar u'\! \equiv \! 0 ~~\tn{modulo}~~E_f\quad \tn{or}\quad\dbar u'\! \overset{E_f}{\equiv} \! 0,
\eEqu
whenever $f'\!\in\! {V}_f$. 

\bRem{rmk:not-Banach}
Recall from Remark~\ref{rmk:Gromov2} that unlike in Section~\ref{sec:canonical-smooth}, the $\ep$-neighborhood $B_\ep(f)$ defined in Page~\ref{l:ball-inequlity} does not have the structure of a Banach manifold; it just wildly embeds inside a Banach space. Therefore, we can not directly apply the Implicit Function Theorem to conclude that $V_f$ is a smooth manifold of the expected dimension. The following theorem gives us an $\aut(f)$-equivariant product (smooth) structure on (\ref{OverAllChart_e}).
\eRem

\begin{theorem}[{Gluing Theorem}]\label{thm:gluing-thm-main}
For $\{\cB_v\}$ and $\ep$ sufficiently small, and $\ve\!<\!\!<\!\!\ep$ and $\de\!\!<\!\!<\!\ve$\footnote{The construction of $\cB$ depends on the parameters $\ve$ and $\de$ and the definition of $B_\ep(f)$ depends additionally on $\ep$.} in the construction of $\cC$, there exists a ``natural'' $\aut(f)$-equivariant continuous one-to-one map
\bEqu{equ:gluing-map}
\mf{gl}\colon V_{f;0}\times \De^{E_\Gamma} \lra V_f
\eEqu
such that the restriction of $\mf{gl}$ to $V_{f;0}\!\subset\!V_f$ is the inclusion map and the action of $\aut(f)\!\subset\!\aut(C)$ on $\De^{\Gamma_E}$ is given by (\ref{equ:gluing-action}).
\eThm

\noindent
In simple words, Theorem~\ref{equ:gluing-map} states that for every $(u',C_a)\!\in\!\wt{V}_{f;0}$ and every deformation $C_{a,\la}$ of the domain in the normal direction to the stratum $\cM_{\Gamma}\!\subset\!\cM_{g,k}$, there exists a systematic way of deforming $(u',C_a)$ to $(u'_\la,C_{a,\la})$ satisfying $\dbar u'_\la\!\stackrel{E_f}{\equiv}\! 0$. Here the word ``natural" means that $u'_\la$ is constructed from $u'$ via some standard procedure that depends on some choices, however modulo those choices, the map $\mf{gl}$ is uniquely determined. More precisely, the construction of gluing map is by constructing first a pre-gluing map  \cite[Definition 5.1]{FOOO-detail3} via a set of appropriate cut-off functions \cite[(4.1)-(4.6)]{FOOO-detail3} and then deforming that by an \textit{alternating method} \cite[Section 5II]{FOOO-detail3} into an actual solution of (\ref{equ:nodal-dbar}). Therefore, the construction of $\mf{gl}$ depends on the choice of holomorphic coordinates (\ref{equ:holo-functions}), cut-off functions, etc..
Nevertheless, the ``weak" smooth structure on $V_f$ given by (\ref{equ:T-weak}) bellow is independent of these choices; see \cite[Section 8]{FOOO-detail3}.

\noindent
By (\ref{equ:translated-space}), (\ref{equ:fp-intersection}), and Theorem~\ref{thm:gluing-thm-main}, the corresponding trivialization $(\wt{\pr}_f\colon\!\wt{U}_f\!\lra\! \wt{V}_f, \wt{s}_f)$ of $(\pr_f\colon\!U_f\!\lra\! V_f, s_f)$ is of the form 
\bEqu{UfVf}
\aligned
&\wt{V}_f=\wt{V}_{f;0}\times \De^{E_\Gamma}, \quad \wt{U}_f= \wt{V}_{f}\times E_f=\wt{V}_{f;0}\times \De^{E_\Gamma}\times E_f= \wt{U}_{f;0}\times \De^{E_\Gamma},\\
&\tn{and}\quad \wt{s}_f(\ze,a,\la)=  \tn{P}_{u,\ze',a}^{-1}(\dbar_{J,\mfj_{\varphi,a}} u')\in E_f  , 
\endaligned
\eEqu
such that $\mf{gl}(\ze,a,\la)\!=\!(u',C_{a,\la})$ with 
$$
u'|_{C_{a,\la}^\circ}=u'|_{C_{a}^\circ}=u_{\ze'}\quad \tn{for some}\quad \ze'\!=\!\ze'(\ze,a,\la)\!\in\!\Gamma^{\ell,p}(\Si^\circ, u^*TX)
$$ 
and $\tn{P}_{u,\ze',a}$ is (the restriction to the core region of)  the parallel translation map in (\ref{family-P_e}). Therefore, Theorem~\ref{thm:gluing-thm-main} defines an $\aut(f)$-equivariant product smooth structure on $\pr_f\colon U_f\lra V_f$.

\noindent
In Definition~\ref{def:kur-chart}.\ref{l:K-map}, we require the Kuranishi map to be smooth. While $s_{f;0}$  in the gluing theorem is smooth, it is not clear whether the extended Kuranishi map $s_f$ remains smooth or not (with respect to the product smooth structure on $V_f$). For this reason, and also to make coordinate change maps smooth, we will equip $V_f$ and $U_f$ with a weaker product smooth structure in the normal direction. 

\noindent
Assuming $\de\!<\!1$, define
\bEqu{equ:T-weak}
T\colon \De\lra \C, \quad T(z)\!=\!\frac{-1}{\tn{log}(r)}\tn{e}^{\mfi \theta} \quad\forall z\!=\!r\tn{e}^{\mfi \theta}\!\in\!\C^*,~~T(0)\!=\!0.
\eEqu
We define the \textbf{weak} smooth  structure $\De_\tn{weak}$ on $\De$ to be the smooth structure given by the coordinate chart $T(z)$. 

\begin{theorem}[{Smoothness Theorem}]\label{thm:gluing-thm-smooth}
With respect to the ``weak'' smooth structures $V_f^{\tn{weak}}$ and $U_f^{\tn{weak}}$ given by $\De_\tn{weak}$ instead of $\De$ in (\ref{UfVf}), the Kuranishi map of 
Theorem~\ref{thm:gluing-thm-main} is smooth.
\eThm

\noindent
We do not go into the proofs of Theorems~\ref{thm:gluing-thm-main} and \ref{thm:gluing-thm-smooth} in this article and refer the interested reader to \cite{FOOO-detail3} which is dedicated to the proof of gluing theorem and its properties (including the case of bordered pseudoholomorphic maps). In the rest of this section, we will only recall  some of the main steps and unique features of the proofs. 

\bRem{rem:smoothness-clarification}
In the upcoming sections, any argument involving a smooth structure on the Kuranishi chart (\ref{OverAllChart_e}), e.g. transversality of certain evaluation maps in Section~\ref{sec:natural}, refers to the weak smooth  structure provided by Theorem~\ref{thm:gluing-thm-smooth}. Since this is the only way we put a smooth structure on (\ref{OverAllChart_e}), we will omit the superscript ``weak" from the notation for simplicity.
\eRem

\noindent
The proof of Theorem~\ref{thm:gluing-thm-main} in \cite{FO, FOOO-detail} and with more details in  \cite{FOOO-detail3} uses weighted Sobolev spaces  \cite[Definition 3.4]{FOOO-detail3} with $\ell\!>\!\!>\!1$ and $p\!=\!2$. The choice of $p\!=\!2$ allows us to define an inner product \cite[(4.11)]{FOOO-detail3} on the source and target of $\tn{D}_u\dbar$. While the overall process is similar to the classical proof of \cite[Chapter 10]{MS2004}, the weighted norms and cut-off functions used in the proof of the gluing theorem in \cite{FOOO-detail3} are different than those of \cite[Section 10.3]{MS2004}; see \cite[Remark 6.17]{FOOO-detail3}. They allow for a better control of the gluing process near the nodes which is essential for establishing Theorem~\ref{thm:gluing-thm-smooth}.
Theorem~\ref{thm:gluing-thm-main}  in \cite{FOOO-detail3} is spread among few statements which we have gathered into one theorem here. The case where the gluing parameter $\la$ and the complex structure of the domain is fixed is covered in \cite[Theorem 3.13]{FOOO-detail3}. The injectivity and surjectivity of the gluing map is discussed in \cite[Section 7]{FOOO-detail3}.  The argument uses the alternative description of the Gromov topology via the notion of $\ep$-close maps or $\ep$-convergence, see Remarks~\ref{rem:varphi-dependence} and~\ref{rmk:Gromov2}. The necessary adjustments for the case where the complex structure of the domain and the gluing parameter are allowed to change are discussed in \cite[Section 8.1]{FOOO-detail3}. Theorem~\ref{thm:gluing-thm-smooth} is a consequence of Exponential Decay Theorem \cite[Theorem 6.4]{FOOO-detail3} in simple words; see \cite[Proposition 8.31]{FOOO-detail3}.

\bRem{rmk:cylindrical-coordinates}
In \cite{FO} and its follow-ups, the authors make use of cylindrical coordinates 
$$
(-s_{e,v}+\mfi \theta_{e,v})\!=\!\log(w_{e,v})\!\in\!(-\infty,\log(\ve)]\!\times\! S^1
$$ 
at the nodes. In these coordinates, Equation~\ref{equ:smoothing-nodes} takes the form
\bEqu{equ:log-smoothing}
s_{e,v_1}\!+\!s_{e,v_2}\!=\!R_{e}, \quad \theta_{e,v_1}\!+\!\theta_{e,v_2}\!=\!\theta_e,\quad \tn{where}~~\log(\la_e)\!=\!-R_e\!+\!\mfi\theta_e.
\eEqu
Geometrically, it corresponds to removing $(-\infty,-R_e)\!\times\! S^1$ from $(-\infty,\log(\ve)]\!\times\! S^1$ at each node and identifying the remaining finite length cylinders $[-R_e,\log(\ve)]\!\times\! S^1$ according to (\ref{equ:log-smoothing}).
It is easier to express the weighted norms, cut-off functions, and the gluing process in terms of the cylindrical coordinates. In cylindrical coordinates, the weak smooth structure is given by $T(\la_e)\!=\!\tn{e}^{\mfi\theta}/R_e$. Since we do not go into the details of the gluing theorem in this article, we will stay with the usual holomorphic coordinates for simplicity. 
\eRem

\bRem{rem:embedding2}
Similarly to Remark~\ref{rem:embedding}, if $E'_f\!\subset\!E_f$, then the (germ of) Kuranishi chart  $\cU'_f$ in (\ref{OverAllChart_e}) corresponding to $E'_f$ embeds in the Kuranishi chart  $\cU_f$ corresponding to $E_f$ (defined via the same trivialization $\varphi$) as in Definition~\ref{def:intersection}. This embedding also satisfies the tangent bundle condition of Definition~\ref{def:Kuranishi-tangent}.
\eRem

\subsection{Case of un-stable domain}\label{sec:natural}

In this section, we extend the construction of Section~\ref{sec:canonical-nodal} to stable maps with un-stable domain.
This is done by adding an invariant set of extra marked points to the domain and intersecting with divisorial slices on the image (at the end) to reduce the dimension back to the expected dimension. In the light of Remarks~\ref{rmk:Gromov2} and~\ref{rem:varphi-dependence}, this stabilization process also allows us to obtain a rather different description of the Gromov topology.

\bDef{def:stabilizer}
Given a stable map $f\!=\!(u,\Si,\mfj,\vec{z})$, an un-ordered $\aut(f)$-invariant set of marked points ${w}\!\subset\! \Si$ disjoint from $\vec{z}$ and the nodes of $\Si$
is called a \textbf{stabilizing} set if  
\bEqu{equ:stab-data}
 \aut (\Si,\mfj,\vec{z},{w})\!=\aut(f)\quad\tn{and}\quad \nd u(w^i)\!\neq\!0\qquad \forall w^i\!\in\!{w},
\eEqu
where $ \aut (\Si,\mfj,\vec{z},{w})$ is the automorphism group defined in (\ref{equ:un-order-group}).
\eDef

\noindent
Since the restriction of $u$ to every unstable component is non-trivial and the set of regular values is open dense, every stable map admits plenty of stabilizing sets. 

\begin{example}\label{exa:1}
Let $(\Si,\mfj)\!\cong\! \P^1\!=\!\C\cup\{\infty\}$, $\vec{z}=\emptyset$, and 
$$
u\colon \P^1 \lra X=\P^1, \quad u(z)\!=\!z^2;
$$
thus, $[f\!=\!(u,\P^1)]\in \cM_{0,0}(\P^1,[2])$. 
Then any quadruple of points ${w}=\{a,b,-a,-b\}$, with $a,b\!\in\! \C^*$ and $a\!\neq\!\pm b$, is a stabilizing set. These are the minimal stabilizing sets.
\end{example}

\noindent
Let $f\!=\!(u,\Si,\mfj,\vec{z})$ be a stable $J$-holomorphic 
map and ${w}$ be a stabilizing set for $f$. 
By the second assumption in (\ref{equ:stab-data}), $\aut(f)$ does not fix any of $w^i$ and ${w}$ decomposes to finitely many $\aut(f)$-orbits $[w^i]$, each of which is isomorphic to $\aut(f)$; i.e. if $w\!\neq\!\emptyset$,  then the first assumption in (\ref{equ:stab-data}) follows from the second one.
For every orbit $[w^i]\!\in\! {w}/\aut(f)$, let $H_{[w^i]}$ be a sufficiently small codimension $2$ submanifold of $X$ which $J$-orthogonally\footnote{i.e. $T_{u(w^i)}H_{[w^i]}$ is the $J$-orthogonal complement of $T_{u(w^i)}u(\Si)$. } intersects $u(\Si)$ only at $u(w^i)$. A set of \textbf{slicing divisors} is a collection of disjoint submanifolds 
$$
H\!=\!\{H_{[w^i]}\}_{[w^i]\in {w}/\aut(f)}
$$ 
with this property. In what follows, by a \textbf{stabilizing pair} or \textbf{stabilizing data} we mean a pair $(w,H)$ consisting of a stabilizing set of points together with a corresponding set of slicing divisors for a fixed $J$-holomorphic $f$.

\noindent
Given a stable $J$-holomorphic map $f\!=\!(u,\Si,\mfj,\vec{z})$ and a stabilizing pair $(w,H)$, we construct a (class of) primary Kuranishi chart centered at $[f]$ in the following way.  In the following argument let $H_{w^i}\!:=\!H_{[w^i]}$. Choose an obstruction space $E_{f}$ for $f$ as in Definition~\ref{def:obs-space-2}. Such $E_f$ can also be thought\footnote{Note that the support of $E_f$ is allowed to have overlap with the stabilizing points. It should only be away from nodal points.} of as an obstruction space for $f^+\!=(u,\Si,\mfj,\vec{z},w)$. Take a (standard) universal family $\cC_+$ around $C^+\!=\!(\Si,\mfj,\vec{z},w)$ and  let 
\bEqu{equ:pre-slice}
\cU^+= (\pr^+\colon \!U^+\!\lra\! V^+, \aut(f^+)\!=\!\aut(f), s^+,\psi^+)
\eEqu
be the resulting primary Kuranishi chart constructed in Section~\ref{sec:canonical-nodal} around $[f^+]\!\in\!\ov\cM_{g,k+\ell}(X,A)$. Recall that if $\Si$ is nodal, the smooth structure on $V^+$ considered in the following argument is the weak smooth structure of Theorem~\ref{thm:gluing-thm-smooth}; see Remark~\ref {rem:smoothness-clarification}. 

\bRem{rmk:orderf+}
The map $f^+$ comes with an un-ordered set of marked points. Therefore, $[f^+]$ is not actually an element of $\ov\cM_{g,k+\ell}(X,A)$, unless we let the notation $\ov\cM_{g,k+\ell}(X,A)$ denote the moduli of genus $g$ maps with $k$ ordered and $\ell$ un-ordered marked points. Fix an ordering $\vec{w}$ on $w$. Then the Kuranishi chart $\cU^+$ is a Kuranishi chart for $[u,\Si,\mfj,\vec{z}\cup \vec{w}]\!\in\!\ov\cM_{g,k+\ell}(X,A)$ with a possibly larger group action. 
\eRem

\noindent
For sufficiently small $V^+$, the evaluation map at $w$-points\footnote{More precisely, we again need to fix an ordering on $w$ but the conclusion is independent of the choice of that.}, 
\bEqu{equ:ev_w}
\ev_w\colon V^+\lra X^w\!\equiv\!\prod_{w^i\in w} X,
\eEqu
remains transverse to $(H_{w^i})_{w^i\in w}\!\subset\!X^w$. By abuse of notation, we will also denote the last tuple by $H$. The choice of weak smooth structure has no effect on this transversality.
In order to cancel out the effect of the added points to the dimension of $\cU^+$, we replace $V^+$ with
$$
V\!=\!\ev_w^{-1}(H).
$$
Then $V$ is $\aut(f)$-invariant and  
\bEqu{equ:sliced-chart}
\cU\equiv (\pr\!=\!\pr^+|_{U}\colon U\!=\!U^+|_V\!\lra\! V, \aut(f), s\!=\!s^+|_V ,\psi),
\eEqu
where $\psi$ is the composition of 
$$
\psi^+|_{(V\cap (s^+)^{-1}(0))}\colon (V\cap (s^+)^{-1}(0))\!\lra\! \ov\cM_{g,k+\ell}(X,A)
$$ 
with the forgetful map of removing $w$-points $\ov\cM_{g,k+\ell}(X,A) \!\lra\! \ov\cM_{g,k}(X,A)$ as in (\ref{equ:pi-I}),
defines a natural primary Kuranishi chart centered at $[f]$. 

\noindent
Since the local slicing divisors $\{H_{w^i}\}$ are intersecting $J$-orthogonally, the normal bundles 
$$
N_{H_{w^i}}X=\frac{TX|_{H_{w^i}}}{TH_{w^i}}\cong TH_{w^i}^{\perp \om}
$$ 
are canonically oriented.  We equip $TV$ with the natural orientation given by the natural orientation of $V^+$ and the exact sequence
$$
0\lra TV \lra TV^+|_{V} \stackrel{[\nd \ev_{w}]}{\xrightarrow{\hspace*{1.5cm}}} \bigoplus_{w^i\in w} N_{H_{w^i}}X \lra 0\;.
$$
Since, by assumption, $u$ is interesting every $H_{w^i}$  positively, the induced linear map
$$
[\nd u]\colon T_{w^i}\Si\lra N_{H_{w^i}}X
$$ 
is an orientation preserving isomorphism. The same holds for every other map in $V$. This implies that our orientation scheme is compatible with the forgetful maps that arise in the following arguments.

\bRem{rem:adding-more-points}
Let $f\!=\!(u,C\!=\!(\Si,\mfj,\vec{z}))$ be a $J$-holomorphic map with smooth domain, and $(w_0,H_0)$ and $(w_1,H_1)$ be two stabilizing sets with $w_0\!\subsetneq\! w_1$. Let 
$$
(\cC_0\lra \cB_0,~~\vec{\mfz},\mf{w}_0\colon \cB_0\lra \cC_0),\quad\tn{and}\quad (\cC_1\lra\cB_1,~~\vec{\mfz},\mf{w}_1\colon \cB_1\lra \cC_1)
$$ 
be universal families around $C_0\!=\!(\Si,\mfj,\vec{z},w_0)$ and $C_1\!=\!(\Si,\mfj,\vec{z},w_1)$, respectively. After possibly restricting to smaller subfamilies, the process of removing the extra marked points  $\mf{w}_1\!\setminus\!\mf{w}_0$ (in other words, extra sections)  gives us a holomorphic $\aut(f)$-equivariant forgetful map 
$$
(\pi^\cC,\pi^\cB)\colon (\cC_1,\cB_1)\lra (\cC_0,\cB_0),
$$
and thus a projection map
$$
\pi_{1,0}\colon\!W^{\ell,p}(\cC_1,X)\lra W^{\ell,p}(\cC_0,X).
$$
Let $\varphi_0$ and $\varphi_1$ be smooth trivializations for $\cC_0$ and $\cC_1$ as in (\ref{equ:smooth-trivial-phi}), respectively.
Then, with notation as in (\ref{equ:Bspace-triv}), the resulting projection map $\pi_{\varphi_1,\varphi_0}$ given by 
$$
\xymatrix{
W^{\ell,p}(C_1,X)\times \cB_1\ar[d]^{\pi_{\varphi_1,\varphi_0}}\ar[rr]^{T_{\varphi_1}}&& W^{\ell,p}(\cC_1,X)  \ar[d]^{\pi_{1,0}} \\
W^{\ell,p}(C_0,X)\times \cB_0 		\ar[rr]^{T_{\varphi_0}}	                              && W^{\ell,p}(\cC_0,X) 
}
$$
can not be identity on the first component. In fact, for each $a\!\in\!\cB_{0}$, 
$$
\pi_{\varphi_1,\varphi_0}\colon W^{\ell,p}(\Si,X)\!\times\! (\pi^\cB)^{-1}(a)\lra W^{\ell,p}(\Si,X)\!\times\! \{a\}
$$
corresponds to twisting with the non-constant $(\pi^\cB)^{-1}(a)$-family of diffeomorphisms 
$$
\rho_{a';a}\colon \varphi_0\circ\pi^\cC\circ\varphi_1^{-1}|_{\Si\times\{a'\}}\colon \Si\times\{a'\}\lra \Si\times \{a\}\quad \forall a'\!\in\! (\pi^\cB)^{-1}(a),
$$ 
each of which fixes $\vec{z}$ and ${w}_0\!\subset\!w_1$ but possibly maps $w_1\!\setminus\!w_0$ to a variation of these points.
In other words, $\pi_{\varphi_1,\varphi_0}$  is of the form 
$$
(u',a')\lra (u'\!\circ\!\rho_{a';a}^{-1},a\!=\!\pi^\cB(a'))\quad \forall (u',a')\!\in\!W^{\ell,p}(\Si,X)\!\times\!\cB_1.
$$
Since the map depends on a non-trivial reparametrization of the domain,  the forgetful map
$$
\pi_{1,0}\colon W^{\ell,p}_{\varphi_1}(\cC_1,X)\lra W^{\ell,p}_{\varphi_0}(\cC_0,X)
$$
is continuous but not smooth. Similarly to (\ref{equ:Frho}), there will be a loss of differentiability caused by the derivative of the reparametrization diffeomorphism. However, if $E_0\!\subset\!E^{\ell-1,p}_{\varphi_0}(\cC_0,X)$ is a ``nice" smooth finite dimensional sub-bundle where each fiber is made of smooth sections (e.g. an obstruction bundle $E_{\varphi_0,f_0}$ as in (\ref{EVphitilde_e})), the pull back bundle $E_1\!=\!\pi_{1,0}^*E_0$ still has the structure of a smooth finite dimensional sub-bundle of $E^{\ell-1,p}_{\varphi_1}(\cC_1,X)$.  We will  extensively use this argument in the next section to define induced Kuranishi charts.
\eRem

\bRem{rmk:No-Comparison}
Proceeding with the notation of Remark~\ref{rem:adding-more-points}, let 
$$f_0\!=\!(u,C\!=\!(\Si,\mfj,\vec{z}, w_0))\quad\tn{and}\quad f_1\!=\!(u,C\!=\!(\Si,\mfj,\vec{z}),w_1).
$$ 
Given an obstruction space $E_f$ for $f$, since $E_f$ can also be thought of as an obstruction space for $f_0$ and $f_1$, let 
$$
E_{f_1,\varphi_1}\!\lra\!B_\ep(f_1)\subset W^{\ell,p}_{\varphi_1}(\cC_1,X)\quad\tn{and}\quad E_{f_0,\varphi_0}\!\lra\!B_\ep(f_0)\subset W^{\ell,p}_{\varphi_0}(\cC_0,X)
$$ 
be the corresponding obstruction bundles of (\ref{EVphitilde_e}). Then in general, because of the reparametrization map $\rho$ in  Remark~\ref{rem:adding-more-points}, $\pi_{1,0}^*E_{f_0,\varphi_0}$ is different from $E_{f_1,\varphi_1}$ and there is no clear relation between the two. Therefore, there exists no natural map between the resulting manifolds $V_{f_1}(\varphi_1)$ and $V_{f_0}(\varphi_0)$ away from the zero set of the Kuranishi maps (i.e. the subset of actual $J$-holomorphic maps); see also Remark~\ref{rem:chart-varphi}. 
Nevertheless, restricted to the subsets $s_{f_1}^{-1}(0)$ and $s_{f_0}^{-1}(0)$ of actual $J$-holomorphic maps, the projection map
$$
\pi_{\varphi_1,\varphi_0}\colon (V_{f_1}(\varphi_1)\cap s_{f_1}^{-1}(0)) \lra (V_{f_0}(\varphi_1)\cap s_{f_0}^{-1}(0))
$$
is a continuous fiber bundle. After interesting with the slicing divisors, the resulting projection map
$$
\pi_{\varphi_1,\varphi_0}\colon (V_{f_1}(\varphi_1)\cap s_{f_1}^{-1}(0) \cap \tn{ev}_{w_1}^{-1}(H_1)) \lra (V_{f_0}(\varphi_1)\cap s_{f_0}^{-1}(0)\cap \tn{ev}_{w_0}^{-1}(H_0))
$$
is a homeomorphism (between sufficiently small neighborhoods  of $f_1$ and $f_0$).
Therefore, with the Kuranishi chart $\cU$ as in (\ref{equ:sliced-chart}), the topology defined via open sets $\psi(s^{-1}(0))$ in a neighborhood of $[f]\!\in\!\ov\cM_{g,k}(X,A)$ is independent of the choice of a stabilizing pair $(H,w)$ used to obtain the Kuranishi chart (\ref{equ:sliced-chart}).
This topology coincides with the Gromov topology and extends the notion of $\ep$-closed topology in Remarks~\ref{rem:varphi-dependence} 
and \ref{rmk:Gromov2} to maps with unstable domain.
\eRem

\noindent
Let us summarize the outcome of our construction so far in the following statement.
For every stable marked $J$-holomorphic map $f\!=\!(u,C\!=\!(\Si,\mfj,\vec{z}))$ with dual graph $\Gamma$, our construction of a primary Kuranishi chart $\cU_f\!=\!(\tn{pr}_f\colon U_f\!\lra\! V_f, G_f,s_f,\psi_f,)$ centered at $f$ depends (at most) on a set of \textbf{auxiliary data} $\tn{Aux}(f)$ consisting of 
\bEnum
\item\label{l:E} an obstruction space $E_f$, 
\item\label{l:w} a set of stabilizing points $w$,
\item\label{l:H} a set of slicing divisors $H$,
\item\label{l:cCv} a set of sufficiently small universal families $\cC^+_v\!\lra\!\cB^+_v$ for each smooth component $C^+_v\!=\!(\Si_v,\mfj_v,\vec{z}_v,w_v)$ of $C^+\!=\!(\Si,\mfj,\vec{z},w)$,
\item\label{l:varphiv} smooth trivializations $\varphi\!=\!\{\varphi_v\}_{v\in V_\Gamma}$ of $\{\cC^+_v\}_{v\in V_\Gamma}$, 
\item\label{l:cC} a standard extension of the universal families and trivializations in \ref{l:cCv} and \ref{l:varphiv} above to a universal family $\cC^+\!\lra\!\cB^+$ around $C^+$ and a trivialization $\varphi$ of the core region $(\cC^+)^\circ$ via holomorphic functions (\ref{equ:holo-functions}) and sufficiently small parameters $\de$ and $\ve$, 
\item\label{l:ep} and a sufficiently small positive number $\ep$ which via the exponentiation identification and trivialization $\varphi$ above determines an $\ep$-close neighborhood $B_\ep(f)$ of $f$ in the corresponding Banach space.
\eEnum
\begin{condition}\label{small-enough_c}
The families and parameters in $(4)(6)(7)$ are taken sufficiently small such that $\cC^+\!\lra\!\cB^+$ is a universal family around every $a\!\in\!\cB^+$, the argument of Remark~\ref{aut-limit_rmk} for every $f'\!\in\!(V_f\cap s_f^{-1}(0))$ holds, $E_{\varphi,f}$ restricted to every such $f'$ is an obstruction space, and every $f'\!\in\!V_f$ intersects $H$ transversely.
\end{condition}

\noindent
Our construction of primary charts readily extends to the $\mc{J}$-family case of Remark~\ref{rmk:family-version}. In this situation, the linearization of Cauchy-Riemann equation would also involve deformations of $J$. We refer to \cite[Section 3.1]{MS2004} for the necessary adjustments.

\subsection{Induced charts}\label{sec:natural-KUR}

Finally, in this section and the next, we construct a cobordism class of natural oriented Kuranishi structures on the moduli space of pseudoholomorphic maps. Every chart of a natural Kuranishi structure depends on a fixed finite collection of primary charts constructed in Sections~\ref{sec:canonical-smooth}, \ref{sec:canonical-nodal}, and \ref{sec:natural}. We call them \textbf{induced} charts. We will discuss coordinate change maps in the next section.
\bThm{thm:Nat-Kur}
Let $(X^{2n},\om)$ be a closed symplectic manifold, $J\!\in\! \cJ(X,\om)$ be an arbitrary compatible (or tame) almost complex structure, $A\!\in\! H_2(X,\Z)$, and $g,k\!\in\! \Z^{\geq 0}$. 
With the dimension $d$ as in (\ref{equ:exp-dim}), the construction of this section provides a class of natural oriented $d$-dimensional Kuranishi structures on $\ov{\cM}_{g,k}(X,A,J)$, every two of which are deformation equivalent as defined after Definition~\ref{def:cobo-Kur}. If $\om_0$ and $\om_1$ are isotopic symplectic structures on $X$, $J_0\!\in\! \cJ(X,\om_0)$, and  $J_1\!\in\! \cJ(X,\om_1)$, then every two of such Kuranishi structures on $\ov{\cM}_{g,k}(X,A,J_0)$ and $\ov{\cM}_{g,k}(X,A,J_1)$ are cobordant in the sense of Definition~\ref{def:cobo-Kur}.
\eThm

\noindent
For every $J$-holomorphic map $[f]\!\in\!\ov\cM_{g,k}(X,A)$, after fixing a choice of auxiliary data $\tn{Aux}(f)$, as listed at the end of last section, let 
$$
\cU_f\!=\! (\pr_f\colon U_f\lra V_f , \aut(f), \psi_f, s_f)
$$
be the resulting primary Kuranishi chart with the footprint $F_f\!\subset\!\ov\cM_{g,k}(X,A)$. For each $f$, fix a relatively compact connected open sub-set $F'_f\!\subset\!F_f$ around $[f]$.
Since $\ov\cM_{g,k}(X,A)$ is compact, there exists a finite set
$$
\cP\!=\!\{f_\alpha\}_{\alpha\in S}\!\subset\!\ov\cM_{g,k}(X,A)
$$
such that $\{F'_{f_\alpha}\}_{\al\in S}$ covers the moduli space; we call such $\cP$ a \textbf{primary collection}.

\noindent
Fix a primary collection $\cP$. For every $f_\al\!\in\! \cP$, let $\Gamma_\al$ denote its dual graph and
\bEqu{equ:aux-data}
\tn{Aux}(f_\al)\!=\!\big( E_\al, w_\al,H_\al,\varphi_\al,\cC_\al\lra\cB_\al, \ep_\al\big)
\eEqu
be the auxiliary data of Items~\ref{l:E}-\ref{l:H} and \ref{l:varphiv}-\ref{l:ep} in page \pageref{l:E} involved in the construction of primary chart $\cU_{f_\al}$.

\noindent
For every arbitrary $[f\!=\!(u,C\!=\!(\Si,\mfj,\vec{z}))]\!\in \!\ov\cM_{g,k}(X,A)$ define
\bEqu{equ:Sf}
S_{f}\!=\!\{\alpha\!\in\! S: [f]\!\in\! \tn{cl}(F'_{f_\alpha})\}.
\eEqu
Then 
\bEqu{equ:Sf'inSf}
S_{f'}\!\subset \!S_f,
\eEqu 
whenever $[f']$ is sufficiently close to $[f]$. Corresponding to every $\al\!\in\!S_f$ we get a unique set of points $w_{\al,f}\!\subset\!\Si$ at which $f$ transversely intersects  the slicing divisors $H_\al$.
For every $[f]$, we fix an arbitrary stabilizing set $(w_0,H_0)$ as in Section~\ref{sec:natural}, disjoint from $(w_{\al,f},H_\al)$, for all $\al\!\in\!S_f$.

\noindent
Before discussing the general case, which involves several components labeled by dual graph, let us consider the simpler case where $\Si$ is smooth.
Let
\bEqu{equ:various-add-ons}
\aligned
&f^{+0}\!=\!(u,C^{+0}\!=\!(\Si,\mfj,\vec{z},w_0)), \quad f^{+\alpha}\!=\!(u,C^{+\alpha}\!=\!(\Si,\mfj,\vec{z},w_{\al,f})),\\
&\tn{and}\quad f^{+0\alpha}\!=\!(u,C^{+0\alpha}\!=\!(\Si,\mfj,\vec{z},w_{0\al}))  ,\qquad\forall \al\!\in\!S_f,
\endaligned
\eEqu
where $w_{0\al}\!=\!w_0\!\cup\! w_{\al,f}$. 

\noindent
By Condition~\ref{small-enough_c}, for every $\al\!\in\!S_f$, the universal family $\cC_\al$ is also a universal family around $C^{+\alpha}$. Fix an $\al\!\in\!S_f$. Let $\cC_{0}$ and $\cC_{0\al}$ be sufficiently small universal families around $C^{+0}$ and $C^{+0\al}$, respectively. The process of removing $w_0$ and $w_{\al,f}$ from $w_{0\al}$ gives us forgetful maps 
\bEqu{equ:forget0al}
\pi_0\colon (\cC_{0\al},\cB_{0\al})\lra (\cC'_{\al},\cB'_{\al}) \quad\tn{and}\quad \pi_\al\colon (\cC_{0\al},\cB_{0\al})\lra (\cC_{0},\cB_0), 
\eEqu
respectively, where $(\cC'_{\al}\!\lra\! \cB_\al')$ is a sufficiently small sub-family\footnote{In which all the curves have smooth domain.} of $(\cC_{\al}\!\lra\!\cB_\al)$. Recall from Remark~\ref{aut-limit_rmk} that $\aut(f)\!\subset\!\aut(f_\al)$ acts on $(\cC'_{\al},\cB'_\al)$, $(\cC_0,\cB_0)$, and $(\cC_{0\al},\cB_{0\al})$.

\noindent
For a choice of 
\bEqu{equ:epsilon-choices}
0\!<\!\ep'_\al\!<\!\!<\ep_\al,\quad    0\!<\!\ep_0\!<\!\!<\ep'_\al,\quad\tn{and} \quad 0\!<\!\ep_{0\al}\!<\!\!<\ep_0,
\eEqu
let 
\bEqu{equ:Banach-0-al-0al}
B_{\ep'_\al}(f^{+\al}),\quad B_{\ep_0}(f^{+0}),\quad\tn{and} \quad B_{\ep_{0\al}}(f^{+0\al})
\eEqu
be the Banach neighborhoods of $(\ell,p)$-smooth maps within $\ep'_\al$, $\ep_0$, and $\ep_{0\al}$ distances from $f^{+\alpha}$, $f^{+0}$, and $f^{+0\alpha}$ as in (\ref{equ:Bepf}), respectively. 
Note that the definition of $\ep$-neighborhoods in (\ref{equ:Banach-0-al-0al}) as open subsets of 
$$
W^{\ell,p}(\cC'_\al,X), \quad W^{\ell,p}(\cC_0,X),\quad \tn{and}\quad W^{\ell,p}(\cC_{0\al},X),
$$
 depends on the choices of smooth trivializations $\varphi'_\al$, $\varphi_0$, and $\varphi_{0\al}$, for $\cC'_\al$, $\cC_0$, and $\cC_{0\al}$ as in (\ref{equ:smooth-trivial-phi}), respectively. However, the final construction of the \textit{induced} obstruction bundle over a neighborhood of $f^{+0}$ is independent of the choice of these trivializations.

\bRem{rem:varphivsvarphi}
Let $\Si_\al$ be the domain of $f_\al$. While $\cC'_\al$ is a subfamily of $\cC_\al$,  the trivializations 
$$
\varphi'_\al\colon\cC'_\al \lra  \Si\!\times\!\cB'_\al \quad\tn{and}\quad \varphi_\al\colon \cC_\al\lra \Si_\al\!\times\!\cB_\al ~(\tn{or}\footnote{If $\Si_\al$ is nodal.}~\Si^\circ_\al\!\times\!\cB_\al)
$$
have often no connection to each other. Even the topological type of $\Si_\al$ and $\Si$ could be different; $\Si_\al$ could be nodal.
\eRem

\noindent
By (\ref{equ:Sf}), $f^{+\al}$ is an element of ${V}_{f_\al}$. In the following, we will use $f^{+0\alpha}$ to relate a neighborhood of $f^{+\al}$ in $B_{\ep_\al}(f_\al)$ to some smaller neighborhood $B_{\ep'_0}(f^{+0})$ of $f^{+0}$ in $B_{\ep_0}(f^{+0})$. This relation will give us a finite dimensional  obstruction bundle 
\bEqu{equ:0al-bundle_e}
E_{0;\al}\lra B_{\ep'_0}(f^{+0}),
\eEqu
whose fiber at each point is obtained from parallel translation of $E_\al$ with respect to $\varphi_\al$ as in (\ref{equ:translated-space})\footnote{So it is independent of the choice of $\varphi'_\al$, $\varphi_0$, and $\varphi_{0\al}$.}. We will add up these obstruction bundles, for all $\al\!\in\!S_f$, and this gives us the total obstruction bundle $E_0$. Then the induced chart $(\pr_f\colon\!{U}_f\!\lra\! {V}_f, {s}_f)$ around $f$ is given by
$$
{V}_f=(\dbar^{-1}(E_0)\subset B_{\ep'_0}(f^{+0}))\cap \tn{ev}_{w_0}^{-1}(H_0),\quad {U}_f=E_0|_{V_f},\quad {s}_f=\dbar.
$$ 

\noindent
By the first assumption in (\ref{equ:epsilon-choices}) and Remark~\ref{rem:varphi-dependence}, $B_{\ep'_\al}(f^{+\al})$ continuously embeds in $B_{\ep_\al}(f_\al)$. Let
$$
E_{\varphi_\al,f^{+\al}}=E_{\varphi_\al,f_\al}|_{B_{\ep'_\al}(f^{+\al})}
$$
be the restriction of the obstruction bundle $E_{\varphi_\al,f_\al}$ defined in (\ref{ObsEvarphif_e}) to $B_{\ep'_\al}(f^{+\al})$. 
This is an $\aut(f)$-equivariant continuous vector bundle; see Remark~\ref{aut-limit_rmk}.
The forgetful maps of (\ref{equ:forget0al}) extend to similarly denoted continuous $\aut(f)$-equivariant forgetful maps
\bEqu{equ:forget0almaps}
\pi_0\colon B_{\ep_{0\al}}(f^{+0\al})\lra B_{\ep'_\al}(f^{+\al}) \quad\tn{and}
\quad \pi_\al\colon B_{\ep_{0\al}}(f^{+0\al})\lra B_{\ep_0}(f^{+0});
\eEqu
see Remark~\ref{rem:adding-more-points}.
Let
\bEqu{equ: cE_0al}
 E_{0\al}\!=\!\pi^*_{0}E_{\varphi_\al,f^{+\al}}\!\subset\!E^{\ell-1,p}(\cC_{0\al},X)|_{B_{\ep_{0\al}}(f^{+0\al})}
\eEqu
be the resulting $\aut(f)$-equivariant pull-back bundle. At this point one can pursue two different paths:
\bEnum
\item we may first consider the  fiber product space
$$
\times_{\substack{\al\in S_f}} \big(B_{\ep_{0\al}}(f^{+0\al})\lra B_{\ep_0}(f^{+0})\big),
$$
build a Kuranishi chart in this fiber product, and then apply all the slicing conditions $H_0$ and $\{H_\al\}_{\al \in S_f}$ at the end;
\item or we may use the slicing conditions $H_\al$ to define an $\aut(f)$-equivariant continuous section 
\bEqu{equ:sal}
s_\al\colon B_{\ep'_0}(f^{+0})\lra B_{\ep_{0\al}}(f^{+0\al})
\eEqu
of $\pi_\al$, for some $\ep'_0\!<\!\!<\!\ep_0$, pull-back  $E_{0;\al}$ to an $\aut(f)$-equivariant vector bundle $E_{0;\al}$ over $B_{\ep_0}(f^{+0})$,
 and add up the pull-back bundles $E_{0;\al}$, for all $\al\!\in\!S_f$, to construct the total obstruction $\aut(f)$-equivariant bundle $E_0$. 
 \eEnum
\noindent
We will pursue the second method. Let 
\bEqu{equ:B0al-slice}
B_{\ep_{0\al}}(f^{+0\al})_{\cap H_\al} \!=\! \ev_{w_\al}^{-1}(H_\al) \!\subset\!  B_{\ep_{0\al}}(f^{+0\al}),
\eEqu
where 
$$
\ev_{w_\al}\colon B_{\ep_{0\al}}(f^{+0\al})\lra X^{w_\al}
$$
is the evaluation map at the $w_\al$-points as in (\ref{equ:ev_w}). Since $\ev_{w_\al}$ remains transverse to $H_\al$ in a neighborhood of $f^{+0\al}$, the left-hand side of (\ref{equ:B0al-slice}) is the image of a contiunous section (\ref{equ:sal}) of the projection map $\pi_\al$ over some smaller neighborhood $B_{\ep'_0}(f^{+0})$; i.e. the section $s_\al$ in  (\ref{equ:sal}) obtained by adding the intersection points of $f^{+0}$ with $H_\al$. Finally, define 
$$
E_{0;\al}=s_\al^*E_{0\al} \!\subset\!E^{\ell-1,p}(\cC_0,X)|_{B_{\ep'_0}(f^{+0})}.
$$

\bLem{lem:smoothness-of-sal}
For any choice of smooth trivialization $\varphi_0$ for $\cC_0$, $E_{0;\al}$ is an $\aut(f)$-equivariant smooth sub-bundle of $E^{\ell-1,p}_{\varphi_0}(\cC_0,X)|_{B_{\ep'_0}(f^{+0})}$.
\eLem

\noindent
The proof this lemma relies on the fact that the obstruction space $E_\al$ is made of smooth sections and thus the change of trivialization map preserves the smoothness of the corresponding obstruction bundles; see Remark~\ref{rem:adding-more-points}. We skip the proof and refer to \cite[Section 34.4]{FOOO-detail} and \cite[Appendix A-B]{FOOO-detail3}; especially see \cite[Figure 14]{FOOO-detail3}. 

\noindent
Let
\bEqu{equ:formal-d-sum}
\pr_0\colon E_{0}\!=\!\bigoplus_{\al\in S_f}  E_{0;\al} \lra B_{\ep'_0}(f^{+0})
\eEqu
be the resulting induced direct sum $\aut(f)$-equivariant obstruction bundle.
The finite rank bundle $E_0$ is a domain-dependent version of the obstruction bundle in (\ref{EVphitilde_e}), with the difference that at each point, the $\al$-summand of the fiber is obtained from the obstruction space $E_\al$ by a parallel-transport procedure from $f_\al$ (instead of $f$) to that point.
\bRem{rmk:general-position}
We say a primary collection $\cP$ is \textbf{in general position} if
$$
E_{0;\al}|_{f^{+0}}\cap E_{0;\beta}|_{f^{+0}} \!=\!\emptyset \qquad \forall [f]\!\in\!\ov\cM_{g,k}(X,A),~ \al,\beta\!\in\!S_f;
$$
i.e. if the formal direct sum in (\ref{equ:formal-d-sum}) is isomorphic to the usual sum.
By \cite[Lemma 18.8]{FOOO-detail}, every $\ov\cM_{g,k}(X,A)$ admits a plethora of primary collections in general position. In fact, it is shown that any given primary collection $\cP$ can be slightly perturbed so that the resulting Kuranishi charts form a primary collection in general position. The proof is somewhat deliberate.
Assuming $\cP$ is in general position, one can simplify the following construction to some extent; see Remark~\ref{rem:GP-simplification} below.
However, for the sake of generality and avoiding the proof of this argument, we continue \textbf{without} this assumption.
\eRem

\noindent
Finally, we construct the {induced} chart at $f$ in the following way. Let
\bEqu{equ-induced-wtV}
V^{+0}\!=\!\big\{ (f',(\eta_\al)_{\al\in S_f})\in E_0 \colon \dbar f'-\sum_{\al \in S_f} \eta_\al=0\big\},\quad U^{+0}\!=\! (\pr_0^* E_0)|_{V^{+0}},
\eEqu
and 
\bEqu{equ:wts0}
 \quad {s}^{+0}\colon {V}^{+0}\lra {U}^{+0}, \quad {s}^{+0}(f',(\eta_\al)_{\al\in S_f})\!=\!\bigoplus_{\al \in S_f} \eta_\al \!\in\!E_0|_{f'}.
\eEqu
In order to clarify the notation in (\ref{equ-induced-wtV}), by $(f',(\eta_\al)_{\al\in S_f})\!\in\! E_0$ we mean the point corresponding to the vector $\oplus_{\al\in S_f} \eta_\al\!\in\!E_0|_{f'}$  in the total space of $E_0$. Then ${U}^{+0}$ is the pull-back of $E_0$ to the total space of $E_0$ restricted to ${V}^{+0}$. Restricted to the zero set of ${s}^{+0}$, we get the obvious footprint map 
\bEqu{equ:foot-plus}
\psi^{+0}\colon (s^{+0})^{-1}(0)\lra \ov\cM_{g,k+|w_0|}(X,A), \quad \psi^{+0}(f',(\eta_\al)_{\al\in S_f}) \!=\![f'];
\eEqu
i.e. $\cU^{+0}=(\pr^{+0}\colon\!U^{+0}\!\lra\! V^{+0},\aut(f),s^{+0},\psi^{+0})$ is a Kuranishi chart around $[f^{+0}]$; see Remark~\ref{rmk:orderf+} for the ordering of extra marked points. We will intersect with the $H_0$ slices to obtain a Kuranishi chart around $[f]$.

\bRem{rem:GP-simplification}
In the case of primary collections in general position, as in Remark~\ref{rmk:general-position}, we have $E_0\!=\!\sum_{\al\in S_f}  E_{0;\al}$ and (\ref{equ-induced-wtV}) simplifies to
$$
\aligned
&{V}^{+0}\!=\!\big\{ f'\in B_{\ep'_0}(f^{+0})\colon \dbar f' \in E_0\big\},\quad {U}^{+0}\!=\!E_0|_{{V}^{+0}},\\
 &\tn{and}\quad {s}^{+0}\colon {V}^{+0}\lra {U}^{+0}, \quad {s}^{+0}(f')\!=\!\dbar f'.
 \endaligned
$$
In this case, the manifold ${V}^{+0}$ in (\ref{equ-induced-wtV}) is the graph of ${s}^{+0}$ above.
\eRem
\noindent
Similarly to the case of primary charts in Section~\ref{sec:canonical-smooth}, by Implicit Function Theorem\footnote{With respect to the smooth structure determined by any choice of smooth trivialization $\varphi_0$ for $\cC_0$ in Lemma~\ref{lem:smoothness-of-sal}.} and Condition~\ref{small-enough_c}, for $\ep_0'$ and $\cC_0$ sufficiently small, 
\bItem
\item $\pr^{+0}\colon {U}^{+0}\lra {V}^{+0}$ s a smooth $\aut(f)$-equivariant vector bundle;
\item the evaluation map 
$$
\ev_{w_0}\colon {V}^{+0}\lra X^{w_0}
$$
is $\aut(f)$-invariant and transverse to $H_0$;
\item and the \textbf{induced chart} 
\bEqu{equ:sliced-chart2}
\cU_f\!\equiv\!(\pr_f\colon U_f\lra V_f , \aut(f), s_f,\psi_f),
\eEqu
where 
$$
V_f\!=\!\ev_{w_0}^{-1}(H_0)\cap {V}^{+0}, \quad U_f\!=\!{U}^{+0}|_{V_f}, \quad {s_f}\!=\!{s}^{+0}|_{V_f},
$$
and $\psi_f$ is the composition of $\psi^{+0}$ with the forgetful map of removing $w_0$-points\footnote{See Remark~\ref{rmk:orderf+}.}
$$\ov\cM_{g,k+|w_0|}(X,A) \!\lra\! \ov\cM_{g,k}(X,A),
$$
is a Kuranishi chart of the expected dimension (\ref{equ:exp-dim}) centered at $[f]$.
\eItem

\noindent 
We now return to the general case where the domain of $f$ is allowed to be nodal. 
Together with the modification above on each smooth component, the construction is similar to that of Section~\ref{sec:canonical-nodal}.
Let $\Gamma$ be the dual graphs of $f$. With notation as in the special case above, dual graphs of $f^{+0}$, $f^{+\al}$, and $f^{+0\al}$ have the same set of vertices and edges; they only have different set of flags. Since the following argument does not involve the flags, by abuse of notation, we will think of $\Gamma$ as the dual graph of all three maps. 

\noindent
For every $v\!\in\!V_{\Gamma}$ let $\cC_{0;v}\!\lra\!\cB_{0;v}$ be a sufficiently small universal family around the smooth component $C^{+0}_v\!=\!(\Si_{v},\mfj_v,\vec{z}_v,w_{0;v})$ of $C^{+0}$. 
In the product family
$$
\cC_0(\Gamma)\!=\!\prod_{v\in V_{\Gamma}} \cC_{0;v} \lra \cB_0(\Gamma)\!=\! \prod_{v\in V_{\Gamma}}\cB_{0;v},
$$
all the curves have the same fixed dual graph $\Gamma$.
Let 
$$
\cC_0 \lra \cB_0\!=\!\cB_0(\Gamma)\times \De^{E_{\Gamma}}
$$ 
be a standard universal family  around $C^{+0}$ obtained from $\cC_0(\Gamma)$ by smoothing the nodes as in (\ref{equ:smoothing-type-family}). Similarly to  Section~\ref{sec:canonical-nodal}, we first construct a fiber-product Kuranishi chart covering a neighborhood of $[f]$ in $\ov\cM_{g,k}(X,A)_{\Gamma}$, and then, via a gluing construction, extend it to a Kuranishi chart for $[f]$ in $\ov\cM_{g,k}(X,A)$.

\noindent
For every $\al\!\in\!S_f$, let 
$$
\cC_{0\al;v}\!\lra\!\cB_{0\al;v}, \quad \cC_{0\al}(\Gamma)\lra \cB_{0\al}(\Gamma), \quad\tn{and}\quad  \cC_{0\al} \lra \cB_{0\al}\!=\!\cB_{0\al}(\Gamma)\times \De^{E_{\Gamma}}
$$
be similarly defined objects for $C^{+0\al}$. Once again, by Condition~\ref{small-enough_c}, the universal family $\cC_\al$ of (\ref{equ:aux-data}) is also a universal family (of the standard form) around $C^{+\alpha}$, and the process of removing $w_0$ and $w_{\al,f}$ from $w_{0\al}$ (as in (\ref{equ:forget0al})) gives us forgetful maps 
\bEqu{equ:forget0al2}
\pi_0\colon (\cC_{0\al},\cB_{0\al})\lra (\cC'_{\al},\cB'_{\al}) \quad\tn{and}\quad \pi_\al\colon (\cC_{0\al},\cB_{0\al})\lra (\cC_{0},\cB_0), 
\eEqu
respectively, where $(\cC'_{\al}\!\lra\! \cB_\al')$ is a sufficiently small sub-family of $(\cC_{\al}\!\lra\!\cB_\al)$. 
For each $\al\!\in\!S_f$, the dual graph $\Gamma_\al$ of $f_\al$ is a degeneration of $\Gamma$ as in (\ref{equ:deg-graphs}); i.e. every neck region of $\cC'_\al$ (and the other two universal families) is also a neck region of the primary family $\cC_\al$.  Let $\cC'_\al(\Gamma)\subset\cC'_\al$ be the sub-family of curves with dual graph $\Gamma$; the forgetful maps $\pi_0$ and $\pi_\al$ take $\cC_{0\al}(\Gamma)$ to $\cC'_{\al}(\Gamma)$ and $\cC_{0}(\Gamma)$, respectively.

\noindent
For a choice of $\ep'_\al$, $\ep_0$, and $\ep_{0\al}$ as in (\ref{equ:epsilon-choices}),
let 
\bEqu{equ:Sets-0-al-0al}
B_{\ep'_\al}(f^{+\al}),\quad B_{\ep_0}(f^{+0}),\quad\tn{and} \quad B_{\ep_{0\al}}(f^{+0\al})
\eEqu
be the set of maps $\ep$-close to $f^{+\al}$, $f^{+0}$, and $f^{+0\al}$, as defined in Page~\pageref{l:ball-inequlity}, respectively.
The same push-pull procedure, as in the case of $f$ with smooth domain above, with (\ref{equ:forget0al2}) and (\ref{equ:Sets-0-al-0al}) in place of (\ref{equ:forget0al2}) and (\ref{equ:Banach-0-al-0al}), respectively, gives us an obstruction bundle 
\bEqu{finalObs_e}
\pr_0\colon E_{0}\!=\!\bigoplus_{\al\in S_f}  E_{0;\al} \lra B_{\ep'_0}(f^{+0})
\eEqu
as in (\ref{equ:formal-d-sum}). Define 
\bEqu{plus-chart-induced_e}
\cU^{+0}\!=\!(\pr^{+0}\colon\!U^{+0}\!\lra\! V^{+0},\aut(f),s^{+0},\psi^{+0})
\eEqu
as in 
(\ref{equ-induced-wtV})-(\ref{equ:foot-plus}). Similarly to Remark~\ref{rmk:not-Banach}, we can not immediately apply the Implicit Function Theorem to conclude that $V^{+0}$ is a smooth manifold of the expected dimension. We again need a gluing theorem to define a smooth structure on $V^{+0}$.

\noindent
Let $B_{\ep'_0}(f^{+0};0)\!\subset\!B_{\ep'_0}(f^{+0})$ as in (\ref{Bepf0_e}) be the subset of maps with dual graph the same as the dual graph of $f^{+0}$. This subspace lies in the fiber-product of Banach manifolds 
$$
W^{\ell,p}(\cC_{0;v},X)\qquad \forall v\!\in\! \Gamma
$$
over the evaluation map at the nodal points as in (\ref{equ:Fiber-Product}). The restriction 
$$
 \pr_{0;\Gamma}\colon E_{0;\Gamma}\!=\! E_0|_{B_{\ep'_0}(f^{+0};0)}\lra B_{\ep'_0}(f^{+0};0)
$$
is a smooth $\aut(f)$-equivariant vector bundle of finite rank as in the case of smooth domain. We define 
\bEqu{equ-induced-wtV0}
{V}^{+0}_\Gamma\!=\!\big\{ (f',(\eta_\al)_{\al\in S_f})\in E_{0;\Gamma} \colon \dbar f'-\sum_{\al \in S_f} \eta_\al=0\big\},
\quad {U}^{+0}_\Gamma\!=\!\pr_{0;\Gamma}^*E_{0;\Gamma}|_{{V}^{+0}_\Gamma}.
\eEqu
and ${s}^{+0}_\Gamma\colon{V}^{+0}_\Gamma\!\lra\!{U}^{+0}_\Gamma $ as in (\ref{equ-induced-wtV}) and (\ref{equ:wts0}). 
Restricted to the zero set of ${s}^{+0}_\Gamma$, we get the obvious footprint map 
$$
\psi^{+0}\colon ({s}^{+0})^{-1}(0)\lra \ov\cM_{g,k+|w_0|}(X,A)_\Gamma, \quad \psi^{+0}(f',(\eta_\al)_{\al\in S_f}) \!=\![f'].
$$
By Implicit Function Theorem,
for $\ep'_0$ and $\cC_0(\Gamma)$ sufficiently small, the projection map 
$\pr^{+0}_\Gamma\colon {U}^{+0}_\Gamma\lra {V}^{+0}_\Gamma$ is a smooth $\aut(f)$-equivariant vector bundle. Moreover, the evaluation map 
$$
\ev_{w_0}\colon {V}^{+0}_\Gamma\lra X^{w_0}
$$
is $\aut(f)$-invariant and transverse to $H_0$. 

\noindent
The fiber-product Kuranishi chart 
$$
\cU^{+0}_\Gamma\!\equiv\!(\pr^{+0}_\Gamma\colon\!{U}^{+0}_\Gamma\!\lra\! {V}^{+0}_\Gamma , \aut(f) , {s}^{+0}_\Gamma, \psi^{+0}_\Gamma),
$$
with footprint a neighborhood of $[f^{+0}]$ in\footnote{See Remark~\ref{rmk:orderf+}.} $\ov\cM_{g,k+|w_0|}(X,A)_\Gamma$ is the analogue of pre-gluing Kuranishi chart in 
(\ref{equ:bundle-0})-(\ref{equ:psif0}). Next, by a gluing procedure similarly to Section~\ref{sec:canonical-nodal},  we extend this chart to the entire $\cC_0$.
The following gluing theorem is the analogue of Theorem~\ref{thm:gluing-thm-main} for the induced charts. The only difference between Theorem~\ref{thm:gluing-thm-induced} and Theorem~\ref{thm:gluing-thm-main} is in how the obstruction bundles are defined; in the former, the obstruction bundle was defined via parallel translate of the obstruction space at the base map $f$; in the latter, the obstruction bundle is a direct sum of several bundles, each of which is obtained from parallel translate of the obstruction spaces $E_\al$ at  the primary maps $f_\al$.

\begin{theorem}[{\cite[Theorem 19.3]{FOOO-detail}}]\label{thm:gluing-thm-induced}
For $\ep'_0$ and $\cC_0$ sufficiently small, there exists a ``natural'' $\aut(f)$-equivariant continuous one-to-one map
$$
\mf{gl}\colon V_{\Gamma}^{+0}\times \De^{E_\Gamma} \lra V^{+0}
$$
such that the restriction of $\mf{gl}$ to $V_{\Gamma}^{+0}\!\subset\! V^{+0}$ is the inclusion map and the action of $\aut(f)\!\subset\!\aut(C)$ on $\De^{\Gamma_E}$ is given by (\ref{equ:gluing-action}).
\eThm

\noindent
Similarly to Theorem~\ref{thm:gluing-thm-smooth}, we will use the product smooth structure given by the weak smooth structure on $\De$ in (\ref{equ:T-weak}).
Finally, similarly to the case of smooth domain before, the evaluation map 
$$
\ev_{w_0}\colon {V}^{+0}\lra X^{w_0}
$$
is $\aut(f)$-invariant and transverse to $H_0$, and the \textbf{induced chart} 
\bEqu{equ:final-induced-chart}
\cU_f\!\equiv\!(\pr_f\colon\!U_f\!\lra\! V_f , \aut(f), s_f, \psi_f),
\eEqu
defined as in (\ref{equ:sliced-chart2})
is a Kuranishi chart of the expected dimension (\ref{equ:exp-dim}) centered at $[f]$. This finishes our construction of induced Kuranishi charts covering the entire moduli space $\ov\cM_{g,k}(X,A)$. We will discuss the existence of coordinate change maps (c.f. Definition~\ref{def:intersection}) in the next section.

\subsection{Coordinate change maps}\label{sec:transition}
\noindent
Assume $f_0\!=\!(u_0,C_0\!=\!(\Si_0,\mfj_0,\vec{z}_0))$ is a $J$-holomorphic map in $\ov\cM_{g,k}(X,A)$, $(w_0,H_0)$ is an arbitrary stabilizing pair for that, and $\cU_0\!=\!\cU_{f_0}(w_0,H_0)$ is the resulting induced Kuranishi chart of (\ref{equ:final-induced-chart}).
Let $f_1\!=\!(u_1,C_1\!=\!(\Si_1,\mfj_1,\vec{z}_1))$ be a $J$-holomorphic map in the footprint of $\cU_0$, $(w_1,H_1)$ be an arbitrary stabilizing pair for $f_1$, and $\cU_1\!=\!\cU_{f_1}(w_1,H_1)$ be the resulting induced Kuranishi chart of (\ref{equ:final-induced-chart}). In particular, if $f_1\!=\!f_0$, we can choose two different stabilizing pairs $(w_0,H_0)$ and $(w_1,H_1)$ which a priori could result in different charts $\cU_0$ and $\cU_1$ around the same point $[f_0]$.
If $f_1$ is sufficiently close to $f_0$, by (\ref{equ:Sf'inSf}), we have $S_{f_1}\!\subset \!S_{f_0}$.
In order to conclude that the induced charts of Section~\ref{sec:natural-KUR} form a Kuranishi structure in the sense of Definition~\ref{def:Kur-structure}, we show that there exists a natural coordinate change map from a sub-chart of $\cU_{1}$ into $\cU_{0}$, as in Definition~\ref{def:intersection}, independent of any choice except the data of primary charts and stabilizing pairs $\{(w_i,H_i)\}_{i=1,2}$. The cocycle condition property of Definition~\ref{def:cocycle} would naturally follow from the construction. Confirming tangent bundle condition is similar to Remark~\ref{rem:embedding}.  We refer to \cite[Section 22-24]{FOOO-detail} for more details on this construction and the properties of coordinate change maps.

\noindent
With notation as above, the points $w_0$, through the family $\cC_0$, give us a unique set of points $w_{1,0}\!\subset\! \Si_1$ at which $f_1$ transversely intersects the slicing divisors $H_0$. The stabilizing pairs  $(w_{1,0},H_0)$ and $(w_1,H_1)$ for $f_1$ could be different.
Let $\cU_{1,0}\!=\!\cU_{f_1}(w_{1,0},H_0)$ be an induced Kuranishi chart as in (\ref{equ:final-induced-chart}) corresponding to the inherited stabilization data $(w_{1,0},H_0)$ at $f_1$. The coordinate change map $\Phi_{10}\colon \!\cU'_{1}\!\lra\! \cU_0$ as we construct below, where $\cU'_1$ is a sufficiently small sub-chart of $\cU_1$ around $f_1$ as in Definition~\ref{def:intersection}, is a composition of 
$$
\Phi_{1,10}\colon \!\cU'_{1}\!\lra\! \cU'_{1,0}\quad \tn{and}\quad \Phi_{10,0}\colon \!\cU'_{1,0}\!\lra\! \cU_{0},
$$
where $\cU'_{1,0}$ is a sufficiently small sub-chart of $\cU'_{1,0}$ around $f_1$, $\Phi_{1,10}$ is a natural isomorphism of Kuranishi charts given by Lemma~\ref{lem:change-stab2} below, and $\Phi_{10,0}$ is a canonical embedding of Kuranishi charts given by Lemma~\ref{lem:same-stab-embedding} below. That these coordinate change maps satisfy the cocycle condition follows from Lemma~\ref{lem:cocyle-stab}.

\bLem{lem:same-stab-embedding}
Assume $f_0$ is a $J$-holomorphic map in $\ov\cM_{g,k}(X,A)$, $(w_0,H_0)$ is an arbitrary stabilizing pair for that, and $\cU_0\!=\!\cU_{f_0}(w_0,H_0)$ is the resulting induced Kuranishi chart of (\ref{equ:final-induced-chart}).
Let $f_1$ be a $J$-holomorphic map in the footprint of $\cU_0$ sufficiently close to $f_0$. If $(w_1,H_1)\!=\!(w_{1,0},H_0)\!$, with $w_{1,0}$ and $\cU_1\!=\!\cU_{f_1}(w_1,H_1)$ defined as above, there exists a canonical coordinate change map centered at $f_1$,
$$
\Phi_{10}\colon \cU'_{1}\lra \cU_0,
$$
from a sub-chart $\cU'_{1}$ of $\cU_{1}$ into $\cU_0$.
\eLem

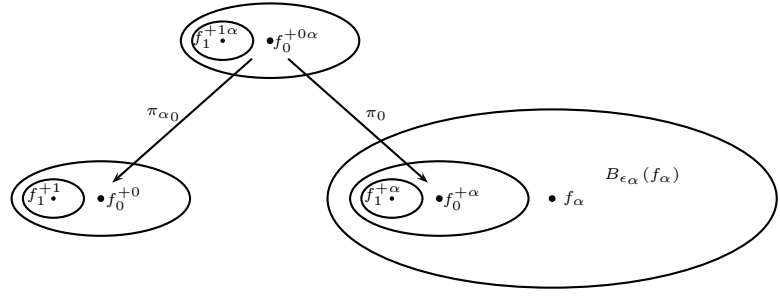
\begin{figure}
\begin{pspicture}(-5,-5.5)(10,-4)
\psset{unit=.3cm}

\psellipse(15,-15)(10,4)\pscircle*(15,-15){.15}\rput(16,-15){\tiny{$f_\al$}}\rput(19,-14){\tiny{$B_{\ep_\al}(f_\al)$}}
\psellipse(10,-15)(4,1.7) \pscircle*(10,-15){.15}\rput(11,-15){\tiny{$f_0^{+\al}$}}
\psellipse(7.9,-15)(1.4,.9)  \pscircle*(7.9,-15){.1} \rput(7.5,-14.8){\tiny{$f_1^{+\al}$}}

\psellipse(2.5,-8)(4,1.7) \pscircle*(2.5,-8){.15}\rput(3.65,-8){\tiny{$f_0^{+0\al}$}}
\psellipse(.4,-8)(1.4,.9)  \pscircle*(.4,-8){.1} \rput(0.1,-7.8){\tiny{$f_1^{+1\al}$}}

\psellipse(-5,-15)(4,1.7) \pscircle*(-5,-15){.15}\rput(-4,-15){\tiny{$f_0^{+0}$}}
\psellipse(-7.1,-15)(1.4,.9)  \pscircle*(-7.1,-15){.1} \rput(-7.5,-14.8){\tiny{$f_1^{+1}$}}

\psline{->}(3.3,-8.8)(9.5,-14.3)\rput(7.2,-11.3){\tiny{$\pi_0$}}
\psline{->}(1.7,-8.8)(-4.5,-14.3)\rput(-2.2,-11.3){\tiny{$\pi_{\al_0}$}}

\end{pspicture}
\caption{The diagram corresponding to the proof Lemma~\ref{lem:same-stab-embedding}.}
\label{fig:inclusions}
\end{figure}

\bProof
Since by assumption $w_1\!=\!w_{1,0}$, with notation as in Section~\ref{sec:natural-KUR},  the germ of universal family $\cC_{1}$ corresponding to $(\Si_1,\mfj_1,\vec{z}_1,w_{1})$ is a subfamily of the universal family $\cC_0$ corresponding to $(\Si_0,\mfj_0,\vec{z}_0,w_{0})$. 
Therefore, for $\ep'_{1}\!<\!\!<\!\ep'_0$, with notation as in Page~\pageref{l:ball-inequlity}, we obtain an embedding $B_{\ep'_{1}}(f_1^{+1})\!\subset\!B_{\ep'_0}(f_0^{+0})$, where 
$$
f_0^{+0}\!=\!(u_0,(\Si_0,\mfj_0,\vec{z}_0,w_{0}))\quad \tn{and}\quad f_1^{+1}\!=\!(u_1,(\Si_1,\mfj_1,\vec{z}_1,w_{1}));
$$ 
this inclusion is illustrated in the left-bottom corner of Figure~\ref{fig:inclusions}.
By assumption $S_{f_1}\!\subset\!S_{f_0}$; therefore, similarly, with notation as in (\ref{equ:forget0almaps}), for every $\al\!\in\!S_{f_1}$, $\ep_{1\al}\!<\!\!<\!\ep_{0\al}$ sufficiently small, and $\ep_{\al_1}\!<\!\!<\!\ep_{\al_0}\!<\!\!<\!\ep_\al$ we have\footnote{$f_i^{+i\al}\!=\!(u_i,(\Si_i,\mfj_i,\vec{z}_i,w_{i}\cup w_{\al,f_i}))$, with $i\!=\!0,1$. }
$$
B_{\ep_{1\al}}(f_1^{+1\al})\!\subset\!B_{\ep_{0\al}}(f_0^{+0\al})\quad \tn{and}\quad B_{\ep_{\al_1}}(f_1^{+\al})\!\subset\!B_{\ep_{\al_0}}(f_0^{+\al})\!\subset\!B_{\ep_\al}(f_\al);
$$ 
These inclusions are illustrated in the top and right-bottom corner of Figure~\ref{fig:inclusions}.
The projection maps 
$$
\pi_0\colon\! B_{\ep_{0\al}}(f^{+0\al}_0)\!\lra\! B_{\ep_{\al_0}}(f^{+\al}_0)\quad\tn{and}\quad 
\pi_1\colon\! B_{\ep_{1\al}}(f^{+1\al}_0)\!\lra\! B_{\ep_{\al_1}}(f^{+\al}_1)$$
and 
$$
\pi_{\al_0}\colon\! B_{\ep_{0\al}}(f^{+0\al}_0)\!\lra\! B_{\ep_0}(f_0^{+0})\quad\tn{and}\quad 
\pi_{\al_1}\colon\! B_{\ep_{1\al}}(f^{+1\al}_0)\!\lra\! B_{\ep_1}(f_1^{+1})
$$
as in (\ref{equ:forget0almaps}) commute with these emebeddings. Therefore,
for every $\al\!\in\!S_{f_1}$, with
$$
E_{0\al}\lra B_{\ep_{0\al}}(f_0^{+0\al}) \quad \tn{and}\quad E_{1\al}\lra B_{\ep_{1\al}}(f_1^{+1\al})
$$
as in (\ref{equ: cE_0al}), we have 
$$
E_{1\al}\!=\!E_{0\al}|_{B_{\ep_{1\al}}(f_1^{+1\al})}\;.
$$
Since we use the same slicing conditions $H_\al$ to define the section
$s_\al$ in (\ref{equ:sal}), we have
$$
E_{1;\al}\!=\!E_{0;\al}|_{B_{\ep'_1}(f_1^{+1})}.
$$
With notation as in (\ref{equ:formal-d-sum}), we conclude that $E_1\!=\!\bigoplus_{\al \in S_1} E_{1;\al}$ embeds in $E_{0}|_{B_{\ep'_1}(f_1^{+1})}$.
Therefore, with notation as in (\ref{plus-chart-induced_e}), the germ of\footnote{Because we are possibly restricting to smaller values of $\ep$.} pre-slicing Kuranishi chart $\cU^{+1}$ corresponding to $E_1$ embeds in $\cU^{+0}$. Since we apply the same slicing condition $H_0$ to both $\cU^{+1}$ and  $\cU^{+0}$ to obtain $\cU_1$ and  $\cU_{0}$, we get a canonical embedding of the germ of $\cU_{1}$ around $f_1$ into  $\cU_{0}$, as defined in (\ref{equ:sliced-chart2}).
\eProof

\bLem{lem:change-stab}
Assume $f$ is a $J$-holomorphic map in $\ov\cM_{g,k}(X,A)$, $(w_1,H_1)\!\subset\!(w_0,H_0)$ are two stabilizing pairs for $f$, and $\cU_0\!=\!\cU_{f}(w_0,H_0)$ and $\cU_1\!=\!\cU_{f}(w_1,H_1)$ are the resulting induced Kuranishi chart of (\ref{equ:final-induced-chart}).
Then the forgetful map $\cU_0\!\lra\!\cU_1$ is well-defined and is a canonical isomorphism. 
\eLem

\bProof
The claim follows from the commutativity of the following diagram:
\bEqu{0to1stab_e}
\xymatrix{
			              & \ar[ld]^{\pi_\al}B_{\ep_{0\al}}(f^{+0\al}) \ar[rdd]^{\pi_0} \ar[d]^{\pi_{0/1}}&  \\
B_{\ep_{0}}(f^{+0})\ar@<3pt>[ru]^{s_\al}\ar[d]^{\pi_{0/1}} & \ar@<3pt>[u]^{s_{0/1}}\ar[ld]^{\pi_\al}B_{\ep_{1\al}}(f^{+1\al}) \ar[rd]^{\pi_1}&  \\
B_{\ep_{1}}(f^{+1})\ar@<3pt>[ru]^{s_\al} \ar@<3pt>[u]^{s_{0/1}}&  				     & B_{\ep_{\al}}(f^{+\al}) \;,
}
\eEqu
where $\pi_0$, $\pi_1$, and $\pi_\al$ are the forgetful maps of (\ref{equ:forget0almaps}), $\pi_{0/1}$ is the forgetful map of removing $w_0/w_1$ stabilizing points, and $s_\al$ and $s_{0/1}$ are the section maps of (\ref{equ:sal}).
\eProof

\noindent
For $i\!=\!1,2$, let $(w_i,H_i)$ be two arbitrary stabilization data for a fixed map $f$. Let $\{\cU_i\!=\!\cU(w_i,H_i)\}_{i=1,2}$ be the resulting induced Kuranishi charts of (\ref{equ:final-induced-chart}). Another arbitrary stabilizing pair $(w_\al,H_\alpha)$ for $f$ disjoint from $\{(w_i,H_i)\}_{i=1,2}$ gives us larger stabilizing pairs
$$
(w_{\al i},H_{\al i})\!=\!(w_i,H_i)\cup (w_\al,H_\al)\qquad \forall i\!=\!1,2,
$$ 
and similarly defined induced Kuranishi charts 
$$
\{\cU(w_{\al i},H_{\al i})\}_{i=1,2}\quad\tn{and} \quad\{\cU(w_{\al},H_{\al})\}_{i=1,2}.
$$
By Lemma~\ref{lem:change-stab}, for $i\!=\!1,2$, after possibly restricting to sub-charts 
$$
\cU'_i\!\subset\! \cU_i, \quad \cU'(w_{\al i},H_{\al i})\!\subset\! \cU(w_{\al i},H_{\al i}),\quad \tn{and}\quad \cU'(w_{\al},H_{\al})\!\subset\!\cU(w_{\al},H_{\al})
$$ 
around $f$, the forgetful maps obtained by removing the $w_\al$ and $w_i$ points, respectively, give us isomorphisms
\bEqu{prep-Phi_e}
\Phi_{\al i,i},\Phi_{\al i,\al}\colon \cU'(w_{\al i},H_{\al i})\lra \cU'_i,\cU'(w_\al,H_\al)
\eEqu
\bLem{lem:change-stab2}
With notation as in (\ref{prep-Phi_e}), the orbibundle isomorphism
\bEqu{equ:natural embedding}
\Phi_{12}\!=\!\Phi_{\al2,2}\circ \Phi_{\al2,\al}^{-1}\circ  \Phi_{\al1,\al}\circ \Phi_{\al1,1}^{-1}\colon \cU'_{1}\lra \cU'_2,
\eEqu
is independent of the choice of $(w_\al,H_\al)$; see Figure~\ref{fig:wedge-map}.
\eLem

\bProof
The claim follows from the commutativity of a diagram similar to (\ref{0to1stab_e}): for two non-overlapping choices of $(w_\al,H_\al)$, say $(w_{\al_1},H_{\al_1})$ and $(w_{\al_2},H_{\al_2})$, we consider the union the claim follows from Lemma~\ref{lem:change-stab}. The general case follows from continuity (or a bigger commutative diagram).
\eProof

\begin{figure}
\begin{pspicture}(-5,-5.5)(10,-4)
\psset{unit=.3cm}

\rput(18,-12){$\cU(w_{\al2},H_{\al2})$}
\rput(-5,-12){$\cU(w_{\al1},H_{\al1})$}
\rput(6,-12){$\cU(w_{\al},H_{\al})$}
\rput(-7,-18){$\cU_1$}
\rput(19,-18){$\cU_2$}
\psline{<-}(-7,-17)(-7,-13)
\psline{->}(-1,-12)(2,-12)
\psline{<-}(9.5,-12)(13,-12)
\psline{->}(19,-13)(19,-17)
\psline[linestyle=dashed,dash=3pt]{->}(-5,-18)(17,-18)
\rput(6,-17){$\Phi_{12}$}
\end{pspicture}
\caption{Diagram of natural embedding (\ref{equ:natural embedding}).}
\label{fig:wedge-map}
\end{figure}

\noindent
Consequently, we call  $\Phi_{12}$ the \textbf{natural isomorphism} between the germ of Kuranishi charts $\cU_1$ and $\cU_2$.

\bLem{lem:cocyle-stab}
For $i\!=\!1,2,3$, let $(w_i,H_i)$ be arbitrary stabilizing pairs for a fixed map $f$. Then the natural isomorphisms $\{\Phi_{ij}\}_{1\leq i<j\leq 3}$ satisfy $\Phi_{13}\!=\!\Phi_{23}\circ \Phi_{12}$. 
\eLem

\begin{proof}
Choose $(\om_\al,H_\al)$ disjoint from each of $(w_i,H_i)_{i=1,2,3}$. Let 
$$
\Phi_{ij}\!=\!\Phi_{\al;ij}\quad \forall 1\!\leq \!i\!<\!j\leq\! 3
$$ 
be the isomorphims of (\ref{equ:natural embedding}). Then
$$
\aligned
&\Phi_{23}\!\circ\! \Phi_{12} =
\big(\Phi_{\al3,3}\!\circ\! \Phi_{\al3,\al}^{-1}\!\circ\! \iota_{23}\!\circ\! \Phi_{\al2,\al}\!\circ\! \Phi_{\al2,2}^{-1}\big)\circ
\big(\Phi_{\al2,2}\!\circ\! \Phi_{\al2,\al}^{-1}\!\circ\! \iota_{12}\!\circ\! \Phi_{\al1,\al}\!\circ\! \Phi_{\al1,1}^{-1}\big)=\\
&
\Phi_{\al3,3}\!\circ\! \Phi_{\al3,\al}^{-1}\!\circ\! (\iota_{23} \circ
 \iota_{12})\!\circ\! \Phi_{\al1,\al}\!\circ\! \Phi_{\al1,1}^{-1} =
 \Phi_{\al3,3}\!\circ\! \Phi_{\al3,\al}^{-1}\!\circ\! \iota_{13}\!\circ\! \Phi_{\al1,\al}\!\circ\! \Phi_{\al1,1}^{-1}=\Phi_{13}.
\endaligned
$$
\end{proof}

\begin{remark}
Recall from Remark~\ref{rmk:No-Comparison} that similar statements do not hold for primary Kuranishi charts; thus, primary charts do not admit coordinate change maps.
This is essentially the main reason for considering the more complicated induced charts to put a Kuranishi structure on $\ov\cM_{g,k}(X,A)$. Since the obstruction bundle at a particular point $[f]\!\in\!\ov\cM_{g,k}(X,A)$ does not depend on the choice of the induced chart covering that point, it enables us to compare different charts corresponding to different stabilization data and other choices.
\end{remark}

\subsection{GW invariants}\label{sec:GW-VFC}

Let $(X^{2n},\om)$ be a closed symplectic manifold, $J\!\in\! \cJ(X,\om)$ be an arbitrary compatible (or tame) almost complex structure, $A\!\in\! H_2(X,\Z)$, and $g,k\!\in\! \Z^{\geq 0}$. With $d$ as in (\ref{equ:exp-dim}), 
the construction of last section gives us a class of natural oriented $d$-dimensional Kuranishi structures on $\ov{\cM}_{g,k}(X,A,J)$. Thus it establishes the first statement of Theorem~\ref{thm:Nat-Kur}. Let 
\bItem
\item $\om_0$ and $\om_1$ be isotopic symplectic structures on $X$, 
\item $J_0\!\in\! \cJ(X,\om_0)$ and  $J_1\!\in\! \cJ(X,\om_1)$, 
\item $\cP_0$ and $\cP_1$ be arbitrary primary collections for $\ov{\cM}_{g,k}(X,A,J_0)$ and $\ov{\cM}_{g,k}(X,A,J_1)$, 
\item and $\cK_0$ and $\cK_1$ be the resulting Kuranishi structures on $\ov{\cM}_{g,k}(X,A,J_0)$ and $\ov{\cM}_{g,k}(X,A,J_1)$.
\eItem

\noindent
Choose a path $J_{[0,1]}\!\equiv \!\{J_t\}_{t\in[0,1]}$ of smooth compatible almost complex structures on $(X,\om)$ connecting $J_0$ to $J_1$. For simplicity, we choose a path $J_{[0,1]}$ such that $J_t$ is equal to constant $J_0$ on $[0,\ep]$ and is equal to $J_1$ on $[1\!-\!\ep,1]$.
Let $\ov{\cM}_{g,k}(X,A,J_{[0,1]})$ be the moduli space of all genus $g$ degree $A$ $k$-marked $J_t$-holomorphic maps (\ref{equ:famlyM}). When equipped with the Gromov topology, this is again a compact metrizable topology space with boundary 
$$
\partial\ov{\cM}_{g,k}(X,A,J_{[0,1]})\!=\!\ov{\cM}_{g,k}(X,A,J_0)\cup \ov{\cM}_{g,k}(X,A,J_1).
$$
Every chart $\cU_{f_\al}$ in $\cP_0$ trivially extends to a Kuranishi chart with boundary $\cU_{f_\al}\!\times\![0,\ep)$ covering  a neighborhood of $f_{\al}\times \{J_0\}$ in $\ov{\cM}_{g,k}(X,A,J_{[0,1]})$. Similarly, $\cP_1$ extends to a primary collection covering a neighborhood of $\ov{\cM}_{g,k}(X,A,J_1)\!\subset\!\ov{\cM}_{g,k}(X,A,J_{[0,1]})$. By adding more primary Kuranishi charts covering a neighborhood of the middle part
$$
\ov{\cM}_{g,k}(X,A,J_{[\ep,1-\ep]}),
$$
we can extend the primary collections $\cP_0$ and $\cP_1$ to a primary collection $\cP$ of Kuranishi charts which covers the entire $\ov{\cM}_{g,k}(X,A,J_{[0,1]})$. The construction of previous section then gives a Kuranishi structure with boundary $\cK$ on $M\!\times\![0,1]$ whose restriction to $M\!\times\!\{0\}$ and $M\!\times\!\{1\}$ coincides with $\cK_0$ and $\cK_1$, respectively. This gives us a cobordism between $\cK_0$ and $\cK_1$ in the sense of Definition~\ref{def:cobo-Kur} and thus establishes the second part of Theorem~\ref{thm:Nat-Kur}.

\noindent
The evaluation maps $\ev$ and forgetful map $\st$ in (\ref{equ:eval}) and (\ref{equ:st}), respectively, readily lift to continuous maps on every one of the Kuranishi charts constructed in the previous sections. These lifts are compatible with coordinate change maps.  We conclude that all the natural Kuranishi structures of Theorem~\ref{thm:Nat-Kur} are indeed $(\ev\!\times\!\st)$-Kuranishi structures in the sense of Definition~\ref{def:rho-Kur-structure}. Together with Theorem~\ref{thm:Nat-Kur} and Proposition~\ref{pro:VFC-map}, we obtain a singular homology class 
$$
[\ov{\cM}_{g,k}(X,A,J)]^\vfc\in H_d(X^k\!\times\! \ov{\cM}_{g,k},\Q)
$$ 
which only depends on the isotopy class of $(X,\om)$. This establishes Theorem~\ref{thm:VFC}. Gromov-Witten invariants with primary insertions are defined by integration of cohomology classes in $X^k$ and $\ov{\cM}_{g,k}$ against this homology class.

\noindent
In order to define GW invariants involving $\psi$ and Hodge classes, defined in Section~\ref{sec:VFC-intro}, we start with a Kuranishi structure $\cK$ provided by Theorem~\ref{thm:Nat-Kur}. The $\psi$ and Hodge classes in the question define a vector bundle $\mf{E}$ on $\cK$ in the sense of Definition~\ref{def:Kur-bundle}. This Kuranishi vector bundle gives us an augmented Kuranishi structure $\cK^{\mf{E}}$ in the sense of  (\ref{equ:Augmented-chart}). The augmented Kuranishi structure $\cK^{\mf{E}}$ is still an $(\ev\!\times\!\st)$-Kuranishi structure. Every two of such $\cK^{\mf{E}}$ are still cobordant. Let 
$$
[\mf{E}\!\to\! \ov{\cM}_{g,k}(X,A,J)]^\vfc:=[\cK^{\mf{E}}]^\vfc \in H_d(X^k\!\times\! \ov{\cM}_{g,k},\Q)
$$
denote the resulting virtual fundamental class of Section~\ref{sec:VFC2}.  Gromov-Witten invariants involving a particular class of type $\mf{E}$ and  primary insertions are defined by integration of cohomology classes in $X^k$ and $\ov{\cM}_{g,k}$ against $[\mf{E}\!\to\! \ov{\cM}_{g,k}(X,A,J)]^\vfc$. They only depend on the isotopy class of $(X,\om)$ and the type of $\mf{E}$. We will not further go into establishing the functorial properties of these virtual fundamental classes corresponding to tensor product, direct sum, etc. of these bundles.

\section{Examples}\label{sec:examples}
In this section, we go over some examples of moduli spaces of pseudoholomorphic maps with pure orbibundle Kuranishi structure.

\subsection{Degree zero maps}\label{subsec:degree0}
Let $(X^{2n},\om)$ be a non-trivial symplectic manifold and $g,k\!\in \!\Z^{\geq 0}$, with $2g\!+\!k\!\geq \!3$.  Every pseudoholomorphic map with trivial homology class $0\!\in\! H_2(X,\Z)$ is a constant map. Therefore, $\ov\cM_{g,k}(X,0)$ (which is independent of the choice of $J$) is homeomorphic to
\bEqu{actual}
\ov\cM_{g,k}(X,0)\!=\!X\times \ov\cM_{g,k}.
\eEqu
In this case, for every $f\!=\!(u,C\!=\!(\Si,\mfj,\vec{z}))\!\in\!\ov\cM_{g,k}$, we have  
$$
\tn{Im}(u)\!\equiv\!p\!\in\!X, \quad u^*TX\!=\!\Si\times T_pX,
$$ 
and 
$$
\tn{D}_u\dbar(\ze)\!=\!\dbar \ze \qquad \forall \ze\colon \Si\lra T_pX\cong \C^n.
$$
In this situation we get 
\bEqu{equ:Obf}
\tn{Obs}(u)\!\cong\!\tn{Obs}(f)\!\cong\!T_pX \otimes_\C H^1(\cT\Si_\mfj).
\eEqu
see Remark~\ref{rem:actual-Obs}. Therefore, if $g\!=\!0$, then $\ov\cM_{g,k}(X,0)$ has the structure of a smooth orbifold given by the right-hand side of (\ref{actual}). Let $\cE_{g}$ be the Hodge orbibundle over $\ov\cM_{g,k}$ defined in (\ref{Hodge-fiber_e}).
For $g\!\geq\! 1$, 
$$
\dim^\vir_\C \ov\cM_{g,k}(X,0)= \dim_\C\!X(1\!-\!g)\!+\! \dim_\C \ov\cM_{g,k}
$$ 
and the right-hand side of (\ref{equ:Obf}) defines an orbibundle 
\bEqu{equ:pure-orbibundle}
TX\boxtimes \cE_{g} \lra X\times \ov\cM_{g,k}
\eEqu
whose fiber at every point, by (\ref{equ:Obf}), is the minimal possible obstruction space; see Remark~\ref{rmk:minimal-obs}. Therefore, (\ref{equ:pure-orbibundle}) with Kuranishi map $s\!\equiv\!0$ defines a pure orbibundle Kuranishi structure on $\ov\cM_{g,k}(X,0)$ whose top chern class is the VFC of Section~\ref{sec:VFC}; see Example~\ref{exa:global-kur}.

\noindent
Note that if $(g,k)\!=\!(1,1)$ or $(2,0)$, then the orbifold structure of $\ov\cM_{g,k}$ is not effective. 
If $n\!>\!0$, we may consider a larger obstruction space to obtain an effective orbifold. 

\noindent
As an example, if $\dim_\C X\!=\!3$, $k\!=\!0$, and $g\!\geq\! 2$, then the Gromov-Witten VFC is a zero-dimensional class, thus a rational number, and it is equal to 
$$
\chi(X) \int_{\ov\cM_{g,0}}c_\top(\mc{E}_{g})=\chi(X) \frac{B_{2g}}{4g(g-1)},
$$
where $B_n$ is the $n$-th Bernoulli number.

\subsection{Elliptic surfaces}\label{subsec:elliptic}

Let $(X,\om,J)$ be a K\"ahler surface which has the structure of an elliptic fiberation $\pi:X\!\lra\! S$ over a smooth curve of genus $g$.
Moreover, assume that $\pi$ has no singular fiber other than multiple-fibers of order $(m_i)_{i\in [\ell]}$ over $(p_i)_{i\in [\ell]}\subset S$.
Let
$$
\rho(m)=\dpst\sum_{a|m} a\qquad \forall m\in \N.
$$
By \cite[Theorem 18.2]{BHPV}, the complex structure of reduced fibers are all the same; say equal to some fixed $\mfj \!\in\! \H/\tn{PSL}(2,\Z)$. 
Let $\cS$ denote the orbifold structure on $S$ given by replacing a neighborhood of each $p_i$ with a cyclic $m_i$-quotient of a disk as in Section~\ref{sec:orbifold}.
The moduli space $\ov\cM_{1,0}(X,dF)$, where $F$ is a generic fiber class and $d\in \Z^+$, has virtual dimension zero. Therefore, the Gromov-Witten virtual fundamental class $[\ov\cM_{1,0}(X,dF)]^\vfc$ is a zero-dimensional cycle in $H_0(X^0\times \cM_{1,0},\Q)\!=\!H_0(\tn{pt},\Q)\!\cong\! \Q$; we denote by $N_d$ to be the degree of this cycle.
In the simplest case of $d=1$, the moduli space $\ov\cM_{1,0}(X,F)$ consists of 
\bItem
\item a unique map over every point of $S\!-\!\{p_i\}_{i \in [\ell]}$,
\item and $\rho(m_i)$ different $m_i$-covering maps $f_{ij}$ of the reduced fiber over $p_i$, for every $i\!\in\! [\ell]$, only one of which, say $f_{i1}$, is in the limit of the first group of maps.
\eItem
Therefore,  $\ov\cM_{1,0}(X,F)$ is homeomorphic to the disjoint union
$$
S \cup \dpst\bigcup_{\substack{i\in [\ell]\\ 2\leq j\leq  \rho(m_i)}} [f_{ij}];
$$
this is an example of a moduli space with connected components of various dimension.
Over the connected component identified with $S$, the obstruction space $\tn{Obs}(f_s)$ at every map $f_s$, with $s\!\in\! S$, has constant rank and is isomorphic to 
$$
H^0(\Om^{0,1}_{\Si_\mfj})\otimes T_sS.
$$ 
The automorphism group of every holomorphic map over $S\!-\!\{p_i\}_{i \in [\ell]}$ is trivial and $f_{ij}$ has an automorhism group of order $m_i$. Therefore, over the component identified with $S$, the moduli space has the structure of a pure orbibundle Kuranishi structure 
$$
H^0(\Om_{\Si_\mfj}^{0,1})\otimes T\cS \cong T\cS \lra \cS
$$ 
with the Kuranishi map $s\!\equiv\! 0$.
On the otherhand, each of the holomorphic maps $f_{ij}$, with $i\!\in\! [\ell]$ and  $2\!\leq\! j\!\leq\!  \rho(m_i)$, is regular (i.e. $\tn{Obs}(f)\!=\!0$). Therefore, the simplest natural Kuranishi structure on every $[f_{ij}]$ is a non-effective orbifold structure on a  single point (with trivial obstruction bundle and) with an isotropy group of order $m_i$. Similarly to before, in order to get an effective orbifold, we may enlarge the obstruction bundle to take care of the non-effective action. Nevertheless, the contribution of every $[f_{ij}]$ to $N_1$ is $1/m_i$. 
We conclude that 
$$
N_1\!=\!\chi(\cS)+\sum_{i=1}^\ell \frac{\rho(m_i)-1}{m_i}=(2-2g-\ell)+\sum_{i=1}^\ell \frac{\rho(m_i)}{m_i}\in \Q.
$$

\subsection{Genus zero maps in quintic}\label{subsec:quintic}

In this section, we consider the example of genus $0$ maps in a quintic Calabi-Yau $3$-fold in $\P^4$. 

\noindent
Let $(X,\om,J)$ be a smooth degree $5$ divisor in $\P^4$ (i.e. $\om$ is the restriction of Fubini-Study K\"ahler form and $J$ is the complex structure). 
By Lefschetz hyperplane theorem, $H_2(X,\Z)\!\cong \!H_2(\P^4,\Z)\!\cong\! \Z$. For every $[d]\!\in\! H_2(X,\Z)$, 
the moduli space $\ov\cM_{0,0}(X,[d])$ has virtual dimension zero. Let
$$
N_d\!=\![\ov\cM_{0,0}(X,[d])]^\vfc\!\in\! H_0(X^0\times \cM_{0,0},\Q)\!\cong\! H_0(\tn{pt},\Q)\!\cong \!\Q\quad \forall d\!\in\! \N
$$
be the corresponding  GW invariants.

\noindent
Unlike previous examples, an explicit description of  the moduli space and of natural Kuranishi structures
is beyond the reach. However, for this (class of) example(s), we can construct an abstract pure orbibundle Kuranishi structure which a priori has no trivial connection to the natural Kuranishi structures of Section~\ref{sec:natural-KUR}. We first recall the construction of this particular abstract pure orbibundle Kuranishi structure and then discuss the known facts about its relation to the natural Kuranishi structures of Section~\ref{sec:natural-KUR}. 

\noindent
Suppose $X$ is the zero set of a degree $5$ homogenous polynomial $F\!\in\!H^0(\cO_{\P^4}(5))$, i.e. we think of $F$ as a section of the line bundle $\cO_{\P^4}(5)$. Every holomorphic map into $X$ can be thought of as a holomorphic map into $\P^4$ (of the same degree and genus). By \cite[Exercise 24.3.4(7)]{mirror}, $\ov\cM_{0,0}(\P^4,[d])$ has the structure of a smooth complex orbifold of complex dimension $5d$. A holomorphic map $f\!=\![u,\Si,\mfj]\!\in\! \ov\cM_{0,0}(\P^4,[d])$ belongs to $\ov\cM_{0,0}(X,[d])$ if and only if the pullback section $u^*F\!\in\! H^0(u^*\cO_{\P^4}(5))$ vanishes.
The  complex vector space $H^0(u^*\cO_{\P^4}(5))$ is of constant rank $5d$. In fact, there exists a rank $5d$ orbibundle 
 \bEqu{E0d}
 \cE_{0,d}\lra \ov\cM_{0,0}(\P^4,[d]),
 \eEqu 
with the fiber $[H^0(u^*\cO_{\P^4}(5))/ \aut(f)]$ over each $f$, and 
$$
s_F\colon\ov\cM_{0,0}(\P^4,[d]) \lra \cE_{0,d},\quad s_F(f)\!=\!u^*F 
 $$
is a section of this bundle; see \cite[Section 26.1.3]{mirror}.
The conclusion is that the orbibundle (\ref{E0d}) together with the orbifold section $s_F$ defines an ``abstract" zero-dimensional pure orbibundle  Kuranishi structure on $\ov\cM_{0,0}(X,[d])$. Let $N_d'$ be the orbifold top-chern number of this orbibundle; see Section~\ref{sec:euler}.

\noindent
Then the main claim is that 
\bEqu{equ:NDND}
N_d\!=\!N_d'\;,
\eEqu 
i.e. the VFC obtained from the abstract pure orbibundle Kuranishi structure above coincides with that of Section~\ref{sec:natural-KUR}. The numbers $N_d'$ can be calculated via localization technique; see \cite[Chapter 29.1]{mirror}. Therefore, the equality (\ref{equ:NDND}) provides a way of calculating genus $0$ GW invariants of quintic threefold (as well as several other similar examples).  This trick does not extend beyond genus $0$; see \cite{Zig1} for the case of genus $1$ invariants and the higher genus case is still a field of active research.

\noindent
We finish this section and thus this manuscript with few words on the proof of (\ref{equ:NDND}). 
In the algebraic case, this is stated as \cite[Theorem 26.1.1]{mirror} without proof. Cox-Katz \cite[Example 7.1.5.1]{CK} outlines a scheme of proving that via the algebraic approach of Li-Tian \cite{LT}. In an argument aimed at the moduli spaces of pseudoholomorphic disks which also applies to this example, \cite[Section 5.3]{PSW} outlines a way of proving (\ref{equ:NDND}) by extending a Kuranishi structure on $\ov\cM_{0,0}(X,[d])$ to a Kuranishi structure on $\ov\cM_{0,0}(\P^4,[d])$ and comparing that with the trivial Kuranishi (orbifold) structure of $\ov\cM_{0,0}(\P^4,[d])$. However, a complete proof in the analytic setting has yet to be written in details.

\vspace{.2in}

\begin{multicols}{2}

\noindent
{\it Mohammad Farajzadeh-Tehrani, \\
mtehrani@scgp.stonybrook.edu}

\noindent
{\it Kenji Fukaya,\\
kfukaya@scgp.stonybrook.edu}
\end{multicols}

\noindent
{\it Simons Center for Geometry and Physics, Stony Brook University}


\printindex
\bibliographystyle{plain}
\bibliography{GWbib}

\end{document}